\documentclass[12pt,a4,oneside]{amsart}

\usepackage{mathtext,amsmath,amsfonts,amssymb,amsthm}
\usepackage{verbatim}
\usepackage{enumerate}
\usepackage{preview}

\usepackage{version,etex}
\usepackage[top=1in, bottom=1in,left=1.5in, right=0.5in]{geometry}
\usepackage[pdftex]{graphicx,color}
\usepackage{epsfig}
\usepackage{hyperref}
\usepackage{float}
\usepackage{subfig}

\usepackage{color}

\reversemarginpar

\catcode`@=11 \catcode`!=11

\expandafter\ifx\csname fiverm\endcsname\relax
  \let\fiverm\fivrm
\fi
  
\let\!latexendpicture=\endpicture 
\let\!latexframe=\frame
\let\!latexlinethickness=\linethickness
\let\!latexmultiput=\multiput
\let\!latexput=\put
 
\def\@picture(#1,#2)(#3,#4){%
  \@picht #2\unitlength
  \setbox\@picbox\hbox to #1\unitlength\bgroup 
  \let\endpicture=\!latexendpicture
  \let\frame=\!latexframe
  \let\linethickness=\!latexlinethickness
  \let\multiput=\!latexmultiput
  \let\put=\!latexput
  \hskip -#3\unitlength \lower #4\unitlength \hbox\bgroup}

\catcode`@=12 \catcode`!=12

\input pictex.tex

\catcode`@=11 \catcode`!=11
  
\let\!pictexendpicture=\endpicture 
\let\!pictexframe=\frame
\let\!pictexlinethickness=\linethickness
\let\!pictexmultiput=\multiput
\let\!pictexput=\put

\def\beginpicture{%
  \setbox\!picbox=\hbox\bgroup%
  \let\endpicture=\!pictexendpicture
  \let\frame=\!pictexframe
  \let\linethickness=\!pictexlinethickness
  \let\multiput=\!pictexmultiput
  \let\put=\!pictexput
  \let\input=\@@input   
  \!xleft=\maxdimen  
  \!xright=-\maxdimen
  \!ybot=\maxdimen
  \!ytop=-\maxdimen}

\let\frame=\!latexframe

\let\pictexframe=\!pictexframe

\let\linethickness=\!latexlinethickness
\let\pictexlinethickness=\!pictexlinethickness

\let\\=\@normalcr
\catcode`@=12 \catcode`!=12

\def\endpfclaim{$\checkmark$}  

\begin{document}
\setcounter{secnumdepth}{4}
\setcounter{tocdepth}{1}
\newcommand{\proclaim}[2]{\medbreak {\bf #1}{\sl #2} \medbreak}

\newcommand{\ntop}[2]{\genfrac{}{}{0pt}{1}{#1}{#2}}

\let\newpf\proof \let\proof\relax \let\endproof\relax
\newenvironment{pf}{\newpf[\proofname]}{\qed\endtrivlist}

\def\basins{\mathcal{B}}
\def\immbasins{\mathcal{B}_{0}}
\def\LL{\mathcal{L}}
\def\comp{\mathrm{Comp}}
\def\N{\mathbb{N}}
\def\I{\boldsymbol{I}}
\def\L{\boldsymbol{L}}
\def\K{\boldsymbol{K}}
\def\J{\boldsymbol{J}}
\def\U{\boldsymbol{U}}
\def\poin{\boldsymbol{Poin}}
\def\Y{\mathcal{Y}}
\def\Comp{\mathrm{Comp}}

\def\cascend{I^{\hat{m}}}
\def\PC{\mathrm{PC}}
\def\crit{\mathrm{Crit}}
\def\R{\mathbb{R}}
\def\C{\mathbb{C}}
\def\dom{\mathrm{Dom}}

\newcommand{\parabolic}{\operatorname{Parabolic}}
\newcommand{\remark}{\noindent\textbf{Remark.}}

\newtheorem{thm}{Theorem}[section]
\newtheorem{cor}[thm]{Corollary}
\newtheorem{lem}[thm]{Lemma}
\newtheorem{prop}[thm]{Proposition}

\numberwithin{equation}{section}
\newcommand{\thmref}[1]{Theorem~\ref{#1}}
\newcommand{\propref}[1]{Proposition~\ref{#1}}
\newcommand{\secref}[1]{\S\ref{#1}}
\newcommand{\lemref}[1]{Lemma~\ref{#1}}
\newcommand{\corref}[1]{Corollary~\ref{#1}} 

\theoremstyle{definition}
\newtheorem{defn}{Definition}[section]
\newtheorem{definition}{Definition}[section]

\newcommand{\diam}{\operatorname{diam}}
\newcommand{\dist}{\operatorname{dist}}
\newcommand{\cl}{\operatorname{cl}}
\newcommand{\inter}{\operatorname{int}}
\renewcommand{\mod}{\operatorname{mod}}

\def\free{free}
\def\Comp{\mathrm{Comp}}
\def\comp{\mathrm{Comp}}
\def\LL{\mathcal{L}}
\def\G{\mathbb G}
\def\R{\mathcal{R}}
\def\E{\mathcal{E}}
\newtheorem*{mainthm}{Main Theorem (Informal Statement)}

\definecolor{patriarch}{rgb}{0.5, 0.0, 0.5}

\def\co{\textcolor{red}}
\def\tco{\textcolor{patriarch}}
\def\sco{\textcolor{blue}}

\newenvironment{proofof}[1]{\medskip
\noindent{\bf Proof of #1.}}{ \hfill\qed\\ }

\thanks{Acknowledgements: This material is based upon work supported by   the ERC AdG grant no 339523 RGDD, 
the National Science Foundation under grant no NSF 1045119 (Trevor Clark), CONACyT grant no 21222 (Sofia Trejo) and FAPESP grant no 2014/09418-0 (Sofia Trejo). Trevor Clark also thanks NSERC for their support.  The authors thank the Mathematics Institute at the University of Warwick and Imperial College London for their hospitality. The authors would also like to thank Oleg Kozlovski and in particular Weixiao Shen for some very useful comments.}

\title{Complex bounds for real maps}

\author{Trevor Clark, Sebastian van Strien and Sofia Trejo}
\address{Imperial College, Imperial College London and Sao Paolo University}
\email{trevorcclark@gmail.com, s.van-strien@imperial.ac.uk, sofia.trejo.a@gmail.com}
\date{December 18, 2016}
\begin{abstract} In this paper we prove complex bounds, also referred
  to as a priori bounds for $C^3,$ and, in particular, for
    analytic maps of the interval. Any $C^3$ mapping of
  the interval has an asymptotically holomorphic extension to a
  neighbourhood of the interval. We associate to such a map, a complex
  box mapping, which provides
a kind of Markov structure for the dynamics. Moreover, we prove universal
geometric bounds on the shape of the domains and on the moduli between
components of the range and domain.
Such bounds show that the first return maps to these domains are well controlled, and consequently 
such bounds form one of the corner stones in many recent results
in one-dimensional dynamics, for example:
renormalization theory, rigidity, density of hyperbolicity, and
local connectivity of Julia sets.
\end{abstract}

\maketitle

\tableofcontents

\section{Introduction and statement of the main results}
The purpose of this paper is to develop a unified technique that
allows one to treat a real analytic or even a
$C^3$ map of the interval as a complex dynamical system, 
where the
domain and range provide a Markov-like structure
for the dynamics.
This problem has a long history that we will discuss later.  
Our results are new for analytic maps, and,
in some cases, even when the mapping is a polynomial.

If the interval mapping is real analytic, the associated complex
mapping is holomorphic, and when the mapping is only smooth, the complex extension is
a quasiregular mapping that is asymptotically holomorphic on its real trace.
Crucially, we will obtain geometric bounds (usually referred to as {\em complex bounds}
or {\em a priori bounds}) for this complex extension. 
In  \hyperref[subsec:applications]{Subsection~\ref{subsec:applications}}
we will explain why complex bounds are useful. 
The following \textit{informal} statement, which we will make more
precise in the next subsection, summarizes our results:

\begin{mainthm}
Let $M\subset \mathbb{R}$ be a compact interval. Assume that $f\colon M\rightarrow M$ 
is a real analytic map (respectively a $C^3$ map with critical points of
integer orders) with no critical points on $\partial M$, and a critical point $c$ such that $f$ is persistently recurrent on $\omega(c)$. 
Then there exist arbitrarily small, combinatorially defined
neighbourhoods $\hat I\subset M$ of $\mathrm{Crit}(f)\cap\omega(c)$
such that we can associate to the (real) first return map to $\hat I$
a complex box mapping (respectively a quasiregular box mapping) with complex bounds.
\end{mainthm}
\begin{center}
\includegraphics[width=11cm]{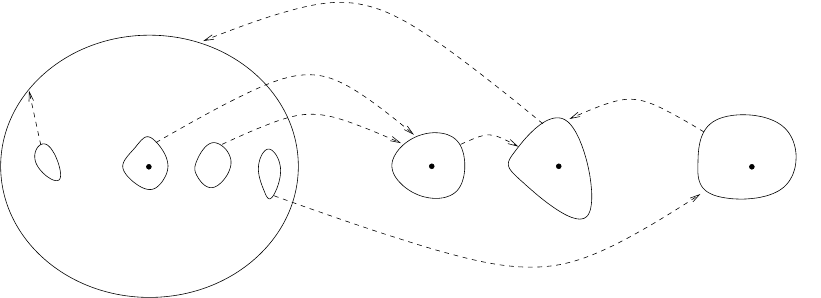}
\end{center}
The notion  \hyperref[persistently recurrent]{\emph{persistently
    recurrent on}} $\omega(c)$ is defined on page
\pageref{persistently recurrent}.

Return maps which are not persistently recurrent on $\omega(c)$ are much easier to work with, 
and complex box mappings will be constructed for such maps in
\cite{CvS}. 

\begin{defn}[Complex Box Mapping]\label{defn:box mapping}
 A mapping $F\colon\mathcal{U}\rightarrow \mathcal{V}$ 
is a complex box mapping\label{def:box mapping} if $F$ is holomorphic and $\mathcal{U}\subset \mathcal{V}$ are open subsets of the complex plane where
\begin{itemize}
\item $\mathcal{V}$ is a union of finitely many pairwise disjoint Jordan disks;
\item every connected component $V$ of $\mathcal{V}$ is either a connected component of $\mathcal{U}$ or $V\cap\mathcal{U}$ is a union of Jordan disks with pairwise disjoint closures that are compactly contained in $V$; 
\item for each component $U$ of $\mathcal{U}$, $F(U)$ is a component
  of $\mathcal{V}$ and $F|U$ is a proper map; 
\end{itemize}
We say that $F\colon\mathcal U\rightarrow\mathcal V$ is a \emph{quasiregular
complex box mapping}, abbreviated \emph{qr box mapping},
if $F$ is a quasiregular mapping that satisfies the remaining
conditions of the definition of a complex box mapping,
see Subsection~\ref{sub-sec:kqcboxS}.
When it will not cause confusion, we will refer to qr box mappings
simply as box mappings.
As usual, a \emph{polynomial-like mapping} is a holomorphic, proper mapping
$F\colon U\rightarrow V$ between topological disks $U\Subset V\subset \C$.
We call a mapping $F\colon U\rightarrow V$,
a \emph{qc polynomial-like map} if $U\Subset V$ are topological disks
in $\C$ and $F\colon U\rightarrow V$ is a proper mapping that can be
expressed as $F=P\circ h$ where $h\colon U\rightarrow U$ is quasiconformal
and $P\colon U\rightarrow V$ is holomorphic.

A complex box mapping (or a qr box mapping) $F\colon \mathcal U\rightarrow\mathcal V$ is called
\emph{real-symmetric} if $\mathcal U$ and $\mathcal V$ are both
real-symmetric, and $\overline{F(\overline z)}=F(z)$. The box mappings
constructed in this paper are real-symmetric, and indeed every construction in this paper
is carried out real-symmetrically.

We say that the map
$F\colon \mathcal{U}\rightarrow \mathcal{V}$  has \emph{complex bounds} (also referred to as \emph{a priori bounds}) if one has estimates on the geometry of $\mathcal{U}$ and $\mathcal V$, see Subsection~\ref{subsec:complex bounds}.
\end{defn}

Throughout this paper, we  only consider first return maps restricted
to the components of their domains that intersect 
$\omega(c)$, for a critical point $c$, and we will implicitly assume this throughout.

Since $f$ is persistently recurrent on $\omega(c)$, 
$\omega(c)$ is a minimal set for $f$, and hence it is compact. Thus,
because each component of $\mathcal{U}$ intersects $\omega(c)$,
$\mathcal{U}$ has at most finitely many components. Complex box
mappings whose domains contain only finitely many components are also
known as \emph{generalized polynomial-like maps}.


\subsection{Some terminology}
Before we can state the main results of this paper, we need some terminology.

\subsubsection{First return and landing maps} Consider a complex box mapping
$F\colon \mathcal{U}\rightarrow\mathcal{V}$.
We will say that $P$ is a \emph{(complex) puzzle piece}
if $P$ is a component of $F^{-n}(V)$, where $V$ is any connected
component of $\mathcal{V}$ and $n\geq 0$. 
To a puzzle piece, we associate three mappings:
the (first) return mapping, the (first) landing map and
the (first) entry map.
Let
$$\mathrm{Dom}(P)=\{z\in\mathcal U:F^k(z)\in P \text{ for some }
k\in\mathbb N\},$$
where we take $\mathbb{N}=\{1,2,3,\dots\}$.
For any $z\in\mathrm{Dom}(P)$, let $k(z)\in\mathbb N$ be minimal so
that $F^{k(z)}(z)\in P$. The \emph{first entry mapping to} $P$ is the
mapping from $\mathrm{Dom}(P)$ to $P$ defined by $z\mapsto
F^{k(z)}(z).$
The \emph{first landing map to} $P$ is defined by $z\mapsto F^{k(z)}(z)$
for $z\in\mathrm{Dom}(P)\setminus P$ and by the identity on $P$.
We define the \emph{first return mapping to} $P$ as the
restriction of the first entry mapping to $P$:
$$R_P\colon \mathrm{Dom}(P)\cap P\rightarrow P,\ \text{where}\ R_P(z)=F^{k(z)}(z).\label{R_I}$$

Suppose that $z\in\mathrm{Dom}(P)$. The connected component of
$\mathrm{Dom}(P)$ that contains $z$ is denoted by $\mathcal{L}_z(P)$ and will be called
a \textit{first entry domain to} $P$ \emph{that contains} $z$. 
We will also call $\mathcal{L}_z(P)$ a \textit{return domain to } $P$, if
$z\in\mathrm{Dom}(P)\cap P.$
\label{entry domain}
We define the \emph{first landing domain to} $P$ containing $z$ by
\begin{equation*}\label{landing domain}
\hat{\mathcal{L}}_{z}(P)= \left\{ 
\begin{array}{rl}
P & \mathrm{if\ } z\in P,\\ 
\mathcal{L}_{z}(P)&\mathrm{if}\  z\notin P.
\end{array}\right.
\end{equation*}
We will also use these definitions in the real case where $P\subset M$
is an interval.


\subsubsection{Complex bounds ($\rho$-nice, $\rho$-free, $\rho$-bounded geometry)}\label{subsec:complex bounds}
Let $\rho>0$.
A puzzle piece $P$ is called $\rho$-\emph{nice}\label{rho-nice domain}
if for any $x\in P\cap\omega(c)$ one has
$\mod(P\setminus\mathcal{L}_x(P))\geq \rho$,
and $\rho$-\emph{free}\label{rho-free domain} if there are puzzle pieces $P^+\supset P\supset P^-$
such that $(P^+\setminus P^-)\cap\omega(c)=\emptyset$,
$\mod(P^+\setminus P)\geq \rho$ and $\mod(P\setminus P^-)\geq\rho$.
We refer to the annulus $P^+\setminus P^-$, which is disjoint from
$\omega(c)$, as \emph{free space}.
We say that a simply connected domain $U$ has 
$\rho$-\emph{bounded geometry with respect to $x\in U$}
if the Euclidian ball \mbox{$B(x,\rho\cdot\mathrm{diam}(U))\subset U$}.
A domain $U$ is said to have 
$\rho$-\emph{bounded geometry}\label{rho-bounded geometry} if there is an $x\in U$ such that $U$ has $\rho$-bounded geometry with respect to 
$x$.

Let $\mathcal{V}$ be the range of a (qc) complex box mapping. We say that
$\mathcal{V}$, respectively, is $\rho$-\emph{nice}, is $\rho$-\emph{free} or has
$\rho$-\emph{bounded geometry} if $\mathcal V$ is the union of
puzzle pieces $V$, such that each $V$, respectively,
is $\rho$-nice, is $\rho$-free or has $\rho$-bounded geometry.

If there exist $\rho>0$, universal, and a neighbourhood
$\mathcal V$ of $\mathrm{Crit}(f)\cap\omega(c)$ such that $R_{\mathcal
  V}\colon \mathrm{Dom}(\mathcal V)\cap\mathcal V\rightarrow\mathcal V$ is a
(qc) complex box mapping
such that one has that $\mathcal{V}$ is $\rho$-nice and has
$\rho$-bounded geometry, then $R_{\mathcal V}$ is said to have \emph{complex bounds.}

\subsubsection{Extendible box mappings}
If $F\colon \cup U_{j}\rightarrow \cup V_{i}$ is a complex box mapping with $b$ critical points, we say that $F$ is {\em $\upsilon$-extendible}\label{defn:extendible} if  there exists $\upsilon>0$ such that for each $i$ there are topological disks $V'_{i}\supset V_{i}$ with $\mod(V'_{i}\setminus\overline{V}_{i})>\upsilon$ such that the following hold:
\begin{enumerate}[(1)]
\item for each $i$, 
if 
$k$ is such that
  $F(V_{i})=V_{k}$, then $F|V_{i}$ extends to a branched covering from
  $V_{i}'$ to $V_{k}'$ and there are no critical points of $F$ in $V'_{i}\setminus V_{i}$;
\item for each component $U_j$ of $\mathcal{U}$, if $k$ is such that $F(U_{j})=V_{k}$, then there exists a topological disk $U_{j}'\supset U_{j}$ so that $F|U_{j}$ extends to a holomorphic map from $U'_{j}$ to $V'_{k}$ with no additional critical points in $U'_{j}\setminus U_{j}$;
\item If $U_{j}\subset V_{k}$, then $U'_{j}\subset V_{k}$.
\end{enumerate}
We will define $\upsilon$-extendible for qr box mappings in the
remarks following the statement of Theorem~\ref{thm:asymptotic box mapping}.


The notions of $\delta$-nice and $\delta$-free have real analogues,
see page~\pageref{interval geometry}.

\medskip


\subsubsection{Remormalizable maps, periodic, central and terminating intervals}
An interval $I$ is a \emph{periodic interval of} $f\colon M\to M$ if there exists $s>1$ such that $I, f(I), \dots, f^{s-1}(I)$ have pairwise 
disjoint interiors and $f^{s}(I)\subset I$ with $f^{s}(\partial I)\subset\partial I$. The integer $s$ is called the \emph{period of} 
$I$. We say that a map $f$ of the interval is \emph{renormalizable at}\label{renormalizable} $x$ if it possesses a  periodic subinterval 
containing $x$, and that it is \emph{infinitely renormalizable at} $x$ if there exist arbitrarily large integers $p>0$ such that $f$ 
has a periodic subinterval containing $x$ with period $p$.

Suppose that $J$ is an interval and $c\in J$ is a recurrent critical
point of $f$. Let $p\geq 1$ be minimal so that $f^{p}(c)\in J.$  We
say that the return to $J$ is \emph{central}\label{central return} if
$f^{p}(c)\in \mathcal{L}_{c}(J).$ Let $\mathcal{L}^{0}_{c}(J)=J$
  and
  $\mathcal{L}^{n}_{c}(J)=\mathcal{L}_{c}(\mathcal{L}^{n-1}_{c}(J))$
  for all $n\geq 1$, then 
we say that $J$ is \emph{terminating}\label{terminating interval} if the returns of $c$ to $\mathcal{L}^{n}_{c}(J)$ are central for all $n$.

If $J$ is terminating, we define $J^{\infty}=\cap_{n\geq
  1}\mathcal{L}^{n}_{c}(J)$; under these circumstances,
$J^{\infty}$ \label{intersection of terminating intervals} is a
periodic interval of period $p$, and $f$ is renormalizable at $c$.

%
%
%
%

Suppose that $\hat I$ and $\hat J$ are unions of intervals.
We say that a mapping $g\colon \hat J\rightarrow\hat I$
\emph{extends to} $G\colon \mathcal U\rightarrow\mathcal V$
if for each connected component $J$ of $\hat J$, there exists
a unique connected component $U$ of $\mathcal U$ such that $J\subset
U$ and $G|J=g|J$. Moreover, we require that each component $U$ of
$\mathcal U$ contains a component $J$ of $\hat J$, that each component $V$ of
$\mathcal V$ contains a unique component $I$ of $\hat I$, and that the maps
 $G$ and $g$ have the same critical points.

\subsection{Complex bounds in the real analytic case} 
We can now state our main theorem for analytic maps 
(see Subsection~\ref{subsec:real puzzle pieces} for the
  definition of a real puzzle piece).

\begin{thm}\label{complex bounds}
Let $M\subset \mathbb{R}$ be a compact interval. Assume that $f\colon M\rightarrow M$ is a real analytic map with a critical point $c$ such that $f$ is persistently recurrent on $\omega(c)$. Then 
there exist $\rho>0$ and combinatorially defined intervals (puzzle pieces) $I\owns c$ of arbitrarily small diameter so that the following holds. Let
$$\hat{I}:=\bigcup_{c'\in\mathrm{Crit}(f)\cap\omega(c)}\hat{\mathcal{L}}_{c'}(I).$$
 \begin{enumerate}[1.]
\item Suppose that $f$ is non-renormalizable. Then the first return map to $\hat{I}$ extends to a complex box mapping 
$$F\colon\mathcal{U}\rightarrow\mathcal{V}\mbox{ so that }\mathcal{V}\cap\mathbb{R}=\hat{I}\mbox{ and }$$
\begin{itemize}
\item for each component $U$ of $\mathcal U$, $F|U$ has at most one
  critical point,
\item each component of $\mathcal{V}$ is $\rho$-nice and $\rho$-free,
\item each component of $\mathcal{V}$ has $\rho$-bounded geometry.
\end{itemize}
\item Suppose that $I$ is a terminating interval for $f$. Then the
  return map to $I^{\infty}$ extends to a
polynomial-like map $F\colon U\rightarrow V$ such that $\mathrm{mod}(V\setminus U)>\rho$.
\end{enumerate} 
\end{thm}

\medskip

\noindent\textbf{Remarks.}
\begin{enumerate}
\item[(a)] If $I$ is sufficiently small, then each component of
  $\hat{\mathcal{L}}_{c'}(I), c'\in\crit(f)\cap\omega(c),$
contains exactly one critical point of $f$, so each
  component $V$ of $\mathcal V$ contains exactly one critical point of
  $F$. 
\item[(b)] The intervals $I$ will be obtained from the {\em
    generalized enhanced nest} defined on page~\pageref{the enhanced nest}.
In the non-renormalizable case this nest coincides with the {\em enhanced nest} from \cite{KSS}.
In Theorems \ref{thm:complex bounds non-ren} and \ref{thm:complex bounds infinitely ren} we restate this theorem making the choice of 
$I$ explicit. 
\item[(c)] When $f$ is infinitely renormalizable, then the assertion holds for any periodic interval $I$ which is 
sufficiently small. If $f$ is at most finitely renormalizable, then as
in the non-renormalizable case, $I$ can be taken to be a sufficiently small 
pullback of some fixed interval $I_{0}$.
\item[(d)] The number $\rho$ is \emph{universal}, or \emph{beau}\label{def:beau}, (a universal bound that
  holds eventually) in the following sense: there exists
$\varepsilon>0,$ which may depend on $f$, such that whenever the
combinatorially defined interval $I$ has $|I|<\varepsilon$,
Theorem~\ref{complex bounds} holds with $\rho>0$ dependant {\em only} on the
number of critical points in $\omega(c)$ and their orders, and not on
$f$.
(It is crucial that we take $I$ to be sufficiently small depending on $f$. 
To get beau bounds, we use the real bounds from Section~\ref{sec:real bounds} and an argument similar to 
 \cite[Theorem IV.B]{dMvS} which states that there exists a beau $\epsilon>0$ so that for 
each $f$ there exists $N$ so that each periodic orbit of period at least $N$ has multiplier  $\ge 1+\epsilon$.) 
\item[(e)] The reason that we do not include $\delta$-free
in the definition of complex bounds is that we require that the free
space be defined by puzzle pieces, and in the infinitely
renormalizable case, we do not prove this. However, in any case,
the mapping $F$ is $\rho$-extendible.
\end{enumerate}

\medskip




The terminology {\lq}beau bounds{\rq} was introduced by Dennis Sullivan.
The underlying concept is
a crucial property for results on renormalization, see Subsection~\ref{subsec:applications}.

\subsection{Complex bounds in the $C^3$ setting}

Let us now explicitly state our main result for $C^3$ interval mappings.
Suppose that $c\in\crit(f)$ is a persistently recurrent critical
point.
Let $c_0\in\omega(c)\cap\crit(f)$ be of even order,
if there is such a critical point in $\omega(c)$, and otherwise choose
$c_0$ arbitrarily. We let $I_0\supset I_1\supset I_2\supset\dots$
denote the \emph{enhanced nest} about $c_0$. See
Section~\ref{sec:enhanced nest} for the definition of this nest.
Whenever $f$ is at most finitely renormalizable at $c_0$, 
we assume that  $I_0$ is contained in the smallest periodic interval
for $f$ containing $c_0$.
If $I\subset\mathbb R$ is an interval 
we let $D_\theta(I)$ the Poincar\'e disk with angle $\theta$
with real trace $I$, see page~\pageref{defn:poincare disk}
for the definition.
 
\begin{thm}[Complex bounds in the $C^3$ case]\label{thm:asymptotic box mapping}
Let $M\subset \mathbb{R}$ be a compact interval. Suppose that $f\colon
M\rightarrow M$ is $C^3$ and its critical points are of integer order
(i.e. $f$ is contained in the class $\mathcal A^3_{\underline b}$ defined in Subsection~\ref{subsubsec:Abk}). 
There exists $\rho>0$ and $C>0$, universal (depending only on the number of 
critical points of $f$ and their order, i.e. on   $\underline b$, but not on $f$), such that for each $n$ sufficiently large the following holds. Suppose that $c_{0}$ is a critical point such that $f$ is persistently recurrent on $\omega(c_{0})$.
\begin{enumerate}
\item Suppose that $f$ is finitely renormalizable. 
 Let $I_{0}\supset I_{1}\supset I_{2}\supset\dots$ be the enhanced
 nest for $f$ about $c_{0}$(constructed so that $I_0$ is contained in
 the smallest periodic interval of $f$ that contains $c_0$).
Let $\hat I_n=\cup_{c\in\omega(c_0)\cap\crit(f)}\hat{\mathcal L}_{c}(I_n)$.
Then for all $n$ sufficiently
  big, the first return map to $\hat I_n$ extends to a $\kappa_n$-qr box mapping
  $F:\mathcal{U}\rightarrow\mathcal{V}$ with $\crit(F)=\crit(f)$ that
  is $\rho$-extendible and such that 
$\mathcal V\cap \mathbb R=\hat I_n$; for each component $U$ of $\mathcal U,$ $F|U$ is
at most unicritical;
the components of $\mathcal{V}$
  are $\rho$-nice, $\rho$-free; the components of $\mathcal{U}$ and
  $\mathcal{V}$ have $\rho$-bounded geometry. Moreover, there exists
  $\theta'\in(0,\pi)$ such that for each component $U$ of  $\mathcal{U}$ or 
  $\mathcal{V}$, there exists an interval $\tilde
  I_U\supset(1+2\rho)(U\cap\mathbb R)$ such that $U$ is contained in
  $D_{\theta'}(\tilde I_U)$. 

\item If $f$ is infinitely renormalizable, then for all $s$ sufficiently big, if $J$ is a periodic interval for $f$ of period $s$, then $f^s\colon J\rightarrow J$ extends to a $\kappa(V)$-qc polynomial-like mapping $F\colon U\rightarrow V$ with $\crit(F)=\crit(f)$, $\mod(V\setminus U)\ge \rho,$ $U$ has $\rho$-bounded geometry and $\mathrm{diam}(V)<C|J|$.
\end{enumerate}
\end{thm}

The number $\kappa_n $ 
depends on $\max_{V}\diam(V),$ where the maximum is taken 
over the component of $\mathcal{V}$,
$\kappa_n,$ and $\kappa(V)$ depends only on
$\diam(V)$ and both $\kappa_n$ and $\kappa(V)$,
tend to 1 as the corresponding diameters tend to $0$.

\bigskip

%
%

\subsection{Complex bounds for induced mappings}
We remark that the proof of our results imply 
the following corollary, which will be useful in applications.

\begin{cor}\label{cor:thm unions}
Let $N$ denote the circle or interval.
Suppose that $f\colon
N\rightarrow N$ is $C^3$ and its critical points are of integer order
(i.e. $f$ is contained in the class $\mathcal A^3_{\underline b}$ defined in Subsection~\ref{subsubsec:Abk}). 
Assume that  $M$ is a union of intervals in $N$
and that $M'$ is a union
of intervals $J'\subset M$ such that for each interval $J',$
there exists $k_{J'}\in\mathbb N$ such that $f^{k_{J'}}(J')\subset M$,
and $f^{k_{J'}}$ does not have a critical point on $\partial J'$.
Define $F\colon M'\rightarrow M$ on each interval $J'\subset M'$ by
$F|J'=f^{k_{J'}}|J'$.
Assume that $c$ is a critical point of $F$ such that $F$ is
persistently recurrent on
$\omega(c)$.
Then the conclusion of 
Theorem~\ref{thm:asymptotic box mapping} holds for F at $c$. 
\end{cor}

For example, 
suppose that $f\colon S^1\rightarrow S^1$ is an analytic mapping of the circle that is
not injective.
Then $f$ has a periodic point $p$, see \cite{CvS}. Let $s$ be the
period of $p$, and let $\mathcal O=\{p,f(p),\dots, f^{s-1}(p)\}$.
Then $\mathcal Y_0=S^1\setminus f^{-1}(\mathcal O)$ is a partition of
$S^1$ by real puzzle pieces, see Subsection~\ref{subsec:real puzzle
  pieces}. Suppose that $f$ has critical point $c_0$ such that $f$ is
persistently recurrent at $c_0$, and
let $Y_0$ be the component of $\mathcal Y_0$ that contains $c_0$. Then
the results of this paper hold for the return mapping to $Y_0$
restricted to the components of the domain that intersect $\omega(c_0)$.

\subsection{Previous results on complex bounds} 
Let us first give some historical background, before discussing in the
next subsection why complex bounds are crucial for results
on renormalization, quasiconformal rigidity and ergodic properties
of one-dimensional dynamical systems.
Complex bounds were first proved by Sullivan for certain infinitely renormalizable unimodal maps \cite{Sullivan}, see also \cite{dMvS}.
They were proved for real unicritical polynomials, in \cite{LS-lc},
\cite{Shen2}, \cite{LY} and \cite{GS}.
Let us now summarize some of the past work for
multimodal analytic maps with all critical points real and of even order:
\begin{itemize}
\item Complex bounds were proved 
for infinitely renormalizable maps with bounded combinatorics in \cite{Sm}.
\item In  \cite{Shen} 
complex bounds were proved for infinitely renormalizable maps. In addition, for at most finitely renormalizable maps,
Shen proves a somewhat weaker version of complex bounds. Namely, the existence of complex box mappings with the property that each domain 
of an iterate of the box mapping is contained in a Poincar\'e disk with real trace of length comparable to the range of the (real) return map, see \cite[Theorem 3']{Shen}. 
\item Complex bounds, analogous to 
those in this paper, 
were proved for at most finitely renormalizable real polynomials with all critical points even and real 
in \cite{KSS}. 
\end{itemize}


Complex bounds for various classes of unicritical analytic maps whose critical point is of odd order, including covering maps of the circle, and certain real polynomial maps and Blaschke products, were obtained in \cite{LS-inflection}. The methods in that paper, and the corresponding paper \cite{Le} in which real bounds are proved, 
do not seem to go through to the case of two or more critical points of odd order. Indeed, also in our proof the presence of odd critical points 
requires us to overcome significant additional difficulties. 

We should note that if $f$ is a non-renormalizable {\em polynomial} (not necessarily real) 
with only hyperbolic periodic orbits in the complex plane,  then the construction of a complex box mapping
follows immediately from the Yoccoz puzzle construction, see \cite{KS}. In fact, if there are neutral periodic points
which are of parabolic type, this construction can be easily made as well: simply consider rays landing on 
repelling periodic points in the boundary of one of the petals of the
periodic point. However, if $f$ has a non-parabolic neutral orbit then in general it may be impossible to find a 
complex box mapping or even a periodic point with two rays landing on it. 
If $f$ is a real polynomial such rays, and therefore a complex box mapping, exists, but nevertheless it seems the only
way to obtain a complex box mapping for which  each puzzle piece contains at most one critical point is
through complex bounds, see  \cite[Section 2.2]{KS}. 

If $f$ is at most finitely renormalizable polynomial with only hyperbolic periodic points, 
then the results in \cite{KS} imply that complex bounds indeed hold for
  $f$. The method there relies
on an important lemma by Kahn-Lyubich, see \cite{KL}, and 
the results in \cite{KS} do not require
the polynomial to be real (and it does not matter whether the 
critical points are of even or odd order).  
Because of this, the proof in \cite{KSS} can be simplified:
one can replace Sections 8 to 11 in \cite{KSS}
by the more general results derived in \cite{KS}.
However the methods in \cite{KS} do not provide complex bounds when $f$ is infinitely renormalizable  or
when  $f$ is $C^3$. 

So combining the puzzle construction and the complex bounds from \cite{KS} shows that if $f$ is a polynomial which  
is at most finitely renormalizable and has only hyperbolic or parabolic periodic  points, then there exists
a complex box mapping so that each puzzle piece contains at most one critical point. 
But if $f$ is a real polynomial which is either infinitely renormalizable (possibly at a  non-real critical point)
or with a non-parabolic periodic point (and a non-real critical point), then \cite{KS} does not provide complex bounds. 

%

Thus, {\em  even for real polynomials our theorem is new:} previous results did not establish {\em beau} complex bounds for (at most) finitely renormalizable real polynomials 
with either non-real critical points or real critical points of odd order, and also not for infinitely renormalizable real polynomials with critical points
of odd order. 

Previous proofs of complex bounds often require dividing the proof into the (essentially) bounded geometry and the big geometry cases. Having big geometry 
simplifies the construction of a complex box mapping at a single level, but when there are no bounds on the scaling factors, it is difficult to transfer estimates to deeper levels, see 
the comment before Proposition~\ref{prop:KSS 11.3}.

The purpose of this paper is to treat all situations in a fully unified manner, dealing with non-renormalizable
and renormalizable maps, with any combinatorics, through essentially the same inductive framework. 
The methods we develop allow us to deal with critical points of any (integer) order. The construction builds on the
one given in Sections 8 to 11 of \cite{KSS}, but encompasses the infinitely renormalizable case,
and overcomes the issues which arise from the presence of odd critical points.
Another important feature of the proof is that it allows one to
associate quasiregular
box mappings to $C^3$ maps, see Theorem \ref{thm:asymptotic box mapping}. This is also an important reason why
we did not aim for short-cuts in the proof in the setting of real analytic maps, 
which do not generalize to the setting of maps with asymptotically holomorphic extensions.
In this paper we will not discuss complex bounds for smooth homeomorphisms of the circle, but merely refer to
\cite{dFdM1,dFdM,Ya}. 

For unicritical, (at most) finitely renormalizable complex polynomials without neutral periodic orbits,
moduli bounds were proved in \cite{KL-lc}.
We should emphasise that complex bounds do not hold in general for
complex maps; there are infinitely renormalizable complex quadratic maps for which 
complex bounds are known to fail, see for example \cite{Le2}. However,  they do hold for
unicritical maps with certain combinatorics,
see \cite{Ka,KL2,KL3}.
For an early work on non-renormalizable mappings with specific combinatorics,
where the shape of puzzle pieces was well
controlled we draw the readers attention to \cite{Sm2}.

\subsection{The usefulness of complex bounds: applications}\label{subsec:applications}
Because of the Koebe Distortion Theorem and the behaviour of the
mapping $z\mapsto z^\ell$, one immediately sees that complex bounds
for $F\colon \mathcal U\rightarrow\mathcal V$
give control on the distortion of diffeomorphic mappings onto
components of $\mathcal V$ and on the shapes of certain puzzle pieces.
However, the implications of complex bounds are much deeper.   
The results in this paper are an important ingredient in results on
the topology  Julia sets, 
renormalization, quasiconformal rigidity and ergodic properties.

\subsubsection{Topological and ergodic properties}
 By methods that are by now standard, see for instance \cite{LS-lf}, \cite{Shen-lf} and \cite{KS}, the results of this paper imply
\begin{thm}
Suppose that $f$ is a real polynomial with real critical points. Then
the Julia set of $f$ is locally connected, and $f$ supports no measurable invariant line field on its Julia set.
\end{thm}
Note that, complex bounds for reluctantly recurrent maps are 
easier to obtain and the proof can be found in
\cite{CvS}.

Before now, in the multicritical case, such a result was only known
when all critical points were of even order. In the unicritical case,
local connectivity was proved in the presence of one even critical
point  in \cite{LS-lc}, \cite{LY}, \cite{GS} and \cite{Hu-Jiang} and
for the Julia sets of certain Blaschke products with a single
  critical point on the unit circle in
\cite{LS-inflection}. For multicritical polynomials with all critical
points real and of even order, local connectivity was proved for
certain infinitely renormalizable maps in \cite{Sm}, in  \cite{Shen}
(in the case of {\lq}bounded geometry{\rq}) and without assumptions on
the geometry in \cite{KSS}.  Absence of invariant line fields was first proved in \cite{McMullen-renorm} and subsequently in \cite{LS-lf},  \cite{LS-inflection}, \cite{Shen-lf} and \cite{KS}.

\subsubsection{Quasisymmetric rigidity}
Complex bounds also play a key role in proving quasisymmetric
rigidity. In particular, complex bounds are a crucial hypothesis in the QC-Criterion
of \cite{KSS}.
The results of this paper immediately extend the results of 
\cite{KSS} to include polynomials with odd critical points:

\begin{thm}
Suppose that $f$ and $\tilde f$ are two real polynomials, with real critical
points. Assume that $f$ and $\tilde f$ are topologically conjugate as
dynamical systems on the real line, that corresponding critical points
for $f$ and $\tilde f$ have the same order and that parabolic points
correspond to parabolic points, then $f$ and $\tilde f$ are
quasiconformally conjugate as dynamical systems on the complex plane.  
\end{thm}

The results of this paper are a vital ingredient in extending this
result for polynomials to all analytic maps and to a broad class of
$C^3$ mappings.

\begin{thm}[Clark - van Strien, \cite{CvS}]\label{thm:main}
Assume that $f,g\colon [0,1]\to \mathbb{R}$ with $f(\{0,1\})\subset
\{0,1\}$  are real analytic and topologically conjugate.
Alternatively, assume that $f,g\colon S^1 \to S^1$ are topologically conjugate
and that $f$ and $g$ each have at least one critical point or at least one periodic point. 
Moreover, assume that the topologically conjugacy is a bijection between
\begin{itemize}
\item the  sets of critical points and the orders of corresponding critical points are the same;
\item   the set of parabolic periodic points.
\end{itemize}
Then the conjugacy between $f$ and $g$ is quasisymmetric.
\end{thm}

For smooth mappings, we have:

\begin{thm}[Clark - van Strien, \cite{CvS}]
Suppose that $f,g\colon [0,1]\rightarrow [0,1]$,
or alternatively that $f,g\colon S^1 \to S^1$ each has least one critical point or at least one periodic point,
are $C^{3}$, each with a finite number of critical points. Suppose that at each $c\in\crit(f),$ one can locally express $f(x)=[\phi(x)]^\ell+f(c)$ where $\phi$ is a $C^3$ diffeomorphism 
with $\phi(c)=0$ and $\ell$ is an integer $\ge 2$, and likewise for $g$. Assume that $f$ and $g$ have only repelling periodic points. Suppose that $f$ and $g$ are topologically conjugate and that the conjugacy is a bijection between $\crit(f)$ and $\crit(g)$ and that the orders of corresponding critical points is the same. Then $f$ and $g$ are quasisymmetrically conjugate. 
\end{thm}
Under some additional smoothness and genericity assumptions, we can remove the condition that all periodic orbits be repelling, as in Theorem \ref{thm:main}.
These theorems extend earlier work for quadratic polynomials
\cite{Lyubich-quadratic-dynamics,GS2,GS3}, and for polynomials without odd critical points \cite{KSS}.
Partial results in this direction are proved in \cite{Shen, LS-inflection}. 

Quasisymmetric rigidity is a crucial ingredient in proving {\em
  density of hyperbolicity}, see
\cite{Lyubich-quadratic-dynamics,GS2,GS3} for quadratic polynomials,
for real polynomials without odd critical points \cite{KSS} and
\cite{AKLS}  for unicritical polynomials. Density of hyperbolicity in
the space of $C^2$ maps of the interval was proved in \cite{Shen}.

Quasisymmetric rigidity  can also be proved for a large class of real transcendental maps, see \cite{RvS1} and \cite{RvS2}.
Another motivation is to extend results about monotonicity of entropy for real polynomials
with only real critical points, see \cite{BvS}, to real polynomials with non-real critical points. This is work 
in progress by the 2nd author joint with Cheraghi. 

\subsubsection{Renormalization results}
In the 1970's, Feigenbaum and Coullet-Tresser \cite{Fe, TC}
observed surprising universal scaling laws in one-dimensional
dynamics.
They noticed that, in the family $$f_{\lambda}\colon x\mapsto x^2+\lambda,$$
the sequence of points, $\lambda_n,$
as $\lambda$ decreases from $1/4$, where the mapping $f_{\lambda_n}$
passes through a period doubling bifurcation has the following
property: the ratio
$$\frac{\lambda_{n-1}-\lambda_{n}}{\lambda_n-\lambda_{n+1}}$$
converges, and even more, the limit is independent of the family
of unimodal maps, provided that family of maps is chosen so that each
map in the family has a non-degenerate critical point.
They also observed similar universality
properties of $\omega(0).$
In the period doubling case, there is a sequence of periodic intervals
$J_n$ with period $2^n$ under $f$,
for which the ratios $|J_n|/|J_{n+1}|$ converge to a fixed
value, which does not depend on the choice of family.
To explain these observations, they introduced the 
period doubling renormalization operator.
They conjectured that this operator has a unique fixed point and that
this fixed point is hyperbolic with a one-dimensional unstable manifold.
Exponential convergence of renormalization has strong implications for the
rigidity of such maps. For example:
\begin{thm}\label{thm:smooth rigidity} Assume that $f,g$ are real analytic infinitely
  renormalizable maps with bounded geometry. Then any topological
  conjugacy between $f$ and $g$ is differentiable at the critical point. 
\end{thm}

The original renormalization conjecture has been extended to cover all unimodal
infinitely renormalizable combinatorial types. Roughly, 
the generalized version of the conjecture states that,
there is a renormalization operator, $\mathcal R$, that acts on an
appropriate space of functions, and has an invariant set, $\mathbb K$, called the full renormalization horseshoe.
Furthermore, at each point $f\in \mathbb K$, $\mathcal R$ has a one-dimensional
unstable manifold, and the stable manifold of $f$ corresponds to its
topological conjugacy class. 
This conjecture was settled for quadratic maps by Lyubich,
\cite{Lyubich-renormalisation,L-ae}.

Partial proofs of such results were first obtained by \cite{CE, La}
using bounds obtained with computer assistance.
Sullivan, \cite{Sullivan}, was the
first to prove convergence of renormalization for real
analytic infinitely 
renormalizable maps of bounded type.
It is precisely for this reason that he derived complex bounds for such maps.
A crucial ingredient in the proof was to have that on a sufficiently small scale, the bounds are independent of the map, i.e. that
these bounds are {\em beau}, see page~\pageref{def:beau}.
This means that there exists a compact class $\mathcal C$ of maps,
so that after renormalizing a map $f$, possibly a large number of
times,
its renormalization is in $\mathcal C$.
Subsequent renormalization results, including {\em exponential} convergence of renormalisation, were obtained by 
\cite{McMullen-renorm, McMullen-3mflds, Lyubich-renormalisation,L-ae}.
The most recent
proof in \cite{AL} of the convergence of renormalization shows that, in fact, the property of complex bounds with beau estimates
is essentially the only ingredient that is required. 
For analytic interval maps with several critical points, results on
renormalization have been proved by Smania, \cite{Smania-phase}, \cite{Smania-shy}.
It is clear from the results just cited that
our results will be a key to extending 
results about the hyperbolicity of 
renormalization to more general settings.

For critical circle maps,
there is a renormalization theory that is closely related to the
theory for unimodal maps, see for example \cite{Ya2,
  KT,dFdM1,dFdM}.

In addition to explaining the universal scaling laws 
in both the phase and parameter spaces
observed by Feigenbaum and Coullet-Tresser,
the hyperbolicity
of renormalization is a vital ingredient in the proof
of the celebrated theorem that in the real family $z\mapsto z^2+c, c\in[-2,1/4],$
almost every map is regular or stochastic, \cite{L-ae}.
The hyperbolicity of renormalization, together with the fact that the leaves of the lamination of the
space of polynomial-like maps by the hybrid classes are analytic
manifolds, 
imply that the set of
infinitely renormalizable maps have measure zero in generic families
of unimodal maps. To complete the proof of the regular or stochastic
theorem a parameter exclusion argument and a geometric
characterization of stochastic mappings are used to show that in the
set of non-regular, non-renormalizable parameters, almost every map is
stochastic, \cite{L-III, MN}. This result has been generalized and
improved: in generic families of analytic unimodal maps almost every map is
regular or Collet-Eckmann, see \cite{AM1, AM2,ALdM, ALS, BSS, C}

\subsection{An outline of the paper and a sketch of the proof}
In this paper, we construct
box mappings with complex bounds
around a critical point $c_0$ with the property that $f$ is
\hyperref[persistently recurrent]{\emph{persistently recurrent}} on
$\omega(c_0)$; recall that either 
$c_0$ is even or all critical points in $\omega(c_0)$ are odd.
When $f$ is not persistently recurrent at $c_0$, one can go with bounded degree from arbitrarily small scales around $c_0$ to a fixed large scale; therefore for such critical points the construction of
complex bounds uses rather different methods which can be found in \cite{KSS2} and \cite{CvS}.
For maps with a persistently recurrent critical point $c_0$, we define a sequence of nested intervals around $c_0$ called the \hyperref[the enhanced nest]{\emph{generalized enhanced nest}}
$$I_{0}\supset I_{1}\supset I_{2}\supset\dots.$$ In the
non-renormalizable setting, this generalized enhanced is identical to the enhanced nest in \cite{KSS};
however, we extend the construction so that it also covers infinitely
renormalizable maps. This nest allows us to 
construct quasi-box mappings, see
Subsection~\ref{subsec:desc quasi-box}, in the infinitely renormalizable and
  non-renormalizable cases simultaneously. This is the crucial step in the
construction of complex box mappings.


The enhanced nest possesses key features that play important roles
throughout the proof. In the non-renormalizable case, each interval
$I_{n+1}$ in the enhanced nest is a pullback of $I_{n}$ with bounded
degree depending only on the vector $\underline b$. 
The combinatorics of the enhanced nest are very well controlled
- if a chain starts in a deep level of the enhanced nest,
and returns to it, we exploit the fact that the chain had to visit
each higher level several times prior to returning. Finally, the enhanced nest provides us with dynamically defined space, free (disjoint) from $\omega(c_{0}),$ even in the (infinitely) renormalizable case. It is worth noticing that the enhanced nest is never a subsequence of the principal nest and that while the combinatorics of the enhanced nest are far more complicated than those of the principal nest, the principal nest does not provide the same geometric control as the enhanced nest.
 
\subsubsection{Sections~\ref{sec:real
    bounds}~-~\ref{sec:pullback along enhanced nest}:
Developing the required tools}
We use the same strategy to prove complex bounds as was used in
  \cite{KSS}, but we extend it to allow for odd critical points
  and for infinitely renormalizable maps. For this reason we follow the notation from
\cite{KSS}, and refer as much as possible to results and proofs from
that paper. We have attempted to indicate to which past results our
statements correspond, even though, the proofs and some of the
statements require modifications.
 It turned out that to overcome the additional difficulties
for dealing with this generality required new ideas. In particular,
many of the results or proofs in Sections~\ref{sec:pulling back disks}
and \ref{sec:pullback along enhanced nest} have no analogues in
\cite{KSS}. 

Our aim is to construct a complex box mapping with complex bounds that
extends $R_{I_n}$ for any $n$ sufficiently large.
Remember that we always restrict return maps to the components that intersect
$\omega(c_0)$.
Our main goal is the construction of a quasi-box
mapping. A quasi-box mapping is similar to a box mapping,
except that the components of its domain need not be compactly contained
in its range and we do not require the components of its domain 
to be pairwise disjoint
(see Definition \ref{defn:quasi-box mapping}).
This construction occupies most of the paper, Sections 3 to 9. 
In the final section, Section 10, we make use of
quasi-box mappings to build box mappings
and we show complex bounds.

To construct quasi-box mappings we use
Poincar\'e disks based on intervals from
the generalizad enhanced nest.
See page~\pageref{defn:poincare disk},
for the definition of the 
Poincar\'e disk $D_{\theta}(I)$, where $I$ is an interval in
$\mathbb R$ and $\theta\in(0,\pi)$.
We let $\mathbb C_I$ denote the
\emph{slit complex plane}
$\mathbb C_I=\mathbb C\setminus(\mathbb R\setminus I)$.
Let us explain the basic idea behind the construction. 
Fix $n$ large and let $\mathcal L_{x}(I_n)$ be the landing domain to
$I_n$ containing  a point $x\in I_n\cap \omega(c_0)$.
Consider the disjoint chain $\{G_j\}^{s}_{j=0}$ given by $G_{s}=I_n$
and $G_0=\mathcal L_{x}(I_n)$ and let
$$U_{s}= D_\theta(I_{n-N})\cap \C_{I_n} \, \,\text{and} \,
\,U_{j}=\comp_{G_j}f^{-1}(U_{j+1}),$$ for $0\leq j <s$ and some fixed
$3<N<n$. Assume there exists a constant $\mu\in(0,1)$ so that 
$$U_{0}\subset D_{\mu\theta}(K)\subset U_{s},$$
where $K$ is an interval well-inside $I_{n-N}$ with $\mathcal
L_{x}(I_n)\subset K$.
Observe that the map $f^{s}: U_{0}\to U_{s}$ is  a (qc) covering map. 
Hence, we can construct a quasi-box mapping $F:\mathcal U\to \mathcal
V$,
by considering $\mathcal U$ equal to the union of all sets $U_{0}$ constructed as above, for points $x\in I_{n}\cap \omega(c_0),$ 
and range $\mathcal V=D_\theta(I_{n-N}).$ The key part of the proof to show the existence of a beau bound on the constant $\mu$, depending only on the number of critical points of $f$ and their orders, and a universal constant $N$ so that the above construction holds for all intervals $I_n$ with $n$ sufficiently large.
To find such constants we first study the geometric properties of the generalized enhanced nest,  and then we study pullbacks of Poincar\'e domains in various circumstances.

\medskip
Let us now survey what is done in each section.

In \hyperref[sec:real bounds]{Section~\ref{sec:real bounds}}, we prove
the necessary ``real bounds.'' 
Before going into further detail, we refer the reader to the beginning of Section~\ref{sec:real bounds} for the definitions of 
$\delta$-nice and $\delta$-free.
It is worth remarking that while some
of the results in this section are very close to those
in \cite{KSS}, the proofs usually differ in significant ways. 
For instance, when all critical points of $f$ are
even as in \cite{KSS}, all intervals $I_n$ are $\delta$-free. However,
in our setting,
this need not be the case.
Very roughly, using the terminology from Subsection~\ref{subsec:complex bounds}, 
 we show the existence of a universal constant $\delta>0$ such that: if $I_{n}$ is a non-terminating level of the generalized
 enhanced nest, then 
\begin{itemize}
\item $I_{n}$ is $\delta$-nice, and 
\item $((1+2\delta)I_{n+1}\setminus
  I_{n+1})\cap\omega(c_{0})=\emptyset.$ 
\end{itemize}

Bounded scaling factors between sufficiently many nearby levels 
in the enhanced nest has important consequences for the
geometry of $\omega(c_0)$.
For example, we show that
if there is a small return domain to $I_{n}$, then $I_{n+1}$ is small
compared to $I_n$ , and
if $I_n$ and $I_{n+1}$ are non-terminating, 
then then $I_{n}$ is $\delta$-free
for some $\delta>0$ depending on $|I_{n}|/|I_{n+1}|$.
Furthermore, we are able to use the generalized enhanced
nest and certain bounded geometry conditions to
control the geometry of the post-critical set for infinitely
renormalizable maps, see Proposition~\ref{delta free}.
Control of the post-critical set is vital throughout this paper
since it makes it possible to control the shape of pullbacks of
Poincar\'e disks in the complex plane, which are the basic pieces we will use to construct quasi-box mappings.

In \hyperref[sec:Poincare disks]{Section~\ref{sec:Poincare disks}} we present some necessary facts about
Poincar\'e disks and show how to pull them back by first return maps. 
When we pull back one of these domains (i.e. we take its preimage) by a unicritical
branched covering, 
knowing that the critical values of the map are not close to the
boundary of the disk, allows us to control the shape of this pullback,
see Lemma~\ref{lem42} and Lemma~\ref{lem:z^l lower bound}.
Of greatest importance to us is that we control the \emph{loss of
  angle}. The pullback of a Poincar\'e disk
  with angle $\theta$ is contained in a Poincar\'e disk of angle
  $\theta'\leq \theta$ and we bound $\theta'$ from below (in terms of $\theta$), 
which gives us some control on the geometry of the pullbacks.

\begin{figure}[htb] \centering \def\svgwidth{170pt}
 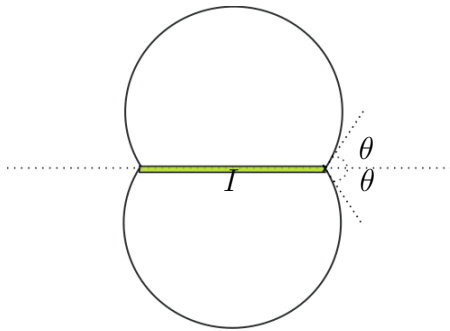
 \caption{The Poincar\'e disk $D_{\theta}(I)$
of angle $\theta$  based on an interval $I$}
 \label{fig:PD}
 \end{figure}
 
 Asymptotically holomorphic extensions were used by Lyubich to prove complex bounds and 
 quasisymmetric rigidity for the quadratic Fibonacci map \cite{Lyubich-Fibonacci2}.
For an application to infinitely renormalizable mappings with bounded combinatorics
see \cite{Sullivan}.
In \cite{GSS2} the theory of asymptotically
holomorphic extensions was developed further, and we use their results in our proof.
In \hyperref[sec: smooth extension]{Section~\ref{sec: smooth
    extension}} 
    we recall these results, and use them to
    develop  tools for dealing with $C^3$ maps
$f\colon M\to M$. 
These maps yield quasiregular extensions of
$C^{3}$ interval maps to a neighbourhood in the 
complex plane, which are asymptotically holomorphic. The key
property of these extensions is that, like analytic maps, they satisfy
an ``Almost Schwarz Inclusion Principle:'' the loss of angle when pulling back a
Poincar\'e disk by a diffeomorphism is small, see Lemma~\ref{GSSProp2}. 
We also show in this section that the results in
Section~\ref{sec:Poincare disks} go through in the smooth setting.
The reader who is only interested in the main theorem in the real
analytic setting can skip this section. 



In \hyperref[sec:pulling back disks]{Section~\ref{sec:pulling back
    disks}} we show how to pull back Poincar\'e disks in a few
different situations: along monotone branches, by maps of bounded
degree, and through long cascades of central returns. All these arguments are
complicated by the fact that our map is not a polynomial; therefore
we will lose angle even when we pull back by a diffeomorphism.

In Propositions~\ref{monotone pullbacks}~and~\ref{good start or
    good deep}, the problem that we face is pulling back a slit
Poincar\'e disk $D_{\theta}(\hat G_s)\cap \mathbb C_{G_s}$ along a chain $\{G_j\}^{s}_{j=0}$ under which the chain $\{\hat G_j\}^{s}_{j=0}$ could have any degree.
To deal with this issue, we make use of 
dynamically defined free space and the existence of fundamental domains.
Of particular importance is Proposition~\ref{monotone pullbacks}
		 (which has no analogue in \cite{KSS}). This proposition allows 
us to pull back along a single monotone branch for as long as we like,
using fundamental domains of definite size to control the loss of
angle.

Proposition~\ref{good start or good deep}
allows us to pull back along long cascades of central returns. It is
similar to \cite[Proposition 11.1]{KSS}; however, it also covers the case of
terminating intervals. Its proof, compared to its analogue in
\cite{KSS},
is complicated by the fact 
each time we pullback along a monotone branch we lose some angle, 
so we have to bound the number of times we switch
between monotone branches.
This argument is new and non-trivial; ideas of 
this proof are also used in Section \ref{sec:combinatorial depth}.

In \hyperref[sec:combinatorial depth]{Section~\ref{sec:combinatorial
    depth}} 
we use the notion of \emph{combinatorial depth} of a chain from
\cite{KSS}. 
For a given critical point the combinatorial depth measures how 
close a chain comes to the critical point $c$ 
in terms of the number central cascades around $c$ that the chain
enters.
Proposition~\ref{combinatorial depth angle control} bounds the loss of angle for the pullback of the Poincar\'e disk along a chain in terms of its combinatorial depth.

In Section 8 we show how to control the loss of
angle as we pull back from one level of the
enhanced nest to the next, under two different circumstances: bounded
scaling factors and big scaling factors. The first case, which
occupies
the majority of the section, is roughly explained as follows.
There exist $s_i<s_{i-1}$, $\lambda\in(0,1)$ and
$\delta>0$, so that
$$\comp_{G_{s_i}}f^{-(s_{i-1}-s_i)}(D_{\theta}(I_{n-i-1})\cap\mathbb C_{G_{s_{i-1}}})\subset
  D_{\lambda\theta}((1+2\delta)^{-1} I_i).$$
This argument and the strategy of the proof are subtle, in particular due the presence
	of terminating intervals. We must control the
        combinatorial depth of (certain segments of) the chain
        $\{G_{j}\}_{j=0}^{s_i},$
which can be guaranteed only when the scaling factor between
$I_{n-i-1}$ and $I_{n-i+1}$ is bounded.
In this part of the proof, the difference between the non-renormalizable and the
renormalizable cases is pronounced. Roughly, in the non-renormalizable
case we are able to decompose the map
$f^{s_{i-1}-s_i}|D_{\theta}(I_{n-i-1})\cap\mathbb C_{G_{s_{i-1}}}$
into a bounded sequence of that maps that we can control,
and use the results of the previous
sections to control the loss of angle at each stage. 
This is impossible to do when there are terminating intervals in the
generalized enhanced nest; if $I_{n-i-1}$ is terminating, $I_{n-i}$ is
not a pullback of $I_{n-i-1}$. In this case, 
we pull back from a terminating interval in the enhanced nest to the
next 
level in two stages: first we pull back from the 
terminating interval to the largest periodic interval contained in it,
and then we pull back from the periodic interval to the next level 
of the generalized enhanced nest.
These strategy is described in Diagrams \ref{diagram1} and \ref{diagram2}.

When there are big scaling factors
between levels of the enhanced nest, 
it is enough to make use of the dynamically defined external free
space to control the loss of angle. 
Furthermore, big scaling factors allow us to construct
of complex box mappings without the aid of quasi-box mappings.

\subsubsection{\hyperref[sec:ext to quasi-box mapping]{Section~\ref{sec:ext to
      quasi-box mapping}: Construction of quasi-box mappings.}
  }\label{subsec:desc quasi-box}
In this section 
we use an inductive argument to construct
a
\hyperref[defn:quasi-box mapping]{\emph{quasi-box mapping}} that extends the
return mapping to any interval $I_n$ of the enhanced nest, provided $n$ is  sufficiently large. The notion of a quasi-box mapping was introduced in
\cite[Remark 12.4]{LS-lc},  used explicitly in \cite{LS-inflection},  
and is defined as follows:

\begin{defn}\label{defn:quasi-box mapping}
Let $U_{i}$ and $V_{j}$ be open Jordan disks in
$\mathbb{C}$, $0\leq j< b$, $i\in \mathcal I$ where $\mathcal I$ is at
most countable. Set $\mathcal{U}=\cup U_i$ and $\mathcal{V}=\cup V_j$.
A mapping $F\colon \mathcal{U}\rightarrow \mathcal{V}$ is a
\emph{holomorphic quasi-box mapping} if for each $i$, there exists $j$ and a
holomorphic branched covering 
$F_i\colon U_{i}\rightarrow V_j$, and the following hold:
\begin{itemize}
\item $V_0,\dots, V_{b-1}$ are pairwise disjoint Jordan disks;
\item every connected component $V_j$ of $\mathcal{V}$ is either a connected component of $\mathcal{U}$ or the intersection of $V_j$ with $\mathcal{U}$ is a union of sets $U_i$ where each of these is  contained in $V_j$, not necessarily compactly;
\item if $U_i\subset V_j$, then $\overline{U}_{i}\setminus V_j\subset \mathbb{R}$.
\end{itemize}
Notice that we do not assume the components of $\mathcal U$ are
disjoint,
so $F\colon \mathcal{U}\rightarrow\mathcal{V}$ may be multi-valued, but
$F|U_i=F_i$ is well-defined as a single valued function on each
$U_i$. If we only require that $F_i$ be quasi-regular on each $U_i$ then
we say that $F$ is a \emph{qc quasi-box mapping}, see
Subsection~\ref{sub-sec:kqcboxS}.
When it will not cause confusion, we will refer to
these mappings as
\emph{quasi-box mappings.}
\end{defn}
\noindent We define the \emph{filled Julia set} of a quasi-box mapping
$F\colon \mathcal U\rightarrow\mathcal V$ as the set 
$$K(F)=\{z\in\mathcal
U:F^k(z)\in\mathcal U\mbox{ for all } k=0,1,2,\dots\}.$$

To construct a qc quasi-box
mapping we only need to control the loss of angle 
when we pull back a Poincar\'e disks under a
certain chain.
The angle control that we have is the same for both the asymptotically
holomorphic extensions we use
as it is for analytic maps, so the proofs
in the smooth and analytic cases are essentially the same.
It is worth remarking that we do not control the dilatation of a qc
quasi-box mapping; 
the presence of
long central cascades seems to make this impossible using our
methods.

 We start the construction of quasi-box mappings with a Poincar\'e disk, based on a
slit domain $$D_{\theta}(I_{n-M})\cap \mathbb{C}_{I_n} \mbox{ where }\mathbb{C}_{I_n}=\mathbb{C}\setminus(\mathbb{R}\setminus I_{n})$$
and pull it back from one level of the enhanced nest to the next using the arguments in
Section~\ref{sec:pullback along enhanced nest}.

If the scaling factor
between the level $I_{n-i}$ and $I_{n-i+1}$ 
is ever big enough, we can easily construct a box mapping for
the return map to a deeper puzzle piece, which gives us a quasi-box
mapping for the return map to deeper levels.
So in what follows we asume the scaling factors
between the level $I_{n-i}$ and $I_{n-i+1}$ are bounded.

The combinatorics of the
enhanced nest make it possible to select times to pull back to in such
way that once we have pulled back to level $I_{n-i}$, the forward orbit
of the chain only visits $I_{n-i}$ a couple of times. 
This bound is used in following way:
if $I_{n-i}$ is comparable to $I_n$, then by the real bounds, the
components of the domain of the return map to $I_{n-i}$ are comparable
to $I_{n-i}$, so the derivative of the return mapping to $I_{n-i}$ is
bounded. Since the number of returns of the chain $G_{s_{n-i}},\dots,
G_s$ to $I_{n-i}$ is bounded, we have that $G_{s_{n-i}}$ is comparable
to $I_n$, and hence to $I_{n-i}$ too.
In the presence
of  terminating intervals, the proof of the existence of suitable times
to pull back to, at each level, involves a combinatorial analysis of the
return maps. Let us explain the pull back argument in more detail.
Let $x\in\omega(c)\cap I_n$ and consider the chain $\{G_j\}_{j=0}^s$ where $G_s=I_n$ and $G_0=\mathcal{L}_x(I_n)$.
From Section~\ref{sec:pullback along enhanced nest}, 
we know there exist $s_{n-M+1}<s,$ $\lambda\in(0,1)$ and
$\delta>0$
such that
$$\comp_{G_{s_{n-M+1}}}f^{-(s-s_{m-M+1})}(D_{\pi/2}(I_{n-M})\cap\mathbb
C_{I_n})\subset D_{\lambda\pi/2}((1+2\delta)^{-1}I_{n-M+1}).$$
In other words, as we pull back a point from one level of the enhanced
nest to the next, we know how much angle we lose. 
If we do not lose any angle; that is, for the set of points in
$$\comp_{G_{s_{n-M+1}}}f^{-(s-s_{m-M+1})}(D_{\pi/2}(I_{n-M})\cap\mathbb
C_{I_n})\cap D_{\pi/2}(I_{n-M+1}),$$
we repeat the argument and control the loss of angle up to a time
  $s_{n-M+2}<s_{n-M+1}$ with a Poincar\'e disk based on $I_{n-M+2}$. 
This argument can be carried on inductively, because of the
combinatorics of the enhanced nest.
If we lose angle when we pull back; that is, for the set of points in
$$\comp_{G_{s_{n-M+1}}}f^{-(s-s_{m-M+1})}(D_{\pi/2}(I_{n-M})\setminus D_{\pi/2}(I_{n-M+1})),$$
then
this is a ``jumping time'' associated to this point (see \cite{LY}),
and also Figure \ref{fig:DI} on page \pageref{fig:DI}. 
In this case,
since $G_{s_{n-M+1}}$ is comparable to $I_{n-M+1}$, as assumed above,
there exists $\lambda'\in(0,1)$ so that
$$\comp_{G_{s_{n-M+1}}}f^{-(s-s_{m-M+1})}(D_{\pi/2}(I_{n-M})\cap\mathbb
C_{I_n})\setminus D_{\pi/2}(I_{n-M+1})\subset D_{\lambda'\pi/2}(G_{s_{n-M+1}}),$$
and we are are able to capture
this set by a Poincar\'e disk based on $G_{s_{n-M+1}}$ without losing
too much angle. 
Since the chain $\{G_i\}^{s-1}_{i=0}$ is disjoint, it is not hard to
control the loss of angle as
we pull back all the way to the start. 
So there are two cases, given  a point
 $z\in D_{\pi/2}(I_{n-M})\cap\mathbb
C_{I_n}$ either the point $f^{-(s-s_{n-i})}(z)\subset
D_{\pi/2}(I_{n-i})$ for all $i$, or there exists some $i$, for which
$s_{n-i}$ gives a jumping time. In either case, choosing $M$ properly,
we show $f^{-s}(z)\subset D_{\pi/2}(I_{n-1}).$ Hence, the return map
to $I_n$ extends to a  quasi-box mapping. 
This is 
Theorem~\ref{quasi-box mapping}. 
It is important to remark that, throughout this section,
 the presence of periodic intervals
significantly complicates the arguments.

\subsubsection{\hyperref[sec:box mapping]{Section 10: }Construction of complex box mappings with complex bounds.}
In this final section, we use of quasi-box mappings 
to construct box mappings. In the non-renormalizable case,
we employ a method that was
first used in \cite{LS-inflection}: we first construct {\lq}by hand{\rq}
a \emph{smooth box mapping} with the desired topological properties:
the domain being compactly contained in the range and the components
of the domain being disjoint, and then intersect it with the quasi-box
mapping to obtain a complex box mapping. In the renormalizable case we
make use of a result of \cite{LY}.
We then prove that complex
bounds hold for these complex box mappings. 
The work in this section is done primarily to deal with the smooth, as
opposed to the analytic, case. 
The complex bounds for at most finitely renormalizable maps
follow immediately from the
``Upper and Lower Bounds,''\label{page:upperlower} see \cite{KSS}:

\begin{prop}[Upper Bounds]\label{prop:upper bounds2}
There exists a constant $\eta>0$ such that for all $n$ sufficiently
large the following hold.
\begin{itemize}
\item $\mathrm{diam}(\I_{n})\leq \eta|I_{n}|$, and;
\item there exists a topological disk $\Omega\supset \I_{n}$ such that
$(\Omega\setminus\I_{n})\cap\omega(c_0)=\emptyset$ and
$$\mod(\Omega\setminus\I_{n})>1/\eta.$$
\end{itemize}
\end{prop}

\begin{prop}[Lower Bounds]\label{prop:lower bounds}
There exist beau constants $\eta>0$ and $\varepsilon>0$ such that for all $n>0$ with $|I_n|<\varepsilon,$
$$B(c_0,\eta|I_n|)\subset\I_n.$$
\end{prop}

For infinitely renormalizable maps, the key estimate is to show that
if $F':U'\rightarrow V'$ is the quasi-box mapping map constructed in
Section~\ref{sec:ext to quasi-box mapping} that extends the
return map to a periodic interval $J$ of sufficiently high period, then
$\mod(V'\setminus K(F'))$ is bounded away from zero. From this we
obtain a polynomial-like map that extends the return map to $J$ with
complex bounds from Lemma~\ref{lem:LY}.


\subsection{Relation to complex bounds for non-renormalizable
  polynomials: \cite{KS}} Up to
Section~\ref{sub-sec:upperbounds}, this paper is concerned with the
construction of a complex box mapping associated to a real map of the
interval. Before this is accomplished, the methods of \cite{KS} do not
seem to apply. 
However, once we have constructed a complex box mapping, 
in the analytic case, it is possible to replace the arguments of this
paper with the arguments of \cite{KS} to establish complex bounds.

In fact, it is possible to generalize \cite[Lemma 9.1 (Small
Distortion of Thin Annuli)]{KS} to the case when
$F\colon \mathcal U\rightarrow\mathcal V$ is quasiregular, and obtain a similar statement.
This gives an alternative approach to establishing complex bounds when
$f\in\mathcal A^3_{\underline b}$ is non-renormalizable. However, with
the preparation that we have already done in this paper, the approach
of \cite{KSS} is more natural, and that is the route we take to prove
complex bounds.

\subsection{The complications of having smooth maps, with critical points of odd order
which are  possibly  infinitely renormalizable.}
Let us highlight some of the main complications in our proof, compared to the proof of complex bounds
in \cite{KSS} in the non-renormalizable case and to \cite{Shen} in the infinitely
renormalizable 
one.

\subsubsection{Diffeomorphic pullbacks} 
In both the analytic and the $C^3$ settings,
we have a loss of angle along diffeomorphic pullbacks of
Poincar\'e disks (see Lemmas~\ref{GSSProp2} and \ref{Almost Schwarz
  Inclusion}).
To bound this loss, we need to ensure that certain 
chains of intervals are disjoint 
or have that the sum of the lengths of their elements are bounded.
This is not a problem in the polynomial case where there is no
loss of angle under diffeomorphic pullbacks.

 \subsubsection{Odd critical points}
If $\omega(c_0)$ contains only critical points of odd order, then the
real bounds require a more elaborate treatment than if at least one
of the critical points is even. This is due to the lack of symmetry 
around critical points of odd order, i.e. an odd critical point
need not be close to the centre of a first return domain.    
As a result, the statement and the proof
of the real bounds become more subtle, 
Indeed, Theorem~\ref{real geometry} 
covers more cases than
the corresponding statement
in \cite{KSS}, namely Proposition~8.1.

When there are any odd critical points, we need
to deal with saddle-cascades occurring in monotone branches.
We do this in Proposition~\ref{monotone pullbacks}.
This is used in Proposition~\ref{good start or good deep},
 which is a key step in the proof, and it is
a generalization of
\cite[Lemma 11.1]{KSS}.
To prove this proposition, 
we need to control the
number of times we switch from one monotone branch to another
since each
time we switch we lose angle.
To make the argument work requires a subtle change to the proof of
Lemma 11.1 in \cite{KSS}.

\subsubsection{Infinitely renormalizable maps}
 If the map is infinitely renormalizable, then the argument used
   to prove complex bounds in \cite{KSS} no longer applies.
For this reason we define
the generalized enhanced nest, see Subsection~\ref{subsec:gen}.
One of the advantages of working with this nest
is that it has better geometric properties than
the principle nest, which was used in \cite{Shen}.
Indeed, we use the control on the geometry of the
  post-critical set, which we obtain under some straightforward
  bounded geometry conditions,
to pull back from one level of the generalized enhanced nest 
to the next. 
This is done in Section~\ref{subsec:renorms}, see also diagrams
in Figures~\ref{diagram1} and \ref{diagram2}, see page~\pageref{diagram1}.


\subsection{Notation and terminology}\label{subsec:notation}
Unless otherwise stated, we adopt the convention that $C>0$ is a large constant, and $\varepsilon>0$ will be a small constant.

\subsubsection{Scaled sets, Poincar\'e domains and components} 
We let $\mathbb{R}$ be the real line. We will always use $I$ to denote an interval in $\mathbb{R}$. 
If $I$ is a bounded interval, then write $I=(a+x,a-x)$ and for $\gamma>0$ define $\gamma I := (a+\gamma x,a-\gamma x)$. \label{delta nhd}
We let $\mathbb{C}$ denote the complex plane. We let $\mathbb{C}_{I}=\mathbb{C}\setminus(\mathbb{R}\setminus I).$
We let $\mathbb{H}$ be the upper-half plane and $\mathbb{H}^{-}$ be the lower-half plane.
If $\theta\in(0,\pi)$, we denote by $D^{+}_{\theta}(I)$ (respectively
$D^{-}_{\theta}(I)$) the region in $\mathbb{H}$ (respectively
$\mathbb{H}^{-}$) bounded by $I$ together with the circle arc
subtending $I$ that meets the real axis with external angle $\theta$
at each boundary point of $I$. We let $D_{\theta}(I)=
D^{+}_{\theta}(I)\cup D^{-}_{\theta}(I)\cup I$\label{defn:poincare disk}.
This set corresponds to the set of points with a fixed distance to $I$ in the Poincar\'e metric
in $\C_I$.
Given a set $K\subset \C$ we let $\comp_{x}(K)$ denote the connected component of $K$ containing $x$. \label{comp}

\subsubsection{Definition of the class of functions}\label{subsubsec:Abk}
We let $\mathcal{A}_{\underline{b}}^{k}$  be the set of $C^{k}$ maps
of the interval $[0,1]$, 
(i.e. are $C^k$ on some small neighbourhood of 
$[0,1]$),
with critical points $(c_{1},\dots, c_{b})$
with integer orders $\underline{b}=(\ell_{1},\dots,\ell_{b}),$ where
$b=|\underline{b}|$ \label{b} is the number of critical points, such
that at each point $c_{i}\in [0,1]$ we can express $f$ locally as
$f(x)=[\phi_i(x-c_{i})]^{\ell_i}+f(c_{i})$ where $\phi_i$ is a local
$C^{k}$ diffeomorphism, $\phi_i(0)=0$ and $\ell_i>0$.  If $\ell_i$ is 
even (odd) we say that the corresponding critical point $c_i$ has even
(odd) order. We will call a critical point of even order a
\emph{turning point}. We let $\mathcal{A}_{\underline{b}}$ denote the
set of such maps that are analytic. 
As usual, we will say that a map is {\em univalent} if it is holomorphic and injective.

\subsubsection{First return maps, pullbacks, periodic intervals and nice intervals}

As usual, 
we let $\omega(x)$ denote the \emph{omega-limit set} of a point $x$.  

Given a set $K\subset \mathbb{C}$ we will only consider branches of  $R_K$ (the first return map to $K$, see  p.\pageref{R_I}) that intersect $\omega(c_{0})$ for a fixed critical point $c_{0}$. Where $c_0$ will be chosen later.

Let $f\colon M\to M$ be a continuous interval map. An interval $J$ is called a \emph{pullback} of an interval $I$ if it is a component of $f^{-s}(I)$ for some $s\in\mathbb{N}$. If $J$ is a pullback of $I$ by $f^{-s}$, we associate to this pullback a \emph{chain} of intervals $\{J_{i}\}_{i=0}^{s}$ with $J_{s}=I$ and $J_{i}=\mathrm{Comp}_{f^{i}(J)}(f^{-1}(J_{i+1}))$ for $i$ satisfying $0\leq i<s$. We say that the order of the chain $\{J_{i}\}_{i=0}^{s}$ is $N$ if precisely $N$ of the intervals $J_{0},\dots, J_{s-1}$ contain a critical point. 

We say that an open interval $J$ is \emph{nice} if $f^{n}(x)\notin\mathrm{int}(J)$ for each $x\in \partial J$ and $n>0$.
This implies that for any $k>n\ge 0$, if a  component $J_{k}$ of $f^{-k}(J)$ intersects a component $J_{n}$ of $f^{-n}(J)$ 
then $J_{k}$ is contained in $J_{n}$; so any two pullbacks of a nice interval are either nested or disjoint. In particular, if the first return time of $x$ to $J$ is equal to $s$ and we consider the chain $\{J_{i}\}_{i=0}^{s}$ with $J_{s}=I$ and $J_{i}=\mathrm{Comp}_{f^{i}(x)}(f^{-1}(J_{i+1}))$ we get the following:  $J_0=\mathcal L_x(J)$,  the intervals $J_0,\ldots J_{s-1}$ are pairwise disjoint and the order of the chain $\{J_{i}\}_{i=0}^{s}$ is bounded by the number of critical points of $f$.

We will say that two intervals $I$ and $J$ have \emph{nested or
  disjoint pullbacks} if for any $m,n\in\mathbb N$, any component $I_1$
of $f^{-m}(I)$ and any component $J_1$ of $f^{-n}(J)$ either $J_1$ and
$I_1$ are nested or $J_1$ and $I_1$ are disjoint.

Give  two nice intervals  $J\subset I$, they are called a {\em nice pair} if 
all iterates of $\partial J$ remain outside the interior of $I$. Under these circumstances, if $J_1,J_2$ are pullbacks of $J$ and $I_1,I_2$ are pullbacks of $I$ with $J_i\subset I_i$ for $i=1,2$ and $I_2\subset I_1$, then either
$$J_2\subset I_2\subset J_1\subset I_1\,\, \mbox{ or }\,\, J_2\subset I_2\subset I_1\setminus J_1.$$ 
\label{nice pair}



\subsubsection{Notation table}
Here is table of some of the notation and terminology used on the paper. 

\bigskip

\noindent
\begin{tabular}{cc}
\begin{tabular}{|l|l|}
\hline
$R_I$ & page \pageref{R_I} \\
\hline
$\mathcal{L}_{x}(I)$ & page \pageref{entry domain} \\
\hline
$\hat{\mathcal{L}}_{x}(V)$ & page \pageref{landing domain}\\
\hline
renormalizable &page \pageref{renormalizable}\\
\hline
central return &page \pageref{central return}\\
\hline
terminating interval &page \pageref{terminating interval}\\
\hline
$J^\infty$ & page \pageref{intersection of terminating intervals}\\
\hline
$\rho$-nice domain & page \pageref{rho-nice domain}\\
\hline
$\rho$-free domain &page \pageref{rho-free domain}\\
\hline
$\rho$-bounded geometry & page \pageref{rho-bounded geometry}\\
\hline
$\gamma \cdot I$ &page \pageref{delta nhd}\\
\hline
$\comp_x(K)$ & page \pageref{comp} \\
\hline
$b, \underline b$ & page \pageref{b} \\
\hline
child & page \pageref{child} \\
\hline
persistently recurrent & page \pageref{persistently recurrent} \\
\hline
successor & page \pageref{successor} \\
\hline
$\Gamma(I)$ & page \pageref{last successor} \\
\hline
central return domain &page \pageref{central} \\
\hline
non-central return domain &page \pageref{non-central} \\
\hline
principal nest &page \pageref{the principal nest} \\
\hline
\end{tabular} &
\begin{tabular}{|l|l|}
\hline
$\hat{m}$ & page \pageref{m-hat} \\
\hline
$\tau$ & page \pageref{tau} \\
\hline
$\mathcal{R}(I)$  & page \pageref{renormalization symbol} \\
\hline
$Y_i$, $\tilde{\mathcal{Y}}_{\gamma}$ & page \pageref{notation: top level}\\
\hline
$\mathcal{A}(I)$, $\mathcal{B}(I)$& page \pageref{A, B in enhanced nest} \\
\hline
$\mathcal{E}(I)$& page \pageref{next level of enhanced nest} \\
\hline
(generalized) enhanced nest & page \pageref{the enhanced nest} \\
\hline
$p_n$ & page \pageref{iterate for pullback of enhanced nest} \\
\hline
$r(I), \hat{r}(I)$ & page \pageref{return times to enhanced nest} \\
\hline
$\rho$-nice interval & page \pageref{interval geometry} \\
\hline
$\rho$-externally free interval & page \pageref{interval geometry} \\
\hline
$\rho$-internally free interval & page \pageref{interval geometry} \\
\hline
$\rho$-free interval & page \pageref{interval geometry} \\
\hline
$\rho$-strongly nice interval & page \pageref{rho-strongly nice} \\
\hline
well-inside, deep-inside & page \pageref{well-inside} \\
\hline
$\mathcal{C}_{c}(I)$ &page \pageref{first non-central return} \\
\hline
$m$&page \pageref{m time} \\
\hline
$\crit(I;J), \crit(I;\{G_j\}_{j=0}^s)$ & page \pageref{crit symbol}\\
\hline
$k(I,\{G_j\}_{j=0}^s), k(I,J), \hat{k}(I,J)$ & page \pageref{comb depth symbols}\\
\hline
$\mathcal T_{\xi}$ & page \pageref{Txi}\\
\hline
\end{tabular}
\end{tabular}


\section{The generalized enhanced nest}\label{sec:enhanced nest}

\subsection{Real puzzle pieces}\label{subsec:real puzzle pieces}

Let us assume that $f\colon M\rightarrow M$ is a $C^{3}$ map of the interval
with $b<\infty$ critical points and let 
$c$ be a recurrent, non-periodic critical point. We say that a set $Z$ is \emph{admissible} if it is a finite forward invariant set, disjoint from the postcritical set of $f$ such that every point of $Z$ is a preimage of a repelling periodic point under an iterate of $f$. 

Given a $Z$ admissible set, we will say that $I$ is a \emph{(real) puzzle piece} of {\em depth n}  (with respect to $Z$)
if it is a component of $f^{-n}(Y)$, where $Y$ is a component 
of $M\setminus f^{-1}(Z)$. We observe that puzzle pieces are
nice intervals, so any two puzzle pieces are either nested or
disjoint.

\subsection{Combinatorics of puzzle pieces}
While the objects in this subsection are defined for real puzzle
pieces, the definitions in this subsection hold whether a puzzle piece
is real or complex.

We say that a puzzle piece is $\omega(c)$-\emph{critical} if it contains a critical point in $\omega(c)$.
Let $P$ be an $\omega(c)$-critical puzzle piece containing the critical point $c'\in\omega(c)$. An $\omega(c)$-critical puzzle piece $Q$ is a called a
\emph{child}\label{child} of $P$ if it is a unicritical pullback of
$P$; that is, there exists
a positive integer $n$ such that $Q$ is a component of $f^{-n}(P)$
containing a critical point in $\omega(c)$, and there exists a puzzle piece
$Q'\supset f(Q)$ such that the map $f^{n-1}\colon Q'\rightarrow P$ is a
diffeomorphism.  A \emph{successor}\label{successor} of $P$ is a puzzle piece of the form $\hat{\mathcal{L}}_{c'}(Q)$, where $Q$ is a child of $\hat{\mathcal{L}}_{c''}(P)$ for some critical point $c''\in\omega(c)$.  By construction,  a successor of $P$ is a pullback of $P$ of order bounded by $2b-1$.

Let $I$ be a puzzle piece containing a recurrent point $x$. We define the {\em principal nest}\label{the principal nest} around $x$ as follows. We set $I^{0}=I$ and inductively define 
 $I^{n+1}=\mathcal{L}_{x}(I^{n}).$ Given a puzzle piece $I^{n}$ the puzzle
 piece $I^{n+1}$ will be called \emph{central},  while\label{central}
 any other return domain to $I^n$ will be called
 \emph{non-central}. \label{non-central}.

Let $c$ be a recurrent, non-periodic, critical point and consider the principal nest
$I^0\supset I^1\supset I^2\supset\dots$ about $c$. We define
$\hat{m}\in\mathbb N\cup\{\infty\}$ to be the smallest number, if it exists, such that a there is a critical point $c'$ of $R_{I^0}|I^{1}$ with $(R_{I^0}|I^{1})(c')\notin I^{\hat{m}}$.\label{def: m hat} If no such integer exists, we set $\hat{m}=\infty$\label{m-hat}. 
 If $\hat{m}<\infty$ we say that $I$ is \emph{non-terminating} and
 otherwise we say that $I$ is \emph{terminating}. In  the terminating case, we let $I^{\infty}=\cap_{i\geq 0} I^{i}$.

\subsection{Persistent recurrence}\label{subsec:persistentrec}
A map $f$ is called \emph{persistently recurrent on $\omega(c)$}\label{persistently recurrent}
if $c$ is recurrent, non-periodic and each $\omega(c)$-critical puzzle
piece has only finitely many children. Under these circumstances, we
will also say that $f$ is \emph{persistently recurrent at} $c$ and that $c$
is a \emph{persistently recurrent critical point.} If $f$ is
persistently recurrent on $\omega(c)$, then $\omega(c)$ is minimal,
but the converse of that statement is false. 
We observe the following, if $f$ is persistently recurrent on $\omega(c)$ each $\omega(c)$-critical puzzle piece $P$ has a smallest successor,
which we denote by $\Gamma(P)$\label{last successor}; and if $Q$ is an entry domain to $P$ intersecting $\omega(c)$, then $\hat{\mathcal{L}}_{c}(Q)$ is a successor of $P$, and thus $P\supset \hat{\mathcal{L}}_{c}(Q)\supset\Gamma(P)$.

From  \cite{KSS}, p.771-772, we know the following: 
\begin{lem}\label{puzzle}
If $f$ is persistently recurrent on $\omega(c)$, then for any
$\varepsilon>0$ there exists an admissible set $Z_\varepsilon$ such that any
real puzzle piece of depth zero that intersects $\omega(c)$ has length less than $\varepsilon$.
\end{lem}

As a direct consequence of the No Wandering Intervals Theorem, see
\cite[page 751]{vSV}, we get that for any $\delta>0$ there exists
$\varepsilon(\delta,f)>0$ such that the length of any pullback of a
critical real puzzle piece with size less than $\varepsilon$ has size
less than $\delta.$ From this fact and  Lemma~\ref{puzzle}, we can
assume that any real puzzle piece that intersects $\omega(c)$ 
has length less than $\varepsilon_0$,
where we can choose $\varepsilon_0>0$ as small as we like.
We will use this observation without further comment.

\medskip
\begin{remark} Suppose that $c\in\crit(f)$ and that $f$ is
  persistently recurrent on $\omega(c).$ 
Then for any sufficiently small nice interval $I\owns c$ the following holds.
\begin{itemize} 
\item The interval $I$ contains no other point from $\crit(f)$. 
\item If $c'\in\crit(f)$ is any critical point such that
there exists a chain $\{G_i\}_{i=0}^{s}$, with $G_s=I$ and
$G_0=\mathcal L_x(I)$ for $x\in \omega(c)\cap I,$
and $c'\in G_j$ for some $0\leq j<s$, then $c'\in \omega(c).$
\end{itemize}
\end{remark}\label{remark SA}

\medskip

Assume $f$ is persistently recurrent on $\omega(c)$.
By Lemma~\ref{puzzle} and the remark above we know that (for $I$ is sufficiently small) all critical values of $f^r|_{I^{1}}$ are contained in $\omega(c)$. 
Since $f$ is persistently recurrent,
$\omega(c)$ is not a periodic orbit. These two facts, along with
\cite[Theorem III.4.1]{dMvS} imply $f^r|_{I^{1}}$ does not have a periodic attractor (otherwise $\omega(c)$ would be a periodic orbit). 
Moreover, since $\omega(c)$ is minimal, and 
the period of all attracting or parabolic cycles is bounded
\cite[Chapter IV, Theorem B]{dMvS}, 
there exists a neighbourhood of $\omega(c)$
which does not intersect any immediate basin of a periodic attractor
or any parabolic cycle.
Thus if $I$ is sufficiently small, for any $x\in\omega(c)$,
the chain $\{G_j\}_{j=0}^s$ with $G_s=I$ and $G_0=\mathcal{L}_x(I)$ avoids a neighbourhood of any immediate basin of a periodic attractor
or any parabolic point.

It will be useful for us to select a critical point $c_0\in\crit(f)$
about which we will focus our construction.
If $\omega(c)$ contains a turning point, we take $c_0$ to be a
turning point, otherwise choose $c_0\in\omega(c)\cap\crit(f)$ arbitrarily.  
Observe that $f$ is persistently recurrent on $\omega(c_0)$, so $c_0$ is recurrent, non-periodic, $\omega(c_0)$ is minimal and $\omega(c)=\omega(c_0)$. 



\subsection{Terminating intervals}
Suppose that $I\owns c_0$ is a terminating interval. Then, since $c$ is recurrent and non-periodic, $R_{I}$ has a critical point of even order,
and hence $c_0$ is of even order too.
Since $c_0$ is a turning point, there exists a neighbourhood $J\supset I$ 
of $c_0$ and an involution $\tau\colon J\rightarrow J$\label{tau} so that $f=f\circ\tau$ on $J$. 
Let $r$ be the return time of $c_{0}$ to $I$. 
Since $I$ is terminating, it follows that 
$ I^{\infty}$ is a periodic interval. More precisely,  
$f^{r}(I^\infty)\subset I^\infty$, and $f^{r}(\partial I^\infty)\subset\partial I^\infty$ and  all of the critical points of $f^{r}|I^{1}$ are contained in $I^\infty$ along with their orbits under $f^r$. We let $\beta$ denote the fixed point of $f^r$ on the boundary of $I^\infty$

Since $c_0\in I^\infty$, a periodic interval, and $c_0$ is recurrent, $f^r|I^\infty$ has at least one repelling orientation reversing fixed point.
Let $\alpha$ be the orientation reversing fixed of $f^r|I^\infty$ closest to $c_0;$  so $(\alpha,\tau(\alpha))\ni c_0$ is the smallest $\tau$-symmetric interval with one repelling fixed point on its boundary. 

Given a terminating interval $I$ we define  $\mathcal{R}(I)$\label{renormalization symbol} as
$$\mathcal{R}(I) := (\alpha,\tau(\alpha)).$$
See Figure \ref{fig:firstreturn} for some examples. If
$\mathcal{R}(I)$ is a periodic interval, then it has period two under $R_I$ and we say $\mathcal{R}(I)$ is \emph{Feigenbaum}.

\begin{figure}[htb]
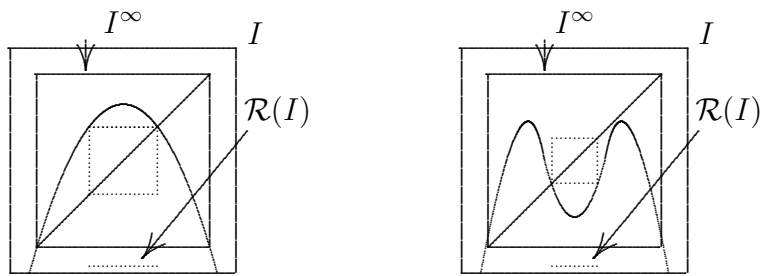

\beginpicture
\dimen0=0.05cm 
\setcoordinatesystem units <\dimen0,\dimen0> point at -120 0 
\setplotarea x from -40 to 40, y from -30 to 30
\setlinear
\put {$I$} at 35 34 
\setsolid
\plot -30 -30 30 -30 30 30 -30 30 -30 -30   /
\plot -23 -23 23 -23 23 23 -23 23 -23 -23   /
\plot -23 -23 23 23 /
\setquadratic 
\plot -25 -30   -15 8   -8 0    0 -15    8 0    15 8  25 -30 /
\setlinear
\setdots <0.6mm>
\plot -6 -6  -6 6  6 6   6 -6 -6 -6 / 
\plot -6 -28 6 -28 /
\setsolid 
\arrow <10pt> [.2,.67] from -8 32 to -8  24
\put {$I^{\infty}$} at 0 37 
\arrow <10pt> [.2,.67] from 33 8 to 5  -26
\put {${\mathcal{R}}(I)$} at 41 13 
\setcoordinatesystem units <\dimen0,\dimen0> point at 0 0 
\setplotarea x from -80 to 40, y from -30 to 30
\setlinear
\put {$I$} at 35 34 
\setsolid
\plot -30 -30 30 -30 30 30 -30 30 -30 -30   /
\plot -23 -23 23 -23 23 23 -23 23 -23 -23   /
\plot -23 -23 23 23 /
\setquadratic 
\plot -25 -30 0 15 25 -30 /
\setlinear
\arrow <10pt> [.2,.67] from -10 32 to -10  24
\put {$I^\infty$} at 0 37 
\setdots <0.6mm>
\plot -9 9 9 9 / 
\plot -9 -9 -9 9 9 9 9 -9 -9 -9 /
\plot -9 -28 9 -28 / 
\setsolid
\arrow <10pt> [.2,.67] from 33 8 to 5  -26
\put {${\mathcal{R}}(I)$} at 41 13
\endpicture
\caption{\label{fig:firstreturn0} The figure shows two cases where $I$ is terminating 
together with the  interval $\mathcal{R}(I)$ marked in dots. }
\end{figure}

\begin{figure}[htb]
    \centering
    \subfloat[Partition of $I^{\infty}$]
{{ \def\svgwidth{214pt}
 \scriptsize
 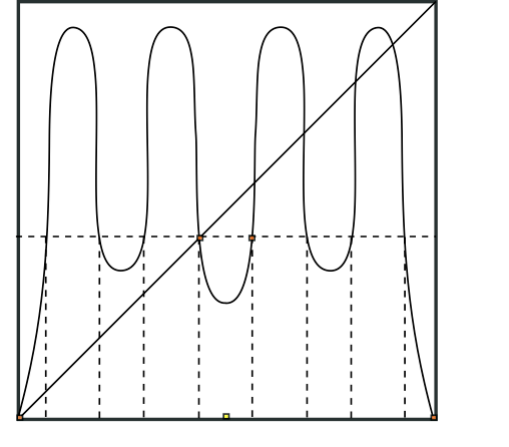
}}%
    \qquad
    \subfloat[Labelling of monotone branches]{{\def\svgwidth{214pt}
 \scriptsize
 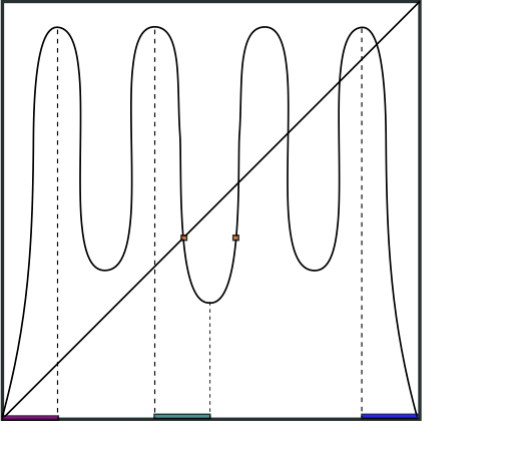
 }}%
    \caption{}%
    \label{fig:Y_i}\label{fig:monotone}%
\end{figure}

We let $Y_{i}$ with  $-a\leq i\leq a$ denote the components of $I^{\infty}\setminus f^{-r}(\alpha)$ labeled as follows: $Y_{0}=\mathcal{R}(I)$, $Y_{-1}\neq Y_0$ is the other component that contains $\alpha$ in its boundary and $Y_{-a}$ and $Y_{a}$ are the components that contain $\beta$ and $\tau(\beta),$ respectively, in their boundaries (see Figure \ref{fig:monotone}).
We let $\tilde{\mathcal{Y}}_{\gamma}$ denote the monotone branch of $f^r| I^{\infty}$ that contains $\gamma$, where $\gamma \in\{\alpha,\beta,\tau(\beta)\}\label{notation: top level}$ (see Figure \ref{fig:monotone}).

%
%

\subsection{The generalized enhanced nest}\label{subsec:gen}
We will extend the construction of the enhanced nest of \cite{KSS} to
cover the renormalizable case.
We will make use of the following combinatorially defined return time.

\begin{lem}[\cite{KSS} Lemma 8.2]\label{KSS82}
Let $I\owns c$ be a $\omega(c_0)$-critical puzzle piece. Then there exists a positive integer $\nu$ with $f^{\nu}(c)\in I$ such that the following holds. Let $U_{0}=\mathrm{Comp}_{c}f^{-\nu}(I)$ and $U_{j}=\mathrm{Comp}_{f^{j}(c)}f^{-(\nu-j)}(I)$ for $0\leq j\leq\nu$. Then
\begin{enumerate}[(1)]
\item $\#\{j:U_{j}\cap\mathrm{Crit}(f)\neq\emptyset, 0\leq j \leq\nu-1\}\leq b^{2}$, and
\item  $U_{0}\cap\omega(c)\subset\mathrm{Comp}_{c}(f^{-\nu}(\mathcal{L}_{f^{\nu}(c)}(I))).$
\end{enumerate} 
\end{lem}

For each  puzzle piece $I\owns c_{0}$ we let $\nu=\nu(I)$ be the smallest positive integer with the properties specified by Lemma \ref{KSS82}. We define
$$\mathcal{A}(I)=\mathrm{Comp}_{c_{0}}f^{-\nu}(\mathcal{L}_{f^{\nu}(c_{0})}(I)),$$
$$\mathcal{B}(I)=\mathrm{Comp}_{c_{0}}f^{-\nu}(I).\label{A, B in enhanced nest}$$

\begin{figure}[H] \centering \def\svgwidth{200pt}
\scriptsize
 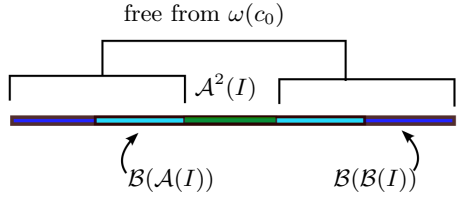
 \caption{The intervals shown are  ${\mathcal A}^2(I)\subset {\mathcal B}(\mathcal A(I))\subset {\mathcal B}^2(I)$. 
 The pair of intervals ${\mathcal B}^2(I)\setminus {\mathcal A}^2(I)$ is disjoint from the postcritical set.}
 \label{AB1}
\end{figure}

By construction
$\mathcal{A}(I)\subset\mathcal{B}(I)$ and
$(\mathcal{B}(I)\setminus\mathcal{A}(I))\cap\omega(c_{0})=\emptyset$,
giving a mechanism for obtaining \emph{free space}, i.e. space disjoint from $\omega(c_0)$, on the outside and inside of the interval $\mathcal{B}(\mathcal{A}(I)).$ See Figure \ref{AB1}.



Next, let $T=5b$, where $b$ is the number of critical points in $\omega(c)$,
and define
$$\mathcal E(I)=\left\{ \begin{array}{ll}
\Gamma^T \mathcal B \mathcal A(I) &\mbox{ if }I \mbox{ is non-terminating,}\\
\mathcal{L}_{c_{0}}(\mathcal{R}(I)) &\mbox{ if }I \mbox{ is terminating.}\end{array}
\right. \label{next level of enhanced nest}$$

The {\em generalized enhanced nest} associated to a nice interval $I\ni c_0$ is inductively defined by 
$$I_0=I\mbox{ and } I_{n+1}=\mathcal E(I_{n})\mbox{ for }n\ge 0 \label{the enhanced nest}.$$
For simplicity, we will refer to the elements from this nest as intervals from the enhanced nest. However it is important to remark this definition differs from the one introduced in \cite{KSS} in the presence of terminating intervals.
If $I_n$ is non-terminating, we let
$p_n>0$\label{iterate for pullback of enhanced nest} be so that $I_{n+1}$ is a component of $f^{-p_n}(I_n)\label{p-n}$. Under these circumstances, Lemma \ref{KSS82} and the definition of successor imply that $I_{n+1}$ is a pullback of $I_n$ of bounded order; $\mathcal{BA}(I_{n})$ is a pullback of 
$I_{n}$ with order bounded from above by $2b^2$,
and each successor $\Gamma^{i+1}(\mathcal{BA}(I_{n}))$ is a pullback of
order at most $2b-1$ of $\Gamma^{i}(\mathcal{BA}(I_{n})).$
Thus
$I_{n+1}$ is a pullback of $I_{n}$ with order at most $2b^2+5b(2b-1)$.
Note that if $c_0$ is a critical point of odd order, $I_n$ is never terminating 
and so we always have that $I_{n+1}=\Gamma^T\mathcal B \mathcal
A(I_n)$. Finally, the interval $I_n$ is called \emph{Feigenbaum} if it
is periodic; necessarily of period two under  $R_{I_{n-1}}.$

For each $n\ge 0$, let $r(I_{n})$ be the minimal return time for $x\in I_{n}\cap \omega(c_0)$ back to $I_n$ 
and $\hat{r}(I_{n})$ be
the maximal entry time for $x\in\omega(c_{0})$ to $I_{n}$.\label{return times to enhanced nest}

Although the following was stated in \cite{KSS}, we provide a proof, since the proof was not given in full in \cite{KSS} and the assumptions here are slightly weaker.
\begin{lem}[cf. \cite{KSS} Lemma 8.3]\label{KSS83}
Assume that $I_n$ and $I_{n+1}$ are both non-terminating, then the following holds:
\begin{enumerate}
\item $3r(I_{n+1})\geq p_{n},$
\item $\hat{r}(I_{n})\leq(1/2^{5b-1})r(I_{n+1}).$
\end{enumerate}
\end{lem}
\begin{pf}
For each $n\geq 0$  and for $0\leq j\leq T-1$ define
 $$L_n=\mathcal{A}(I_n),   \,  \,\,  \,  \,  \, \,  M_{n,0}= K_n= \mathcal{B}(L_n)   \,  \,\,  \,  \, \, \,  \text{and}    \,  \, \, M_{n,j+1}= \Gamma(M_{n,j}).$$ 
 
Let $s_n$ and $t_n$ be such that $L_n$ is a pullback of $I_n$ under $f^{s_n}$ and $K_n$ is a pullback of $L_n$ under $f^{t_n}$. For each $0\leq j \leq T-1$ let $q_{n,j}$ be such that $M_{n,j+1}$ is a pullback of $M_{n,j}$ under $f^{q_{n,j}}$. Finally define $p_n= s_n+t_n+q_{n,1}+\ldots +q_{n,T}.$
Using the same proof given in Lemma 8.2 on \cite{KSS} we get 
\begin{itemize}
\item[(i)] $2b^2r(I_{n+1}) \geq s_n\geq r(I_n),$
\item[(ii)] $b^2r(K_n)\geq t_n\geq r(L_n).$
\end{itemize}
For each $ j$ consider the chain $\{G_{j_i}\}_{i=0}^{q_{n,j}}$, with $G_{j_{q_{n,j}}}= M_{n,j-1}$ and $G_{j_0}= M_{n,j}$. By definition of $M_{n,j}$ we know that  the interval $G_{j_i}$  does not contain $c_0$ if $0<i< q_{n,j}$. So we conclude that $r(M_{n,j})\geq q_{n,j}$.  The fact that $M_{n,j-1}$ is non-terminating for every $0\leq j\leq T$ implies that $R_{M_{n,j-1}}(M_{n,j})\cap M_{n_j}= \emptyset$. Therefore,  $q_{n,j}\geq 2r(M_{n,j-1})$. Putting these two inequalities together we get
\begin{itemize}
\item[(iii)] $r(M_{n,j})\geq q_{n,j}\geq 2r(M_{n,j-1}).$
\end{itemize}
Since $M_{n,1}$ is the smallest successor of $M_{n,0}$ and $M_{n,0}\subset I_n$ we have $\hat r(I_n)\leq \hat r(M_{n,0})\leq q_{n,1}.$ Using this fact and equation (iii) we get
$$\hat r(I_n)\leq q_{n,1}\leq \frac{1}{2}q_{n,2}\leq \ldots \leq \frac{1}{2^{T-1}}r(I_{n+1}).$$
The previous inequality together with equations (i) and (iii) gives us
\begin{eqnarray*}
   p_n  &=&   s_n+t_n+q_{n,1}+\ldots +q_{n,T} \nonumber \\
   &\leq & 2b^2r(L_n)+  b^2r(K_n)+ q_{n,1}+\ldots +q_{n,T}\nonumber \\
   &\le& 2b^2r(K_n)+  b^2r(K_n)+ q_{n,1}+\ldots +q_{n,T}\nonumber \\
   &\le& 2b^2q_{n,1}+  b^2rq_{n,1}/2+ q_{n,1}+\ldots +q_{n,T}\nonumber \\
 &\le& \frac{b^2}{2^{T-1}}r(I_{n+1})+  \frac{b^2}{2^{T-2}}r(I_{n+1})+ \left( \frac{1}{2^{T-1}}+\ldots+ \frac{1}{2}\right)r(I_{n+1}) \nonumber \\
 &\le& \frac{3b^2}{2^{T-2}}r(I_{n+1})+ \left( \frac{1}{2^{T-1}}+\ldots+ \frac{1}{2}\right)r(I_{n+1}) \nonumber \\
 &\leq& 3r(I_{n+1}).
\end{eqnarray*}
\end{pf}

\section{Real bounds associated to the enhanced nest}\label{sec:real bounds}


In this section we will consider maps $f\in
\mathcal{A}^{3}_{\underline{b}}$ and show geometric bounds for the
intervals of the generalized enhanced nest
Before stating the main result we need to introduce some terminology.



\medskip
Given a constant $\rho>0,$ a nice interval $I$ is called:
\begin{equation*}\begin{array}{ll}\label{interval geometry}
\mbox{$\rho$-\emph{nice}}&\mbox{ if for each $x\in I\cap\omega(c_{0}),$ }  (1+2\rho)\mathcal{L}_{x}(I)\subset I;\\
\mbox{$\rho$-\emph{free}} &\mbox{ if }((1+2\rho)I\setminus(1+2\rho)^{-1}I)\cap\omega(c_0)=\emptyset; \\
\mbox{$\rho$-\emph{externally free}} &\mbox{ if there exists a nice interval $J\supset (1+2\rho)I,$ so that $J\supset I$ is a nice pair}\\
&\mbox{ and }J\cap \omega(c_{0})\subset I;\\
\mbox{$\rho$-\emph{internally free}} &\mbox{ if there exists a nice interval $J'$, so that  $(1+2\rho)J'\subset I,$}\\
&\mbox{ $J'\subset I$ is a nice pair and $I\cap\omega(c_0)\subset J'$.}\\
\end{array}\end{equation*}
We say that $I$ is $\rho$-\emph{strongly nice}\label{rho-strongly nice} if it is $\rho$-nice and if
for each $x,y\in I\cap\omega(c_{0})$ either $\mathcal{L}_{x}(I)=\mathcal{L}_{y}(I)$,
$(1+2\rho)\mathcal{L}_{x}(I)\cap \mathcal{L}_{y}(I)=\emptyset$ or  $(1+2\rho)\mathcal{L}_{y}(I)\cap \mathcal{L}_{x}(I)=\emptyset$. 
We will say that an interval $J$ is \emph{well-inside}\label{well-inside} an interval $I$ if $(1+2\rho)J\subset I$ for $\rho$ universal. If $\rho$ is large, then we say that $J$ is \emph{deep-inside} $I$.


The purpose of this section is to prove the following  theorem. 
\begin{thm}[Real geometry of the enhanced nest]\label{real geometry}
Suppose that $f\in \mathcal{A}^{3}_{\underline{b}}$. 
There exists $\varepsilon_f>0$ such that the following holds.
Assume that $c_0$ is a critical point at which $f$ is persistently
recurrent, and that either $c_0$ has even order or
that every critical point in $\omega(c_0)$ has odd order.
Suppose that $I_0\owns c_0$ is a nice interval with $|I_0|<\varepsilon_f.$ 
Let $I_{0}\supset I_{1}\supset\dots$ be the generalized enhanced nest for $f$ at  $c_0$. Then the following hold:
\begin{enumerate}[(a)]
\item There exists $\rho>0$ such that if $I_n$ is non-terminating, then $I_n$ is $\rho$-nice. In addition, if $I_{n-1}$ is non-terminating then $I_n$ is $\rho$-externally and $\rho$-internally free, where the externally free space is given by an interval $J\supset (1+2\rho)I_n$ and the internal free space is given by an interval $J'\subset (1+2\rho)J'\subset I_n.$  Moreover, if $c_0$ is even $|J'|\geq\rho|I_n|$; if $c_0$ is odd then for each $\nu>0$ there exists $\rho'>0$ so that if $|I_{n-1}|/|I_n|<\nu$ then $|J'|\geq\rho'|I_n|$.

\item Suppose that $I_{n-1}$ is terminating and $I_n$ is non-terminating. Then for each $\nu>0$ there exists $\rho'>0$ so that if $|I_{n-1}|/|I_n|, |I_n|/|I_{n+1}|<\nu$, then $I_n$ is $\rho'$-free.

\item For any $C'>0$ there exists $C>0$ such that if  $I_n$ is non-terminating and there is some $x\in\omega(c_0)
\cap I_{n}$ with $(1+2C)\mathcal{L}_{x}(I_{n})\subset I_{n}$  
then $(1+2C')I_{n+1}\subset  I_n$.
\item For any $C'>0$, there exists $C>0$ such that if $n\ge 0$, $I_n$ is non-terminating and
$C$-nice, then $I_{n+1}$ is $C'$-nice and $C'$-externally  free.

\item Suppose that $I_{n}$ is non-terminating, then for each $\varepsilon'>0$ there exists $\varepsilon>0$ such that if there is some $x\in\omega(c_0)\cap I_{n}$ with $|\mathcal{L}_{x}(I_{n})|\leq\varepsilon|I_{n}|$, then $|I_{n+1}|\le \varepsilon'|I_n|$.

\item For any $C'>0$ there exists $\varepsilon>0$ so that if $I_{n}$ and $I_{n+1}$ are both non-terminating and 
$|I_{n+1}|\le \varepsilon |I_{n}|$, then $I_{n+2}$ is a $C'$-strongly nice and $C'$-externally free.
\item For any $C'>0$, there exists $C>0$ such that if $I_n$ is non-terminating, $I_{n+1}$
is terminating and $I_{n}\supset(1+2C)I_{n+1}$, then
$I_{n+1}^{\infty}$ is $C'$-externally free.
\end{enumerate}
\end{thm}

\begin{remark}
\begin{enumerate}
\item In the statements, the constants $C,C'$ are large, while
  $\varepsilon, \varepsilon'$ are small. This leads to intervals which
  are deep-inside, large free space, large ratios between intervals
  etc. For this reason, the bounds with $C,C'$ and
$\varepsilon, \varepsilon'$ will be referred to as  ``big bounds.''

\item If $c_0$ is even, because of symmetry about the critical point,
in Part (a) we have that $I_{n+1}$ is $\rho$-nice and $\rho$-free. Also, in Part (a) notice that if $|F'|\geq\rho'|I_n|$, then the components of $I_n\setminus F'$ have size comparable to $I_n$. If $c_0$ is odd, then critical points can be close to the boundary of the return domains. This causes the map to lose free space inside the intervals. Because of this, the statement in Part (c) is no longer sufficient for our purposes, so we supplement it with (e).
\item  To get internal free space, in general, we require that the lengths of the intervals in few consecutive levels are comparable, see Corollary \ref{cor:delta free}.

\item In Part (a) the free space around $I_{n}$ is given directly by the construction of the enhanced nest, since return domains to $I_{n-1}$ are well-inside $I_{n-1}$. The return domains to $I_{n-1}$ can be close to the boundary, but when this happens they are very small. In case $I_{n-1}$ is terminating, the free space, when we have it, is a consequence of certain real bounds, rather than the construction of the enhanced nest.
\end{enumerate}
\end{remark}

Except for $\varepsilon_f$, which depends on $f$, the above bounds are {\it universal}, or in Sullivan's terminology {\em beau}.
In other words, one can choose the constants in the above theorem so that they do not depend on $f$ but {\em only}
on the vector $\underline b=(\ell_{1},\dots, \ell_{b}),$  provided we begin the construction with an interval $I_0$ small enough.  

Theorem~\ref{real geometry} generalizes Proposition 8.1 of \cite{KSS}, 
to maps $f$ that are $C^{3}$ with critical points of any order (i.e. not necessarily even order) and allowing for maps which are possibly (infinitely) renormalizable.
Proposition~\ref{delta free} is an addendum to Theorem~\ref{real geometry} that will also be used later on.  In the case of terminating intervals, some specific details of the proof of Proposition~\ref{delta free} will be frequently used throughout the paper.

An important difference between this work and \cite{KSS} is that we no longer get that $$(I_n\setminus(1+2\rho)^{-1}I_n)\cap\omega(c_0)=\emptyset$$ in the case where all critical points
in $\omega(c)$ have odd order, which means that there can be small
return domains to $I_n$ that intersect $\omega(c_0)$ close to the
boundary of $I_n$. This complicates the proof significantly.

\medskip

The proof of Theorem \ref{real geometry} occupies most of this section.

\subsection{Distortion bounds and their consequences}

\begin{thm}[\cite{vSV}, Theorem C and \cite{Li-Shen}, Theorem A]\label{real Koebe}
Suppose that $f\in\mathcal{A}^{3}_{\underline{b}}$. Then one has the following properties:
\begin{enumerate}[(1)]
\item \textit{Improved Macroscopic Koebe Principle.} For each $\xi>0$,
  there exists $\xi'>0$ such that if $I$ is a nice interval, $V$ is a
  nice interval that is $\xi$-well-inside $I$, $x\in I$ and $f^{k}(x)\in V$ with $k\geq 1$ not necessarily minimal, then the pullback of $V$ along $\{x, f(x),\dots f^{k}(x)\}$ is $\xi'$-well-inside the return domain to $I$ containing $x$.
\item \textit{Improved Koebe Principle.} For each $\tau>0$,
there exists $K>0$ and $\xi=\xi(f,\tau)>0$ 
satisfying the following. Let $T\subset M$ be an open interval, 
and let
$J$ be a closed subinterval of $T$ such that the following hold:
\begin{itemize}
\item $J\cap\omega(c_0)\neq\emptyset$,
\item $f^{n}|T$ a diffeomorphism,
\item $|f^n(T)|<\xi$,
\item $f^{n}(J)$ $\tau$-well-inside $f^{n}(T)$.
\end{itemize}
Then $f^{n}|J$ has bounded distortion; that is, for any $x,y\in J,$
$$\frac{|Df^{n}(x)|}{|Df^{n}(y)|}\leq K.$$
Furthermore, $K\rightarrow 1$ as $\tau\rightarrow\infty$.
\item Negative Schwarzian derivative. For each critical point $c$ that is not in the basin of a periodic attractor, there exists a neighbourhood $U$ of $c$ such that whenever $f^{n}(x)\in U$ for some $x\in I$ and $n\geq 0$, the Schwarzian derivative of $f^{n+1}$ at $x$ is negative.
\end{enumerate}
\end{thm}

\medskip
Theorem~\ref{real Koebe}(2) follows from \cite[Theorem A]{Li-Shen}. There the
result is stated for mappings with all periodic
orbits repelling. However, the Theorem holds, without any change,  if we consider $f$ restricted to puzzle pieces intersecting $\omega(c_0);$  since all periodic points contained in those pieces are repelling.

\medskip

From now on we assume the following
\medskip

\noindent\textbf{Standing Assumptions:}\label{standing assumptions}
\begin{itemize}
\item $f\in\mathcal A_{\underline b}^3$ is persistently recurrent on
  $\omega(c)$, where $c\in\crit(f)$. 
\item Either $c_0\in\omega(c)$ is even or $c_0$ is odd and every critical point in $\omega(c)$ is odd.
\item Any nice interval $I$ that intersects $\omega(c_0)$ is so small
  that for any chain $\{G_j\}_{j=0}^s$ with $G_s=I$ and
  $G_0\cap\omega(c)\neq\emptyset$ each
  $G_j$ avoids a neighbourhood of any immediate basin of attraction of
  any periodic attractor or any parabolic cycle.
\item If $I$ is a nice interval containing a critical point
  $c'\in\omega(c)$, then $I$ is so small that the remark on page \pageref{remark SA} holds for $I$.
\item Suppose  $I$ and $J$ are nice intervals with $(1+2\delta)J\subset I$ with $\delta>0$
    universal. We assume $|I|<\xi,$ where
    $\xi=\min\{\xi(f,\tau):\delta\leq \tau\leq 10^{10}\},$ and $\xi(f,\tau)$ is the constant given by Theorem \ref{real Koebe} (2). 
\end{itemize}
We can assume these since, as we have  seen in Section
\ref{subsec:persistentrec}, we can guarantee that puzzle pieces
intersecting $\omega(c_0)$ are arbitrarily small. It is worth
observing that if $(1+2C)J\subset
I$ for $C$ large, then taking $J'=(1+2(C-1/2))J$, we can apply the
Improved Koebe Principal with $\xi=\xi(f,1/2)>0$, so that 
when $C$ is large our control on the distortion given by Theorem \ref{real Koebe}
remains bounded. 
\medskip

Fact 9.1 of \cite{KSS} does not generalize to our present setting (the proof in \cite{KSS} relies on the maps having negative Schwarzian derivative). However, the following analogue holds.

\begin{lem}\cite[Fact 9.1]{KSS}\label{KSSFact91}
For each $N\in\mathbb{N}$ and  $\rho>0$ there exists $\rho'>0$ such that the
following holds. Let $\{G_{j}\}_{j=0}^{s}$ and $\{G'_{j}\}_{j=0}^{s}$ be chains such that $G_{j}\subset G'_{j}$ for all $j$, $0\leq j\leq s$ and $G_0\cap\omega(c_0)\neq\emptyset.$ Assume that the order of $\{G'_{j}\}_{j=0}^{s}$ is at most $N$ and that $(1+2\rho)G_{s}\subset G'_{s}$. Then $(1+2\rho')G_{0}\subset G_{0}'.$ Furthermore, for fixed $N$, $\rho'\rightarrow\infty$ as $\rho\rightarrow\infty.$
\end{lem}
\begin{pf}
Let $\hat G_s=(1+\rho)G_s$  and consider the chain $\{\hat G_j\}^s_{j=0}$ with $G_j\subset \hat G_j.$ Observe that $\hat G_s$ is $\rho/2$-well-inside of $G_s'.$ Let $0<s_{1}<\dots<s$ be the times $j$ so that $G'_{j}$ contains a critical point. Then we can decompose $f^{s}\colon G'_{0}\rightarrow G'_{s}$ into at most $N$ maps of the form $f\colon  G'_{s_{j-1}}\rightarrow G'_{s_{j-1}+1}$ followed by $f^{s_{j}-s_{j-1}-1}\colon G'_{s_{j-1}+1}\rightarrow G'_{s_{j}}$. It follows from Theorem~\ref{real Koebe}(2) that each of the factors $f^{s_{j}-s_{j-1}+1}|_{\hat G_{s_{j-1}+1}}$ has bounded distortion. This, along with the fact that $f$ has non-flat critical points, implies that existence of $\rho'>0$ so that $(1+2\rho')G_0\subset \hat G_0\subset G_0'.$
\end{pf}

%


Two immediate consequences of the previous two results are:
\begin{cor}\label{nice}
For each $\rho>0$ there exists $\rho'>0$ so that the following holds. If $I$ is a nice interval, $J=\mathcal L_x(I)$ for some $x\in\omega(c_0)$ and $(1+2\rho)J\subset I$ then, $J$ is $\rho'$-nice. Furthermore, if $I$ is $\rho$-nice then $\rho'\to \infty$ as $\rho\to \infty$. 
\end{cor}

\begin{cor}\label{cor:free-space}
For each $\rho>0$ and each integer $N\geq 0$, there exists $\rho'>0$
so that the following holds. Let $I$ and $J$ be nice intervals with
$J$ a pullback of $I$ of order bounded by $N$ that intersects
$\omega(c_0)$. Then if $I$ is $\rho$-free, $J$ is $\rho'$-free. Furthermore, for fixed $N,$ we get that $\rho'\to \infty$ as $\rho\to \infty$. 
\end{cor}
Observe the above implies that if $I$ is $\rho$-externally free then $J$ is $\rho'$- externally free; however, we do {\em not} claim that if there exists $\rho>0$ such that $(I\setminus (1+2\rho)^{-1}I)\cap\omega(c_0)=\emptyset,$ then there exists $\rho'>0$ such that $(J\setminus(1+2\rho')^{-1}J)\cap\omega(c_0)=\emptyset$.
This statement may be wrong, because the internal free space that we pull back may be small 
compared to the total interval $J$, unless we also have external free
space to control the distortion of the mapping.

\medskip

The next two lemmas are closely related and are useful when a return
domain to $I$ is not well-inside $I$.
\begin{lem}\cite[Lemmas 2 and 3]{vSV}\label{KSSLemma94}
There exists a constant $\rho=\rho(\underline{b})>0$ with the following property. Let $I$ be a nice interval,
 consider a point  $x\in I$ which returns to $I$ in time $s$ and let
$I^{1}=\mathcal{L}_{x}(I)$.  Then either of the following holds:
\begin{enumerate}[(1)]
\item $(1+2\rho)I^{1}\subset I,$
\item there exists an interval $G_s\supset (1+2\rho)I^1$ with $I\subset G_s$, so that the chain $\{G_{j}\}_{j=0}^{s}$ with $G_{0}\owns x$ has intersection multiplicity bounded from above by a constant $N$ that depends only of the modality of $f$. 
\end{enumerate} 
\end{lem}

\subsection{The existence of suitable fundamental domains}


Let $J\subset J'$ be intervals and assume that $g\colon J\rightarrow J'$ is monotone. An interval $J^{*}$ is called a 
\emph{fundamental domain for $g$} if it is of the form $(x,g(x))$ (or $(g(x),x)$). We will repeatedly use
the property that the pullbacks of a fundamental domain under iterates of $g$ are disjoint.

\begin{lem}\label{fundamental domains}
For each $\sigma\in (0,1)$ there exist $\kappa>0$ and for each 
$f\in \mathcal A^3_{\underline b}$ there exists an integer $N$ such that the following holds. Let $I$ be a nice interval and let $J$ be a first return 
domain of $I$ with return time $n>N.$ Assume that $f^n\colon J\to I$ 
is an orientation preserving diffeomorphism 
with fixed point $p$ and let $J'$ be a component of $J\setminus \{p\}$. Then 
for each $x\in J'$ with $|x-p|\ge \sigma |J'|$
there exists a fundamental domain $F$ containing $x$ with  $|F|\ge \kappa \cdot d(F,p)$.
\end{lem}
\begin{pf} 
By Theorem IV.B in \cite{dMvS} there exists $\delta_1>0$ so that for each 
$f\in \mathcal A^3_{\underline b}$ there exists $N$
so that $Df^n(p)\ge 1+\delta_1$ provided the period $n$ of $p$ is at least  $N$.
For simplicity assume that $p=0$ and that $[0,x]\subset J'$. Write $l=\{0\}$, $j=[0,\sigma x]$, $r=[\sigma x,x]$, 
$t=[0,x]$ and  $g=f^n$. 
Define the cross-ratio distortion
$$B(g,t,j) = \frac{|g(j)||g(t)|}{|g(l)||g(r)|}\cdot\frac{|l||r|}{|t||j|}.$$
Take $y=\sigma x$ and take $\delta>0$ so that $g(\sigma x)/(\sigma x)=1+\delta$. Then 
$$ B(g,t,j)= \frac{1}{g'(0)}  \frac{g(\sigma x)}{(\sigma x)}\frac{g(x)}{x}\frac{x-\sigma x}{g(x)-g(\sigma x)}
\le \frac{1+\delta}{1+\delta_1}\frac{1-\sigma }{1-\sigma (1+\delta),}
$$
since $g(x)\ge x$. We would like to show that $\delta$ uniformly bounded from below for all $x\in J'$. If this is not the case, taking $\delta>0$ small we get $B(g,t,j)< 1$ contradicting Theorem 2.1 of Chapter 4 of \cite{dMvS}; which states that $B(g,t,j)\ge 1$.
\end{pf}

\subsection{The principal nest}
Recall that if $I$ is a nice interval that contains a recurrent point $x$, then
the principal nest about $x$ is defined inductively by
$I^0=I$ and $I^{n+1}=\mathcal{L}_x(I^n)$.
The next lemma is an important bound on the geometry of intervals in
the principal nest.

\begin{lem}\cite[Theorem A]{vSV}\label{KSSLemma91}
There exists $\delta>0$ such that if $I$ is a nice interval, $x\in I$ is recurrent and $R_{I}(x)\notin\mathcal{L}_{x}(I)$, then for each $d\in\mathbb{N}$ if $R_{\mathcal{L}^{d}_{x}(I)}(x)\notin\mathcal{L}^{d+1}_{x}(I)$, then $$(1+2\delta)\mathcal{L}^{d+2}_{x}(I)\subset\mathcal{L}^{d+1}_{x}(I).$$
\end{lem}

\begin{lem}\cite[Lemma 5.5]{Shen}\label{SLemma55}
There exists a constant  $\rho=\rho(\underline{b})>0$ with the following property. Let $I$ be a periodic interval
of sufficiently large period $s$.  Then $(1+2\rho)I$ does not 
contain $f^i(I)$ for $i=1,\dots,s-1$. Moreover, the chain $\{G_{j}\}_{j=0}^{s}$ with $G_{s}=(1+2\rho)I$
and $G_0\supset I$ has the following two properties
\begin{itemize}
\item $\{G_{j}\}_{j=0}^{s}$ has intersection multiplicity at most  four
and 
\item the map $f^s$ does not have a critical point in $G_0\setminus I$.
\end{itemize}
\end{lem}

We say that $I:=I^{0}\supset I^{1}\supset \dots\supset I^{m}$  is
\emph{a central cascade} if $m\geq 2$ and the return time of
  $I^m$ to the intervals $I^{0}, I^{1},\dots,I^{m-1}$ is always the
  same.

\begin{lem}[cf. \cite{KSS} Lemma 9.5]
\label{bounds1}
For any $\delta>0$ there exist $\kappa>0$ and $C>0$ with the following properties. Let $I$ be a nice interval (as usual, assumed to be sufficiently small),
having a central cascade $I:=I^{0}\supset I^{1}\supset \dots\supset I^{m}$ with $m\geq 2$ and let $r$ be the return time of $I^{1}$ to $I^{0}$. 
\begin{enumerate}[(1)]
\item
If $|I^{2}|\geq \delta |I^{0}|$, then for any critical point $c$ of the map $R_{I}|I^{2}$ we have
$$|f^{r}(c)-c|\geq\kappa |I^{0}|\; and \; |Df^{r}(x)|\leq C\ for\ all\ x\in I^{2}.$$
\item 
If $|I^{1}|\geq \delta |I^{0}|$ and we let $\tilde I=(1+2\delta)I$ and
$\tilde I^1=\comp_{I^1}(f^{-r}(\tilde{I}))$ the following
holds. Suppose that $f^r$ extends to a map from $\tilde I^1$ to
$\tilde I,$ with the same set of critical points as $f^r|I^1,$ that can be decomposed into a finite composition of maps with bounded distortion and polynomials. Then for any critical point $c$ of the map $R_{I}|I^{1}$ we have
$$|f^{r}(c)-c|\geq\kappa |I^{0}|,\; and \; |Df^{r}(x)|\leq C\ for\ all\ x\in I^{1}.$$
\end{enumerate}
\end{lem}
\begin{pf} As before, let $\underline b=(\ell_1,\dots,\ell_b)$. 
By Lemma~\ref{KSSLemma94},  Theorem \ref{real Koebe} and since $c\in I^2$, there exist beau constants $\rho'>0$, $K<\infty$ and an integer $N$ (only depending on $\underline b$)  and an interval $J$ with either $J=I^2$ or with $J$ of the form $(c-\delta,c+\delta)$ so that
$|J|\ge \rho' |I^2|$,  $f^r(J)\subset (1+\rho')I_0$ and so that $f^{r}\colon J\to f^r(J)$ 
can be written as  a composition of at most $N$ maps whose distortion is bounded by 
 $K$ and polynomials $z\mapsto z^{\ell_{i}}$. Taking $J=I^2$ the 2nd inequality in (1) follows. To prove the 1st inequality, 
 take $J$ as above of the form $(c-\delta,c+\delta)$.  So there exists a beau constant $\hat K<\infty$ so that for all $\epsilon\in (0,1)$ one has 
$|Df(x)|\le \epsilon \hat K |f(J)|/|J|$ for all $x\in (c-\epsilon \delta,c+\epsilon \delta)$ and $|Df^{r-1}(y)|\le \hat K^{N-1} |f^r(J)|/|J|$ for all $y\in f(J)$. It follows that $|Df^r(x)|\le \epsilon \hat K^N |f^r(J)|/|J|$  for all $x\in (c-\epsilon \delta,c+\epsilon \delta)$ and  $|f^r(J)|/|J| \le ((1+\rho')/\rho') |I_0|/|I_2|\le ((1+\rho')/(\delta \rho'):=K^*$.  
If we take $\epsilon >0$ so that $\epsilon \hat K^N  K^* <1/2$ then $|Df^r(x)|\le (1/2)$ for all  $x\in (c-\epsilon\delta,c+\epsilon\delta)$. 
 Now choose $\kappa = \epsilon  /(16\delta\rho')$ and assume by contradiction that $|f^r(c)-c|\le \kappa |I^0|$. 
Then $|f^r(c)-c|\le (\epsilon /(16\delta\rho') |I^0|\le \epsilon/(16\rho')|I^2| \le (\epsilon /16)|J|=\epsilon\delta/8$. Thus we get 
that $f^r$ maps $(c-\epsilon\delta,c+\epsilon\delta)$ into itself and $|Df^r(x)|\le 1/2$ on this interval. Hence $f^r$ has an attracting fixed point, 
which contradicts that $\omega(c)$ is minimal. Thus we have proved the 1st inequality in (1) by contradiction.  The inequalities in (2) follow as in (1). 
\end{pf}

\subsection{Geometry of pullbacks}
\begin{lem}\label{hatI} Let $I$ be a nice interval.  Assume $z'\notin I$ has first entry time to $I$ equal to $k>0$ and that there exists $l>k$ such that $f^l(z')\in I$. Letting $z=f^k(z')$,
$I^1=\mathcal{L}_z (I)$, $\hat I=\mathcal{L}_{z'}(I)$ and 
$\hat K=\mathcal{L}_{z'}(\hat I),$ we have that  $\hat K\subset \mathcal{L}_{z'}(I^1)$.
\end{lem}
\begin{pf} 
Since $\hat I$ is a pullback of $I$ and $\hat K$ is a pullback of $\hat I,$ 
$\hat K$ is a pullback of $I$, by say $f^s$.
Since $k$ is the first entry time of $z'$ to $I$ we have that $k<s$. Let $K'=\comp_zf^{-(s-k)}(I)$.
Since $z'\in \hat K$,  $z\in K'$ 
and $I^1$ is a pullback of $I$ containing $z$, we have that $K'\subset I^1.$
Hence $\hat K\subset \mathcal{L}_{z'}(I^1)$.
\end{pf} 


\begin{lem}[cf. \cite{Shen} Proposition 4.1] \label{KSSLemma93}
For any $C>0$ and any $d\in\mathbb{N}$, there exists $C'>0$ so that the
following holds. Assume $I\supset J$ are nice intervals with nested or
disjoint pullbacks
and $J\supset \mathcal{L}^{d}_{x}(I)$ for some $x\in\omega(c_0).$ If $(1+2C')J\subset I$,
then $(1+2C)\mathcal{L}_{y}(J)\subset\mathcal{L}_{y}(I)$ for any $y\in\omega(c_0).$

If $c_{0}$ is odd, then for any $\delta>0$ and $C>0$ there is exists $C'>0$ 
so that  the following holds. Let $I$
be a nice interval with $((1+2\delta)I\setminus
I)\cap\omega(c_{0})=\emptyset.$ Suppose that $J\subset I$ is a nice
interval with $J\supset \mathcal{L}_{x}(I)$ for some
$x\in\omega(c_0)$.
If $|J|\le | I|/(1+2C')$ then
$|\mathcal{L}_{y}(J)|\le |\mathcal{L}_{y}(I)|/(1+2C)$ for any $y\in\omega(c_0)$.
\end{lem}



\begin{pf}

We will prove the first part of the lemma by induction on $d$. Let us begin with the case $d=1$. In \cite[Theorem B2]{vSV}, the result is stated for $J=\mathcal{L}_{c}(I)$, but the proof holds
for $J\supset \mathcal{L}_{c}(I)$ provided $J,$ $I$ are a nice pair (definition p.\pageref{nice pair}). The fact that $d=1$ implies that $I$ and $J$ are a nice pair so the result follows. Before we prove the general case, we observe that the constant $C'$ will depend on $d$.
Assume the result holds for all $d'<d$ and let us show it for $d$. 
Since $y\in\omega(c_0),$ there exists $s>0,$ minimal, so that $z=f^s(y)\in J$. Consider the chains  $\{G_{j}\}_{j=0}^{s}$ with $G_s=J$ and $G_0=\mathcal L_y(J)$ and $\{G'_{j}\}_{j=0}^{s}$ with $G'_s=I$ and $G'_0=\comp_{y}f^{-s}(I).$
Let $0<s'<s$ be maximal so that $f^{s'}(y)\in I$. If there exists no such $s'$, then $z$ is the first entry time of $y$ to both $I$ and $J$, so the chain $\{G'_{j}\}_{j=1}^{s}$ is disjoint, therefore it has order bounded by $b$.
In this case the result follows from Lemma \ref{KSSFact91}.  If $s'$
is defined, let $z'= f^{s'+1}(y)$ and $\hat I= G'_{s'+1}$. 
Observe $\hat I=\mathcal L_{z'}(I).$
Since $I$ is small and it intersects $\omega(c_0)$ the entry time of
$f^{s'}(y)$ to $I$ is bigger than one, so $z'\notin I$. Let $J'=
G_{s'+1}$. 
By the definition of $s'$, the chain  $\{G'_{j}\}_{j=s'+1}^{s}$ is
disjoint, so Lemma \ref{KSSFact91} implies that $J'$ is deep-inside
$\hat I$. Observe that $\mathcal L^i_{z}(I)=\mathcal L^i_{x}(I)$ for
$i<d$ since $z\in J$. From Lemma \ref{hatI} we conclude that
$J'\supset \mathcal L^{d-1}_{z'}(\hat I)$, and the result follows from
the induction hypothesis.  

The proof of the 2nd part of the lemma is fairly involved and will be given in the next subsection.\end{pf}

\subsubsection{Proof of the 2nd part of Lemma~\ref{KSSLemma93}}\label{subsec:prooflemma93}

If $I$ is a nice interval that contains a recurrent point $x$, then define
\begin{equation*}
\mathcal{C}_{x}(I):= \left\{ 
\begin{array}{rl}\label{first non-central return}
I^{m} & \mathrm{if\ } I\ \mbox{is\ non-terminating\ and}\\ 
I^{\infty}\ &\mathrm{otherwise},
\end{array}\right.
\end{equation*}
where $m$\label{m time} is minimal such that $R_I(x)\notin I^m$
We will sometimes omit the $x$ from the notation when it will not cause confusion. 

By \cite[Lemma 15]{vSV}, if $U\subset I$ is a nice interval and $U\cap\mathcal{L}_x(I)=\emptyset$, then
any pullback of $U$ that contains $x$ is contained in $\mathcal{C}_x(I)$.

\begin{lem}[cf. \cite{vSV} Lemma 16]\label{small domains lemma}
Let $\delta>0$ and suppose that $c_{0}$ is odd (so every critical point in $\omega(c_{0})$ is odd). There exists a function
$\rho_{1}:\mathbb{R}^{+}\rightarrow\mathbb{R}^{+}$ such that
$\rho_{1}(\varepsilon)\rightarrow 0$ as $\varepsilon\rightarrow 0$
with the following properties. Let $I^{0}$ be a nice interval
containing  a point $c\in \omega(c_0)$ with $((1+2\delta)I^0\setminus I^0)\cap\omega(c_0)=\emptyset$.
Let $V$ be an interval such that
$|V|/|I^{0}|<\varepsilon$ and such that $V\subset I^{0}$ is a nice pair.
Suppose that $f^{s}(x)\in V$ for some $s>0$ and some $x\in I^{0}\cap\omega(c_0).$ Let
$G_s=V,\hat{G}_{s}=I^{0}$ and let $\{G_{j}\}_{j=0}^{s},
\{\hat{G}_{j}\}_{j=0}^s$ be the chains so that
$$\hat{G}_{i}=\comp_{f^{i}(x)}(f^{-(s-i)}(\hat{G}_{s}))\  \mathrm{and}\
G_{i}=\comp_{f^{i}(x)}(f^{-(s-i)}(G_{s})).$$
Then there exist $t$, $0\leq t<s$ and an interval $\hat{G}^{1}_{t}$
with $G_{t}\subset \hat{G}^{1}_{t}\subset \mathcal{L}_{f^t(x)}(I^{0})$
such that $$\frac{|G_{t}|}{|\hat{G}^{1}_t|}<\rho_1(\varepsilon).$$
If $t>0,$ then:
\begin{enumerate}[(1)]
\item $G_{t}\subset \hat{G}^{1}_{t}$ is a nice pair and
\item each pullback of $\hat{G}^{1}_{t}$ that intersects $c$ is
contained in $\mathcal{C}_{c}(I^{0})$.
\end{enumerate}
\end{lem}
\begin{pf}
Let $I^1=\mathcal{L}_c(I^0)$.
Let \mbox{$\phi:I^{1}\rightarrow I^{0}$}
be the first return map of $I^1$ to $I^0$.
Since $((1+2\delta)I^0\setminus I^0)\cap\omega(c_0)=\emptyset$, by the Improved Koebe Principle, Theorem \ref{real Koebe}(2), any first entry map to $I^{0}$ decomposes into at most $b$ maps of the form $z\mapsto z^{\ell},$ where $\ell$ is odd, followed by a diffeomorphism with bounded distortion.
Hence, there exists a function
$\hat{\rho}:\mathbb{R}^{+}\rightarrow\mathbb{R}^{+}$ with
$\hat{\rho}(\varepsilon)\rightarrow 0$ as $\varepsilon\rightarrow 0$
so that any pullback, $V',$ of $V$ by a first entry map to $I_{0}$
satisfies $$\frac{|V'|}{|\mathcal{L}_{V'}(I^{0})|}<\hat{\rho}(\varepsilon).$$
It is worth noticing that this is where we use the fact that all critical points in
$\omega(c_0)$ are odd.

Let $0=:t_{0}<t_{1}<\dots<t_{k}:=s$ be the integers such that
$f^{t_j}(x)\in I^{0}$. By the preceding comment,
$$\frac{|G_{t_{k-1}}|}{|\hat{G}_{t_{k-1}}|}=\frac{|G_{t_{k-1}}|}{|\mathcal{L}_{f^{t_{k-1}}(x)}(I^{0})|}<\hat{\rho}(\varepsilon).$$
Hence, if $f^{t_{k-1}}(x)\notin I^{1}$, then all required properties
hold for $t=t_{k-1}$ taking $\hat G^{1}_{t}=\mathcal{L}_{f^{t}(x)}(I^{0})$
and $\rho_{1}=\hat{\rho}$, so from now on we assume that
$f^{t_{k-1}}(z)\in I^1$. Now let $k'$ be minimal such that
$$f^{t_{k'}}(x), f^{t_{k'+1}}(x),\dots, f^{t_{k-1}}(x)\in I^1.$$
Since $R_{I^0}|I^1$ is monotone, we have that $\hat{G}_{t_{k'}}\subset
\hat{G}_{t_{k'+1}}\subset\dots\subset \hat{G}_{t_{k-1}}=I^{1}\subset
\hat{G}_{t_k}=I^0.$  We can assume that $\frac{|I^1|}{|I^0|}$ is not
small: otherwise, if $k'>0$, the required properties hold for
$t=t_{k'-1}$ and $\hat G^{1}_{t}=\mathcal{L}_{f^t(x)}(I^0)$, and if $k'=0,$
we set $t=0$ and $\hat{G}^{1}_{t}=I^{1}.$ This means that we can
additionally assume that $V\subset I^0\setminus I^1$: since $V\subset
I^0$ is a nice pair, if this was not the case, we would have
$I^1\subset V,$ but then $I^1$ would be very small compared to
$I^{0}$.

\noindent\textbf{Claim:} There exists a function $\rho_{1}$, as above
such that $$\frac{|G_{k'}|}{|I^1|}\leq\rho_{1}(\varepsilon).$$

\noindent\textit{Proof of claim.}
We will assume that $\phi\colon I^{1}\rightarrow I^{0}$ is orientation
preserving (if it is not, we replace it by its second iterate). Let
$p$ be a fixed point of $\phi\colon I^{1}\rightarrow I^{0}$, let
$\hat{G}_{t_{i},\pm}$ denote the components of
$\hat{G}_{t_i}\setminus\{p\}$ where $\hat{G}_{t_i,+} $ is on the same
side of $p$ as $V$.

\noindent\textit{Case 1}. If the union of two adjacent fundamental
domains, $\hat{G}_{t_k,+}\setminus \hat{G}_{t_{k-2},+}$ is much
bigger than $V$, then we are done since the pullback of
$\hat{G}_{t_k,+}\setminus \hat{G}_{t_{k-2},+}$ under the first return map to $I^1$
has
intersection multiplicity bounded by three.

\noindent\textit{Case 2}. The condition of Case 1 does not hold, but
$|\hat{G}_{t_{k-1},+}\setminus \hat{G}_{t_{k-2},+}|$ is much bigger
than  $|\hat{G}_{t_{k-2}}|$. Then $G_{t_{i}}$ is much smaller than $I_{0}$
for all $i=k',\dots, k$. Using this statement for $i=k'+1$ and
pulling back once more we are done.

\noindent\textit{Case 3}.
$|\hat{G}_{t_{k-1},+}\setminus \hat{G}_{t_{k-2},+}|$ is much smaller
than $|\hat{G}_{t_{k-2},+}|$. By the decomposition of the return map
into maps with bounded distortion followed by a polynomial, the derivative of 
$\phi\colon \hat{G}_{t_{k-2}}\rightarrow
I^{1}$ is
bounded. Hence, since $G_{t_{k-1}}$ is very
small in $I^1
$ and $G_{t_{k-1}}\subset \hat{G}_{t_{k-1}}\setminus \hat{G}_{t_{k-3}}
,$ we have that $|G_{t_{k-1}}|$ is very small compared to
$|\hat{G}_{t_{k-1}}\setminus \hat{G}_{t_{k-3}}|$. The proof now
follows as in Case 1.
\endpfclaim

To complete the proof of the lemma, if $k'>0$, take $t=t_{k'-1}$ and
$\hat G^{1}_{t}=\mathcal{L}_{f^{t}(x)}(I^0)$, and if $k'=0$, take $t=0$ and
$\hat G^{1}_{t}=\mathcal{L}_x(I^1)$.
\end{pf}

We now prove the second part of Lemma~\ref{KSSLemma93}.

\begin{lem}[cf. \cite{vSV} Proposition 5]\label{small return domains have small return domains}
Suppose $c_{0}$ is odd. Assume that
$I\supset J$ are nice intervals such that their pullbacks are either nested or disjoint, with $((1+2\delta)I\setminus I)\cap\omega(c_{0})=\emptyset$. Let $y\in\omega(c_0)\cap I$. Then 
for any $\varepsilon>0$, there exists $\varepsilon'>0$ such that if $J\supset \mathcal{L}_{y}(I)$ and
$$\frac{|J|}{|I|}< \varepsilon',$$ then for any $x\in\omega(c_0)$, $$\frac{|\mathcal{L}_{x}(J)|}{|\mathcal{L}_{x}(I)|}<\varepsilon.$$
\end{lem}

\begin{pf}
Notice that if $I\supset J\supset\mathcal{L}_{x}(I)$ and $I$ and $J$
are nice intervals with nested or disjoint pullbacks,
then $I\supset J$ is a nice pair.

Let $s$ be the entry time of $x$ to $J$.  
Let $\hat{G}^{0}_{s}=I$ and $G_{s}=J$, and
define chains $\{\hat{G}_{i}^{0}\}_{i=0}^{s}$ and $\{G_{i}\}_{i=0}^{s}$ so
that  $G_{i}\subset\hat{G}^{0}_{i}$ are the pullbacks of
$G_{s}\subset\hat{G}^{0}_{s}$ containing $f^{i}(x).$
Let $s_{1}$ be maximal so that $\hat{G}^{0}_{s_{1}}$ contains a
critical point of $f$. If no such $s_{1}$ exists, we set $s_{1}=0.$ Then by
the Improved Koebe Principle, there exists a function
$\hat{\rho}\colon \mathbb{R}^{+}\rightarrow\mathbb{R}^{+}$ with
$\hat{\rho}(\varepsilon)\rightarrow 0$ as $\varepsilon\rightarrow 0$
such that
$$\frac{|G_{s_1}|}{|\hat{G}^{0}_{s_1}|}<\hat{\rho}(\varepsilon).$$
If $s_{1}=0$, then the proof is finished. If $s_1>0,$ let $c_1$ be the
critical point in $\hat{G}^{0}_{s_{1}}$. Note that
$\hat{G}^{0}_{s_{1}}\supset G_{s_{1}}$ is a nice pair. Now let $t_1$
and $\hat{G}^{1}_{t_{1}}$ be the time and the interval given by the
previous lemma. Then
$$\frac{|G_{t_{1}}|}{|\hat{G}^{1}_{t_{1}}|}\leq\rho_1\circ\hat{\rho}(\varepsilon).$$
If $t_1=0$, this completes the proof. Otherwise, $G_{t_{1}}\subset
\hat{G}^{1}_{t_{1}}$ is a nice pair, and moreover, any pullback of
$\hat{G}^{1}_{t_{1}}$ that intersects $c_1$ is contained in
$\mathcal{C}_{c_1}(\hat{G}^{0}_{s_1})$.

Repeating this construction inductively, we obtain a sequence of times
$t_0\geq s_1>t_1>\dots >s_{\kappa-1}>t_{\kappa-1}\geq s_{\kappa}\geq
t_{\kappa}=0$ and a sequence of nice pairs
$G_{t_i}\subset\hat{G}^{i}_{t_i} $, $i=0,1,\dots,\kappa-1$ such that
\begin{enumerate}[(1)]
\item for $i=1,\dots,\kappa$, $|G_{s_i}|/|\hat{G}^{i}_{s_i}| <
\hat{\rho}((\rho_1\circ\hat{\rho})^{i-1}(\varepsilon));$
\item for $i=0,1,\dots,\kappa$,
$|G_{t_i}|/|\hat{G}^{i}_{t_i}|<(\rho_{1}\circ\hat{\rho})^{i}(\varepsilon);$
\item if $\hat{G}^{i}_{t_{i}}, \hat{G}^{j}_{t_j}$, with $i<j$, both
contain the same critical point $c$, then $\hat{G}^{j}_{t_j}\subset
\mathcal{C}_{c}(\hat{G}^{i}_{t_i})$.
\end{enumerate}
If $\kappa$ is not large, the proposition follows immediately. If
$\kappa$ is large, we use a different argument. By the last property,
there is a critical point $c$ and a sequence $i(1)<i(2)<\dots<i(r)$,
with $r>\kappa/(b-1)$, such that
$$\hat{G}^{i(1)}_{n_{i(1)}}\supset\mathcal{C}_c(\hat{G}^{i(1)}_{n_{i(1)}})\supset
\hat{G}^{i(2)}_{n_{i(2)}}\supset\mathcal{C}_c(\hat{G}^{i(2)}_{n_{i(2)}})\supset \ldots \supset
\hat{G}^{i(r)}_{n_{i(r)}}\supset G_{n_{i(r)}}.$$
By Lemma \ref{KSSLemma91}, it follows that
$\mathcal{C}_{c}(\hat{G}^{i(j+1)}_{n_{i(j+1)}})$ is
$\delta$-well-inside $\mathcal{C}_c(\hat{G}^{i(j)}_{n_{i(j)}})$, so
that $\mathcal{L}_{x}(\mathcal{C}_c(\hat{G}^{i(j+1)}_{n_{i(j+1)})})$
is $\delta'$-well-inside
$\mathcal{L}_{x}(\mathcal{C}_c(\hat{G}^{i(j)}_{n_{i(j)}}))$.  Since
this holds for $j=1,\dots, r-1$, 
$\mathcal{L}_{x}(\hat G^{i(1)}_{n_{i(1)}})$ contains
a $(1+2\delta')^{r-2}$-scaled neighbourhood of
$\mathcal{L}_{x}(\hat G^{i(r)}_{n_{i(r)}})$. 
Since $\mathcal L_x(J)=G_{0}\subset
\mathcal{L}_{x}(\hat G^{i(r)}_{n_{i(r)}})$ 
and 
$\mathcal{L}_{x}(\hat{G}^{i(1)}_{n_{i(1)}})\subset \mathcal{L}_x(I)$.
Hence $\mathcal{L}_x(I)$ is a $(1+2\delta')^{r-2}$-scaled
neighbourhood of $\mathcal{L}_{x}(J)$. This completes the argument if
$\kappa$ and hence $r$ is large.
\end{pf}
This concludes the proof of Lemma~\ref{KSSLemma93}. 

\begin{lem}[cf. \cite{KSS} Lemma 9.6]\label{KSSLemma96}
For any $C>0$ and $\delta>0$ there exists $C'>0$ such that the
following holds.
Let $I$ be a nice interval and $J$ a pullback of $I$
with
$I\supset (1+2C')J$ and $ c_0\in J$.
Assume that either
\begin{itemize}
\item $c_0$ is a critical point of even order or
\item  $I$ is $\delta$-nice. 
\end{itemize}
Then for any $x\in\omega(c_0),$ we have that $(1+2C)\mathcal{L}_{x}(J)\subset\mathcal{L}_{x}(I).$
\end{lem}

\begin{remark} The proof of this lemma in the case when $c_0$ is a
  critical point of odd order is significantly more difficult than
  when it has even order.  The reason for this is that if $I\supset
  K\supset J\ni c_0$ are intervals so that $J$ deep-inside $I$, then
  unless $I$, $K$ and $J$ are symmetric with respect to $c_0$ it is
  quite possible that  $K$ not well-inside $I$ and also that $J$ is
  not well-inside $K$. A similar issue arises in the proof of Lemma~\ref{KSSLemma97}.
\end{remark}

\begin{pf}
Assume that $x\in \omega(c_0).$ Observe that since $f$ is
persistently recurrent on $\omega(c_0)$ and $c_0\in J$, there exists $k>0$ such that
$f^k(x)\in J$.
Let $I^{n}=\mathcal{L}^{n}_{c_0}(I)$ and define $m(0)=0$ and
$m(1)<m(2)<\dots$ as  the positive integers such that
$R_{I^{m(i)-1}}(c_0)\notin I^{m(i)}$. Let $k_{0}$ be maximal such that
$J\subset I^{m(k_{0})}$. Notice that $J$ could be equal to
$I^{m(k_0)}$.

\medskip
\noindent\textit{Claim 1. We can assume that $k_0$ is uniformly bounded.}

\medskip
\noindent\textit{Proof of Claim 1.} By Lemma \ref{KSSLemma91}, for any $i$ with $1< i\leq
k_{0}-1$ the interval $I^{m(i)}$ contains a definite neighbourhood of
$I^{m(i)+1}$. By the Theorem \ref{real Koebe} (1), we know that
$\mathcal{L}_x(I^{m(i)})$ contains a definite neighbourhood of
$\mathcal{L}_{x}(I^{m(i)+1}).$ As $\mathcal{L}_{x}(J)\subset
\mathcal{L}_{x}(\mathcal{L}_{c_0}(I^{m(k_{0}-1)+1}))$ and
$\mathcal{L}_{x}(I^{m(1)})\subset\mathcal{L}_{x}(I^0)$, the lemma
follows if $k_{0}$ is sufficiently large.
\endpfclaim

\medskip

\noindent\textit{Case 1. Assume that $c_0$ is even.} Suppose first
that $k_0=0$. If additionally, $m(1)<4$ or if $I^j\subset J$ for $j\in\{1,2,3\}$ we get that $J\supset I^4=\mathcal L^4_{c_0}(I)$, and the result follows from Lemma \ref{KSSLemma93}. So we assume $J\subset I^2$ and $m(1)\geq 4$. 
If $|I^{j}|/|I^{j+1}|>C_1/3$  is large for $j\in\{0,1\}$ and some
$C_1>1,$ then $(1+2C_1)I^{j+1}\subset I^{j}$. From Lemma
\ref{KSSLemma93}, and making $C_1$ larger if necessary, we know that $(1+2C)\mathcal L_x(I^{j+1})\subset \mathcal L_x(I^j)$. Since $\mathcal L_x(J)\subset \mathcal L_x(I^{j+1})\subset \mathcal L_x(I^j)\subset \mathcal L_x(I^0),$ the result follows. On the other hand if $|I^0|/|I^2|$ is not large,  Lemma \ref{bounds1} implies that $|I^{m(1)-1}|$ is comparable to $|I^{0}|$, so $J$ is deep-inside $I^{m(1)-1}$ and the lemma follows from Lemma \ref{KSSLemma93}.

Let us now assume that $k_0>0$. Then, there are
two possibilities: either 
\begin{itemize}
\item $|I^{m(j)}|/|I^{m(j+1)}|$ is close to one
for all $0\leq j\leq k_0-1$ or 
\item there exist $j_0$ with  $0\leq j_0\leq
k_0-1$ so that $|I^{m(j_0)}|/|I^{m(j_0+1)}|$ is large.
\end{itemize} 
In either case the result follows from the case $k_0=0$ applied to
intervals $I'$ and $J'$ with $J\subset J'\subset I'\subset I^0$ as
follows. In the first case, we consider $I'=I^{m(k_0-1)}$ and $J'=J$,
and in the second case we consider $I'=I^{m(j_0)}$ and
$J'=I^{m(j_0+1)}$. This concludes the proof when $c_0$ is even.

\medskip

\noindent\textit{Case 2. Assume $c_0$ is odd.} Then all critical points in $\omega(c_0)$ are odd. Let $m'\in \N$ be so that $I^{m'}\subset J \subset I^{m'-1}$. If $m'<6$ the result follows from Lemma \ref{KSSLemma93}, so from now on we will assume $J\subset I^5.$  Since $I$ is $\delta$-nice by Corollary \ref{nice} there exits $\delta'>0$ so that $I^1$ and $I^2$ are $\delta'$-nice.

\medskip
\noindent\textit{Claim 2. We can assume that 
that $|I^1|$ is comparable to
  $|I^3|$ and that there exists a large constant $C_2=C_2(C,\delta)>1$ such that
  $(1+2C_2)J\subset I^2$.}

\medskip
\noindent\textit{Proof of Claim 2.}
We can assume that $|I^2|$ is comparable to $|I^1|$. Indeed, if
$|I^2|$ is very small compared $|I^1|,$ then because $I^1$ is
well-inside $I^0$ we get that $I^2$ is deep-inside $I^0$ and the lemma
follows from Lemma \ref{KSSLemma93}. Similarly, we can assume that
$|I^3|$ is comparable to $|I^2|$. Therefore, from now on, we can and
will assume that $|I^3|$ is comparable to $|I^1|$. 
If $|J|$ is comparable to
$|I^1|,$ then, since $J$ is deep-inside of $I^0$, we have that $I^1$ is
deep-inside $I^0$, and the result follows from Lemma
\ref{KSSLemma93}. So we can assume $|J|$
is small compared to $|I^1|$ and indeed also to $|I^3|$.
 Furthermore, since $J\subset I^3,$
$(1+2\delta')I^3\subset I^2$ and $|J|$ is small compared to $|I^3|$ we
get that $J$ is deep-inside $I^2$. 
\endpfclaim

\medskip

Let us first study the cases when $k_0=0$ or when $k_0=1$ and
$I^{m(1)+1}\subset J\subset I^{m(1)}.$ Observe that since $m'>6$ we
must have that $m(1)>5$. Let $r$ be the return time of $c_0$ into $I^0$
and let $f^r\colon I^1\to I^0$ be the return map. Since $c_0$ is odd,
$f^r|I^1$ is a homeomorphism. 
Let $L'=(1+\delta')I^1$ and $L=\comp_{c_0}f^{-r}(L)\supset I^2$. By
the choice of $\delta'$ we know that $L\subset I^0$ and Theorem
\ref{real Koebe} (2) implies that $f^r\colon  L\to L'$ is a (finite)
composition of polynomials and  maps with bounded distortion. From
Claim 2, we can assume that $|I^1|$ is comparable to $|I^2|$, so Lemma~\ref{bounds1}(2) implies that $|Df^r|$ is bounded from above on $I^2$ and  that   $|c -f^{r}(c)|/|I_2|$ is bounded from below for any $c\in \crit (f^r|_{I^2}).$ 
Even more, if we let $g=(f^r|_{I^1})^{-1}$ we get that
$|c_0-g(c_0)|/|I_2|$  is bounded from below. The map $f^r|_{I^1}$ is
monotone so by definition of $I^{m(1)}$ the points $c_0$, $g(c_0)$ and
$f^r(c_0)$ are contained in $I^{m(1)-1}$, and $c_0$ lies
between $g(c_0)$ and $f^r(c_0)$. 
Since $|c_0-g(c_0)|/|I_2|$ and $|c_0-f^{r}(c_0)|/|I_2|$ are bounded from below,  
$\{g(c_0),c_0,f^r(c_0)\}\subset I^{m(1)-1}, c_0\in J,$ 
and $J$ is deep-inside $I^2$, 
we have that $J$ is deep-inside $I^{m(1)-1}.$  
We are assuming $J\supset I^{m(1)+1},$ so in this case the lemma follows from Lemma \ref{KSSLemma93}.

\medskip
Next, let us assume that $I^{m(2)}\subset J \subset I^{m(1)+1}$.  If
$|I^{m(1)+1}|$ is comparable to $|J|,$ then $I^{m(1)+1}$ is
deep-inside $I^2$.  
Applying the case $k_0=1$ to the intervals $I'=I^2$ and $J'=I^{m(1)+1},$ 
we get that  $(1+2C)\mathcal L_x(J')\subset \mathcal L_x(I^2)$. Since $\mathcal L_x(J)\subset \mathcal L_x(J')$ we are done in this case.
Assume that $|J|$ is small compared to $|I^{m(1)+1}|$. Let $k>0$ be minimal so that $f^k(c_0)\subset I^0\setminus I^1$ and let $K=\mathcal L_{f^{k}(c_0)}(I)$. Since $ J \subset I^{m(1)+1}$ and $k$ is minimal we have that $c_0$ visits $K$ before returning to $I^{m(1)}$.  
Since $I$ is $\delta$-nice and $K$ is a return domain to $I^0$,
we can apply Theorem \ref{real Koebe} (1) to find $\delta'>0$ so that
$K'=(1+2\delta')\mathcal L_{f^{k}(c_0)}(I^{m(1)})\subset K$. 
Let $\{G'_i\}_{i=0}^k$ be the chain given by $K=G'_k$ and $G'_0=\mathcal L_{c_0}(K)$ and $\{G_i\}_{i=0}^k$ be the chain given by $K'=G_k$ and $G_0=\mathcal L_{c_0}(K')$.
Since $k$ is the first entry time of $c_0$ to $K$ each of these chains has
order bounded by $b$.  
By Lemma \ref{KSSFact91} there exists  $\delta_2>0$ so that
$(1+2\delta_2)\mathcal L_{c_0}(K')\subset \mathcal L_{c_0}(K)\subset
I^{m(1)}$. Since $ I^{m(1)+1}\subset \mathcal
L_{c_0}(K')$ and $J$ is small compared to
$I^{m(1)+1},$ we have that $J$ is deep-inside $I^{m(1)}.$ 
Since
$I^{m(2)}\subset J \subset I^{m(1)+1}$,
taking $I'=I^{m(1)}$ and $J'=J$
the result follows from the
case $k_0=0$.

\medskip

Finally, let us now consider the case when $k_0\geq 2$. Recall, $k_{0}$ is maximal
such that $J\subset I^{m(k_{0})}$.  If $J$ is deep-inside
$I^{m(k_0)},$ the result follows from the case $k_0=0$, so assume this
is not the case. There are two possibilities, either 
\begin{itemize}
\item[(1)] $J$ is small compared to
$I^{m(k_0)}$ and $d(\partial I^{m(k_0)}, J)/ |J|$ is bounded from
above or 
\item[(2)] $|J|$ is comparable to $|I^{m(k_0)}|$. 
\end{itemize}
Assume that (1) holds. From the definition of the intervals $I^{m(i)}$ and the
fact that $k_0\geq 2,$ since $I^{m(k_0)}\subset \mathcal
L_{c_0}(I^{m(k_0-1)})$, we can apply Lemma \ref{KSSLemma91} to find
$\delta_3>0$ so that $(1+2\delta_3)I^{m(k_0)}\subset
I^{m(k_0-1)}$. This implies that $J$ is deep-inside of
$I^{m(k_0-1)}$. 
Taking $I'=I ^{m(k_0-1)}$ and $J'=J,$
the result follows either from the case $I^{m(1)+1}\subset J\subset I^{m(1)}$ or from the case $I^{m(2)}\subset J \subset I^{m(1)+1}$. 
So we can assume that (2) holds; that is, $|J|$ is comparable to $|I^{m(k_0)}|$. Then
$I^{m(k_0)}$ is small compared to $I^0$. 
Since, by Claim 1, $k_0$ is bounded from above, we must have that one
of the ratios $|I^{m(j)}|/|I^{m(j+1)}|$ is big for some
$j\in\{0,\ldots k_0-1\}$. Let $j_0$ be maximal with this property. 
If $j_0=0,$ then $|I^{m(1)}|$ is comparable to $|J|$, and $|I^{m(j)}|$
is comparable to $|I^{m(j+1)}|$, for any other $j\in{1,\ldots k_0-1}$.
This gives us that $|I^{m(1)}|$ is comparable to $|J|$, so that
$J'=I^{m(1)}$ 
must be deep-inside $I'=I^0$, and we argue as in the case when $k_0=0$.
Suppose $j_0=1.$ 
Since $I$ is $\delta$-nice, $I^{m(1)+1}$ is well-inside $I^{m(1)}$,
so $I^{m(2)}$ is deep-inside $I^{m(1)}.$ Setting $I'=I^{m(1)}$ and
$J'=I^{m(2)}.$ 
If
$j_0>1,$ then by Lemma \ref{KSSLemma91} and the definition of the
intervals $I^{m(j)}$, we have that $I^{m(j_0+1)}$ is well-inside
$I^{m(j_0)}$, so  $I^{m(j_0+1)}$ is deep-inside $I^{m(j_0-1)}.$ 
We take $I'=I^{m(j_0+1)}$ and $J'=I^{m(j_0-1)}$.
In either of these last two cases, we conclude the proof as we did when (1) holds.
\end{pf}

\begin{lem}[cf. \cite{KSS} Lemma 9.7]\label{KSSLemma97}
For each $\delta, \rho>0$ and each integer $N\geq 0$, there exist $\rho',\delta'>0$ so that the following holds. Let $I$ and $J$ be nice intervals with $I\supset  (1+2\rho)J\supset J\owns c_0$ and $J$ is a pullback of $I$ with order 
bounded by $N$. Let $\{G_{j}\}_{j=0}^{s}$ be the chain associated to the pullback and assume one of the following conditions holds:
\begin{itemize}
\item $c_0$ is a critical point of even order. 
\item $c_0$ is a critical point of odd order, $I$ is $\delta$-nice and $G_{i}\cap J=\emptyset$ for $i=1,2,\dots, s-1$.
\end{itemize}
Then $J$ is $\rho'$-nice, where $\rho'$ is defined by a function $\rho'=\rho'(\rho)>0$ depending on $N$ and $\delta$ such 
that $\rho'\rightarrow\infty$ as $\rho\rightarrow\infty$.
\end{lem}
\begin{pf}
Let $\{G_{j}\}_{j=0}^{s}$ be the chain with $G_{s}=I$ and $G_{0}=J$. Let us first assume that none of the intervals from the chain $\{G_i\}_{i=1}^{s-1}$ intersect $J$. Given $x\in J\cap \omega(c_0)$ so that $r$ is the return time of $x$ to $J$ we must have that $r\geq s$. Let $y=f^s(x)$. By the previous lemma we know that there exists $\rho_1>0$, $\rho_1\to \infty$ as $\rho\to \infty$, so that $(1+2\rho_1)\mathcal L_y(J)\subset \mathcal L_y(I)$. Lemma \ref{KSSFact91}  implies the existence of $\rho'>0$, $\rho'\to \infty$ as $\rho_1\to \infty$, so that $(1+2 \rho')\comp_x f^{-s}(\mathcal L_y(J))\subset J$, which means that  $(1+2 \rho')\mathcal L_x(J)\subset J$.  

Assume there exists an element of $\{G_i\}_{i=1}^{s-1}$ that intersects $J$. Let $s'<s$ be maximal so that $G_{s'}\cap J\neq \emptyset$. By assumption we must have that $c_0$ is a critical point of even order, so by symmetry we have that  for  $\rho_2=\rho/2$ either $(1+2\rho_2)G_{s'}\subset I$ or  $(1+2\rho_2)J\subset G_{s'}$. Assume $(1+2\rho_2)G_{s'}\subset I$. The maximality of $s'$ implies that none of the elements of the chain $\{G_i\}_{i=s'+1}^{s-1}$ intersect $J$. So we can apply the previous argument to find $\rho_3>0$ such that $G_{s'}$ is $\rho_3-$nice. This means that $(1+2\rho_3)J\subset G_{s'}$. If we are in this case or in the case that  $(1+2\rho_2)J\subset G_{s'}$ the proof follows by applying the previous argument at most $N-1$ times, since $J$ is a pullback of  $G_{s'}$ of order less than $N$.

\medskip

\noindent \textbf{Remark:} the first part of the proof also shows the following. Assume that $c_0$ is a critical point of odd order and that there exists exactly one $s'\in\{1,2,\dots, s-1\}$ such that $G_{s'}\cap J\neq\emptyset$. Then, if  $x\in J$ is such that $f^k(x)\notin J$, where $k>0$ is minimal so that $f^k(x)\in G_{s'}$, we have that $(1+2\rho')\mathcal{L}_x(J)\subset J$.

\end{pf}

Before we continue, recall that $\Gamma(I)$ is the smallest successor
of $I$; so $\Gamma(I)$ is a pullback of $I$ of order
bounded by $2b-1$. 

\begin{lem}[cf. \cite{KSS} Lemma 9.8]\label{KSSLemma98}
There exists a universal constant $\delta>0$ such that if $I$ is a nice interval containing $c_0$ and both $I$ and $\mathcal{C}_{c_0}(I)$ are non-terminating, then $(1+2\delta)\Gamma^{2}(I)\subset I$.
\end{lem}
\begin{pf} 
Let $I$ be a nice interval containing a point $c_0$. If $U$ is a nice
interval in $I$ that is disjoint from $\mathcal{L}_{c_0}(I)$, then any
pullback of $U$ that contains $c_0$ is contained in
$\mathcal{C}_{c_0}(I)$, see \cite[Lemma 15]{vSV}.  

Since $I$ is non-terminating, $\omega(c_0)$ intersects a non-central domain of the return map  to $I$
and therefore, by the above statement,  $\Gamma(I)\subset \mathcal{C}_{c_0}(I)$. Since $\mathcal{C}_{c_0}(I)$ is non-terminating, $\omega(c_0)$ intersects a non-central domain of the return map to $\mathcal{C}_{c_0}(I)$, so again we have that any pullback of this domain is contained in $\mathcal{C}^2_{c_0}(I)$.
Hence $\Gamma^2(I)\subset \mathcal{C}^2_{c_0}(I)$. By Lemma~\ref{KSSLemma91}, it follows that  $(1+2\delta)\Gamma^2(I)\subset \mathcal{C}_{c_0}(I)$. 
\end{pf}

\begin{lem}\label{lem:rhonice}
Suppose that $I\owns c_0$ is a nice interval and additionally if $c_0$
is odd, that $I$ is $\delta$-nice for some $\delta>0.$
Then, there exists a constant $\hat \delta>0,$ that is universal if $c_0$ is even, and $\hat\delta=\hat\delta(\delta)$ if $c_0$ is odd, so that 
if  $\Gamma^2(I)$  is non-terminating, then $\Gamma^3(I)$ is $\hat
\delta$-nice.
Moreover, for each $\rho>0$ there exists $\rho'>0$ with $\rho'\to \infty$ as $\rho\to \infty$, 
 so that if $(1+2\rho)\Gamma^2(I)\subset I$, then  $(1+2\rho')\Gamma^3(I)\subset \Gamma^2(I)$
and $\Gamma^3(I)$ is $\rho'$-nice. 
\end{lem}
\begin{pf} Assume $c_0$ is even. Since $\Gamma^2(I)$ is non-terminating, both $\Gamma(I)$ and $C_{c_0}(I)$ are non-terminating. So the result follows from Lemma \ref{KSSLemma98} and Lemma~\ref{KSSLemma97}. 
If $c_0$ is odd then we argue as follows. 
Let $J$ be a successor of $I$, and let  $\{J_{j}\}_{j=0}^{k}$ be the chain with $J_{0}=J$ and $J_{k}=I.$ Then, it is not hard to see that the only elements of the chain containing the point $c_0$ are $G_0$ and $G_k$. This means that is we consider the chain $\{G_{j}\}_{j=0}^{s}$  with $G_{0}=\Gamma^2(I)$ and $G_s=I$  there is only one $s'$, $0<s'<s$, with $G_{s'}\cap\Gamma^2(I)\neq\emptyset$. 
Using the remark at the end of the proof of
Lemma~\ref{KSSLemma97} it follows that all but one of the domains of the first return map to
$\Gamma^2(I)$ is $\rho$-well-inside $\Gamma^2(I)$. Since $\Gamma^2(I)$ is non-terminating,
the orbit of $\Gamma^3(I)$ passes through a non-central return domain of $\Gamma^2(I)$.
It follows that $(1+2\rho)\Gamma^3(I)\subset \Gamma^2(I)$. From 
Lemma~\ref{KSSLemma97} it follows that $\Gamma^3(I)$ is   $\rho'(\rho)$-nice. 
\end{pf}

\medskip

\begin{lem}[cf. \cite{KSS} Lemma 9.9]\label{KSSLemma99}
For any $\rho>0$ there exists $\rho'>0$, with $\rho'\rightarrow\infty$ as $\rho\rightarrow\infty$, such that if $I$ is a $\rho$-nice interval containing $c_0$, then $(1+2\rho')\mathcal{A}(I)\setminus\mathcal{A}(I)$ 
is  disjoint from $\omega(c_0)$ and contained in $\mathcal{B}(I)$.
Moreover, if $c_{0}$ is of even order, then
 $\mathcal{B}(I)\setminus (1+2\rho')^{-1}\mathcal{B}(I)$ 
is disjoint from $\omega(c_0)$. If $c_{0}$ is of odd order, then
 for each $\nu>0$ and each $\rho>0$ there exists
$\rho''>0$, so that if $|\mathcal{A}(I)|/|\mathcal{B}(I)|\ge \nu$, then 
 $\mathcal{B}(I)\setminus (1+2\rho'')^{-1}\mathcal{B}(I)$ 
is  disjoint from $\omega(c_0)$.
\end{lem}
\begin{pf} 
By definition we know that $\mathcal{B}(I)$ and $\mathcal{A}(I)$ are both pullbacks of $I$ of bounded order (where the bound depends only on $b$). Since $I$ is $\rho$-nice we can apply Lemma \ref{KSSFact91} to find $\rho'>0$ so that $(1+2\rho')\mathcal A(I)\subset \mathcal B(I)$.  By definition $\mathcal{B}(I)\setminus \mathcal{A}(I)$ is disjoint from $\omega(c_0).$ So $(1+2\rho')\mathcal{A}(I)\setminus\mathcal{A}(I)$ is disjoint from $\omega(c_0).$  If $|\mathcal A(I)|$ is comparable to $|\mathcal B(I)|$ then the above implies the existence of $\rho''>0$ so that $\mathcal{B}(I)\setminus (1+2\rho'')^{-1}\mathcal{B}(I)\subset \mathcal{B}(I)\setminus \mathcal{A}(I)$ and the result follows. So the only case left to consider is when $c_0$ is even and $|\mathcal A(I)|$ is small compared to $|\mathcal B(I)|.$ Since $c_0$ is even, this implies $\mathcal A(I)$ is deep-inside $\mathcal B(I),$ so the components of $\mathcal B(I)\setminus \mathcal A(I)$ have size comparable to $\mathcal B(I)$ and the result follows.
\end{pf}

\subsection{Geometry of periodic intervals}
In this subsection, we study the geometry of terminating intervals. Before we begin, let us refer the reader to Section \ref{sec:enhanced nest} for the relevant notation. 
Recall that given $I_i,$ a terminating
interval from generalized enhanced nest, $\mathcal{R}(I_{i-1})$ is
Feigenbaum if and only if it is a periodic interval, necessarily of
period two.

\begin{prop}\label{delta free}
For each $\nu>0$ there exists $\rho>0$ such that the following holds. Assume that $I_{i-1}$ is a terminating interval from the enhanced nest. 
\begin{enumerate}
\item If $\mathcal{R}(I_{i-1})$ is terminating and not periodic, then $\mathcal{R}(I_{i-1})$ is $\rho$-free.
\item If $\mathcal{R}(I_{i-1})$ and $\mathcal{R}^2(I_{i-1})$ are Feigenbaum, then $\mathcal{R}(I_{i-1})$ is $\rho$-free.
\item If $\mathcal{R}(I_{i-1})$ is non-terminating and $|\mathcal{R}(I_{i-1})| / |I_{i+1}|< \nu$, then $\mathcal{R}(I_{i-1})$ is $\rho$-free.
\item If $\mathcal{R}(I_{i-1})$ is Feigenbaum, $\mathcal{R}(I_{i})$ is not periodic and $|I_{i}|/ |I_{i+2}|< \nu$, then $\mathcal{R}(I_{i-1})$ is $\rho$-free.
\end{enumerate}
 \end{prop}
\begin{pf}
Let $r>0$ be minimal so that $f^r(c_0)\in I_{i-1}.$
Since $I_{i-1}$ is terminating, 
$I_{i-1}^\infty=\cap_{n\geq 0} I_{i-1}^n,$
$f^r(c_0)\in I_{i-1}^\infty,$
and
$f^r|I_{i-1}^\infty$ maps $I_{i-1}^\infty$ into itself.
Thus, we have that
$\mathcal{L}_{c_0}^n(I_{i-1}^\infty)=I_{i-1}^\infty,$ for all
$n\in\mathbb N,$ is an infinite central cascade.
Since $I_{i-1}$ is terminating, 
we know that  $I_{i}=\mathcal{L}_{c_{0}}(\mathcal{R}(I_{i-1}))$.

Apply Lemma \ref{bounds1} with $\delta=1$, fixed,
to obtain constants $C,\kappa>0$ such that:
\begin{enumerate}
\item[(a)] If $c\in \crit f^r|_{ I_{i-1}^\infty}$ we have $|f^{r}(c)-c|\geq\kappa | I_{i-1}^\infty|;$
\item[(b)] $|Df^{r}(x)|\leq C$ for all $x\in  I_{i-1}^\infty;$
\item[(c)] Given $c\in \crit f^r|_{ I_{i-1}^\infty}$ and $p=f^r(p)$ both in $ I_{i-1}^\infty,$ we have $|p-c|\geq\kappa | I_{i-1}^\infty|,$
\end{enumerate}
where the third statement follows from the first two. We will refer to
these estimates throughout the proof.

Let $V_0$ be the closure of the landing domain to
 $\mathcal{R}(I_{i-1})$ containing $\alpha$ in its boundary and
 $\hat{V}_0 =\Comp_\alpha f^{-r}(V_0) \cup V_{0}$.  


\medskip
\noindent\textit{Proof of 1.} 
Since $\mathcal{R}(I_{i-1})$ is terminating and not Feigenbaum,
$f^r(V_0)= \mathcal R(I_{i-1}),$
and $V_0$ is disjoint from $\omega(c_{0}).$ 
This implies that $\hat{V}_0$ is disjoint from $\omega(c_{0})$.
By (a) and (b) there exists $\rho>0$ such that $|V_0|> \rho |\mathcal{R}(I_{i-1})|,$  
and $|\hat{V}_0|$ 
is comparable to $|I_{i-1}|.$ 
Since $c_0$ is even, the result follows by symmetry.

\medskip
\noindent\textit{Proof of 2.}
Let us assume $\mathcal{R}(I_{i-1})$ is not $\rho$-free.
Then there exist $m\in\N$ such that $f^{2rm}(c_0)$ is close to the
boundary of $\mathcal{R}(I_{i-1})$.
Regardless of  whether $f^{2rm}(c_0)$  is close to $\alpha$ or
$\tau(\alpha)$,
we must have that  $f^{2r(m+1)}(c_0)$ is close to $\alpha$, 
which is a contradiction. To see this, observe that
if $R^2(I_{i-1})$ is Feigenbaum, 
either $f^{2rm}(c_0)$ or $f^{2r(m+1)}(c_0)$ must be contained in
$R^2(I_{i-1})$, which is well-inside $\mathcal{R}(I_{i-1})$ by (b) and
(c). 
So it cannot happen that both  $f^{2rm}(c_0)$ and $f^{2r(m+1)}(c_0)$ 
are close to the boundary of $\mathcal{R}(I_{i-1})$.

\medskip
\noindent\textit{Proof of 3.}
Assume $\mathcal{R}(I_{i-1})$ is non-terminating and 
$|\mathcal{R}(I_{i-1})| / |I_{i+1}|< \nu$.
The return map to $I_i$ restricted to $V_0$ is equal to $f^{2r}$. 
If $V_0$ does not intersect $\omega(c_0)$, we are done, 
so let us assume $V_0$ intersects $\omega(c_0)$. 
Let $V_{j+1}$ be the pullback of $V_{j}$ 
under $f^{2r}$ containing $\alpha$ for $j\in \N$. 
Observe that (b) implies that $|V_1|$ is comparable $|\mathcal{R}(I_{i-1})|,$
so if $\omega(c_0)$ does not intersect $V_1$ the proposition follows.
Let us assume $V_1$ contains some point in $\omega(c_0)$. 
We will show that $f^{2r}|_{V_1}$ is monotone. 
Let $c'$ be the turning point of $f^{2r}|_{V_0}$  closest to $\alpha$. By construction $c'$ is an absolute maximum (or minimum if $\tau(\alpha)< \alpha$). If $f^{2r}(c')$ is not in $V_0$, then $f^{2r}|_{V_1}$  is monotone. If $f^{2r}(c')\in V_0$, then $f^{2r}(V_0)\subset V_0$, which contradicts the fact that $\omega(c_0)$ is minimal. Since $\alpha$ is an orientation reversing fixed point of $f^r$, we have that $f^{2r}|_{V_1}$ is an orientation preserving monotone map with at most $b$ critical points. Define  $\hat{V}_1= \comp_{\alpha}f^{-r}(V_1)\cup V_1$ and  $\hat{V}_{j+1}= \comp_\alpha f^{-r}(\hat{V}_j)$ for $j\in \N$.
Let $c''\in\crit (f^{2r})\cap V_1$ be the critical point closest to
$\alpha$. If $c''$ is not defined, then let $c''=c'$ and
$n_0=1$. Otherwise, let $n_0$ be maximal such that $c''\in V_{n_0}.$ 
Statement (c) implies that there exists $\kappa'>0$ so that  $|V_{n_0}|> \kappa'|\mathcal{R}(I_{i-1})|$
The proposition follows if $V_{n_0}$ does not intersect $\omega(c_0)$, so we will assume that $V_{n_0}$ intersects $\omega(c_0)$.  Notice that for every $n>n_0$  the map $f^{(n_0-n)r}\colon  \hat{V}_{n}\to \hat{V}_{n_0}$ is a diffeomorphism.  
Assume $V_{n_0+2}=[\alpha, y]$. Then by Lemma \ref{fundamental domains}, we know that there exists $M>0$ so that $|y-f^{2r}(y)|>M|\mathcal R (I_{i-1})|$. Since the derivative of $f^{r}$ is bounded on $I_{i-1}^\infty$, we can find $\rho>0$ such that 
$$(1+2\rho)\hat V_{n_0+2}\subset \hat V_{n_0+1}.$$
By Theorem \ref{real Koebe} (1) we know that there exists
$\delta_2> 0$ such that for any $n\in \N$ with $n>2$
\begin{equation}\label{eqn: delta free}
(1+2\delta_2)\hat{V}_{n_0+n+1}\subset \hat{V}_{n_0+ n}.
\end{equation}

The minimality of $\omega(c_0)$ implies that there exist $m =n_0+2+m'$ maximal such that $\hat{V}_{m}\cap \omega(c_0)\neq \emptyset$. Let $k_0\in \N$ be minimal with the property that $f^{k_0r}(c_0)\in \hat{V}_{m}$.
Define $k< k_0$  maximal such that $f^{kr}(c_0)\in \mathcal{R}(I_{i-1})$. By the definition of $k,$ we have that $f^{k_0r}(c_0)$ is the first entry  of $f^{kr}(c_0)$ to  both $\hat{V}_{m}$ and $\hat{V}_1$.
For each $j\geq n_0+2$ let $V'_j$ be the pullback of $\hat{V}_{j}$
along the orbit $\{f^{k_0r}(c_0),f^{k_0r-1 }(c_0), \ldots,
f^{kr}(c_0)\}$ containing $f^{kr}(c_0)$. Notice that the intervals
$V'_j$ are nice intervals. Applying Lemma \ref{KSSFact91} we can find
$\delta_3> 0$ such that  $(1+2\delta_3)V'_{j+1}\subset V'_j.$  Using
Theorem \ref{real Koebe} (1) we can find $\delta_4>0$ such that $(1+2\delta_4)\mathcal{L}_{c_0}(V'_{j+1})\subset \mathcal{L}_{c_0}(V'_j)$. So  $$(1+2m'\delta_4)\mathcal{L}_{c_0}(V'_{m})\subset \mathcal{L}_{c_0}(V'_1).$$

From this we can conclude that given any $\nu>0$ there exists $m'\in \N$ such that if $m=n_0+2+m'$ and  $V_{m}\cap\omega(c_0)\neq \emptyset$, then $|\mathcal{R}(I_{i-1})| / |I_{i+1}|>\nu$. This is because $|I_{i-1}^\infty|$  is always comparable to $|\mathcal{R}(I_{i-1})|$ and
\begin{equation}\label{eqn:free3}
I_{i+1}\subset \Gamma^2(\mathcal{R}(I_{i-1}))\subset
\mathcal{L}_{c_0}(V_{m})\subset \mathcal{L}_{c_0}(\hat{V}_{m}) \subset
I_{i-1}^\infty.
\end{equation}

It remains to show that $\Gamma^2(\mathcal{R}(I_{i-1}))\subset
\mathcal{\mathcal L}_{c_0}(V_{m})$. To this end, let $x$ be $f^{k_0}(c_0)$  the image of $c_0$ under the first landing map of $c_0$ to $V_m$. Then $\mathcal L_{c_0} \mathcal L_x \Gamma (\mathcal R (I_{i-1}))$ contains $\Gamma^2 \mathcal R (I_{i-1})$, since $\mathcal L_{c_0} \mathcal L_x \Gamma (\mathcal R (I_{i-1}))$ is the pullback of $\Gamma (\mathcal R (I_{i-1}))$ by two landing maps, and so it has to be bigger than the last child of $\Gamma (\mathcal R (I_{i-1})).$ But now $\mathcal L_{c_0}\mathcal  L_x \Gamma (\mathcal R (I_{i-1}))$ is contained in $\mathcal L_{c_0}(V_m).$

\medskip
\noindent\textit{Proof of 4.} 
Assume $\mathcal{R}(I_{i-1})$ is Feigenbaum, 
$\mathcal{R}(I_{i})$ is not Feigenbaum and 
$|I_{i}| / |I_{i+2}|< \nu$. 
Since $\mathcal R(I_{i-1})$ is Feigenbaum, 
we know that  $R(I_{i-1}) = I_i$. 
 Since $\mathcal{R}(I_{i-1})$ is Feigenbaum, and the return time of
 $c_0$ to $I_{i-1}$ is $r$, 
we have that the return time of $c_0$ to $I_i$ is equal to $2r$. 
Let $\alpha$ be the fixed point of $f^{r}$ contained in the boundary
of $I_i$. 
By statement (c), the distance between $c_0$ and $\alpha$ is comparable to
$|I_{i-1}^\infty|$,
so $|I_i|$ and  $|I_{i-1}^\infty|$ are comparable.

Let $\alpha'$ be the orientation reversing fixed point of $f^{2r}$ contained in the boundary of $\mathcal{R}(I_i)$. Let $W_1$ be the closure of the component of  $I_i\setminus\{ f^{-2r}(\alpha')\}$  that contains $\alpha$ in its boundary. 
Let $\hat{W}_1$ be the interior of $W_1 \cup \comp_{\alpha} f^{-r}(W_1)$  and $\hat{W}_{j+1}= \comp_\alpha f^{-r}(\hat{W}_j)$ for $j\in \N$. Observe that $\hat{W}_j\subset I_{i-1}^\infty$ for all $j\in \N$. 
Since $I_i$ is Feigenbaum, (a) and (b) imply that $|W_2|$ is
comparable to $|I_i|$. 
 If $\omega(c_0)$ does not intersect $W_2$, $I_i$ is $\delta$-free, so assume $\omega(c_0)$ intersects $W_2$. Observe that  $f^{2r}|_{W_1}$ is monotone. Using the same argument as in the proof of $3$ we can find $\delta_2>0$ and $n_0\in\N$ such that given any $j>3$

$$(1+2\delta_2)\hat{W}_{n_0+j+1}\subset \hat{W}_{n_0+ j}.$$

If $I_i$ is not $\delta$-free there exists, $m'$ big such that the orbit of $c_0$ enters $\hat{W}_{m}$, with $m=n_0+2+m'$. Let $k\in \N$ be minimal so that $f^{rk}(c_0)\in \hat{W}_{m}.$ Let $k_0<k$ be minimal so that $f^{rk_0}(c_0)\in I_i$.
By Theorem \ref{real Koebe} (1) there exists $\delta_4>0$ such that
 $(1+2m'\delta_4)\mathcal{L}_{c_0}(
 \mathcal{L}_{f^{rk}(c_0)}(\hat{W}_{m}))\subset I_{i-1}^\infty .$ So
 we have that $I_{i+2}$ is deep-inside $ I_{i-1}^\infty $, 
since
\begin{equation}\label{eqn:free4}
I_{i+2}\subset \Gamma(\mathcal{R}(I_{i}))\subset \mathcal{L}_{c_0}(
\mathcal{L}_{f^{rk_0}(c_0)}(\mathcal{R}(I_{i})))\subset
\mathcal{L}_{c_0}( \mathcal{L}_{f^{rk_0}(c_0)}(\hat{W}_{m})).
\end{equation}
We know that $|I_i|$ and  $|I_{i-1}^\infty|$ are comparable, so $I_{i+2}$ being deep-inside $ I_{i-1}^\infty $  contradicts $|I_i|/|I_{i+2}|<\nu$.
\end{pf}

\begin{cor}\label{first landing}
There exists $\delta>0$ and 
for each $\rho>0$ there exists $\rho'>0$ such that the following
holds. Let $I_{i-1}$ be a terminating interval from the enhanced
nest. 

\begin{itemize}
\item Assume $\mathcal{R}(I_{i-1})$ is non-periodic.
Let $V_0$ be the component of $f^{-2r}(\mathcal{R}(I_{i-1}))$ that
contains $\alpha$ in its boundary and $W_0$ be the component of
$f^{-2r}(\mathcal{R}(I_{i-1}))$ that contains $\tau(\alpha)$ in its
boundary. 
If  $x\in \mathcal{R}(I_{i-1})$ and $\mathcal
L_x(\mathcal{R}(I_{i-1}))$ is different from $V_0$ and $W_0$, then
$\mathcal {L}_x(\mathcal{R}(I_{i-1}))$ is $\delta$-nice.
\item
If $\mathcal{R}(I_{i-1})$ is $\rho$-free, then
 $\mathcal {L}_x(\mathcal{R}(I_{i-1}))$ is
$\rho'$-free.
\end{itemize}
\end{cor}

\begin{pf}
Let $V_1$ be the component of $f^{-4r}(\mathcal{R}(I_{i-1}))$ that contains $\alpha$ in its boundary and $W_1$ be the component of $f^{-4r}(\mathcal{R}(I_{i-1}))$ that contains $\tau(\alpha)$ in its boundary.
From the proof of Proposition~\ref{delta free} we know that there
exists $\rho_1>0$ such that  $|V_1|\geq \rho_1
|\mathcal{R}(I_{i-1})|$.  Since $W_1$ is the symmetric component, with
respect to $c_0$, corresponding to $V_1$ the same bounds hold for
$W_1$. The bounds on the size of $V_1$ and $W_1$ imply that $\mathcal
L_x(\mathcal{R}(I_{i-1}))$ is $\rho_1$-well-inside
$\mathcal{R}(I_{i-1})$ if $\overline{\mathcal
  L_x(\mathcal{R}(I_{i-1}))}\cap \{\alpha,
\tau(\alpha)\}=\emptyset$. This fact, along with Theorem \ref{real
  Koebe} (1) imply the existence of $\delta>0$ such that $\mathcal
L_x(\mathcal{R}(I_{i-1}))$ is $\delta$-nice. 
If $\mathcal R(I_{i-1})$ is $\rho$-free, we can apply Corollary  \ref{cor:free-space} to get that  $\mathcal {L}_x(\mathcal{R}(I_{i-1}))$ is $\rho'$-free, for some $\rho'>0$.

\end{pf}

\medskip
Below we will make use of the sets $Y_{j}$ and $\tilde{\mathcal{Y}}_{\gamma}$ introduced on page \pageref{notation: top level}.

 \begin{lem}\label{rho nice landing domains}
There exists $\rho>0$ with the following property. Let $I_{i}$ be a terminating interval of the enhanced nest and $Y_j$ be a component of  $I_{i}^\infty\setminus (f^r|I^\infty_{i})^{-1}(\alpha)$. Given a point $x\in Y_j$  the following holds.

\begin{enumerate}
\item If $Y_j$ does not intersect  $\tilde{ \mathcal{ Y}}_\alpha$ and $-a<j<a$, then $Y_j$ is $\rho$-nice.
\item Assume $j>-a$ and the first return time of $x$ to $Y_j$ equal to $f^{kr}$. If there exist $0<j_0\leq k $ such that  $f^{rj_0}(x)\notin \tilde {\mathcal{Y}}_\alpha$, then $\mathcal L_x (Y_j)$ is  $\rho$-nice.
\item If $j=-a$ and the return time of $x$ to $Y_j$ is bigger than $f^r,$ then $\mathcal L_x (Y_j)$ is  $\rho$-nice.
\end{enumerate}
 \end{lem}
 \begin{pf}
Let $\tilde Y_0= Y_0\cap \tilde {\mathcal{Y}}_\alpha$ and $\tilde Y_{-1}=Y_{-1}\cap \tilde {\mathcal{Y}_\alpha}.$
By Lemma \ref{bounds1}, we can find $\kappa>0$
such than $|\tilde Y_0|, |\tilde Y_{-1}|> \kappa |I_{i}^\infty|$.
It follows immediately that there exists 
$\rho'>0$ such that the following holds.  Given $Y_j$ with  $j>-a$  and $J=\comp
f^{-r}(Y_j)\subset Y_\ell$ for some $\ell\in \{-a,\dots,a\}$, 
if $J\cap\tilde{ \mathcal{Y}}_\alpha= \emptyset$
then $(1+2\rho') J\subset Y_\ell$.  
This implies (1).
To prove (2), take $x\in Y_j$  with first return time to $Y_j$ 
equal to $f^{kr}$ and such that 
$f^{rj_0}(x)\notin \tilde{ \mathcal{Y}}_\alpha$ for some 
$0<j_0\leq k$. Then,  $(1+2\rho' )\mathcal L_x (Y_j)\subset Y_j$. 
By Theorem \ref{real Koebe},  we know that here exists $\rho>0$ such that $\mathcal L_x (Y_j)$ is $\rho$-nice. One can prove (3) in a similar way.
 \end{pf}

\subsection{Proof of Theorem \ref{real geometry}.}

We will prove Theorem \ref{real geometry} in two separate stages. First we will show parts (a) to (e).

\medskip
\noindent\textbf{Part (a).} 
If $I_{n-1}$ is non-terminating, $I_n=\Gamma^2(\Gamma^{T-2}(\mathcal B(\mathcal A(I_{n-1}))))$, and by Lemma \ref{lem:rhonice} we get that $I_{n}$ is $\rho_1$-nice for some $\rho_1>0$ which depends only on $\underline b.$ If $I_{n-1}$ is terminating the result follows from Corollary \ref{first landing} and the observation that $I_n$ being non-terminating implies  that $I_n=\mathcal L_x(\mathcal{R}(I_{n-1}))$ is different from $V_0$ and $W_0.$ In conclusion, if $I_n$ is non-terminating we have that $I_n$ is $\rho_1$-nice.

Now assume $I_{n-1}$ is non-terminating. Let $\nu_1\in \N$ be so that $\mathcal B\mathcal A(I_{n-1})= \comp_{c_0}f^{-\nu_1}(\mathcal A(I_{n-1})),$ and $\mathcal A^2(I_{n-1})= \comp_{c_0}f^{-\nu_1}(\mathcal L_{f^{\nu_1}(c_0)}\mathcal A(I_{n-1})).$ Since $\mathcal B(I_{n-1})\setminus \mathcal A(I_{n-1})$ is free from $\omega(c_0),$ we get that $\mathcal B^2(I_{n-1})=\comp_{c_0}f^{-\nu_1}(\mathcal B(I_{n-1})).$ 

By construction, we know that $\mathcal B^2(I_{n-1})\setminus \mathcal
B\mathcal A(I_{n-1})$ is free from $\omega(c_0)$ and that the
poscritical set of $ \mathcal B\mathcal A(I_{n-1})$ is contained in
$\mathcal A^2(I_{n-1}).$ 
Moreover, the pairs $\mathcal B^2(I_{n-1})\supset\mathcal{BA}(I_{n-1})$ and
$\mathcal{BA}(I_{n-1})\supset \mathcal{A}^2(I_{n-1})$ are both
nice pairs.
By the above, we have that $I_{n-1}$ is $\rho_1$-nice, Lemma  \ref{KSSLemma99} implies the existence of $\rho_2>0$ such that $(1+2\rho_2)\mathcal A(I_{n-1})\subset \mathcal B(I_{n-1}).$ Since $\mathcal B^2(I_{n-1})$ is a pullback of bounded order (depending only on $b$) of $\mathcal B(I_{n-1})$, Lemma \ref{KSSFact91} implies the existence of $\rho_3>0$ so that $(1+2\rho_3)\mathcal B\mathcal A(I_{n-1})\subset \mathcal B^2(I_{n-1}).$ By an analogous argument we get that $(1+2\rho_3)\mathcal A^2(I_{n-1})\subset \mathcal B\mathcal A(I_{n-1}).$ So $\mathcal B\mathcal A(I_{n-1})$ is internally and externally free.
Since $I_n$ is a pullback of bounded order of $\mathcal B\mathcal
A(I_{n-1}),$ the result follows from Lemma~\ref{KSSLemma99}. Observe
that the external free space is given by $J_n=\Gamma^T(\mathcal
B^2(I_{n-1}))$, the internal free space is given by
$J_n'=\Gamma^T(\mathcal A^2(I_{n-1}))$, and both 
$J_n\supset I_n$ and $I_n\supset J_n'$ are nice pairs.

Finally, if $c_0$ is even the components of $I_n\setminus J_n'$ are comparable to $|I_n|$. If $c_0$ is odd, then by the 2nd part of Lemma  \ref{KSSLemma99} there exists $\rho_4(\nu)>0$ so that if $|I_{n-1}|/|I_n|<\nu$, $|J_n'|\geq \rho_4|I_n|$. 

\medskip
\noindent
\noindent\textbf{Part (b).} Since $I_n$ is non-terminating, we know
that $\mathcal R(I_{n-1})$ is non-terminating. By
Proposition~\ref{delta free} (3) there exists $\rho>0$ so that $\mathcal R(I_{n-1})$ is $\rho$-free. The result follows from the definition of $I_n$ and Corollary  \ref{cor:free-space}.

\medskip
\noindent
\noindent\textbf{Part (c).} By Lemma \ref{KSSLemma93}, $(1+2C')\mathcal{L}_{c_{0}}(\mathcal{L}_{x}(I_{n}))\subset \mathcal{L}_{c_{0}}(I_{n}).$ Since $\Gamma(I_{n})\subset \mathcal{L}_{c_{0}}(\mathcal{L}_{x}(I_{n}))$, the result follows.

\medskip
\noindent 
\noindent\textbf{Part (d).} The first part follows directly from
Corollary \ref{nice}. To show the second part, observe that if $I_{n}$
is $C$-nice, then by Lemma \ref{KSSFact91} there exists $C''>0,$ $C''\to \infty$ as $C\to\infty,$ so that  $(1+2C'')\mathcal{A}(I)\subset \mathcal{B}(I)$. So the result follows from Corollary \ref{cor:free-space}.

\medskip
\noindent
\noindent\textbf{Part (e).}   If $|I_{n+1}|<\varepsilon' / 2|I_n|$ we are done, so let us assume this is not the case. Let us first assume that $I_{n-1}$ is terminating. For each $\varepsilon''>0$ there exists $C>0$ so that if $|I_n|<\varepsilon'' |\mathcal R(I_{n-1})|$, then $(1+2C)I_n\subset \mathcal R(I_{n-1})$.  By Lemma \ref{KSSLemma93} we know that $I_n$ is $C'$-nice, for some $C'(C)>0$ ($C'\to\infty$ as $C\to\infty$). The result follows if $\varepsilon''$ is sufficiently small, so let us assume this is not the case. So there exists $\nu>0$ so that $|R(I_{n-1})|/|I_{n}|, |I_n|/|I_{n+1}| < \nu$. By Proposition~\ref{delta free}  and the definition of $I_n,$ there exists $\rho>0$ such that $I_{n}$ is $\rho$-free. Then since $|\mathcal{L}_{x}(I_{n})|<\varepsilon |I_{n}|,$ there exists $C''>0$ ($C''\rightarrow\infty$ as $\varepsilon\rightarrow 0$) such that $(1+2C'')\mathcal{L}_{x}(I_{n})\subset I_{n}$ and the statement follows by Part (c).  A similar argument works when $I_{n-1}$ is non-terminating and $c_0$ is even. Assume $I_{n-1}$ is non-terminating and $c_0$ is odd, then by Part (a) $I_n$ is $\rho$-nice and $\rho$-externally free. The second part of Lemma \ref{KSSLemma93} implies that taking $\varepsilon$ sufficiently small $|\mathcal L_{c_0}(\mathcal L_x(I_n))|\leq\epsilon' |\mathcal L_{c_0}(I_n)|$, and the result follows.

\medskip
In order to prove the remaining part of the theorem we will need  two extra lemmas, and the following definition.
We say that an interval $J\subset I$ is $\varepsilon$-\emph{small in} $I$ if $|J|\le \varepsilon |I|$.

\begin{lem}\label{lem:small return domains}
Given $\varepsilon\in(0,1)$, $ \rho>0$ and $N\in \N$ exists $\varepsilon'>0$ so that the following holds. Assume that $c_0$ is odd and that $I$ is a $\rho$-nice interval with  $((1+2\rho)I\setminus I)\cap\omega(c_0)=\emptyset.$ Let $J\subset I$ be a pullback of $I$ of order $N$ with $c_0\in J.$ If
$$\frac{|J|}{|I|}<\varepsilon',$$
then for each $x\in J\cap\omega(c_0)$ we have that  $\mathcal L_x(J)$ is $\varepsilon$-small compared to $J.$ 
\end{lem}
\begin{pf}
Let us begin by showing the following claim.

\medskip
\noindent\textit{Claim.
For each $\nu>0$ there exists $\nu'>0$ such that if $J$ is $\nu'$-small in $I$ then for any $y\in\omega(c_0),$ $\mathcal{L}_y(J)$ is $\nu$-small in $\mathcal{L}_y(I).$}

\medskip
\noindent\textit{Proof of claim.}
Let $y\in\omega(c_0). $ There are two cases.  Let us first assume
$J\supset \LL_{z}(I)$ for some $z\in\omega(c_0)$.  Then, the claim follows directly form Lemma \ref{small return domains have small return domains}.   So assume there exists $z\in\omega(c_0)$ so that $J\subset  \LL_{z}(I).$ If $J$ is comparable to $\LL_{z}(I)$, then  $\LL_{z}(I)$ is small in $I$ and the claim follows from Lemma \ref{small return domains have small return domains} applied to $\LL_{z}(I)$ and $I$, since $\LL_y(J)\subset \LL_y(\LL_z(I))$, and the landing domains to $\LL_z(I)$ are small compared to the landing domains to $I$.  Now, if $\LL_{z}(I)$ is much larger than $J$, we get that $J$ is deep inside of $I,$ and the claim follows from Lemma \ref{KSSLemma96}.\endpfclaim

\medskip

Let $\{G_j\}_{j=0}^s$ be the chain with $G_s=I$ and $G_0=J$. Let $x\in J\cap\omega(c_0)$ be fixed, and let $k$ be its first return time to $J$. Let us first assume  $k\geq s$ and let $y=f^{s}(x).$ By the claim, $\LL_y(J)$ is $\varepsilon''$- small compared to $\LL_y(I),$ for some $\varepsilon''>0$ with  $\varepsilon'' \to 0$ as   $\varepsilon' \to 0.$ Since $I$ is $\rho-$nice this implies $\LL_y(J)$ is deep-inside $I,$ and the result follows from Lemma \ref{KSSFact91}.

Assume that $k<s.$ Observe that by definition of $k$, $G_k$ contains
$J$. Consider the times $k\leq s_i<s$ with $s_0=k$ and $s_{i}<s_{i+1}$
such that $J\subset G_{s_i}.$ Recall there are at most $n<N$ such
times. Our aim is to show that $J$ is small in $G_k$, so the proof follows form the argument above. We will do this by showing that $J$ is small compared to each $G_{s_i}$, for $i=0,\ldots n$. If $J$ is small compared to $G_{s_n}$ we are done for $n$, so let us assume $G_{s_n}$ is comparable to $J$. By the argument above we know that for any $z\in \omega(c_0)\cap G_{s_n}$ we have that $\LL_z(G_{s_n})$ is $\varepsilon_n$-small compared to $G_{s_n}$, with $\varepsilon_n\to 0$ as $\varepsilon'\to 0$. In particular, $G_{s_{n-1}}$ is small compared to $G_{s_n}$, so $J$ is small compared to $G_{s_n}$, contradiction. So $J$ is small compared to $G_{s_n}.$ Applying this argument $n-k$ times we show $J$ is $\varepsilon_k$-small in $G_k,$ with $\varepsilon_k\to 0$ as $\varepsilon'\to 0$, and the result follows.
\end{pf}

\begin{lem}\label{big scaling factors real bounds}
Suppose that $I_{n-1}$ and $I_n$ are non-terminating. Then
for any $C>0$, there exists $\varepsilon>0$ such that if
$$\frac{|I_{n}|}{|I_{n-1}|}<\varepsilon,$$  then 
 \begin{enumerate}[(1)]
 \item $\Gamma(I_{n}^{1})$ is $C$-nice;
 \item $I_{n+1}$ is $C$-strongly nice and $C$-externally free.
 \end{enumerate}

\end{lem}
\begin{pf}

 First we will show that \textit{(1)} holds.
 If $c_0$ is a critical point of even order, then for any $C'>0$, there exists $\varepsilon>0$ such that 
 if $|I_{n}|/|I_{n-1}|<\varepsilon$, then $(1+2C')I_{n}\subset I_{n-1}$.
 By Lemma \ref{KSSLemma97}, if $C'$ is chosen large enough, then $I_{n}$ will be 
 $C''$-nice. Finally, since $\Gamma(I_n^1)$ is a pullback of $I_n$ of bounded order (only depending on $\underline b$) we can apply  Lemma \ref{KSSLemma97} to prove (1). So from now on we will assume that $c_0$ is of odd order. 
 Since $c_0$ is odd all intervals from the enhanced nest are
 non-terminating. In particular, by Theorem~\ref{real geometry}~(a) we
 have that there exists $\rho>0$ such that that $I_{n-1}$ and $I_n$ are both $\rho$-nice and $\rho$-externally free.

 In particular, $$I_{n}\subset(1+2\rho)I_{n}\subset\Gamma^{T}\mathcal{B}^2(I_{n-1})=:F_{n}.$$ By the proof Part a) of Theorem \ref{real geometry},  we know that $F_{n}\setminus I_{n}$ is disjoint from $\omega(c_{0}),$ and $I_{n}\subset F_{n}$ is a nice pair. 
 By Lemma~\ref{lem:small return domains}, every component of the return
 map to $I_{n}$ is $\varepsilon'$-small in $I_{n}$, so we have that every component of the domain of the return map to $I_{n}$ that intersects $\omega(c_{0})$, in particular $I_{n}^{1}$, is $C'$-well-inside $F_{n}$ with $C'\rightarrow\infty$ as $\varepsilon'\rightarrow 0$.

 To see that $\Gamma(I_{n}^{1})$ is $C$-nice, let $x\in I_{n}^{1}\cap\omega(c_0)$. Note that by 
 Lemma \ref{small return domains have small return domains} we have that the
 return domains to $I_{n}^{1}$ are small compared to $I_{n}^{1}$. The first return map from $I_{n}^{1}$ to $I_{n}$ extends to a map from a domain $\tilde{J}\supset I_{n}^{1}$ to $F_{n}$. 
 Since $F_{n}\setminus I_{n}$ is disjoint from $\omega(c_{0})$, $(\tilde{J}\setminus I_{n}^{1})\cap\omega(c_0)=\emptyset$, the first return map from $\mathcal{L}_{x}(I_{n}^{1})$ to $I_{n}^{1}$ extends to a map from $J'\supset\mathcal{L}_x(I_n^1)$ to $\tilde J$, and there exists a constant $C_1$, $C_1\rightarrow\infty$ as $\varepsilon\rightarrow 0$ such that $(1+2C_1)\mathcal{L}_{x}(I_{n}^{1})\subset J'$. Since $I_{n}^{1}$ is non-terminating, we have that $J'\subset I_{n}^{1}$, for all but one component: the component of the domain of the return map to $I_n^1$ that has the same return time to $I_n^1$ as $I_n^1$ has to $I_n$. For this one component, $J'=\tilde{J}$. As in Lemma \ref{lem:rhonice}, there exists $C$, $C\rightarrow\infty$ as $C_1\rightarrow\infty$, such that $\Gamma(I_n^1)$ is $C$-nice.

Now we show that (2) holds. Observe that
the map $f^{\nu}\colon \mathcal{B}(I_{n})\rightarrow I_{n}$ extends to a domain $\tilde{F_{n}}$, so that $f^{\nu}\colon \tilde{F}_{n}\rightarrow F_{n}$ has the same critical points as $f^{v}|\mathcal{B}(I_{n}).$

 By Lemma \ref{lem:small return domains}, we know that the landing domains to $I_n$ are small compared to $I_n,$ which implies that  $\mathcal{A}(I_{n})$ is $C'$-well-inside $\tilde{F}_n$, for some $C'>0$ with $C'\rightarrow\infty$ as $\varepsilon\rightarrow 0.$
In addition,  $\tilde{F}_{n}\setminus \mathcal{A}(I_{n})$ is disjoint from $\omega(c_{0})$. If we let $H_{n+1}:=\Gamma^{T}\mathcal{B}(\tilde{F}_{n}),$ then $H_{n+1}\setminus I_{n+1}$ is disjoint from $\omega(c_{0})$. By Lemma \ref{KSSFact91} there exists $C''>0$ so that $(1+2 C'')I_{n+1}\subset H_{n+1},$ where $C''\to \infty$ as $ C'\to \infty.$ So $I_{n+1}$ is $C''$-externally free. Observe that the return time of any $x\in \omega(c_0)\cap I_{n+1}$ to $I_{n+1}$ is the same as its return time to $H_{n+1}$. 
Thus, there exists a constant $C>0$, $C\rightarrow\infty$ as $C''\rightarrow\infty$, such that for each $x\in\omega(c_0)\cap H_{n+1}$ we have that $$\mathcal{L}_{x}(I_{n+1})\subset(1+2C)\mathcal{L}_{x}(I_{n+1})\subset\mathcal{L}_{x}(H_{n+1})\subset I_{n+1}.$$ Moreover, each component of the domain of the return map to $H_n$ that intersects $\omega(c_{0})$ contains a unique component of the domain of the return map to $I_{n+1}$ that intersects $\omega(c_{0})$. Hence $I_{n+1}$ is $C$-strongly nice, where $C\rightarrow\infty$ as $\varepsilon\rightarrow 0$.
\end{pf}

\bigskip

Now we can prove the remaining part of Theorem \ref{real geometry}.

\medskip
\noindent\textbf{Part (f).} 
The fact that $I_{n+2}$ is $C'$-strongly nice and $C'$-externally free follows from Lemma \ref{big scaling factors real bounds}.

\medskip
\noindent
\noindent\textbf{Part (g).} By Lemma~\ref{KSSLemma96}, there exists $C''>0$ ($C''\rightarrow\infty$ as $C\rightarrow\infty$) such that if $(1+2C)I_{n+1}\subset I_{n}$, then $(1+2C'')I_{n+1}^1\subset I_{n+1}$. The result follows from this and the fact that $\omega(c_0)\cap I_{n+1}\subset I_{n+1}^{\infty}$.

\subsection{Existence of free space}
As a corollary to Theorem \ref{real geometry} and  Proposition~\ref{delta free}, we have the following:

\begin{cor}\label{cor:delta free}
For each $\nu>0$ there exists $\rho>0$ such that if $I_{i}$ is an interval from the enhanced nest, then
\begin{enumerate}
\item if $I_{i-1}$ is terminating and $|I_{i-1}|/|I_{i+2}|<\nu$, then $I_{i}$ is $\rho$-free;
\item if $I_{i-1}$ and $I_{i}$ are non-terminating, and $|I_{i-1}|/|I_{i}|<\nu$, then $I_{i}$ is $\rho$-free;
\item if $I_{i-1}$ is non-terminating and $I_{i}$ is terminating, then $I_{i}$ is $\rho$-free.
\end{enumerate}
\end{cor}
\begin{pf}
The first part follows from Proposition~\ref{delta free} and Lemma \ref{KSSFact91}, and the second and third from Theorem \ref{real geometry}. 
\end{pf}

\section{Poincar\'e disks  and tools for dealing with analytic maps}\label{sec:Poincare disks}
To construct box mappings we will make use of Poincar\'e
disks. In this section, we study the basic properties of Poincar\'e
disks and their pullbacks under analytic maps, so throughout this
section we consider maps $f\in\mathcal A_{\underline b}.$ An analytic
map $f$ on the interval extends uniquely to a holomorphic mapping on a neighbourhood of the
interval, and we will also use $f$ to denote this extension. We recall that given an interval $I,$ the Poincar\'e disk with angle $\theta$ based on $I$ is denoted by $D_{\theta}(I).$  For a definition, see Subsection \ref{subsec:notation}, page~\pageref{defn:poincare disk}.

\begin{lem}[Almost Schwarz Inclusion \cite{dFdM}]\label{Almost Schwarz Inclusion}
There exist $K<\infty$, $a_{0}>0$  and a function $\theta:(0,a_{0})\rightarrow(0,\infty)$ satisfying $\theta(a)\rightarrow 0$ and $a/\theta(a)\rightarrow 0$ as $a\rightarrow 0$ such that the following holds. Let $F\colon \mathbb{D}\rightarrow\mathbb{C}$ be univalent and real-symmetric, and assume that $I\subset\mathbb{R}$ is an interval containing 0 with $|I|<a\in(0,a_{0})$. Let $I'=F(I)$. Then
\begin{enumerate}[(a)]
\item for all $\theta\geq\theta(|I|)$, we have $$F(D_{\theta}(I))\subset D_{(1-K|I|^{1+\delta})\theta}(I'),$$ where $0< \delta<1$ is a universal constant;
\item for all $\theta\in(\pi/2,\pi)$ we have $$F(D_{\pi-\theta}(I))\subset D_{\pi-K|I|\theta}(I').$$
\end{enumerate}
\end{lem}

	\begin{lem}[cf. \cite{Shen} Lemma 7.4]\label{lem42}
		Let $\ell\geq 2$ be an integer and consider $P(z)=z^{\ell}.$ Then, for each  $\theta\in(0,\pi)$ the following holds.
		\begin{itemize}
		\item Suppose $\ell$ is even and let $K\geq 1$.
		Then there exists $\lambda=\lambda(K,\ell)\in(0,1)$ such that
				$$P^{-1}(D_{\theta}(-K,1))\subset D_{\lambda\theta}(-1,1)).$$
			
			\item Suppose that $\ell$ is an odd integer. Let $K>0$. Then
			there exists $\lambda=\lambda(\ell)\in(0,1)$ such that 
				$$P^{-1}(D_{\theta}(-K^{\ell},1))\subset D_{\lambda\theta}((-K,1)).$$
		\end{itemize}
	\end{lem}
	\begin{pf}
We will give the proof in the case that $\ell$ is odd; the proof in the even case is similar.
Since the problem is invariant under scaling, we can assume that $K\geq 1$.
We will argue by contradiction. 
If the lemma is false, then for every $\lambda>0$,
there exist $\theta_{1/\lambda}\in(0,\pi)$
and $z_{1/\lambda}^{\ell}\in D_{\theta}(-K^{\ell},1)$
such that $z_{1/\lambda}\notin D_{\lambda\theta}(-K,1)$.
Note that we may as well assume that $\theta\in(0,\pi/2)$, so that is what we will do. 
For each integer $n\geq 1$, let $\lambda_{n}=1/n$
and let $\theta_{n}$ and $z_{n}$ be the corresponding $\theta_{1/\lambda}$ and $z_{1/\lambda}$. 

Since $z_{n}\notin D_{\theta_{n}/n}(-K,1)$, 
$$\arg\frac{z_{n}-1}{z_{n}+K}\leq\frac{\theta_{n}}{n},$$
so that 
$$\tan\arg\frac{z_{n}-1}{z_{n}+K}\leq\tan\frac{\theta_{n}}{n}.$$
Writing $z_{n}=r_{n}e^{it_{n}}$, we have
$$\arg\frac{z_{n}-1}{z_{n}+K}=\arg(z_{n}-1)-\arg(z_{n}+K)=\arctan\Big(\frac{r_{n}\sin t_{n}}{r_{n}\cos t_{n}-1}\Big)-\arctan\Big(\frac{r_{n}\sin t_{n}}{r_{n}\cos t_{n} +K}\Big).$$
So, taking the tangent gives us
$$
\frac{r_{n}\sin t_{n}(K+1)}{r_{n}^{2}+r_{n}(-1+K)\cos t_{n}-K} = 
\frac{\frac{r_{n}\sin t_{n}}{r_{n}\cos t_{n}-1}-\frac{r_{n}\sin t_{n}}{r_{n}\cos t_{n} +K}}{1+\Big(\frac{r_{n}\sin t_{n}}{r\cos t_{n}-1}\Big)\Big(\frac{r_{n}\sin t_{n}}{r_{n}\cos t_{n}+K}\Big)}
=
\tan\arg\frac{z_{n}-1}{z_{n}+K}
\leq\tan\frac{\theta_{n}}{n}.$$
 Similarly, letting $w_{n}=z_{n}^{\ell}=r_{n}^{\ell}e^{i\ell t_{n}}$, we have
$$\frac{r_{n}^{\ell}\sin \ell t_{n}(K^{\ell}+1)}{r_{n}^{2\ell}+r_{n}^{\ell}(-1+K^{\ell})\cos \ell t_{n}-K^{\ell}}=
\tan\arg\frac{w_{n}-1}{w_{n}+K^{\ell}}\geq\tan(\theta_{n}).$$
Dividing the second inequality by the first one, we have that
\begin{equation}\label{eqn:fraction}
\frac{r_{n}^{\ell}\sin \ell t_{n}(K^{\ell}+1)(r_{n}^{2}+r_{n}(-1+K)\cos t_{n}-K)}{(r_{n}^{2\ell}+r_{n}^{\ell}(-1+K^{\ell})\cos \ell t_{n} -K^{\ell})r_{n}\sin t_{n}(K+1)}\geq\frac{\tan\theta_{n}}{\tan(\theta_{n}/n)}
\end{equation}
Observe this expression becomes very large as $n\rightarrow\infty.$ 
Since $\sin \ell t_n/\sin t_n$ is bounded from above, and cancelling an $r_n$
in the numerator and denominator, this is equivalent to 
\begin{equation}\label{eqn:zd}
\frac{r_{n}^{\ell-1}(K^{\ell}+1)(r_{n}^{2}+r_{n}(-1+K)\cos t_{n}-K)}{(r_{n}^{2\ell}+r_{n}^{\ell}(-1+K^{\ell})\cos \ell t_{n} -K^{\ell})(K+1)}
\asymp 
\frac{r_{n}^{\ell+1} +r_n^\ell(K-1)\cos t_n-r_n^{\ell-1}K }{r_{n}^{2\ell}+r_{n}^{\ell}(-1+K^{\ell})\cos \ell t_{n} -K^{\ell}}
\end{equation}
becoming very large as $n$ tends to infinity. We will show that this is impossible.

 Let $R_n=r_n/K$. Equation \ref{eqn:zd} is comparable to

$$\frac{r_{n}^{\ell+1}K^{\ell-1}+r^{\ell}_{n} K^{\ell} \cos t_{n} -r_n^\ell K^{\ell-1}\cos t_n -r_n^{\ell-1}K^{\ell}}{r_{n}^{2\ell}+r_{n}^{\ell}K^{\ell}\cos \ell t_{n}  -r_n^\ell\cos \ell t_n  -K^{\ell}},$$

which equals
$$\frac{\frac{r_n^{\ell+1}}{K^{\ell+1}} +\frac{r_n^\ell}{K^\ell}\cos t_n-\frac{r_n^\ell\cos t_n}{K^{\ell+1}}-\frac{r_n^{\ell-1}}{K^\ell}}{\frac{r_n^{2\ell}}{K^{2\ell}}+\frac{r_n^\ell}{K^\ell}\cos\ell t_n-\frac{r_n^\ell\cos \ell t_n}{K^{2\ell}}-\frac{1}{K^\ell}}
=\frac{R_n^{\ell+1}+R_n^\ell\cos t_n -\frac{R_n^\ell}{K}\cos t_n  -\frac{R_n^{\ell-1}}{K}}{R_{n}^{2\ell}+R_n^\ell\cos\ell t_n-\frac{R_n^\ell}{K^\ell}\cos\ell t_n-\frac{1}{K^\ell}}$$
in the scaled coordinates. 
Since this expression becomes very large as $n\rightarrow\infty$, by comparing powers of $R_n$,
we see that $R_n$ is bounded.
Since $z_{n}\notin D_{\theta_{n}/n}(-K,1),$ we have that $t_n\rightarrow 0$ or $t_n\rightarrow\pi$.
First consider the case that $t_n\rightarrow 0$.
Then as $n\to \infty$,
$$\arg\Big(\frac{z_{n}-1}{z_{n}+K}\Big) \asymp \arg (z_{n}-1)
\mbox{ and }
\arg\Big(\frac{z^{\ell}_{n}-1}{z^{\ell}_{n}+K^\ell}\Big) \asymp \arg (z^{\ell}_{n}-1).$$
However, 
since 
$z_{n}\notin D_{\theta_{n}/n}(-K,1)$
and
$z_{n}^\ell \in D_{\theta_{n}}(-K^\ell,1)$,
combining these estimates with (\ref{eqn:fraction}) we have that
$$\frac{\arg(z_{n}^\ell-1)}{\arg(z_{n}-1)}\rightarrow\infty.$$
We will show that this is impossible.
Suppose that $r_n$ does not converge to 1.
Then, since $t_n\rightarrow 0$, neither $r_n^\ell\cos(\ell t_n) -1$ nor $r_n\cos(t_n) -1$
converges to 0.
But now we see that
$$\frac{\arg(z_{n}^\ell-1)}{\arg(z_{n}-1)}=\frac{\arctan\frac {r_n^\ell\sin(\ell t_n)}{r_n^\ell\cos(\ell t_n) -1}}{\arctan\frac{r_n\sin(t_n)}{r_n\cos(t_n) -1}}\asymp \frac{\frac {r_n^\ell\sin(\ell t_n)}{r_n^\ell\cos(\ell t_n) -1}}{\frac{r_n\sin(t_n)}{r_n\cos(t_n) -1}}\asymp \frac{r_n^\ell\sin(\ell t_n)}{r_n\sin(t_n)}$$
is bounded. So we can assume that $r_n\rightarrow 1.$

Then we have that 
\begin{equation}\label{eqn:argument}
\frac{\arg(z_{n}^\ell-1)}{\arg(z_{n}-1)}
=\frac{\arg((z_n-1)(z_{n}^{\ell-1}+z_{n}^{\ell-2}+\dots+1))}{\arg(z_n-1)}
\end{equation}
$$=\frac{\arg((z_n-1)+\arg(z_{n}^{\ell-1}+z_{n}^{\ell-1}+\dots+1)}{\arg(z_n-1)}
=1+\frac{\arg(z_{n}^{\ell-1}+z_{n}^{\ell-1}+\dots+1)}{\arg(z_n-1)},
$$
which is bounded,
since there exists a constant $C>0$ such that
$$\arg(z_{n}^{\ell-1}+z_{n}^{\ell-1}+\dots+1)\leq C\arg(z_n-1)$$
as $z_n\rightarrow 1$.

Now, suppose that $t_n\rightarrow \pi$.
Then as $n\to \infty$,
$$\arg\Big(\frac{z_{n}-1}{z_{n}+K}\Big) \asymp \pi- \arg (z_{n}+K),\mbox{ and }
\arg\Big(\frac{z^{\ell}_{n}-1}{z^{\ell}_{n}+K^\ell}\Big)\asymp \pi- \arg (z^{\ell}_{n}+K^\ell).$$
But now we have that as $n\rightarrow\infty$
$$\frac{\pi-\arg(z_{n}^\ell+K^\ell)}{\pi-\arg(z_{n}+K)}\rightarrow\infty,$$
but arguing just as we did in the case when $t_n\rightarrow 0$ we see that this is impossible. Suppose that $r_n$ does not converge to $-K$.
Then, since $t_n\rightarrow \pi$, neither $K^\ell- r_n^\ell\cos(\ell t_n)$ nor $K-r_n\cos(t_n)$
converges to 0.
But now we see that
$$\frac{\pi- \arg(z_{n}^\ell+K^\ell)}{\pi -\arg(z_{n}+K)}=
\frac{\arctan\frac {r_n^\ell\sin(\ell t_n)}{K^\ell - r_n^\ell\cos(\ell t_n)}}{\arctan\frac{r_n\sin(t_n)}{K-r_n \cos(t_n)}}\asymp \frac{\frac {r_n^\ell\sin(\ell t_n)}{K^\ell - r_n^\ell\cos(\ell t_n)}}{\frac{r_n\sin(t_n)}{K-r_n\cos(t_n)}}\asymp \frac{r_n^\ell\sin(\ell t_n)}{r_n\sin(t_n)}$$
is bounded. So we can assume that $r_n\rightarrow -K.$
Now, a similar calculation to the one in (\ref{eqn:argument}) shows that
$$\frac{\pi-\arg(z_{n}^\ell+K^\ell)}{\pi-\arg(z_{n}+K)}$$ must be bounded,
which again yields a contradiction.


The proof when $\ell$ is even is similar and is a natural generalization of 
\cite[Lemma 7.4]{Shen}, see
\cite[Lemma 13.2]{KSS}.
\end{pf}

The next lemma is straight forward.
\begin{lem}[cf. \cite{KSS} Lemma 13.3]\label{lem:z^l lower bound}
Let $\ell\geq 2$ be an integer and consider $P(z)=z^{\ell}.$ For any $A>0$ and any $\theta\in(0,\pi)$, there exists $\theta'\in(0,\pi)$ such that the following holds.
\begin{itemize}
\item If $\ell$ is even, then $$P^{-1}(D_{\theta}((-A,1)))\supset D_{\theta'}((-1,1)).$$
\item If $\ell$ is odd, then there exists $\theta'$ such that 
$$P^{-1}(D_{\theta}((-A^\ell,1)))\supset D_{\theta'}((-A,1)).$$
\end{itemize}
\end{lem}

Although we state the next lemma for any integer $\ell\ge 2$, it will only be
used to deal with inflection points.

\begin{lem}\label{different branch angle control}
Let $\ell\geq 2$ be an integer and consider $P(z)=z^{\ell}.$  For each $C<\infty$, there exists $\lambda\in(0,1),$ depending on  
$\ell$ and $C,$ so that for all $\theta\in (0,\pi/2)$ the following holds. 
Let $I_{1}\subset [0,1)$ and $I_{2}\subset (-1,0]$ be intervals and assume that 
$$|I_{1}|\geq 1/C,$$
and that $w\in D_{\theta}(I_{2})$. Let $z \in P^{-1}(w)$ be so that $|\mathrm{arg}(z)|<\pi/\ell$. Then
$$z\in D_{\lambda\theta}(I_{1})\cup A$$ 
where $A=\emptyset$ if $0\in I_1$ and otherwise $A$ is the component bounded by the lines
$\mathrm{arg}(z)=\pm \pi/\ell$ and the boundary of $D_{\lambda\theta}(I_{1})$. 
\end{lem}
\begin{pf}
Fix $\theta<\pi/4$.  
%
%
%
Let $I_1=(\alpha,\beta)$ and take $z\in \partial D_{\lambda\theta}(\alpha,\beta)$. 
Setting $z=re^{it}$ and arguing as in the proof of Lemma \ref{lem42}, we have that
$$\tan\arg\frac{z-\beta}{z-\alpha}=\frac{r(\beta-\alpha)\sin t}{r^{2}-r(\beta+\alpha)\cos t+\beta\alpha}=\tan\lambda\theta.$$
Hence,
$$r^{2}+(-(\beta+\alpha)\cos t-\frac{(\beta-\alpha)}{\tan \lambda\theta}\sin t)r +\alpha\beta=0.$$
This equation has two solutions,

$$r_{\pm}=\frac{1}{2}\Bigg(\Big((\alpha+\beta)\cos t+\frac{(\beta-\alpha)\sin t}{\tan \lambda\theta}\Big)\pm\sqrt{\Big((\alpha+\beta)\cos t+\frac{(\beta-\alpha)\sin t}{\tan \lambda\theta}\Big)^{2}-4\alpha\beta}\Bigg).$$
Given $t\in (0,\pi/\ell)$ we denote by $z_{\pm}=r_{\pm}e^{it}$. To prove the lemma, it suffices to show that we can take $\lambda>0$ small so that for all $t\in (0,\pi/\ell)$
we have $z_{+}^\ell\notin D_\theta(-1,0)$. 
This is equivalent to showing that
\begin{equation}
\tan\arg\frac{z_{+}^{\ell}}{z_{+}^{\ell}+1}\leq\tan\theta, \mbox{ i.e. }
r_{+}^{\ell}\ge \frac{\sin(\ell t)}{\tan\theta} - \cos(\ell t).\label{eq41}
\end{equation}
From the equation for $r_+$ we have
$$r_{+}>\frac{1}{2}(\beta+\alpha)\cos t+\frac{\beta-\alpha}{2\tan\lambda \theta}\sin t.$$
Hence 
\begin{equation}
r_{+}^{\ell}\geq  \frac{1}{2^{\ell-1}}(\beta+\alpha)^{\ell-1}\cos^{\ell-1}t\frac{\beta-\alpha}{2\tan\lambda\theta}\sin t
+  \frac{(\beta-\alpha)^\ell}{2^\ell\tan^\ell \lambda\theta}\sin^\ell t . \label{eq42}
\end{equation}
Provided we take $\lambda>0$ sufficiently small (here $\lambda$
depends on  $\beta-\alpha$ and $\ell$), for each $t\in (0,\pi/\ell)$
the right hand side of this expression dominates the r.h.s. of (\ref{eq41}). 
Indeed, if $\cos(t)$ is not small (which holds automatically if $\ell>2$, because $t\in (0,\pi/\ell)$), then 
we  the first and second term in the r.h.s. of (\ref{eq42}) dominate
respectively  the first and second term in the r.h.s. of (\ref{eq41}) (provided we choose $\lambda>0$
small). If $\cos(t)$ is small, then $\sin(t)$ is not small, and the 2nd term in (\ref{eq42}) dominates
the r.h.s. of  (\ref{eq41}),  provided we choose $\lambda>0$ small. 
\end{pf}

\begin{figure}[htb] \centering \def\svgwidth{200pt}
\scriptsize

 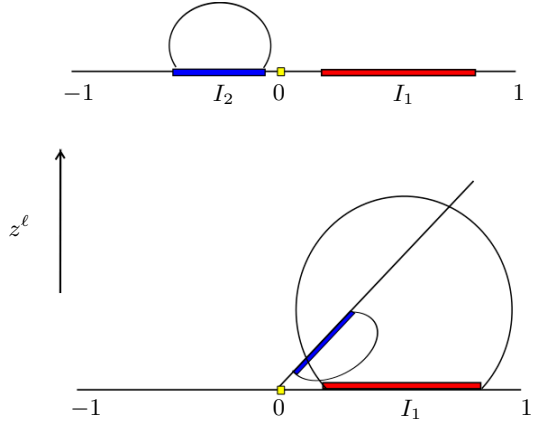
 \label{fig:uno2}
 \caption{The pullback of Poincar\'e disks as in Lemma \ref{different branch angle control}.}
 \end{figure}

\begin{lem}\cite[Lemma 13.4]{KSS}\label{Poincare disks}
One can compare Poincar\'e disks in the following ways.
\begin{enumerate}[(a)]
\item There exists $\theta_{0}$ such that for each $A>1$ and each $\theta\in(0,\pi)$,
$$D_{\theta}([-1,1])\subset D_{\mathrm{min}\{\theta A/2,\theta_{0}\}}([-A,A]).$$
\item For each $\lambda\in (0,1)$ and $\delta>0$ there exists $\lambda'\in (0,1)$ such that for each $\theta\in(0,\pi)$,
$$D_{\lambda\theta}([-1,1])\setminus D_{\theta}(-1-\delta,1+\delta)\subset D_{\lambda'\theta}([-1-\delta,-1]).$$
\item For each $\lambda\in (0,1)$ and $\delta>0$ there exists $\lambda'\in (0,1)$ such that for each $\theta\in(0,\pi)$
$$D_{\lambda\theta}([-1,1])\setminus D_{\theta}([-1-\delta,1+\delta])\subset D_{\lambda'\theta}([1,1+\delta]).$$
\item For each $\lambda\in (0,1)$ and $\delta>0$ there exists $\lambda'>0$ such that for each interval $J\subset[-1,1]$ and each $\theta>0$,
$$D_{\lambda\theta}([-1,1])\setminus D_{\theta}(-1-\delta,1+\delta])\subset D_{\lambda'\theta|J|}(J).$$
\end{enumerate}
\end{lem}

\subsection{Pulling back Poincar\'e disks under disjoint orbits}\label{subsec:pullbackdisjoint}
In this subsection we will show how to pull back Poincar\'e disks
along a chain of intervals under a map $f\in \mathcal A_{\underline
  b}.$ We remind the reader that we always assume that the
Standing Assumptions on page \pageref{standing assumptions} hold.
 

\begin{lem}[cf. \cite{KSS} Lemma 13.5]\label{chain pullback}
For any $\delta>0$ there exists $\lambda>0$ and for each $\theta\in(0,\pi)$ there exists $\varepsilon>0$  so that  following holds. 
Let $\{H_{j}\}_{j=0}^{s}$ and $\{H'_{j}\}_{j=0}^{s}$ be two chains with $H_{j}\subset H'_{j}$ for all $j=0,\dots, s$ and $H_0\cap\omega(c_0)\neq\emptyset$. Assume that the following hold:
\begin{enumerate}[(a)]
\item $|H'_{s}|\leq\varepsilon$,
\item $H'_{s}\supset(1+2\delta)H_{s}$,
\item the chain $\{H'_{j}\}_{j=0}^{s-1}$ is disjoint,
\item $H_{i}\setminus (1+2\delta)^{-1}H_{i}$, $i=0,\dots,s-1$ does not contain a critical point.
\end{enumerate}
Let $U_s=D_{\theta}(H_{s}),$ and $U_0=\comp_{H_{0}}(f^{-s}(U_s))$.  Then
 $$U\subset D_{\lambda\theta}(H_{0}).$$
\end{lem}

\begin{pf}
The Lebesgue measure of $\omega(c_0)$ is zero (see \cite[Theorem E]{vSV}), so if $H_s'$ is sufficiently small the set $\cup_{i=0}^s H_{i}'$ is a small neighbourhood of a subset of $\omega(c_0)$, and we can assume that $\sum_{i=0}^s |H_i'|$ is small.

First suppose that only $H_0'$ contains a critical point. By the Improved Koebe Principle, Theorem~\ref{real Koebe}~(2), we can find
$\rho_1,\delta_1>0$ such that $H'_{j}\supset(1+2\rho_1)H_{j}$ and
$H_{j}\setminus (1+2\delta_1)^{-1}H_{j}$ contains no critical points of $f$ for $1\leq j\leq s$. By Lemma \ref{Almost Schwarz Inclusion} there exists $\lambda_1>0$ such that 
$$U_1\subset D_{\lambda_1\theta}(H_1).$$
Since $H_{1}\setminus (1+2\delta_1)^{-1}H_{1}$ contains no critical values of $f$ we can apply Lemma \ref{lem42} to find $\lambda_2>0$ such that
$$U_{0}\subset D_{\lambda_2\theta}(H_{0}).$$
Since the critical points of $f$ are non-flat, there exists $\delta_2>0$ such that $H_0'\supset (1+2\delta_2)H_0$, and
$(H_0\setminus (1+2\delta_2)^{-1}H_0)\cap\crit(f)=\emptyset.$
Thus the lemma follows by applying the previous argument at most $b$ times.
\end{pf}

This lemma implies the following corollaries.

\begin{cor}\label{free space well inside}  Suppose that instead of (d) in the above lemma there exists $\delta'\in (0,\delta)$
so that the sets $(1+2\delta'){H}_{j}\setminus {H}_{j}$ are free from critical points for all $j\in\{0,\ldots s\}$. Then
there exist an interval $H'$, $\lambda'\in(0,1)$ and $\delta''>0$ such that
$$\hat{U}=\mathrm{Comp}_{{H}_{0}}f^{-s}(D_{\theta}({H}_{s}))\subset D_{\lambda'\theta}(H')$$
where ${H}_{0}\subset H'\subset (1+2\delta'')H'\subset H'_{0}$.
\end{cor}

%
\begin{cor}\label{chain space}
For each $\delta>0$ there exist $\delta'>0$ and $\lambda\in(0,1),$ 
and for each $\theta\in(0,\pi)$ there exists $\varepsilon'>0$ 
such that the following holds for all  $\varepsilon\in(0,\varepsilon')$. Let $I$ be a nice interval with $|I|<\varepsilon$ and let be $J$ a domain of the first entry map to $I$ with first entry time equal to $s>0.$ Let $\{H_{j}\}_{j=0}^{s}$ be a chain with $(1+2\delta)H_{s}\subset I$, 
$H_{0}\subset J$ and $H_0\cap\omega(c_0)\neq\emptyset$. Then there exists an interval $H'$ with $H_{0}\subset H'\subset (1+2\delta')H'\subset J$ so that
$$\mathrm{Comp}_{H_{0}}f^{-s}(D_{\theta}(H_{s}))\subset D_{\lambda\theta}(H').$$
\end{cor}

\section{Tools for dealing with $C^3$ maps}\label{sec: smooth extension}
In this section, we will  develop the additional tools required for dealing  with  $C^3$ maps. A reader who is primarily interested in the real analytic case can skip this section. 


\subsection{Quasiconformal and quasiregular mappings}
We will find it convenient to make use of the analytic definition of
quasiconformal and quasiregular mappings, since this approach facilitates the
definition of quasiregular mappings. We refer the reader to
\cite{AIM} for additional background on these classes of mappings.

We define
$$\partial_{\alpha}f(z)=\cos(\alpha)f_x(z)+\sin(\alpha)f_y(z),\quad\alpha\in[0,2\pi).$$
Suppose that $U$ and $V$ are domains in $\mathbb C$.
A mapping $f\colon U\rightarrow V$ is called $\kappa$-\emph{quasiregular},
abbreviated $\kappa$-\emph{qr}, if it is orientation-preserving, $f$ is in the
Sobolev space
$W^{1,2}_{loc}(U)$, and the directional derivatives satisfy
$$\max_\alpha|\partial_\alpha f(z)|\leq
\kappa\min_{\alpha}|\partial_{\alpha} f(z)|$$
for almost every $z\in\Omega$.
If, in addition $f$ is a homeomorphism, then $f$ is called
$\kappa$-\emph{quasiconformal}, abbreviated $\kappa$-\emph{qc}.

We call the constant $\kappa$ the \emph{quasiconformal distortion of} $f$.

This definition of  quasiconformal mappings is equivalent to the following geometric definition. A mapping $f\colon U\rightarrow V$ is $\kappa$-quasiconformal if and only if 
for any annulus $A\subset U$,
$$\frac{1}{\kappa}\mod(A)\leq\mod(f(A))\leq \kappa\mod(A).$$

We say that a mapping $f\colon U\rightarrow V$ is \emph{quasiregular
  (quasiconformal)}
if it is $\kappa$-quasiregular ($\kappa$-quasiconformal) for some $\kappa$.

We assume that $f\in W^{1,2}_{loc}(U)$
since it implies that $\mathrm{Jac}(f)=|f_z|^2-|f_{\bar z}|^2$ is locally-integrable.
For the reader who is familiar with qc mappings, but not qr mappings,
let us remark that for homeomorphisms $f\in W^{1,1}_{loc}(U)$, implies
$f\in W^{1,2}_{loc}(U)$; however, for quasiregular mappings it is not
sufficient to assume only that $f\in W^{1,1}_{loc}(U)$.



Recall that
$$f_z=\partial f=\frac{\partial f}{\partial
  z}=\frac{1}{2}(\frac{\partial f}{\partial
  x}-i\frac{\partial f}{\partial y}),\quad\mbox{and}\quad
f_{\bar z}=\bar\partial f=\frac{\partial f}{\partial\bar{z}} =
\frac{1}{2}(\frac{\partial f}{\partial x}+i\frac{\partial f}{\partial y}).$$

\begin{thm}
Suppose that $f\colon U\rightarrow V$ is a homeomorphism in
$W^{1,2}_{loc}(U)$ then $f$ is $\kappa$-quasiconformal if and only if
$$\frac{\partial f}{\partial \bar z}(z)=\mu_f(z)\frac{\partial
  f}{\partial z}(z)\quad \mbox{for almost every }z\in U,$$
where $\mu_f$ is a bounded measurable function satisfying
$$\|\mu_f\|_{\infty}\leq\frac{\kappa-1}{\kappa+1}<1.$$
\end{thm}
The function $\mu_f$ is called the \emph{Beltrami coefficient} of $f$.

If $f$ is a $\kappa$-quasiconformal mapping with Beltrami coefficient
$\mu_f$, let $\kappa(f)$ be the minimal $\kappa$ so that
the quasiconformal distortion of $f$ is bounded by $\kappa$. We have that
$$\|\mu_f\|_{\infty}=\frac{\kappa(f)-1}{\kappa(f)+1}\quad\mbox{and}\quad
\kappa(f)=\frac{1+\|\mu_f\|_{\infty}}{1-\|\mu_f\|_{\infty}}.$$
We will drop the subscript on the Beltrami differential when it is
clear to what mapping we are referring.

By Weyl's Lemma a 1-quasiconformal mapping is conformal and a
1-quasiregular mapping is holomorphic (see \cite{AIM}, Lemma A.6.10).
The following are some basic properties of quasiconformal mappings.

\begin{prop}
Let $f\colon U\rightarrow V$ be a $\kappa$-qc mapping onto $V$ and let
$g\colon V\rightarrow\mathbb C$ be a $\kappa'$-qc mapping. Then
\begin{itemize}
\item $f$ is differentiable a.e. on $U$.
\item $f^{-1}\colon V\rightarrow U$ is $\kappa$-qc.
\item $g\circ f\colon U\rightarrow\mathbb C$ is $\kappa'\kappa$-qc.
\end{itemize}
\end{prop}

Quasiregular mappings can be factored into a qc mapping followed by a
holomorphic mapping.
\begin{prop}[Stoilow Factorization]
Suppose that $f$ is quasiregular and defined on a simply connected
domain $\Omega$. Then $f=h\circ g$. where $g\colon \Omega\rightarrow\Omega$
is quasiconformal and $h$ is holomorphic on $\Omega$.
\end{prop}

It follows that if $f\colon U\rightarrow\mathbb C$
is a quasiregular mapping, then $f$ is open and discrete. 
Moreover, if $f\colon U\rightarrow\mathbb C$ is quasiregular, then there
exists $k\in[0,1)$ such that $|f_{\bar z}|\leq k|f_z|$ almost
everywhere in $U$, and we have 
$$\frac{\partial f}{\partial \bar z}=\mu(z)\frac{\partial f}{\partial
  z},\quad \mu_f:=\frac{f_{\bar z}}{f_z}\ \mbox{with}\
\|\mu_f\|_{\infty}\leq k.$$ 


\subsection{Asymptotically holomorphic extensions of a $C^3$ mapping}\label{sub-sec:asympholoS}
Let $\mathcal K\neq\emptyset$ be a compact subset of $\mathbb{R}^{2}$, $U$ an open neighbourhood of $\mathcal K$ and $f\colon U\rightarrow\mathbb{C}$ a $C^{1}$ map. We say that $f$ is $asymptotically\ holomorphic$ of order $t$, $t\geq 1,$ on $\mathcal K\subset\mathbb{R}^{2}$ if for every $(x,y)\in \mathcal K$ 
$$\frac{\partial}{\partial\bar{z}}f(x,y)=0,$$
and 
$$ \frac{\frac{\partial}{\partial\bar{z}}f(x,y)}{d((x,y),\mathcal K)^{t-1}}\rightarrow 0$$
uniformly as $(x,y)\rightarrow \mathcal K$ for $(x,y)\in U\setminus \mathcal K$. 

In our applications of asymptotically holomorphic extensions, $\mathcal K$ will be an interval contained in the real line.

\begin{lem}\cite[Lemma 2.1]{GSS2}\label{lem:GSS2 2.1}
Suppose that $f\colon[0,1]\rightarrow[0,1]$
is a $C^3$ mapping.
Let
$0=u_0<v_0=u_1<v_1=u_2<\dots=u_n<v_n=1$
be a partition of the interval  such that
on each interval $(u_i,v_i)$,
$f$ is either a diffeomorphism , $h_i$, with $h_i'(u_i),h_i'(v_i)\neq 0$, or
$f$ can be expressed in the form 
$f(x)=f(\zeta_i)(1-h_i(x)^\ell)$ for
$\zeta_i\in(u_i,v_i)\cap\crit(f)$ and $h_i$
is a diffeomorphism such that $h_i'(u_i),h_i'(v_i)\neq 0$, and $h_i(\zeta_i)=0$.
In either case, the diffeomorphism $h_i$
has a $C^3$-extension to a diffeomorphism $H_i$
defined on a neighbourhood $U$ of the interval $[u_i,v_i]$
in the complex plane and $H_i$ is asymptotically holomorphic of order
3 on $[u_i,v_i]$ in $U$. Moreover, the neighbourhood $U$ and the
extension $H_i$ are symmetric with respect to the real axis.
\end{lem}

We will also use $f$ to denote the asymptotically holomorphic
extension of order three of a map $f$, and whenever we refer to an
asymptotically holomorphic extension we mean the asymptotically
holomorphic extension of order three given by Lemma~\ref{lem:GSS2
  2.1}.
We refer the reader to \cite{GSS2} for
background on asymptotically holomorphic extensions of real maps.
Let us remark that in general $\kappa$-quasiregular mappings are only
differentiable almost everywhere; however, in our setting,
they are $C^3$.

\medskip
An equivalent definition of $\kappa$-quasiconformal is the following:
$H\colon U\rightarrow V$ is $\kappa$-quasiconformal iff $H\in W^{1,2}_{loc}(U)$
and 
\begin{equation}\label{def:jac}
|DH(z)|^2<\kappa J(H,z),
\end{equation}
where $DH=H_{z}+H_{\bar z}$ and $J(H,z)=|H_{z}|^2-|H_{\bar z}|^2$.
Direct calculation shows that if $H$ is an asymptotically holomorphic
extension given by Lemma~\ref{lem:GSS2 2.1} at a point $x\in \mathbb
R$, then $J(H,x)=|h'(x)|^2.$ Since each $h$ is a diffeomorphism, we
can assume that this quantity is bounded away from 0 by $\mu$. Calculating, 
$$H_{\bar z}(x+iy)=\frac{1}{2}\Big(\frac{\partial H}{\partial
  x}(x+iy)+i\frac{\partial H}{\partial y}(x+iy)\Big)<C|y|^2.$$
This implies that 
$$|H_z(x+iy)|=\frac{1}{2}\Big|\frac{\partial H}{\partial
  x}(x+iy)-i\frac{\partial H}{\partial y}(x+iy)\Big|\asymp |h'(x)|.$$

Thus, if $|y|<\eta,$ $H$ is
$1+\eta^2$-quasiconformal, close to the real line we have that $J(H,x+iy)$ is
close to both $|H_{z}(x+iy)|^2$ and $|h(x)|^2$.
Thus
$$1+\eta^2 > \frac{|DH(z)|^2}{J(H,z)}\asymp\frac{|H_z(z)|^2}{|H_z(z)|^2}+2\frac{|H_z(z)||H_{\bar
    z}(z)|}{|H_z(z)|^2}+\frac{|H_{\bar z}(z)|^2}{|H_{z}(z)|^2}=(1+\mu_H(z))^2.$$

It immediately follows that
$$\mu_H(z)=O(|\mathrm{Im}(z)|^{2}).$$

 \begin{prop}[\cite{GSS2} Proposition 2]\label{GSSProp2}
Let $f\colon I\rightarrow\mathbb{R}$ be a $C^{3}$ diffeomorphism from a compact interval $I$ with non-empty interior into the real line.
There exists $\zeta>0$ and $\delta>0$ such that if $J$ is an interval contained in $I$, $0<\alpha<\pi$ and $\diam(D_{\alpha}(J))<\delta,$ then
$$f(D_{\alpha}(J))\subset D_{\tilde{\alpha}}(f(J)),$$
where $\tilde{\alpha}=\alpha-\zeta|J|\diam D_{\alpha}(J)$, and $\tilde{\alpha}<\pi$.
\end{prop}

This proposition gives an easy generalization and improvement of the Almost Schwarz
Inclusion Principle, Lemma~\ref{Almost Schwarz Inclusion}. 

\begin{cor}[Almost Schwarz Inclusion in the asymptotically holomorphic case]\label{cor:GSSProp2}
For every small $a>0$, there exists $\alpha(a)>0$ satisfying $\alpha(a)\rightarrow 0$ and $a/\alpha(a)\rightarrow 0$ as $\alpha\rightarrow 0,$ such that the following holds. Let $f\colon I\rightarrow\mathbb{R}$ be a $C^{3}$ diffeomorphism from a compact interval $I\supset\{ 0,a\}$ into the real line, with $f(0)=0$ and $f(a)=a$. Let $f$ be a $C^{3}$ extension of $f$ to a complex neighbourhood of $I$, with $f$ asymptotically holomorphic of order 3 on $I$. Then there exists $\zeta>0$ and $\delta>0$ such that if $\alpha(a)<\alpha<\pi$ and $\diam(D_{\alpha}([0,a]))<\delta,$ then
$$f(D_{\alpha}([0,a]))\subset D_{(1-\zeta a^{1+\kappa})\alpha}([0,a]),$$
where $\kappa\in(0,1)$ can be taken arbitrarily close to $1$.
\end{cor}
\begin{pf}
Fix $\alpha_0\in(0,\pi/4)$.
Suppose first that $\alpha<\alpha_{0}<\pi/2.$ Proposition~\ref{GSSProp2}, we have that $$f(D_{\alpha}([0,a]))\subset D_{\tilde{a}}([0,a]),$$
where $\tilde{\alpha}=\alpha-\zeta a \diam(D_{\alpha}([0,a]))=\alpha-\zeta a \frac{a}{\sin \alpha}.$
Then since $\alpha$ is small, we have that $\sin\alpha\asymp \alpha$, so that
$$\tilde{\alpha}=\alpha\Big(1-\zeta\frac{a^{2}}{\alpha^{2}}\Big).$$

Suppose that $\delta\in(0,1)$, and let $\alpha(a)=a^{1/2-\delta/2}$. Notice that as $a\rightarrow 0$, $\alpha(a)\rightarrow 0$ and $a/\alpha(a) = a^{1/2+\delta/2}\rightarrow 0$. Now for any $\alpha>\alpha(a)$, we have that 
$$\tilde{\alpha}\asymp
\alpha\Big(1-\zeta\frac{a^{2}}{\alpha^{2}}\Big)\geq\alpha\Big(1-\zeta\frac{a^{2}}{a^{1-\delta}}\Big)=\alpha(1-\zeta
a^{1+\delta}).$$
Hence
$$f(D_{\alpha}([0,a]))\subset D_{(1-\zeta a^{1+\delta})\alpha}([0,a]).$$
Since this holds for all $\delta\in(0,1)$, we are done.

On the other hand, suppose that $\alpha>\alpha_{0},$ then
$$\tilde{\alpha}=\alpha-\zeta a \diam(D_{\alpha}([0,a]))\asymp \alpha -
\zeta a^{2}=\alpha\Big(1-\zeta\frac{a^{2}}{\alpha}\Big).$$
Take any $\delta\in(0,1)$ and let $\alpha(a)=a^{\delta}$. Note that as $a\rightarrow 0$ both $\alpha(a)$ and $a/\alpha(a)\rightarrow 0$. Moreover, if $\alpha>\alpha(a)$, then we have
$$\tilde{\alpha}\asymp \alpha\Big(1-\zeta \frac{a^{2}}{\alpha}\Big)\geq\alpha(1-\zeta a^{2-\delta}).$$
Hence 
$$f(D_{\alpha}([0,a]))\subset D_{(1-\zeta a^{2-\delta})\alpha}([0,a]).$$
\end{pf}

\begin{remark}
The previous corollary implies that Lemma~\ref{chain
  pullback} and its corollaries 
in Subsection~\ref{subsec:pullbackdisjoint} also hold in the
asymptotically holomorphic case if we consider maps $C^3$ provided
that the Standing Assumptions on page \pageref{standing assumptions} hold.
We will use them in this generality without further comment.
\end{remark}
%

\medskip
The same argument used to prove Proposition~\ref{GSSProp2}
can be used to prove: 
 \begin{cor}\label{cor:diffeo lower bounds}
Let $f\colon I\rightarrow\mathbb{R}$ be a $C^{3}$ diffeomorphism from a compact interval $I$ with non-empty interior into the real line.
There exists $\zeta>0$ and $\delta>0$ such that if $J$ is an interval contained in $I$, $0<\alpha<\pi$ and $\diam(D_{\alpha}(J))<\delta,$ then
$$f(D_{\alpha}(J))\supset D_{\tilde{\alpha}}(f(J)),$$
where $\tilde{\alpha}=\alpha+\zeta|J|\diam D_{\alpha}(J)$.
\end{cor}

\subsection{$\kappa$-qr box mappings}\label{sub-sec:kqcboxS}
Let $U$ and $V$ be open Jordan disks in $\mathbb{C}$. We say that a map $F\colon U\rightarrow V$ is a \emph{$\kappa$-qc branched covering}
if it can be decomposed as $F=G\circ H$ where $H\colon U\rightarrow U$ is $\kappa$-qc homeomorphism and $G\colon U\rightarrow V$
is a holomorphic branched covering.

A mapping $F\colon \mathcal{U}\rightarrow \mathcal{V}$ is \emph{qc quasi-box mapping} if the following holds. The domain $\mathcal U$, is a possibly countable, union of open Jordan disks
$U_i$, $\mathcal V$ is a union of finitely many open Jordan disks
$V_j$ with $j=0,\dots, n-1$ and:
\begin{itemize}
\item $V_0,\dots, V_{n-1}$ are pairwise disjoint;
\item every connected component $V_j\in \mathcal{V}$ is either a connected component of $\mathcal{U}$ or the intersection of $V_j$ with $\mathcal{U}$ is a union of $U_i$'s, each of these contained in $V_j$ (not necessarily compactly);
\item if $U_i\subset V_j$, then $\overline{U}_{i}\setminus V_j\subset \mathbb{R}$.
\item for each $i$, there exists $j$
such that $F|U_i\colon U_i\rightarrow V_j$ is a qc branched
covering.
\end{itemize}
  A \emph{$\kappa$-qr box mapping} is defined
analogously: we modify the definition of a complex box mapping
on page~\pageref{def:box mapping} to only require that on each
component $U$ of $\mathcal U$ that $F|U$ is a $\kappa$-quasiregular
covering map onto a component of $\mathcal V$. These mappings are
$\kappa$-qc branched coverings because of the Stoilow Factorization
Theorem, see for example Corollary 5.5.2 in \cite{AIM}.



\subsection{Additional tools for asymptotically holomorphic maps}
To prove complex bounds 
for asymptotically holomorphic maps, we will need to control
their quasiconformal distortion. First, we have the following useful estimate:

\begin{lem}[Theorem B \cite{Li-Shen}]\label{lem:squared lengths}
		Let $f$ be $C^3$ with all periodic orbits hyperbolic repelling. Then for any $\alpha$ there exists $\eta=\eta(\alpha)$ such that
		for any interval $T$ and any
		$s\in\mathbb{N}$, if $f^s\colon  T\rightarrow f^s(T)$
		is a diffeomorphism, then
		$$\sum_{i=0}^{s}|f^i(T)|^{1+\alpha}<\eta.$$
\end{lem}

In the next lemma, we prove that when we pull back a Poincar\'e disk
$D_{\theta}(J_s)$ by a diffeomorphism, $f^s\colon J_0\rightarrow J_s$,
the total loss of angle is small; so we can estimate the
diameters of the pullbacks by the lengths of their real traces. Since
the extension of $f$ is asymptotically holomprphic of order three, this
will give a bound on the quasiconformal distortion of the extension of $f^s$ in terms of $\sum_{j=0}^{s-1}|J_j|^2$.

Given an interval $I,$ let $$\mu(I)=\max_{J}|J|,$$ where the
maximum is taken over all pullbacks $J$ of $I.$
Let $\mu'(I)=\max_J|J|,$ where the maximum is taken over all
diffeomorphic pullbacks $J$ of $I$.
By \cite[Lemma 5.2]{Kozlovski-Schwarzian}, if $I\cap\omega(c_0)\neq\emptyset,$ there exists a constant
$\tau_1\geq 1$ such that $\mu'(I)\leq \tau_1|I|$.
Consequently, if
$I\cap\omega(c_0)\neq\emptyset$, then $\mu(I)$ tends to zero as $|I|$ tends to zero. 

\medskip

If $f$ is persistently recurrent on $c_0$ and the Standing Assumptions
on page \pageref{standing assumptions} hold we obtain the following two results.

\begin{lem}\label{lem:angle control along diffeos}
For any $\theta\in(0,\pi/2),$
there exist $\varepsilon'>0$ and $\tilde\theta\in(0,\pi/2)$ such that
for any $\varepsilon\in(0,\varepsilon')$
the following holds.
Let $f$ be an asymptotically holomorphic extension of a $C^3$ mapping $f$.
Suppose that $|J_s|<\varepsilon$,  $J_0\cap\omega(c_0)\neq\emptyset$ and $f^s\colon J_0\rightarrow J_s$ is a diffeomorphism. Let
$\{J_{j}\}_{j=0}^s$ be the chain such that
$J_{j}=\comp_{f^{j}(J_0)}f^{-(s-j)}(J_s)$.
Let $U_s=D_{\theta}(J_s)$, and set
$$U_{j}=\comp_{J_j}(f^{-(s-j)}(U_s))\text{   for } j=0,\dots,s.$$
Then $U_j\subset D_{\tilde \theta}(J_j)$.
Moreover, we can make the difference $\theta-\tilde\theta$ as small as we like by taking $\varepsilon$ sufficiently small.
\end{lem}
\begin{pf}
Since $f$ is $C^3$, $f$ has an extension to an asymptotically
holomorphic map of order 3 on $M$. Since we are under the Standing Assumptions of Section \ref{sec:enhanced nest}, we can apply Lemma \ref{lem:squared lengths}
to find $\eta=\eta(1/2)$ such that
$\sum_{i=0}^s|J_i|^2<\eta\max_{0\leq i\leq k} |J_i|^{1/2}$.
There exists a constant $\eta_1>0$ so that for any $i$, $0\leq i\leq s$, and $\alpha\in(0,\pi/2)$, we have
$\diam(D_{\alpha}(J_i))\leq \eta_1|J_i|/\sin \alpha.$
Let $K$ be equal to  the constant $\zeta$ from Proposition~\ref{GSSProp2}.
Let $\tilde \theta\in(0,\theta)$ and define 
$$\theta_s=\tilde\theta+(\eta_1 K/\sin \tilde\theta)\sum_{i=0}^{s}|J_i|^2
\leq \tilde\theta+(\eta_2/\sin \tilde\theta)\max_{0\leq i\leq k} |J_i|^{1/2},$$
where $\eta_2= \eta\eta_1K$. We can assume that the
difference $\theta-\tilde\theta$ is as small as we like.
Provided that $\varepsilon>0$ is small enough, we have that $0<\tilde\theta<\theta_s<\theta,$
and that
for each $i$, $D_{\tilde \theta}(J_i)$ is contained in the domain of the asymptotically
holomorphic extension of $f$.

Now, set
$$\theta_j=\theta_{j+1}-K|J_{j+1}|\diam(D_{\theta_{j+1}}(J_{j+1})).$$

We prove that $$\theta_j\geq \tilde\theta+(\eta_1 K/\sin \tilde\theta)\sum_{i=0}^{j}|J_{i}|^2,$$
for $j=0,1,\dots,s$ by induction.
It holds for $\theta_s$ by definition. Suppose that it holds for $\theta_{j+1}.$
Since 
$$\theta_{j+1}\geq \tilde\theta+(\eta_1K/\sin\tilde\theta)\sum_{i=0}^{j+1}|J_i|^2,$$
we obviously have that $\theta_{j+1}>\tilde\theta.$
So that
$$\theta_j=\theta_{j+1}-K|J_{j+1}|\diam(D_{\theta_{j+1}}(J_{j+1}))\geq
\theta_{j+1}-\eta_1 K|J_{j+1}|^2/\sin \theta_{j+1}$$
$$ \geq
\tilde\theta+(\eta_1 K/\sin \tilde\theta)\sum_{i=0}^{j+1}|J_i|^2-\eta_1 K|J_{j+1}|^2/\sin \theta_{j+1}$$
$$\geq \tilde\theta+(\eta_1 K/\sin \tilde\theta)\sum_{i=0}^{j+1}|J_i|^2-\eta_1 K|J_{j+1}|^2/\sin \tilde\theta$$
$$ \geq \tilde\theta+(\eta_1 K/\sin \tilde\theta)\sum_{i=0}^{j}|J_i|^2\geq \tilde\theta.$$
\end{pf}

\begin{remark}
As we can see from Lemma \ref{lem:angle control along diffeos} the
loss of angle when we pull back Poincar\'e disks under asymptotic
holomorphic extensions decreases at small scales. In the rest of the paper, the loss of angle will be calculated for intervals of a given size, so the same bounds will hold at smaller scales. This fact will be implicitly assumed from now on.
\end{remark}

\medskip



\begin{lem}\label{lem:quasiregular pullbacks}
For any $\theta\in(0,\pi)$, $\delta>0$ and $N\in \N$ there exists
$\eta>0$ and $\varepsilon>0$ such that the following holds.
Assume that $J_s\owns c_0$ is a $\delta$-nice real puzzle piece with $|J_s|<\varepsilon$ and that either
\begin{enumerate}
\item $c_0$ is even and $J_s$ is $\delta$-free or 
\item that $c_0$ is odd and all critical points in $\omega(c_0)$ are odd.
\end{enumerate}
Suppose that the chain $\{J_i\}_{i=0}^s$ has order bounded by $N$ and $J_0\cap\omega(c_0)\neq\emptyset$.
Then the asymptotically holomorphic extension
$$f^s\colon \comp_{J_0}f^{-s}(D_\theta(J_s))\rightarrow D_{\theta}(J_s)$$
is $(1+\eta|\mu(J_s)|^{1/2})$-quasiregular.
\end{lem}
\begin{pf}
Let $\mathbf U_s=D_{\theta}(J_s)$ and set $\mathbf
U_j=\comp_{J_j}f^{-(s-j)}(\mathbf U_s).$
Observe that $f^s|\mathbf U_0$ is a proper map of bounded degree. 
Let $k_1<s$ be maximal such that
$J_{k_1}$ contains a critical value of $f$.
Notice that if $c_0$ is even, there exists a constant $\delta'>0$ such that $J_{k_1}$ is
$\delta'$-free.
Let $k_1'\geq k_1$ be minimal so that $J_{k_1'}\subset J_s$.
Since $f^{s-k'_1}|J_{k'_{1}}$ is a diffeomorphism, the sum
$$\sum_{k=k'_1}^s |J_k|^2\leq
\eta_1\max_{k=k'_1,\dots,s-1}|J_k|^{1/2}$$
where $\eta_1=\eta_1(1/2)$ is the constant from Lemma \ref{lem:squared lengths}.
Then by Lemma \ref{lem:angle control along diffeos}, there exists a
constant $\eta'>0$ such that for any point $z_k\in \mathbf U_k, k=k'_1,k'_1+1,\dots,s-1$,
$$|\mu(z_k)|=\frac{\overline\partial f}{\partial f}(z_k)\leq \eta'(\mathrm{Im}(z_k)^2)\leq
\eta'(\diam(\mathbf U_k))^2\leq \eta'|J_{k}|^2.$$
Where the first inequality follows from the fact that $f$ is
asymptotically holomorphic of order three.

So we get that there exists a constant $\eta''>0$ such that
$$\sum_{k=k_1'}^s|\mu(z_k )|\leq \eta''\max_{k=k'_1,\dots,s-1}|J_k|^{1/2}.$$

Let $s_1<k_1$ be maximal so that 
$J_{s_1}\subset J_s.$
The domains $J_{s_1}, J_{s_1+1},\dots,J_{k_1'}$ are disjoint and small when 
$J_{k_1'}$ is small, so, 
since either all critical points in $\omega(c_0)$ are odd or $c_0$
is $\delta'$-free by Lemma~\ref{lem42}, 
$\diam f^{-1}(\mathbf U_{k_1})\asymp|J_{k_1-1}|$.
Thus, estimating as before,
if $z\in\mathbf U_{k_1'}$,
$$\sum_{t=s_1}^{k'_1-1}|\mu(z_t)|\leq \eta''\max_{k=s_1,s_1+1+\dots k'_1-1}|J_k|^{1/2}.$$

We repeat this argument letting $k_2<s_1$ be maximal so that $J_{k_2}$
contains a critical value of $f$. If no such $k_2$ exists, set
$k_2=0.$
Let
$k_2'\geq k_2$ be minimal so that
$J_{k_2'}\subset J_{s_1}$ and $s_2<k_2'$ maximal so that
$J_{s_2}\subset J_{s_1}$, and if no such $s_2$ exists, set $s_2=0$.
Proceeding inductively, let $p$ be maximal so that $s_p\neq 0$. 

Let $z\in \mathbf J_0$
Then $$\sum_{k=0}^s|\mu(f^k(z))|=\sum_{k=0}^{s_p}|\mu(f^k(z))|+\sum_{j=1}^p\sum_{k=s_j}^{s_{j-1}}|\mu(f^{k}(z))|
\leq \eta''\sum_{j=1}^p\max_{k=s_j,\dots s_{j-1}-1}|J_k|^{1/2}.$$ 

Since $\max_{k=s_j,\dots s_{j-1}-1}|J_k|^{1/2}$ decays
exponentially (every $J_k$ is $\delta'=\delta'(\delta, N)>0$ nice), the result follows.
\end{pf}

\section{Pulling back a Poincar\'e disc through high iterates of first return maps}
\label{sec:pulling back disks}

The results in this section apply to both maps $f\in\mathcal A_{\underline b}$ and maps $f\in \mathcal A_{\underline b}^3$ for which the asymptotically holomorphic extensions are as in Subsection~\ref{sub-sec:asympholoS}. Recall, we will always work under the Standing Assumptions on page \pageref{standing assumptions}.

\subsection{Pullbacks through a monotone branch}

\begin{prop}\label{monotone pullbacks}
For each $\delta>0$ and $C>0$ there exist $\delta'>0$  and
$\lambda'\in(0,1),$  and for each $\theta\in(0,\pi)$ there exists $\varepsilon>0$ 
such that the following holds.
Let $I^0$ be an interval with $|I^0|<\varepsilon$ and let  $I^1\subset I^0$ be a first return domain to $I^0.$ Let $r$ be so that  $R_{I^0}|_{I^1}=f^r.$ Assume $f^r|_{I^1}$ is  monotone and can be decomposed into at most $b$ maps of the form $g\circ p_\ell\circ h$, where $p_\ell=x^\ell$ for some odd integer $\ell>0$ and $g$ and $h$ are diffeomorphisms with  bounded distortion. Assume $$|f^{r}(c)-c|\geq\delta |I^{1}|\; and \; |Df^{r}(x)|\leq C$$
for all critical points $c$ of $f^r |_{I^1}$ and all $x \in
I^{1}.$ Let $\{G_{j}\}_{j=0}^{pr}$ be a disjoint chain with
$(f^r|_{I^1})(G_{jr})=G_{(j+1)r}$ for  $ j=0\ldots  p-1,$ $G_{pr}\subset I^0$
and $G_0\cap \omega(c_0)\neq \emptyset$. Assume there exist disjoint  intervals
$K, F^1,F^2\subset I^0$ 
with the following properties.
\begin{itemize} 
\item $G_{pr}\subset K,$ 
\item $F^1$ and $F^2$ are contained in fundamental domains for
$f^r|_{I^1}$
\item  $|F^i|>\delta|I^0|$ for $i=1,2$, and
\item $K$ is contained in the convex hull of $F^1$ and $F^2.$
\end{itemize} 
Let $$V=D_{\theta}(K)\cap\mathbb{C}_{G_{pr}}\
and\ U_{i}=\mathrm{Comp}_{G_{i}}f^{-(pr-i)}(V), \mbox{ for }i=0,\dots,pr.$$
Then for each $z\in U_{0}$ there exists an interval $K'$ such that
$$z\in D_{\lambda'\theta}(K'),$$ where $G_{0}\subset
K'\subset(1+2\delta')K'\subset I^1.$
\end{prop}

\begin{pf}
Without loss of generality we will assume the map $f^r|_{I^1}$ is orientation preserving and that its fixed point $p$ is equal to $0$. The orientation reversing case will follow from the orientation preserving case by considering by $f^{2r}$ instead of $f^r$ and applying Lemma \ref{chain pullback} once.

Since $f^r|_{I^1}$ has no periodic attractors we have that for each of its critical
points $c$ there exists a maximal $m(c)>0$ so that
$f^{rj}(c)\in I^1$ of all $0\leq j\leq m(c)$. Consider $\Omega=\cup
\{c,\ldots, f^{m(c)}(c)\},$ where the union is taken over all critical
points of $f^r|_{I^1}.$  Since each $F^k$ is contained in a
fundamental domain,
$F^k\cap \Omega$ contains at most $b$ points for $k=1,2$. So by slightly shrinking $F^k,$ if necessary, we can assume $ \Omega\cap F^k=\emptyset$  and that the convex hull of $F^1$ and $F^2$ is $\delta/2$-well inside $I^0.$ Observe this means that no pullback of $F^k$ under $f^r |_{I^1}$ contains a critical point for $k=1,2.$
\medskip

We will first prove the proposition for the unicritical case.
Unless is specify otherwise, we will
assume that all pullbacks of Poincar\'{e} domains under $f$ are taken
under the branch of $f^{-1}$ that maps $G_{i+1}$ to $G_i$  for
$i\in\{0,\ldots pr-1\}$. Observe that the inverse branches of $f$
extend to the complex plane and we are able to express which inverse
branch we are using by choosing the domain : when we write $f^{-1}\colon \mathbb C_{G_{j+1}}\rightarrow \mathbb C$, we always mean the choice of inverse branch that maps $G_{j+1}$ to $G_j$. Consider the chain $\{H_j\}^r_{j=0}$, with $H_r=I^0$ and $H_0=I^1.$

Let $c$ be the critical point of $f^r|_{I^1}$. We will denote  by
$c_j=f^j(c)$ for $j\in \mathbb{Z},$ where the inverse images of $c$
will  be taken under the pullbacks defined above. Since $c$ is a
critical point of $f^r$, there exists  $i$ (maximal)
so that  $0 \leq i< r$ and $c_{i}=c'$ for $c'\in \crit(f)$. 
To orient ourselves we will assume $c>0$. From the fact that $f^r$ has no periodic attractors, and that $f^r(0)=0,$ we have that  $0<c_{-r}<c<c_r$. The case $c<0$ is analogous to this one, taking into consideration that in this case $c_{r}<c<c_{-r}<0$. 

\medskip
We can assume that $c_r$ is not contained in $G_{pr}$. Otherwise, let
$J_r$ be the convex hull of $F^1$ and $F^2$. Since $J_r\setminus
  (1+\delta)^{-1}J_r$ is disjoint from $\Omega$, we can apply
 Corollary \ref{free space well inside} to find $\hat \lambda>0$ so that $U_{(p-1)r}\subset D_{\hat \lambda \theta}(J_0)\cap \mathbb C_{G_{(p-1)r}}$,where $f^{r}(J_0)=J_r.$ Since $G_{(p-1)r}$ does not contain $c_r$ we can start the proof at time $(p-1)r$, instead of $pr$.

\medskip

We will first consider the case when $G_{pr}$ is in the right hand
component of $\mathbb R\setminus\{c_r\}$, $G_0$ is in the left
one, and $F^1$ and $F^2$ are in opposite components.
Without loss of generality assume $F^2$ and $G_{pr}$ are on the same
component of $\mathbb R\setminus\{c_r\}$. Note that the relative
position of $F^1$ and $F^2$ with respect to $0$ has not been
specified; see Figure \ref{fig:uno} for a possible configuration.

\begin{figure}[h] \centering \def\svgwidth{200pt}
\scriptsize
 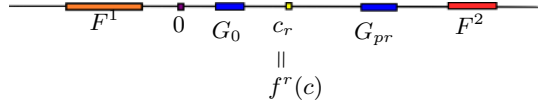
 \caption{\small Configuration w.r.t. $c_r$ ( the critical value of $f^r$).} 
 \label{fig:uno}
 \end{figure}

Let $R_{pr}=F^2$ and $L_{pr}=F^1$ and consider the chain $\{R_j\}^{pr}_{j=(p-1)r}$  given by the pullbacks under 
the branches of $f^{-1}$ which map $G_{j+1}$ to $G_j.$ 
 Observe that all intervals $R_j$ are contained in $\mathbb R$ for $(p-1)r\leq j\leq pr$. On the other hand, if we consider the pullbacks of $L_{pr}$ under the same maps, we get that $L_{(p-1)r+i+1}$ is not in the same component of $\mathbb R\setminus \{c'\}$ as $G_{(p-1)r+i+1}$  (see Figure ~\ref{fig:dos}). 

 \begin{figure}[h] \centering \def\svgwidth{200pt}
\scriptsize
 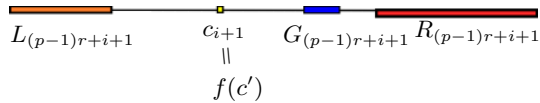
\caption{\small Pullback  at time {i+1}.} 
 \label{fig:dos}
 \end{figure}
\noindent

This means that the pullback of  $L_{(p-1)r+i+1}$ under the branch of
$f^{-1}$ which maps  $G_{(p-1)r+i+1}$ to  $G_{(p-1)r+i}$ does not lie
in the real line.

By Lemma \ref{Poincare disks}, there exists $\lambda\in(0,1)$ so that $$V\subset D_{\lambda \theta}(L_{pr})\cup D_{\lambda\theta}(R_{pr})\cup A,$$
where $A$ is the bounded component of the complement of
$D_{\lambda\theta}(L_{pr})\cup D_{\lambda\theta}(R_{pr})$. In order to control the pullbacks of points in $V$ we will, separately, control the pullbacks of points in $D_{\lambda\theta}(L_{pr}),$ $D_{\lambda\theta}(R_{pr})$ and $A.$ 
%
\medskip

\medskip
Let $q$ be maximal with $0\leq q< p$ and so that  $G_{qr}$ and $G_{pr}$ are in the same component of $\mathbb R\setminus \{c_r\}.$
From the choice of $q$, the intervals $R_{jr}$ and $G_{jr}$ are in the same component of $\mathbb R\setminus \{c_r\}$ for $j=q,\dots,p.$

\medskip
\noindent \textbf{To pullback $ D_{\lambda\theta}(R_{pr})$.}  Let $z^{pr}\in D_{\lambda\theta}(R_{pr})$ and consider the sequence of points $\{z^j\}^{pr}_{j=0}$, where $z^j$ is preimage of $z^{j+1}$ under the branch of $f^{-1}$ specified above.

\begin{itemize}

\item[(a)] \textit{Assume that $q=0$ or that  $q>0$ and the intervals $G_{(q-1)r}$ and $R_{(q-1)r}$ lie in the same component of $\mathbb
R\setminus\{c_r\}.$ }
Then, since the chain $\{R_j\}^{pr}_{j=(q-1)r}$ is disjoint we can apply Lemma \ref{lem:angle control along diffeos} to find $\lambda_0\in(0,1)$ so that 
$$z^{0}\in D_{\lambda_0\theta}(R_{0}).$$
Observe that $R_0$ lies between $0$ and $R_{pr}$ and that the distance between $R_0$ and $0$ is comparable to $|R_0|.$ So there exists $\delta_0>0$ so that $(1+2\delta_0)R_0\subset I^1$, and the proposition follows.

\item[(b)] \textit{Assume that $q>0$ and the intervals $G_{(q-1)r}$ and $R_{(q-1)r}$ lie in opposite components of $\mathbb R\setminus\{c_r\}.$ }
Apply Lemma \ref{lem:angle control along diffeos}  to find $\lambda_1\in(0,1)$ so that 
$$z^{(q-1)r}\in D_{\lambda_1\theta}(R_{(q-1)r}).$$
 If $q=1$ then $(q-1)r=0$ and the result follows as in (a). 
If  $q>1,$ then by Lemma \ref{lem:angle control along diffeos}, we can choose $\lambda_2>0$ so that $$z^{(q-2)r+i+1}\in D_{\lambda_2\theta}(R_{(q-2)r+i+1}).$$
Observe that  $R_{(q-2)r+i+1}$ and  $G_{(q-2)r +i+1}$   are in different components of $\mathbb R\setminus \{c_{i+1}\}$.  By the assumptions on $f^r,$ we know that the interval  $(c_{-r-j},c_{j})$ has size comparable to $|H_{j}|$ for $j\in \{i,i+1\}$. So 
we can apply Lemma \ref{different branch angle control} to find  $\lambda_3>0$ so that image of 
 $D_{\lambda_2\theta}(R_{(q-2)r+i+1})$ under
the branch that maps $G_{(q-2)r+r+i}$
 to $G_{(q-2)r+r+i-1}$
is contained in $D_{\lambda_3\theta}(c_{-r-i}, c_i)$. See Figure \ref{fig:tres}. 
 
 \begin{figure}[htb] \centering \def\svgwidth{200pt}
\scriptsize
 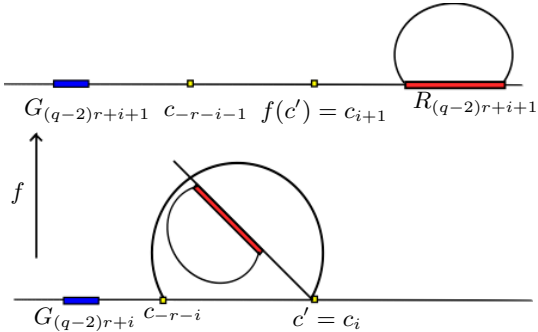
 \caption{\small Controlling  the loss of angle  after $f(c')$ ($c'\in\crit(f)$).} 
\label{fig:tres}
\end{figure}
 
In particular $$z^{(q-2)r+i}\in D_{\lambda_3\theta}(c_{-r-i},c_i).$$  
Finally, if we let $R'_{(q-2)r+i}=(c_{-r-i},c_i)$
and we consider the chain 
$\{R'_j\}^{(q-2)r+i}_{j=0}$ the result follows as in Case (a).

\end{itemize}

In conclusion, exist an interval $\tilde K\subset I^1$ and constants $\tilde \lambda \in(0,1)$ and $\tilde \delta>0$ so that $(1+2\tilde \delta)\tilde K\subset I^1$  and for all $z^{pr}\in D_{\tilde{\lambda}\theta}(R_{pr})$ we have that
$$z^0\in D_{\tilde{\lambda} \theta}(\tilde K).$$

\medskip
\noindent \textbf{To pullback $ D_{\lambda\theta}(L_{pr})$.}  
Pick $z^{pr}\in D_{\lambda\theta}(L_{pr})$ and consider the sequence of points $\{z^j\}^{pr}_{j=0}$, where $z^j$  
is the preimage of $z^{j+1}$ under the branch of $f^{-1}$
  specified above. 

\medskip
By Lemma \ref{lem:angle control along diffeos}, there exists $\lambda'_1>0$ so that $$z^{(p-1)r+i+1}\in D_{\lambda'_1\theta}(L_{(p-1)r+i+1}).$$ 
By assumption $L_{(p-1)r+i+1}$ and $G_{(p-1)r+i+1}$ are in different components of $\mathbb R\setminus \{c_{i+1}\}$. 
Let $\hat L=(f^{-1}|_{\mathbb C_{G_{i+1}}})(D_{\lambda'_1\theta}(L_{(p-1)r+i+1})).$
By Lemma \ref{different branch angle control}, there exists $\lambda'_2>0$ so that 
 $$\hat L\subset D_{\theta}(c_{i}, w_{i})\cup D_{\lambda'_2\theta}(R_{(p-1)r+i})\cup A'$$
where $A'$ is the bounded component of the complement of  $D_{\theta}(c_{i}, w_{i})\cup D_{\lambda'_2\theta}(R_{(p-1)r+i})$ 
 and $w_{i}$ is defined as follows. If $c_r\in I^1$, let
 $w_{i}=c_{2r-i}$. If $c_r\notin I^1$, let $w_{i}$ be the boundary
 point of $H_{i}$ that lies on the same side of $c_{i}$ as $G_{(p-1)r+i}$. In either case,
 the size of the interval $(c_{i}, w_{i})$ is comparable to $|R_{(p-1)r+i}|$ and $|H_{i}|.$ See Figure \ref{fig:51}.
 
 \begin{figure}[htb] \centering \def\svgwidth{200pt}
\scriptsize
 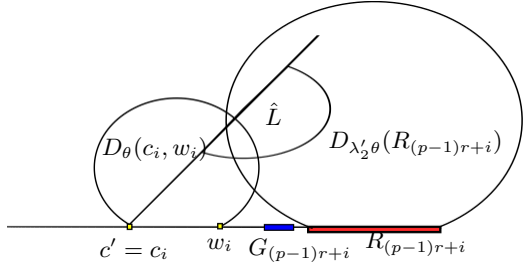
 \caption{\small Use two domains to control the loss of angle at $c_i=c'\in\crit(f)$} 
  \label{fig:51}

\end{figure}
 
\noindent 
Consider the chain $\{B_j\}^{r-1}_{j=0}$ with $B_{r-1}=(c_{i}, w_{i})$ and $B_0$ the pullback of $B_{r-1}$ contained in $H_{i+1}.$ 
By Lemma  \ref{lem:angle control along diffeos}, we can choose $\lambda'_3>0$ so that $$\comp_{B_0} (f^{-(r-1)}D_\theta(B_{r-1}))\subset  D_{\lambda'_3\theta}(B_0).$$
By Lemma \ref{different branch angle control}, we can choose  $\lambda'_4>0$ so that 
the component of the inverse image of $D_{\lambda'_3\theta}(B_0)$ under $f$ that lies on the same side of $c_{i}$ as $G_{(p-1)r+ i}$, denoted by $\hat B$, satisfies the following  
 $$\hat B \subset D_{\theta}(c_{i}, w_{i})\cup D_{\lambda'_4\theta}(R_{(p-1)r+i})\cup C'$$
 where $C'$ is the bounded component of the complement to $D_{\theta}(c_{i}, w_{i})\cup D_{\lambda'_4\theta}(R_{(p-1)r+i})$. See Figure \ref{fig:10}.

 \begin{figure}[htb] \centering \def\svgwidth{200pt}
\scriptsize
 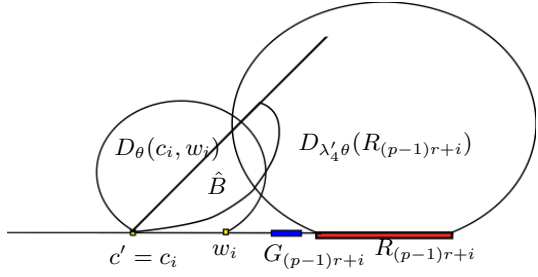
 \caption{Adjust the angle at $c_i=c'\in\crit(f)$} 
 \label{fig:10}
 \end{figure}

Thus,  $$z^{(p-1)r+i}\in D_{\theta}(c_{i}, w_{i})\cup D_{\lambda'_4\theta}(R_{(p-1)r+i})\cup C'.$$



By definition of  $\lambda'_4$ we have that one of the following holds:

\begin{itemize}
\item[(i)] $z^{kr+i}\in D_{\theta}(c_{i}, w_{i})$ for all $k\in \{q-1,\ldots p-1\}.$ In particular, $z^{(q-1)r+i}\in D_{\theta}(c_{i}, w_{i})$. If we let $R'_{(q-1)r+i}=(c_{i}, w_{i}),$ the proposition follows as in (a); considering the chain $\{R'_j\}^{(q-1)r+i}_{j=0},$ instead of the the chain $\{R_j\}^{(q-1)r+i}_{j=0}.$

\item[(ii)] There exists $k\in \{q-1,\ldots p-1\}$ so that $z^{kr+i}\in
  D_{\lambda_4\theta}(R_{(p-1)r+i})$. Then define
  $R'_{kr+i}=R_{(p-1)r+i}$ and apply arguments (a) or (b) above.

\item[iii)] $z^{kr+i}\in C'$ for all $k\in \{q-2,\ldots p-1\}.$ In
  particular, $z^{(q-2)r+i}\in C'.$ From Cases (i) and (ii) we know that there exist $\delta'_5>0,$ $\lambda'_5\in(0,1)$ and intervals $K_1$ and $K_2$ which are $\delta'_5$-well inside $I^1$ so that the following holds. If  $z^{(q-2)r+i}\in
 D_{\theta}(c_i,w_i)$ then $z^0\in D_{\lambda'_5\theta}(K_1)$. And if $z^{(q-2)r+i}\in D_{\lambda'_4\theta}(R_{(p-1)r+i})$ then $z^0\in D_{\lambda'_5\theta}(K_2).$ So by the continuity of $f$, if $z^{(q-2)+i}\in C'$ we have that $z^0$ is contained in the bounded component of the complement of
 $D_{\lambda'_5\theta}(K_1)\cup D_{\lambda'_5\theta}(K_2).$  If we let $K_3$ denote the convex hull of $K_1$ and $K_2$ we have that $z^0\in D_{\lambda'_5\theta}(K_3).$
Finally, the existence of $\lambda'_5$ implies that $|K_1|$ is
comparable to $|K_3|$. Since $K_1$ is a fundamental domain, or contained in one, the distance between $K_1$ and $0$ is comparable to $|K_1|.$ So there exists $\delta'_6>0$ so that $K_3$ is $\delta'_6$-well inside $I^1,$ and the proposition follows.

\end{itemize}

\medskip
\noindent \textbf{To pullback $A$.} 
Making $\lambda'_4$ smaller, if necessary, we can make sure that the following holds. Given $z^{pr}\in D_{\lambda\theta}(R_{pr})$ then $z^{(p-1)r+i}\in D_{\lambda'_4\theta}(R_{(p-1)r+i}).$ Observe this implies that given $z^{pr}\in A$ then
$$z^{(p-1)r+i}\in D_{\theta}(c_{i}, w_{i})\cup D_{\lambda'_4\theta}(R_{(p-1)r+i})\cup C',$$
and the proposition follows from one of the previous cases.


\medskip
This completes the proof of the proposition in this case. The arguments given in this proof 
also deal with the case when $G_{pr}$ and $G_0$ lie on the same side of $c_r$ and when $F^1$ and $F^2$ lie on the same side of $0$. Thus the unicritical case follows.

If $f^r|_{I^1}$ is multicritical it will have at most $b$ critical points. It is clear that, by the monotonicity of $f^r|_{I^1},$ the proposition in this case follows by applying the unicritical case at most $b+1$ times.

\end{pf}

\noindent\textbf{Remark 1.} From the proof of Proposition~\ref{monotone pullbacks} we get the  extra information that $K\subset [x,y]$ where $x,y\in I^0$ are of one the two following types:
\begin{itemize}
\item $x\in \{c^j_{-r},c^j, c^j_{r}\}$, for $j\in\{1,m\}$ and $y\in \partial I^0$;
\item $x\in \{c^j_{-r}, c^j,c^j_{r}\}$ and $y\in\{c^{j+1}_{-r}, c^{j+1}, c^{j+1}_{r}\}$, for $j\in \{1,\ldots m-1\}.$
\end{itemize}

\medskip

\noindent\textbf{Remark 2.} We cannot take real pullbacks of the fundamental domains at each step when there are critical points since we cannot control the order of the pullback of the larger interval $(F^1, F^2)$. 

\noindent\textbf{Remark 3.} If $f^{r}$ is a diffeomorphism, we do not require the bound on the derivative.

\subsection{Pullbacks for which the modality is bounded}

\begin{prop}\label{lem:N-modalpb}    
For any $\delta>0$ and $N\in\mathbb{N}$ there exists $\lambda>0,$ and for any $\theta\in(0,\pi)$ there exists $\varepsilon>0$ so that the following holds. 
Let $I$ be a nice interval with $|I|<\varepsilon$.
Suppose that either:
\begin{itemize}
\item[(1)] $I$ is  $\delta$-free or 
\item[(2)] $c_0$ is odd.
\end{itemize}
Let $J$ be an $N$-modal pullback of $I$ by $f^{t}$ for some $t>0$ with $J\cap\omega(c_{0})\neq\emptyset$. 
Let $V=D_{\theta}(I)$ and  $U_{i}:=\mathrm{Comp}_{f^{i}(J)}f^{-(t-i)}(V)$. Then $U_{0}\subset D_{\lambda\theta}(J)$.
\end{prop}
\begin{pf}
If we are in case (1), then we have the following.
Let $\{G_{j}\}_{j=0}^{t}$ be the chain with $G_{t}=I$ and $G_{0}=J$.
Since $I$ is $\delta$-internally free, $F'\subset I$ is a nice pair, $\omega(c_{0})\cap I\subset F'$
and $F'\subset (1+2\delta)^{-1} I$. Let $\{H_j\}_{j=0}^t$
be the chain corresponding to the pullbacks of $F',$ with $H_j\subset G_j$.  
Let $F_j^{\pm}$ be the components of $G_j\setminus H_j$. Observe that $G_t,H_t$ is a nice pair, so their pullbacks are also nice pairs.
Since each of the intervals $G_j$ intersects $\omega(c_{0})$ whereas
$G_j\cap H_j$ does not intersect $\omega(c_{0})$ it follows that
the intervals $F^\pm_j$ are all disjoint. Moreover, since
$H_j\subset (1+2\delta')^{-1}G_j$, it follows that the intervals
$F_j^+$ and $F_j^-$ are not small compared to $G_j$.

Decompose $f^t\colon J\rightarrow I$ as $f^t|J=D_{N+1}\circ P_{\ell_N}
\circ\dots\circ P_{\ell_2}\circ D_2\circ P_{\ell_1}\circ D_1$, where
each $D_i$ is a diffeomorphism and $P_{\ell_i}\colon z\mapsto z^{\ell_i}$.
The result follows by using Lemma~\ref{lem:angle control along
  diffeos}
to control the loss of angle when we pull back by a $D_i$ and 
Lemma~\ref{lem42}
when we pull back by a polynomial.
The first paragraph of this proof ensures that
whenever we pullback by an even polynomial,
the critical value is not close to the boundary. 

\end{pf}

\medskip

\subsection{Pulling back through a non-monotone branch}
In this subsection we will control pullbacks of Poincar\'e disks through a long
cascade of central returns in the principal nest.

\medskip
Let $I$ be a nice interval containing a critical point $c\in\omega(c_0)$. Recall the 
definition of the principal nest: $I^{0}=I$ and if $I^{n}$ is defined, we set $I^{n+1}=\mathcal{L}_{c}(I^{n})$ for $n\in \N$. Let $r$ be so that $R_I|_{I^1}=f^r.$ If $f^r|_{I^1}$ has a turning point, define 
$\tilde{m}$ to be minimal with the property that there is a turning point $c'$ of $f^r|_{I^1}$ 
such that $f^r(c')\notin I^{\tilde{m}}$. Note that  if $I$ is periodic, $I=I^1=I^2=\dots$ and $\tilde m=\infty$. In this case $f^r(I^{\tilde m})\subset I^{\tilde m}$
and $f^r(\partial I^{\tilde m})\subset \partial I^{\tilde m}$, so $f$ is renormalizable and $I^{\tilde m}=I^\infty$ is  periodic. Under the above assumptions we have the following.

\begin{prop}[cf. \cite{KSS} Lemma 11.1]\label{good start or good deep}
There exists a constant $\hat{\rho}>0$ and for each $\delta>0$ there exist $\delta'>0$ and  
$\lambda\in(0,1)$  so that for each $\theta\in(0,\pi)$ there exists $\varepsilon>0$ 
with the following properties. Let $I^0=I$, $I^1=\mathcal{L}_{c}(I)$ and assume that
$|I|<\varepsilon$ and $(1+2\delta)I^1\subset I$. Let $r$ be so that $R_{I^{0}}|_{I^{1}}=f^r$. Additionally, if $c$ is odd assume $f^r|_{I^1}$ is monotone. 
Let $\{G_{j}\}_{j=0}^{pr}$ be a disjoint chain 
with $G_{jr}\subset I^{1}$ for all $0\leq j\leq p-1,$ so that pullbacks of $G_{pr}$ and $I$ are nested or disjoint. Let $\{\hat{G}_{j}\}_{j=0}^{pr}$ be a chain with 
$G_{pr}\subset \hat{G}_{pr}\subset (1+2\delta)\hat{G}_{pr}\subset I$ and $G_{0}\subset \hat{G}_{0}$ and define
$$V=D_{\theta}(\hat{G}_{pr})\cap\mathbb{C}_{G_{pr}} \mbox{ and } U=\mathrm{Comp}_{G_{0}}f^{-pr}(V).$$
If $c$ is even, then for each $z\in U$, there exists an interval $K$ such that either
\begin{equation}\label{good at start}
z\in D_{\lambda\theta}(K)\mbox{ and } G_{0}\subset K\subset (1+2\delta')K\subset I
\end{equation}
or there exists $0\leq p'\leq p$ and   intervals $K\subset I'\subset I^{\tilde{m}}$ with
\begin{equation}\label{good deep}
f^{p'r}(z)\in D_{\lambda\theta}(K)\mbox{ and }  G_{p'r}\subset K\subset (1+2\delta')K\subset I'
\end{equation}
where $I'$ is $\hat{\rho}$-nice if $\tilde m<\infty$.

If $c$ is odd, then for each $z\in U$, there exists an interval $K$ such that
(\ref{good at start}) holds.
\end{prop}

\begin{pf} Let $z\in U$. If $p \leq 10$ or if $\tilde{m}=1$, then we can apply Corollary \ref{chain space} finitely many times to complete the proof. If $p\le 10$, then (\ref{good at start}) holds and if $\tilde{m}=1$, then (\ref{good deep}) holds for $p'=p-1,$ so we will assume that $p>10$ and $\tilde{m}\geq 2$.

\medskip
Let us first assume that the $R_I$ restricted to $I^1$ is
monotone. Apply Corollary \ref{chain space} twice to find
$\lambda_1>0$ so that $f^{(p-2)r}(z)\subset
D_{\lambda_1\theta}(I^2)$. If the size of $I^2$ is comparable to the
size of $I^1$ we can apply Lemma \ref{bounds1} to $f^r|_{I^2}$. Since
$(1+2\delta)I^1\subset I^0$ there exists $\delta_1>0$ so that $I^2$ is
$\delta_1$ well-inside $I^1$. Using the connected components of
$I^1\setminus(1+\delta_1)I^2$ we can find suitable fundamental domains and apply
Proposition~\ref{monotone pullbacks} to show (\ref{good at start}).

If $|I^2|/|I^1|<\varepsilon$, for $\epsilon>0$ small, we get that
$I^2$ is deep-inside $J=(1+\delta)I^1$. Since $J$ is well-inside $I$
and $f^r$ is monotone there exists $C=C(\varepsilon)>0$ so that $\hat
J=\Comp_{c_0}f^{-r}(J)\subset (1+2C)\hat J\subset J$. Note that $C\to
\infty$ as $\epsilon\to 0$. 
Since $\tilde m>1,$ there are no critical points or critical values of
$f^r|I^1$ contained in $I\setminus I^1$, so we can apply
Lemma~\ref{chain pullback}
to find $\lambda_2\in(0,1)$ so that
$\comp_c f^{-r}D_{\lambda_1\theta}\subset D_{\lambda_2\theta}(J')$.
Making $\epsilon$ small enough we get the following. Given $x\in \hat J$

$$\Comp_c f^{-r} D_{\lambda_1\theta}(J)\subset D_{\lambda_2\theta}(J')\subset D_{\lambda_1\theta}(J),$$
so we get (\ref{good at start}). 

From now on we will assume that $c$ is even.

\medskip
\noindent\textit{Claim 1.} There exists a (universal) constant $\nu\in (0,1)$ such that if $|I^{2}|/|I^{0}|\leq\nu$ 
then (\ref{good at start}) holds. 

\medskip
\noindent\textit{Proof of Claim 1.}
Assume $|I^2|/|I^0|\leq\nu$. Since $c$ is even, there exists
$C=C(\nu)$ so that $(1+2C)I^2\subset I^0$, which in turn implies that
there exists a constant $C'=C'(\nu)$ so that $(1+2C')I^2\subset I^1$. 
Recall that
$f^{(p-2)r}(z)\in D_{\lambda_1\theta}(I^{2}).$ 
By Corollary \ref{chain space},
there exists $\lambda_3\in(0,1)$ such that for each $x\in I^2$,
$$\comp_{x}f^{-r}(D_{\lambda_1\theta}(0.5 I^1))\subset D_{\lambda_3\theta}(I^2).$$
If $\nu$ is sufficiently small, depending only on the choice of 0.5, 
then $$D_{\lambda_3\theta}(I^2)\subset D_{\lambda_1\theta}(0.5I^1),$$
from which it follows that
(\ref{good at start}) holds. \endpfclaim


\medskip


From now on, we will assume that $|I^{2}|/|I^{0}|>\nu$. 

\medskip

 Let $E_{1}$ and $E_{2}$ be the outer monotone branches of $f^{r}|I^{1}$. 
In order to be definite, let $E_1$ be so that $f^r|E_1$ is orientation preserving. 
Let $c_{l}$ be the turning point in $\partial E_{1}$ and $c_{r}$ the turning point in $\partial E_{2}$

\medskip
\noindent \textit{Claim 2.} There exist universal constants $\kappa_{1},\kappa_2>0$ and $C>1$ such that
\begin{itemize}
\item[(1)] for any critical point $c$ of $f^{r}|I^{2}$, $\dist(f^{r}(c),c)\geq\kappa_{1}|I^{0}|$;
\item[(2)] $|(f^{r})'(x)|\leq C$ for any $x\in I^{2}$;
\item[(3)] either  $f^{r}(I^{\tilde{m}})\cap I^{\tilde{m}}=\emptyset$ or $|E_{i}\cap I^{\tilde{m}}|\geq \kappa_{2}|I^{0}|,$ for $i=1,2.$
\end{itemize} 

\medskip
\noindent\textit{Proof of Claim 2.}
From Claim 1 we have that $|I^{2}|/|I^{0}|>\nu$,  so (1) and (2) follow from Lemma \ref{bounds1}. Statement (3) follows from Statements (1) and (2).\endpfclaim

\medskip

\noindent
\textit{Claim 3.} It is enough to prove the proposition for the case
that $G_{pr}\subset I^{\tilde m}$ and  $f^{pr}(z)\in D_{\hat\lambda\theta}(J)\setminus  D_\theta(I^1)$
for some constant $\hat\lambda\in(0,1)$ and $\hat \delta>0$ (which do not depend on 
$\theta)$ and some interval $I^{2}\subset J\subset (1+2\hat \delta)J\subset I^{1}$. 
In particular, we can assume 
that $f^{r}(I^{\tilde{m}})\cap I^{\tilde{m}}\neq \emptyset$ (because otherwise (\ref{good at start}) holds). 

\medskip
\noindent
\textit{Proof of Claim 3.}  
If necessary we can apply Corollary \ref{chain space}  to obtain
$\lambda_4 \in(0,1)$
so that  $f^{(p-1)r}(z)\in D_{\lambda_4\theta}(I^{1})$. Replacing $p$ by $p-1$ and $\theta$ by $\lambda_4\theta$ we may assume that $\hat{G}_{pr}\subset I^{1}$
and $f^{pr}(z)\in D_{\theta}(I^{1}).$ If $z\in D_{\theta}(I^{1})$, then (\ref{good at start}) holds and  the proof is completed.
So we may assume that there exists a maximal $q$ with $0\leq q<p$ such that $f^{qr}(z)\notin D_{\theta}(I^{1})$. Since $f^{(q+1)r}(z)\in D_{\theta}(I^{1})$, by Corollary \ref{chain space} there exist $\delta_{2}>0$ and $\lambda_{5}\in(0,1)$ such that
\begin{equation}f^{qr}(z)\in D_{\lambda_{5}\theta}(J)\setminus D_{\theta}(I^{1}),
\label{eq:claim2prop52}
\end{equation}
where $J$ is an interval with $I^{2}\subset J\subset (1+2\delta_{2})J\subset I^{1}$ (see Figure \ref{escape}).

\begin{figure}[htp]
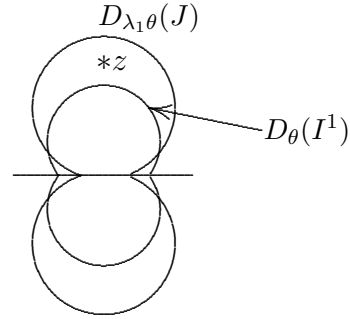
 \hfil
\beginpicture
\dimen0=0.3cm
\setcoordinatesystem units <\dimen0,\dimen0>  point at -15 0 
\setplotarea x from -4 to 0, y from -6 to 6
\setlinear
\plot -8 0 0 0 /
\circulararc -320 degrees from -5 0  center at  -4 3 
\circulararc 320 degrees from -5 0  center at  -4 -3 
\put {{\small $D_{\lambda_1\theta}(J)$}} at -2 7 
\circulararc -252 degrees from -6 0  center at  -4 1.5 
\circulararc 252 degrees from -6 0  center at  -4 -1.5 
\arrow <3mm> [0.2,0.67] from 3 2  to -2 3 
\put {{\small $D_{\theta}(I^1)$}} at 5 2 
\put {{$*$}} at -4 5 
\put {{$z$}} at -3.3 5
\endpicture
\caption{\label{fig:DI}The sets $D_{\theta}(I^1)$ and $D_{\lambda_1 \theta}(J)$.
The point $z$ is {\lq}jumping{\rq} and we will capture it by $J$.}
\label{escape}
\end{figure}

Let us first see what happens if $f^{r}(I^{\tilde{m}})\cap I^{\tilde{m}}= \emptyset.$  In this case the only way that an interval $G_j$ could have a pullback under $f^r$ is if it is a pullback under $E_1$. Assume this is the case. Claim 1 and Lemma~\ref{Poincare disks} allow us to find  $\delta_3, \lambda_6>0$ and  $K$ an interval $\delta_3-$well-inside $E_1$ such that $f^{qr}(z)\in D_{\lambda_6\theta}(K)$. 
We can find fundamental domains of size comparable to $I^{1}$ for the map $f^{r}|E_{1}\colon  E_{1}\rightarrow f^{r}(E_{1})$ as follows. On one side the fundamental domain is given by an interval in $I^{0}\setminus I^{1}$ and on the other side it is given by the $c_{l}$ and  $f^{r}(c_{l})$. Applying Proposition~\ref{monotone pullbacks} we get  (\ref{good at start}). From now on assume
 $f^{r}(I^{\tilde{m}})\cap I^{\tilde{m}}\neq \emptyset.$

If $q\leq2$, then by applying Corollary \ref{chain space} at most twice, we obtain (\ref{good at start}).  So we may assume that $q\geq 3.$ If $G_{qr}\subset I^{\tilde{m}}$, then the claim follows from (\ref{eq:claim2prop52}). So we will assume that $G_{qr}\not\subset I^{\tilde{m}}$. From this assumption and the fact that the intervals $G_i$ are disjoint we get that if $G_{qr}\not\subset E_1$, then  $G_{qr}\cap E_1=\emptyset$. In this case we apply Corollary \ref{chain space} once and the argument used define $q$ and (\ref{eq:claim2prop52}) to prove the claim.  If $G_{qr}\subset E_1$ let $q'\ge 0$ be minimal so that for all $i$
with $q'\leq i\leq q$, $G_{ir}\subset E_1.$
Note that $I^1$ is well-inside $I^0$ (and therefore $I^2$ is well-inside $I^1$)
and since $\hat m>1$ all critical points of $f^r|I^1$ are contained in $I^2$.
Observe that since $I^1$ is well-inside $I^0$ and by (1) and (2) of Claim 2 that $(f^r| E_1)^{-1}(E_1)$ is well-inside $E_1$.
We can apply Proposition~\ref{monotone pullbacks} to $f^r|E^1$ and therefore we obtain $\lambda_7>0$ 
so that $f^{q'r}(z)\in D_{\lambda_{7}\theta}(E_1\cap I^1)$. If $q'=0,1$, then applying Corollary \ref{chain space} gives (\ref{good at start}).
By the choice of $q'$, \, $G_{ir}\subset I^{\tilde m}$ for all $0\le i\le q'-2$.  If $f^{(q'-2)r}(z)\in D_{\theta}(I^1)$, then we can repeat the above argument choosing an integer  $q<q'-2$ as above. Hence the claim follows.\endpfclaim

\medskip

Since $f^r$ is a composition of folding maps, $f^{r}(I^{\tilde{m}}\cap
E_{1})=f^{r}(I^{\tilde{m}}\cap E_{2})=f^{r}(I^{\tilde{m}})$. 
Let $J_{1}, J_{2}$ be the outermost connected components of
$(f|I^{\tilde m})^{-r}(I^{\tilde m})$ that intersect $E_{1},
E_{2}$ respectively. Note that $c_l$ and $c_r$ are the turning points
in the boundaries of $E_1$ and $E_2$ respectively. 
Let $J_{1,1}$ and $J_{1,2}$ (if it exists) be the two outermost components of
$(f|J_1)^{-r}(J_1)$ such that $f^r|J_{1,j}, j=1,2$ is monotone.
See Figure~\ref{fig:firstreturn}.


\begin{figure}[ht]
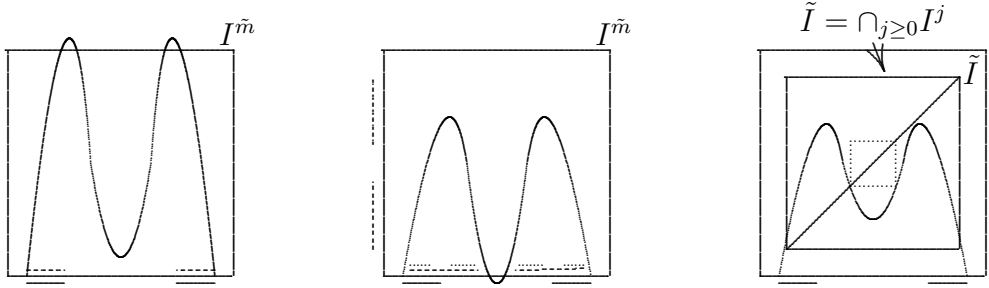

\beginpicture
\dimen0=0.05cm 
\setcoordinatesystem units <\dimen0,\dimen0> point at 80 0 
\setplotarea x from -40 to 30, y from -30 to 30
\setlinear
\setsolid
\put {$I^{\tilde m}$} at 31 35  
\plot -30 -30 30 -30 30 30 -30 30 -30 -30   /
\setquadratic 
\plot -25 -30   -15 32   -8 0    0 -25    8 0    15 32   25 -30 /
\setlinear
\plot  -25 -32 -15 -32 /
\plot 25 -32 15 -32 /
\setdashes 
\setdashes <0.5mm>
\plot  -25 -28.5 -15 -28.5 /
\plot 25 -28.5 15 -28.5 /
\setsolid
\setcoordinatesystem units <\dimen0,\dimen0> point at -20 0
\setplotarea x from -40 to 40, y from -30 to 30
\setlinear
\setsolid
\put {$I^{\tilde m}$} at 31 35 
\plot -30 -30 30 -30 30 30 -30 30 -30 -30   /
\setquadratic 
\plot -25 -30   -15 10   -8 0    0 -32    8 0    15 10   25 -30 /
\setlinear
\plot  -25 -32 -15 -32 /
\plot 25 -32 15 -32 /
\setdashes 
\setdashes <0.5mm>
\plot -23 -28.5 -5 -28.5 /
\plot 23 -28 5 -28.5 /
\plot -33 -23 -33 -5 /
\plot -33 22 -33 5 /
\setdots <0.5mm>
\plot -23 -27 -18 -27 /
\plot -12 -27 -6 -27 /
\plot 23 -27 18 -27 /
\plot 12 -27 6 -27 /
\setcoordinatesystem units <\dimen0,\dimen0> point at -120 0 
\setplotarea x from -40 to 40, y from -30 to 30
\setlinear
\put {$\tilde I$} at 26 25 
\setsolid
\plot -30 -30 30 -30 30 30 -30 30 -30 -30   /
\plot -23 -23 23 -23 23 23 -23 23 -23 -23   /
\plot -23 -23 23 23 /
\setquadratic 
\plot -25 -30   -15 8   -8 0    0 -15    8 0    15 8  25 -30 /
\setlinear
\plot  -25 -32 -15 -32 /
\plot 25 -32 15 -32 /
\setdots <0.6mm>
\plot -6 -6 -6 6 6 6 6 -6 -6 -6 /
\setsolid 
\arrow <10pt> [.2,.67] from 0 32 to 3   24.4
\put {$\tilde I=\cap_{j\ge 0} I^j$} at 0 38 
\endpicture
\caption{\label{fig:firstreturn} Two examples of cases where $\tilde m<\infty$  
and one with $\tilde m=\infty$ and $I^{\tilde m}=\cap_{m\ge 0}I^m$.  The intervals $E_1,E_2$ are marked in solid lines;
the dashed lines refer to the intervals $J_1,J_2$ and the dotted lines are $J_{ij}$.}
\end{figure}

\medskip

\noindent
{\it Claim 4.} 
If  $I^{\tilde m}$ is well-inside $I^{\tilde m-1}$
then  (\ref{good at start}) or (\ref{good deep}) holds for $p'=p-1$
and $I'=I^{\tilde m}$.

\medskip

Let us remark that when we say well-inside, we mean that 
$(1+2\eta)I^{\tilde m}\subset I^{\tilde m-1}$ for some universal
$\eta>0,$
one is welcome to think of $\eta$ as being the constant $\delta>0$
from Lemma~\ref{KSSLemma91}.

\medskip
\noindent
{\it Proof of Claim 4.} 
By
Claim 2 w have that $|I^{\tilde m}|$ is comparable to $|I^0|$, and by
 Claim 3 we have that $f^{pr}(z)\in D_{\hat\lambda\theta}(J)\setminus D_\theta(I^1)$.
Hence by Lemma \ref{Poincare disks} there exists a constant $\lambda_8\in(0,1)$ so 
 that $f^{pr}(z)\in D_{\lambda_8\theta}(I^{\tilde m})$.  
Since $I^{\tilde m}$ is well-inside $I^{\tilde m-1}$, by Corollary \ref{chain space}
we obtain constant $\lambda_9\in(0,1)$, so that $f^{(p-1)r}(z)\in
D_{\lambda_9\theta}(K')$ for some interval $K'$ which is well-inside
$I^{\tilde m}$ by Lemma~\ref{KSSFact91}.
\endpfclaim

\medskip

Let us describe the current situation.
The interval $I^{\tilde m}$ is not well-inside $I^{\tilde m-1}$,
and $|I^{\tilde m}|$ is comparable to $|I^0|$.
Moreover, there exists a universal constant $\kappa_3\in(0,1)$ such that the intervals
$J_{i}, J_{1,i}, i=1,2$ all have length at least
$\kappa_3|I^{\tilde{m}}|,$
the intervals $J_1,$ $J_2$ and $J_{1,2}$ are not necessarily
well-inside $I^{\tilde m}$,
and $J_{1,2}$ is not necessarily well-inside $J_1$.
Moreover, we can assume that the
properties from Claim 3 hold (and in particular that $G_{pr}\subset I^{\tilde m}$).

\medskip

\noindent
{\it Claim 5.} There exists a constant 
$\kappa_4\in (0,1)$ so that each critical point $c'$ of $f^r|I^{\tilde m}$
is contained in $(1+2\kappa_4)^{-1} I^{\tilde m}$ and 
$f^r(J_{1,1})\setminus J_{1,1}$ contains an interval of size $\kappa_4|J_{1,1}|$.

\medskip
\noindent
\textit{Proof of Claim 5.}
The first statement follows from the Claim 4, from Lemma \ref{bounds1} and
since $f^r|I^{\tilde m}$ has no attracting fixed points. The second
statement follows from the first one and Claim 2.
\endpfclaim

\medskip

\noindent
{\it Claim 6.} There exist $\delta_{0}>0$, $\lambda_{0}\in(0,1)$ depending on $\delta$ and for each $z\in U$ there exists an interval $K_{0}$ such that either (\ref{good at start}) holds with $K=K_{0}$ and $\delta'=\delta_{0}$ or there exists $p_{0}<p$ such that
\begin{equation}\label{space deep}
f^{p_{0}r}(z)\in D_{\lambda_{0}\theta}(K_{0})\, \mbox{ and } \, G_{p_{0}r}\subset K_{0}\subset (1+2\delta)K_{0}\subset I^{\tilde{m}}.
\end{equation}

\medskip
\noindent
{\it Proof of Claim 6.} 
If $G_{(p-1)r}\cap J_{1,1}=\emptyset$, then by Claim 5, 
$G_{(p-2)r}$ is well-inside $I^{\tilde m}$ and
by applying Lemma \ref{chain pullback} twice we obtain equation (\ref{space deep}) for $p_0=p-2$
for some interval $K_0\supset G_{(p-2)r}$ which is well-inside $I^{\tilde m}$. 
So let $p_1<p$ be minimal so that
$G_{ir}\cap J_{1,1}\ne \emptyset$ for $i=p_1,\dots,p-1$.
Suppose first that $p-p_{1}\geq 2$. We pull back twice using Corollary \ref{chain space}, to obtain an interval $K_{1}$ which is well-inside $J_{1,1}$.
One component of $f^{r}(J_{1,1})\setminus J_{1,1}$ is a fundamental domain to one side of $K_{1}$ and one component of $I^{0}\setminus I^{1}$ is a fundamental domain to the other side.  Claims 4 and 5 imply that the fundamental domains and the gap between them all have length comparable to $I^{\tilde{m}}$.
So we can 
apply Proposition~\ref{monotone pullbacks} to $f^r| J_{1,1}\to f^r(J_{1,1})\subset I^{\tilde m}$ and thus
obtain an interval $K'$ which is well-inside $J_{1,1}$ and a constant $\lambda_{10}>0$
so that $f^{p_{1}r}(z)\in D_{\lambda_{10}\theta}(K')$. If $p_1=0$, then (\ref{good at start}) holds 
and if $p_1>0$ (including the case $p-p_{1}<2$), then
$G_{(p_1-1)r}\cap J_{1,1}=\emptyset,$ and we argue as in beginning of the proof
of this claim. \endpfclaim

\medskip
\noindent
{\it Claim 7.} If  $\tilde m<\infty$, then there exist $0\leq p'\leq
p_{0},$ $\delta'>0$ and $\lambda'\in(0,1)$ so that
\begin{equation}\label{conc prop}
f^{p'r}(z)\in D_{\lambda'\theta}(K') \mbox{ and } G_{p'r}\subset K'\subset (1+2\delta')K'\subset I^{\tilde{m}},
\end{equation}
where either $p'=0$ or where $p'>0$ and $(1+2\delta')K'\subset I'$,
where $I'$ is a $\hat\rho$-nice interval in $I^{\tilde m}$ which is
equal to $\mathcal{L}_{G_{p'r}}(I^{\tilde m})$ or to $\mathcal{L}_{G_{p'r}}(\mathcal{L}_{G_{(p'-1)r}}(I^{\tilde m}))$.

\medskip
\noindent
{\it Proof of Claim 7.} Let $A_k$ be the component of $f^{-r}(I^{\tilde{m}})$ containing $G_{kr}$. 
As before we can assume that $p_0\ge 1$. Let $p_{1}\le p_0$ be minimal so that for all $i$ with $p_{1}\leq i< p_{0}$,  
$G_{ir}$ is a pullback of $G_{(i+1)r}$ through $J_{1,1}$. If $p_{1}$ is less than 3, then (\ref{good at start}) holds by Corollary \ref{chain space}, so we will assume that $p_{1}>3$.
Because of Claim 5, as in Claim 6, we 
can apply Proposition~\ref{monotone pullbacks} to $f^r\colon J_{1,1}\to f^r(J_{1,1})\subset I^{\tilde m}$.
If $p_1<p_0$, then this implies that there exist  $\delta_1'>0$, $\lambda_1'\in(0,1)$ and an interval $K'_1$ so that
$$f^{p_1 r}(z)\subset D_{\lambda'_1\theta}(K_1'),$$ where $G_{p_1 r}\subset K_1'\subset (1+2\delta_1')K_1'\subset J_{1,1}=A_{p_1}$.
So in any case we get that there exist  $\delta'_2>0$, $\lambda'_2\in(0,1)$ and an interval $K'_2$ so that
$$f^{(p_1-1)r}(z)\subset D_{\lambda'_2\theta}(K_2'),$$ where $G_{p_1 r}\subset K_2'\subset (1+2\delta_2')K_2'\subset A_{p_1-1}$
(here we use that $G_{p_1 r}\subset K'_1\subset (1+2\delta'_1)K'_1$). 
To finish the proof, we show that we can pull back until we arrive in
a return domain to a nice interval that is $\rho$-well-inside that nice
interval,for a universal $\rho$, so that the return domain is 
$\hat{\rho}$-nice, for $\hat \rho>0$, universal, by Corollary~\ref{nice}.

In the remainder of the proof, when we say that an interval $J$
well-inside an interval $I$, we mean that $(1+2\eta)J\subset I$, 
where $\eta$ depends only on the universal constants from Claims 2,3
and 5.
Note that $A_{p_1-1}$ is well-inside $I^{\tilde m}$ 
unless $A_{p_1-1}=J_{1}$ or $J_2$. So unless $A_{p_1-1}=J_{1}$ or $J_2$ the claim follows taking $I'=A_{p_1-1}$.
If $A_{p_{1}-1}=J_{1}$, then we are done if $G_{(p_{1}-1)r}$  is
contained in a landing domain to $J_{1}$ that is well-inside
$J_{1}$. The only way that this does not happen is if
$G_{(p_{1}-1)r}\subset J_{1,2}$. So suppose that this is the case. Under
one more pullback $G_{(p_{1}-2)r}$ is either inside a return domain to
$I^{\tilde{m}}$ that is well-inside $I^{\tilde{m}}$ or inside a return
domain to $J_{1}$ that is well-inside $J_{1}$ or it is contained in
$J_{2}$. We are done except in the last case, but then we do the following.
If $A_{p_1-1}=J_2$, we pull back once more. If $A_{p_{1}-2}$ is not
$J_{1}$ or $J_{2}$ we are done. If $A_{p_{1}-2}=J_{2}$ and $J_{2}$ is
not monotone, then we are well-inside a landing domain to $J_{2}$ that
is well-inside $J_{2}$. Suppose $J_{2}$ is monotone. Pulling back once
more, we are done unless $A_{p_{1}-3}=J_{1}$. In this case we are done
after one more pullback.
So the claim and the proposition follow from Corollary~\ref{chain space}.
\end{pf}

\section{Pulling back a Poincar\'e disc along a chain with bounded combinatorial depth}\label{sec:combinatorial depth}

The results in this section apply to maps $f\in\mathcal A_{\underline
  b}$ and  to asymptotically holomorphic extensions of maps $f\in\mathcal
A^3_{\underline b}$, as in Subsection~\ref{sub-sec:asympholoS}. Once
again, we recall that we always assume that the Standing Assumptions on page
\pageref{standing assumptions} hold.


\medskip
Let $I$ be a nice interval containing the critical point $c_{0},$ and let $m\geq 1$ be minimal so that $R_{I}(c_{0})\notin I^{m}$. Notice that $m\geq\hat{m}$  
(defined on page \pageref{def: m hat}). Recall,
\begin{equation*}
\mathcal{C}(I):= \left\{ 
\begin{array}{rl}
I^{m} & \mathrm{if\ } I\ \mbox{is\ non-terminating\ and}\\ 
I^{\infty}\ &\mathrm{otherwise}.
\end{array}\right.
\end{equation*}

If $J$ is a return domain to an arbitrary nice interval $I$, and $\{G_{i}\}_{i=0}^{r}$ is the chain with $G_{r}=I$ and $G_{0}=J$ where $r$ is the return time of  $J$ to $I$, we define
$$\mathrm{Crit}(I;J)=\Big(\bigcup_{i=0}^{r-1}G_{i}\Big)\cap\mathrm{Crit}(f).$$
Similarly, if $\mathbb{G}=\{G_{j}\}_{j=0}^{s}$ is an arbitrary chain such that the pullbacks of $G_s$ and $I$ are either nested or disjoint, 
$G_{0}\subset I$ and $0=n_{0}<n_{1}<\dots<n_{p}=s$ are the integers with $G_{n_{i}}\subset I$, we define
$$\mathrm{Crit}(I;\mathbb{G})=\bigcup_{i=0}^{p-1}\mathrm{Crit}(I;\mathcal{L}_{G_{n_{i}}}(I)).$$\label{crit symbol}

For any nice interval $I$ and any critical point $c$ we define 
$$k_{c}(I,\mathbb{G})=\inf\{k_{c}\ge 0:\mbox{ there exists no }j=1,2,\dots,s-1\mbox{ with }G_{j} \subset  \mathcal{C}^{k_{c}}(\hat{\mathcal{L}}_{c}(I))
\},$$
where we take  $\mathcal{C}^0(\mathcal{\hat L}_{c}(I))= \hat{\mathcal{L}}_{c}(I)$
and
$$k(I,\mathbb{G})=\sum_{c\in\mathrm{Crit}(I;\mathbb{G})}k_{c}(I,\mathbb{G}).$$
The \emph{combinatorial depth of the chain $\mathbb{G}$ with respect to $I$} is defined to be $k(I;\mathbb{G})$. Note that $k(I,\mathbb{G})$ is well-defined even if $I$ does not contain a critical point. If $J$ and $I$ are nice intervals
with nested or disjoint pullbacks and $c\in J\subset I$, then we define
$$k(I,J)=\min\{k\ge 0:\mathcal{C}^{k}(I)\subset J\}$$\label{comb depth symbols}
and take
$$\hat{k}(I,J)=\sum_{c\in\mathrm{Crit}(f)}k(\hat{\mathcal{L}}_{c'}(I),\hat{\mathcal{L}}_{c'}(J)).$$ \label{comb depth symbols 2}

\noindent
\textbf{Remark.} In the presence of terminating intervals it is possible for the combinatorial depth to be infinite. Let $\{G_{j}\}_{j=0}^s$ be a chain and suppose that $\mathcal{C}_{c}^n(I)$ is terminating.
Then for all $k>n$,
$\mathcal{C}^k_{c}(I)=\mathcal{C}_{c}^{n+1}(I)$, so if some interval $G_{k_0}$ is contained in $\mathcal{C}_{c}^{n+1}(I)$ we get that $k(I,\{G_{j}\}_{j=0}^s)$ is infinite.

\medskip

\begin{prop}[cf. \cite{KSS}, Proposition 11.2]\label{combinatorial depth angle control}
For each $\delta>0, k\geq 0$ and $N\geq 0$ there exist $\mu(k,N,\delta)\in(0,1)$ and $\delta'>0,$ and for each $\theta\in(0,\pi)$ 
there exists $\varepsilon>0$ so that 
the following holds. Let $I$ be a $\delta$-nice interval with $|I|<\varepsilon$.
Suppose that $\mathbb{G}:=\{G_i\}_{i=0}^s$ is a chain
such that $G_0, G_s$ are nice intervals contained in $I$,
the pullbacks of $G_s$ and $I$ are nested or disjoint, the intervals 
$G_0,\dots,G_{s-1}$ are pairwise disjoint
and $G_0\cap\omega(c_0)\neq\emptyset.$
Assume that 
$$k(I,\mathbb{G})\leq k\;\; and\; \;\#\mathrm{Crit}(I;\mathbb{G})\leq N.$$
Let $\hat{G}_{s}$ be an interval with $G_{s}\subset \hat{G}_{s}\subset (1+2\delta)\hat{G}_{s}\subset I$. Let 
$V=D_{\theta}(\hat{G}_{s})\cap\mathbb{C}_{G_{s}}$ and $U_{i}=\mathrm{Comp}_{G_{i}}f^{-(s-i)}(V)$ for $i=0,\dots,s$.
Then, there exists an interval $\hat{I}\supset G_{0}$ with $(1+2\delta')\hat{I}\subset I$ such that 
$$U_{0}\subset D_{\mu(k,N,\delta)\theta}(\hat{I}).$$ 
\end{prop}
\begin{pf}
The proof of this proposition goes by induction on $(N,k)$ with the
lexicographical ordering. If $N=0$, then all the branches are
diffeomorphisms, so the proposition follows because the sum of the lengths of the intervals $\hat G_0,\dots,\hat G_{s-1}$
is uniformly bounded. This can  be seen as follows: let $F\colon
\mathcal J\to I$
be equal to the first return map to $I$
restricted to its  diffeomorphic branches.
Let $J_n$ be a domain of $F^n$ and $J_{n+1}\subset J_n$
a domain of $F^{n+1}$.
Since $F^n\colon J_n\to I$ is a diffeomorphism,
by Theorem~\ref{real Koebe}~(2) $F^n|J_n$ has bounded distortion. Since
each component of $\mathcal J$ is $\delta$-well-inside $I$, 
it follows that there exists $\kappa\in(0,1)$ depending only on $\delta$ such
that $|J_{n+1}|\le \kappa |J_n|$.
It follows that the sum of the lengths
of the intervals $\hat G_i$ contained in  $I$ is universally bounded,
and now Theorem~\ref{real Koebe}~(2) implies that 
the sum of the lengths of  the intervals $\hat G_i$ is universally bounded. 

\medskip
 

Assume now that the statement holds for all $(N',k')$ for which either $N'<N$ or $N'=N$ and $k'<k$. We will prove that the 
statement holds for $(N,k)$.
Let $\mathbb{G}=\{G_{j}\}_{j=0}^{s}$ be a chain as above with $k(I,\mathbb{G})=k$ and 
$\#\mathrm{Crit}(I;\mathbb{G})=N.$ 

\medskip
Without loss of generality we can assume that $I\cap \crit(f)\neq
\emptyset$. If this is not the case we can define $t< s$ to  be
maximal so that $G_{t}\subset \mathcal L_{c}(I),$ for a critical point
$c\in \crit(f)$. 
Suppose that for no $t',$ $t<t'<s$, we have that $G_{t'}\subset I$,
then $G_s$ is the first entry of $G_t$ into $I$, so we can pull back
to time $t$ using Lemma~\ref{chain space}. 
Otherwise, let
$t<t'<s$ be minimal so that $G_{t'}\subset I.$ 
Using the argument for $N=0,$ we pull back from time $s$ to time $t'$
and
using Lemma~\ref{chain space}
to pull back from time $t'$ to time $t$. 
In either case, we obtain
$$U_{t}\subset D_{\lambda_0 \theta}(H_{t})\cap\mathbb{C}_{G_{t}},$$ where
$G_{t}\subset \hat G_t \subset H_{t}\subset (1+2\delta_0) H_t\subset
\mathcal L_{c}(I),$ where $\lambda_0\in(0,1)$ and $\delta_0>0$ depend
on $\delta$.
Since $I$ is $\delta$-nice, Lemma~\ref{KSSFact91} implies that  
after shrinking $\delta_0$, if necessary, $\mathcal L_c(I)$ is
$\delta_0$-nice. 

From now on we will assume that $I$ contains a critical point $c$.
Let $I^1=\mathcal{L}_c(I)$. 
Let $s'<s$ be maximal such that $G_{s'}\subset I$.
Applying Lemma \ref{chain space} we get that $$U_{s'}\subset D_{\lambda_1 \theta}(H_{s'})\cap\mathbb{C}_{G_{s'}},$$ where
$G_{s'}\subset \hat G_{s'} \subset H_{s'}\subset (1+ 2\delta_1)
H_{s'} \subset \mathcal L_{G_{s'}}(I)$, where  
$\lambda_1\in(0,1)$ and $\delta_1>0$ depend on $\delta$.  
Depending on the position of $G_{s'}$ we have two cases.

\medskip
\noindent\textbf{Case 1.}
$G_{s'}\subset  I \setminus I^1$.  Let $J=\mathcal L_{G_{s'}}(I).$
By Corollary~\ref{nice}, there exists $\delta''>0$ so that $J$ is $\delta''$-nice.
Define $s_1<s'$ minimal such that $G_{s_1}\subset J$. If there exists
no such $s_1$, then $G_{s'}$ is the first entry of $G_0$ to $J$ and
the proposition follows from Lemma~\ref{chain space}. 

Let us assume $s_1$ is defined. The structure of the real puzzle and
the fact that $c$ escapes $I^{1}$ imply that
$\mathcal{L}_{c}(J)\subset \mathcal{C}_{c}(I).$ 
It follows easily from this fact and the  definition of $k$ that if any of the intervals $G_{s_1},\dots,G_{s'}$ enters $\mathcal{C}_{c}(I)$, then
$$k(J,\{G_{j}\}_{j=s_1}^{s'})\leq k(I,\{G_{j}\}_{j=0}^{s})-1=k-1,$$ and if
none of the  intervals $G_{s_1},\dots,G_{s'}$  enters
$\mathcal{C}_{c}(I)$, then
$\mathrm{Crit}(J,\{G_{j}\}_{j=s_1}^{s'})<N$. In either 
case, replacing $I$ by $J$ we may apply the induction hypothesis to the chain $\{G_{j}\}_{j=s_{1}}^{s'}$ to obtain
$$U_{s_1}\subset D_{\lambda_2 \theta}(H_{s_1})\cap\mathbb{C}_{G_{s_1}},$$
where $G_{s_{1}}\subset \hat G_{s_{1}}\subset H_{s_1}\subset (1+2\delta_2) H_{s_1}\subset J$,  where
$\lambda_2\in(0,1)$ and $\delta_2>0$ depend on $\delta$, $N$ and $k$.
By the minimality of $s_{1}$ and Lemma~\ref{chain space} we have that
$$U_{0}\subset D_{\lambda_3\theta}(H_{0})\cap\mathbb{C}_{G_{0}},$$
where $G_{0}\subset \hat{G}_{0}\subset H_{0}\subset(1+2\delta_3) H_0\subset I,$
with $\delta_3>0$ and $\lambda_3\in(0,1)$ depending only on $\delta_2$.
This completes the proof in this case.

\medskip
\noindent\textbf{Case 2.} $G_{s'}\subset I^{1}$. There are two
possibilities: either the return map to $I$ restricted to $I^1$ is
monotone or it has a turning point. 

If the return map is monotone, let $0\leq t''\leq s'$ be minimal so
that 
$G_{t''}\subset I^{1}$ and $G_{j}\cap(I\setminus
I^{1})=\emptyset$ for all $j=t'',\dots,s'$. We apply Proposition~\ref{monotone pullbacks}  to find $\lambda''\in(0,1)$ so that $$U_{t''}\subset D_{\lambda''\theta}(I^2).$$
If $t''=0$, we have proved the proposition. If $t''>0$, there exists $0\leq q<t''$ maximal such that $G_q\subset I$. By definition of $q$ we have that $G_{q}\subset I\setminus I^{1}$ and repeating Case 1, we obtain the induction statement.

Assume that the return map to $I$ restricted to $I^1$ is not
monotone. Let $0\leq s''<s'$ be maximal such that $G_{s''}\subset
I$. If 
$s''=0$ the proposition follows from Lemma \ref{chain space},
so suppose $0<s''$. 
If $G_{s''}\subset I\setminus I^1$, then the proposition follows from
Case 1,
so let us assume $G_{s''}\subset I^1$.
Let $r_1\leq s'$ be maximal such that $G_{r_1}\subset \mathcal{ \hat
  L}_{c'}(I)$, where $c'\in \crit(f)$ is a turning point.
Note that $ s'' <r_1\leq s'$. Let $J= \mathcal {L}_{c'}(I)$.
By Lemma~\ref{chain space}, we have that  $$U_{r_1}\subset D_{\lambda'_1 \theta}(H_{r_1})\cap\mathbb{C}_{G_{r_1}},$$ where
$G_{r_1}\subset\hat G_{r_1}\subset H_{r_1}\subset (1+2\delta'_1)
H_{r_1}\subset J,$ and  $\lambda'_1\in(0,1)$ and $\delta'_1>0$ both
depend only on $\delta$. Once again, making $\delta'_1$ smaller, if
necessary, we may assume that $J$ is a $\delta'_1$- nice. Let $q'\leq r_1$ be minimal such that $G_{q}\subset J$ and let  $r_1'<r_1$ be maximal such that $G_{r_1'}\subset J$.
If $r_1'=q'$, the proposition follows from Lemma \ref{chain space},
so we assume $q'<r_1'$. Then we have two possibilities: either
$G_{r_1'}\subset J\setminus J^1$ or $G_{r_1'}\subset J^1$. If the
first holds, then we are in Case 1 and we obtain the induction
statement. 

Assume $G_{r_1'}\subset J^1$. Let $t_1\leq r_1'$ be minimal so that
$G_{t_1}\subset J^{1}$ and $G_{j}\cap(J\setminus J^{1})=\emptyset$ for
all $t_1\leq j\leq r_1'$ with $G_j\subset J.$
Note that if $J\owns c$ is terminating, since $G_{t_1}\cap\omega(c_0)\neq\emptyset$, then $G_{t_1}\subset J^\infty$, so $k(\mathcal{L}_c(I), \mathbb{G})$ is infinite. Hence we have that $J$ is non-terminating.
By Proposition~\ref{good start or good deep} we obtain
$\alpha, \hat\rho>0$, $\lambda_2'\in(0,1)$ and for each $z\in U_{t_1}$ an interval $K$ so that either
\begin{enumerate}[(a)]
\item $f^{t_1}(z)\in D_{\lambda_2'\theta}(K)$ and $G_{t_1}\subset K\subset (1+2\alpha)K\subset J$ or
\item there exist $r_{2}$ with $t_1\leq r_{2}\leq r_{1}'$ and a $\hat{\rho}$-nice interval $J'\subset J^{\hat{m}}$, so that $f^{r_{2}}(z)\in D_{\lambda_2'\theta}(K)$ and $G_{r_{2}}\subset K\subset (1+2\alpha)K\subset J'.$
\end{enumerate}
If (a) holds, then the same argument used for $t''$ in the monotone case proves the proposition.
If (b) holds, then
 by the definition of $\tilde{m}$, there exists
 $c''\in\mathrm{Crit}(J,\mathbb{G})$ which enters $J\setminus J^{1}$
 before it enters $J^{\tilde{m}}\supset J'$.
Let $q'$ be minimal such that for every $q'\leq j \leq r_2$ with $G_j\subset J$ we have that $G_j\subset J^{\tilde m}$. If we take  $q'\leq r\leq r_{2}$ minimal so that $G_{r}\subset J',$ then $c''\notin \mathrm{Crit}(J',\{G_{j}\}_{j=r}^{r_{2}}),$ and therefore
$\#\mathrm{Crit}(J', \{G_{j}\}_{j=r}^{r_{2}})<N$.
Replacing $I$ with $J'$, and applying the induction hypothesis to the
chain
$\{G_{j}\}_{j=r}^{r_{2}}$, we have that $f^{r}(z)\in
D_{\lambda'_3\theta}(\hat{J}')$, where $\hat{J}'$ is  an interval such
that $G_{r}\subset\hat{J}'\subset(1+2\alpha')\hat{J}'\subset J'$ and
$\lambda_3'\in(0,1)$ depend only on $\alpha$, $N$ and $k$.
If $G_{0},\dots, G_{r-1}$ avoid $J'$, then Lemma~\ref{chain space} completes the proof of the proposition. Otherwise,
there exists $r_{3}<q'$ so that $G_{r_{3}}\subset J\setminus J^{1}$.
Let $r_{3}$ be maximal with respect to this property. Then $G_{r_{3}},\dots,G_{r-1}$ avoid $J'$ and so
applying Lemma~\ref{chain space} we obtain $f^{r_{3}}(z)\in D_{\lambda'' \theta}(\hat{J}'')$
for some interval $\hat{J}''$ satisfying $G_{r_{3}}\subset\hat{J}''\subset(1+2\alpha'')\hat{J}''\subset J$
where $\alpha''>0$ depends on $\alpha'$, and $\lambda''\in(0,1)$
depends on $\alpha', N$ and $k$. Applying Case 1, replacing $s$ with $r_{3}$ completes the proof in this case.
\end{pf}

\section{Pulling back a Poincar\'e disc  one step along the enhanced nest}\label{sec:pullback along enhanced nest}

The results in this section apply to maps $f\in\mathcal A_{\underline
  b}$ and  to asymptotic holomorphic extensions of maps $f\in\mathcal
A^3_{\underline b}$, as in Subsection~\ref{sub-sec:asympholoS}. We
recall that we always assume that the Standing Assumptions on page
\pageref{standing assumptions} hold.

First, we state and prove a lemma which relates a bound on the scaling factors between two levelswith some combinatorial information. 

\begin{lem}\label{lem61}
For each $\nu\geq 1$, there exists $K\in\mathbb{N}$ with the following properties. 
Assume that $I_i$ and $I_{i+1}$ are non-terminating and that both
$$|I_{i}|/|I_{i+1}|\;and\;|I_{i+1}|/|I_{i+2}|\leq\nu . $$
Then
$$k(I_{i},I_{i+1})\;and\;
k(\mathcal{L}_{c}(I_{i}),\mathcal{L}_{c}(I_{i+1}))\le K \mbox{ for
  each } c\in \crit(f)\cap \omega(c_0).$$
\end{lem}
\begin{pf}
If $k(I_{i},I_{i+1})$ is large, then by Lemma \ref{KSSLemma91},  $|I_{i}|/|I_{i+1}|$ is large. If $c$ is a critical point and $k(\mathcal{L}_{c}(I_{i}),\mathcal{L}_{c}(I_{i+1}))$ is large, then by Lemma \ref{KSSLemma96} $|\mathcal{L}_{c_{0}}\mathcal{L}_{c}(I_{i})|/|\mathcal{L}_{c_{0}}\mathcal{L}_{c}(I_{i+1})|$ is large too. By construction, we have that
$$I_{i}\supset\mathcal{L}_{c_{0}}\mathcal{L}_{c}(I_{i})\supset\Gamma(I_{i})\supset I_{i+1}$$
and
$$I_{i+1}\supset\mathcal{L}_{c_{0}}\mathcal{L}_{c}(I_{i+1})\supset\Gamma(I_{i+1})\supset I_{i+2}.$$
So, if $|\mathcal{L}_{c_{0}}\mathcal{L}_{c}(I_{i})|/|\mathcal{L}_{c_{0}}\mathcal{L}_{c}(I_{i+1})|$ is large, 
and since $I_i\supset \mathcal{L}_{c_{0}}\mathcal{L}_{c}(I_{i}) \supset I_{i+1}\supset 
\mathcal{L}_{c_{0}}\mathcal{L}_{c}(I_{i+1})\supset I_{i+2}$, 
either $|I_{i}|/|I_{i+1}|$ or $|I_{i+1}|/|I_{i+2}|$ is large, which contradicts our assumption. 
\end{pf}

\subsection{The renormalizable case}\label{subsec:renorms}
In this subsection, we deal with pullbacks along consecutive intervals
of the enhanced nest in the presence of terminating intervals.

\begin{figure}[H] \centering \def\svgwidth{370pt}
 \scriptsize
 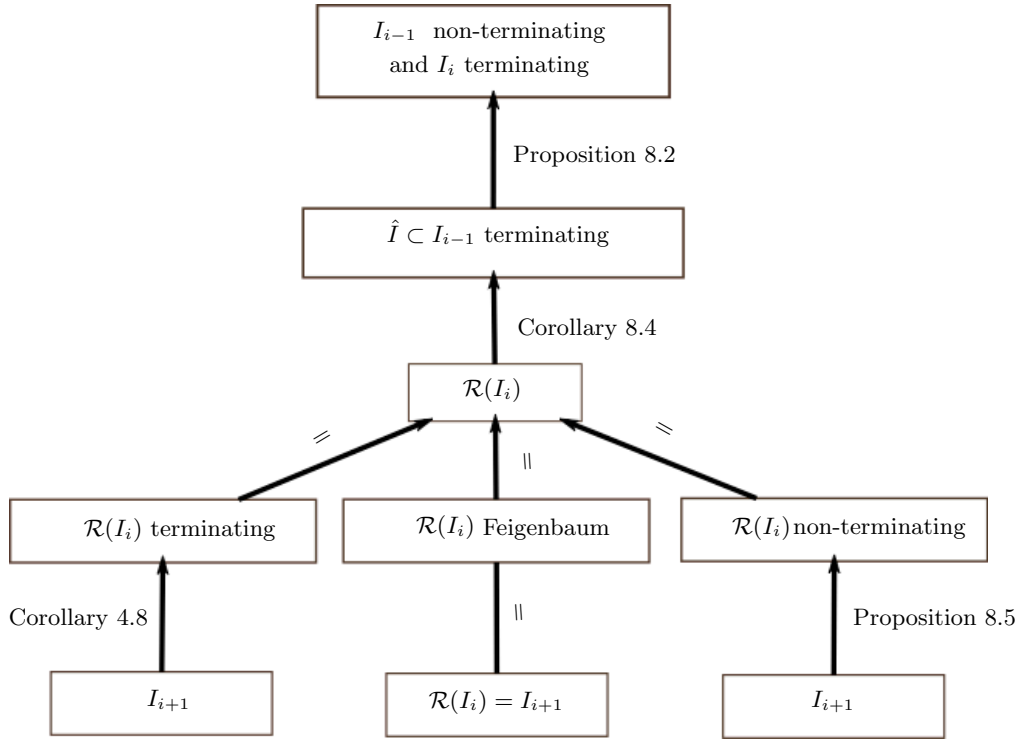
 \caption{\small Pulling back a Poincar\'e disk based on $I_{i-1}$ to
   one based on $I_{i+1}$ when $I_{i-1}$ is non-terminating and $I_i$
   is terminating. We start at the top of the diagram and the arrows
   indicate pulling back. The $=$ symbol indicates that there is
   nothing to do.\label{diagram1}}
 \end{figure}

\begin{figure}[H] \centering \def\svgwidth{350pt}
  \scriptsize
 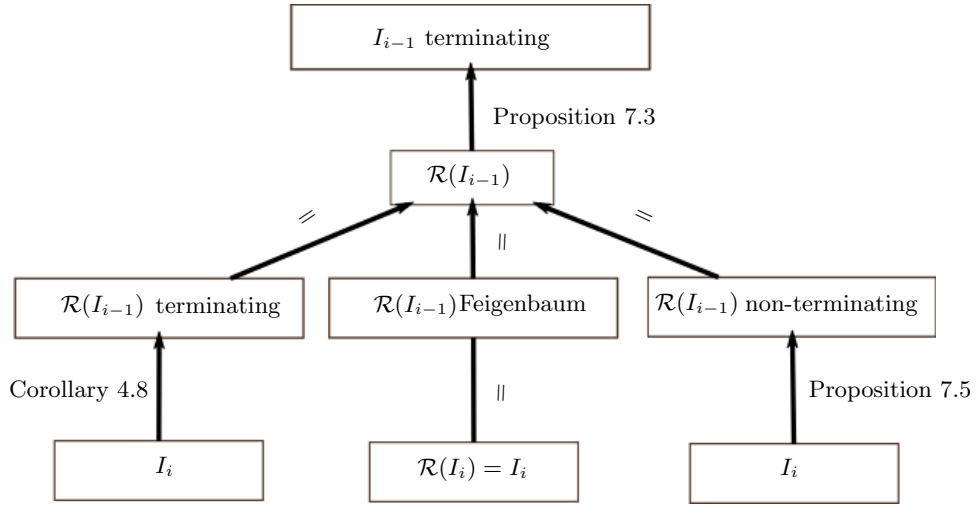
 \caption{\small Pulling back a Poincar\'e disk based on $I_{i-1}$ to
   one based on $I_{i}$ when $I_{i-1}$ is terminating. We start at the
   top of the diagram and the arrows indicate pulling back.
The $=$ symbol indicates that there is
   nothing to do.\label{diagram2}}
 \end{figure}

Figures \ref{diagram1} and \ref{diagram2} show how these pullbacks are structured.

In the next proposition we will deal with the first step of the pullback when $I_{i-1}$ is non-terminating and $I_{i}$ is terminating, see Figure \ref{diagram1}.

\begin{prop}\label{folding lemma}
For any $\nu\geq 1$ there exists $\lambda\in(0,1),$ 
and for each $\theta\in(0,\pi)$ there exists $\varepsilon>0$ 
such that the following holds. 
 Suppose that $I_{i-1}$ is non-terminating, $I_{i}$ is terminating, $|I_{i-1}|<\varepsilon$ and
$$\frac{|I_{i-2}|}{|I_{i+1}|}<\nu.$$
Let $r$ be such that $f^{r}|_{I_i^{\infty}}=R_{I_{i}^\infty}.$
Let $\{G_{j}\}_{j=0}^{s}$ be a disjoint chain with 
$G_{0}\subset I_i^\infty$ and $G_{0}\cap\omega(c_0)\neq\emptyset$. Assume that $G_{s_{i-1}}\subset I_{i-1}$ for some $s_{i-1}> 2r.$ Then there exists $s_{i}'<s_{i-1}$ with
$G_{s_i'}\subset I_i^\infty$ such that if
$$U_{s_{i-1}}=D_{\theta}(I_{i-1})\cap\mathbb{C}_{G_{s_{i-1}}}\ \mathrm{and}\ U_{0}=\mathrm{Comp}_{G_{0}}f^{-s_{i-1}}(U_{s_{i-1}}),$$ 
then there is an interval $K$ which is well-inside  $\hat I$ such that $$U_{s_{i}'}\subset D_{\lambda\theta}(K),$$
where $\hat I\supset I^{\infty}_{i-1}$ is the largest terminating interval in the principal nest 
$$I_{i-1}=I_{i-1}^0\supset I_{i-1}^1\supset \dots.$$ 
Moreover, there are at most two integers $j$ with $ s_{i}' \leq j\leq s_{i-1}$ and $G_{j}\cap I_{i}\neq\emptyset.$
\end{prop}

\medskip
\noindent\textbf{Remark.} The time $s_i'$ is chosen so that the base
of the Poincar\'e domain containing $U_{s_{i}'}$ is contained in $\hat
I$. Pulling back to time $s_{i}'$ is the first step needed to
eventually get a Poincar\'e disk based on $I_{i+1}$ (see  Figure
\ref{diagram1}).

\begin{pf}
Let $I_{i-1}=I^{0}\supset I^1\supset I^2\supset I^3\supset\dots$ be the principal nest about $c_0$. Let $m(0)=0$ and let $m(0)<m(1)<m(2)< m(3)\dots< m(l)$ be the levels the of principal nest such that $R_{I^{m(i)-1}}(c_0)\notin I^{m(i)},$ $i=1,2,\dots, l$. 
Note that $l<\infty$ since $I_{i}$ is terminating. By definition $I^{m(l)}$ is terminating,  so $I^{m(l)}\supset I_{i}^\infty$ and the return time
of $I^{m(l)+1}$ to $I^{m(l)}$ is $r$.
We will begin by defining a time $j_0$ 
from which we will be able to
pull back to a time $s_i',$ as in the statement of the proposition, along a bounded
number of monotone branches in the principal nest about $c_0$ and
first return maps.
We will first pick a sequence of times $k_i \leq k_i'$ 
that give us good control of the orbit of
$c_0$ in the principal nest, $I^{0}\supset I^1\supset I^2\supset\dots\supset I^{m(l)}.$

Let $k_l< s_{i-1}$ be maximal such that $G_{k_l}\subset I_{i}^{\infty}$.
Let $$k_l< k_{l-1}\leq  k_{l-1}'< k_{l-2}\leq \dots < k_{0}\leq k_0' < k^*\leq j_0\leq s_{i-1}$$ be defined as follows.
See Figure \ref{fig:t_0} for an illustration of this sequence.
 Let $k_{l-1}>k_{l}$ be minimal such that $f^{k_{l-1}}(c_0)\in I^{m(l)-1}$.
 Since $G_{k_l}\subset I_{i}^\infty$, $k_{l-1}$ coincides with the 
 minimal time $> k_l$ such that $G_{k_{l-1}}\subset  I^{m(l)-1}$
 and $f^{k_{l-1}}(c_0)\notin I^{m(l)}$. 
Assuming that $k_{j}$ is defined, let $k'_{j}$ be maximal such that for all $k$, $k_{j} \leq  k \leq k_{j}'$, if $f^{k}(c_0)\in I^{m(j)}$, then $f^k(c_0)\in \mathcal{L}_{c_0}(I^{m(j)})$, and define 
$k_{j-1}=k_{j}'+n_{j}$ where $n_{j}$ is the return time of 
$\mathcal{L}_{f^{k'_{j}}(c_0)}(I^{m(j)})=\mathcal{L}_{c_0}(I^{m(j)})$ to $I^{m(j)-1}$. Note that the return of $f^{k'_{j}}(c_0)$ to $I^{m(j)-1}$ is non-central, i.e. $f^{k_{j-1}(c_0)}$ is not contained in $I^{m(j)}$. 
With $k_0'$ defined, we set $k_*=k_{0}'+q_0$, where $q_0$ is the return time
of  $\mathcal{L}_{f^{k'_{0}}(c_0)}(I^{m(0)})=I^{m(0)+1}$ to $I^{m(0)}$.
Notice that $f^{k_*}(c_0)\notin I^{m(0)+1}$.
 Let $j_0$, $s_{i-1}-r<j_0 \leq s_{i-1}$ be such that $G_{j_0}\subset \mathcal{L}_{f^{k_*}(c_0)}(I_{i}^\infty)$;  remember $r$ is the period of $I_i^\infty$.  Since $s_{i-1}>2r$, we know that $j_0\geq k_*$. 

The pullback of $G_{j_0}$ will follow the same path as the pullback of $f^{k_*}(c_0)$ 
since they are both contained in the orbit of $I_{i-1}^\infty$. We will use this property below Claim 2, see 
Figure~ \ref{fig:to m-hat}.

 \begin{figure}[htb!] \centering \def\svgwidth{400pt}
 \scriptsize
 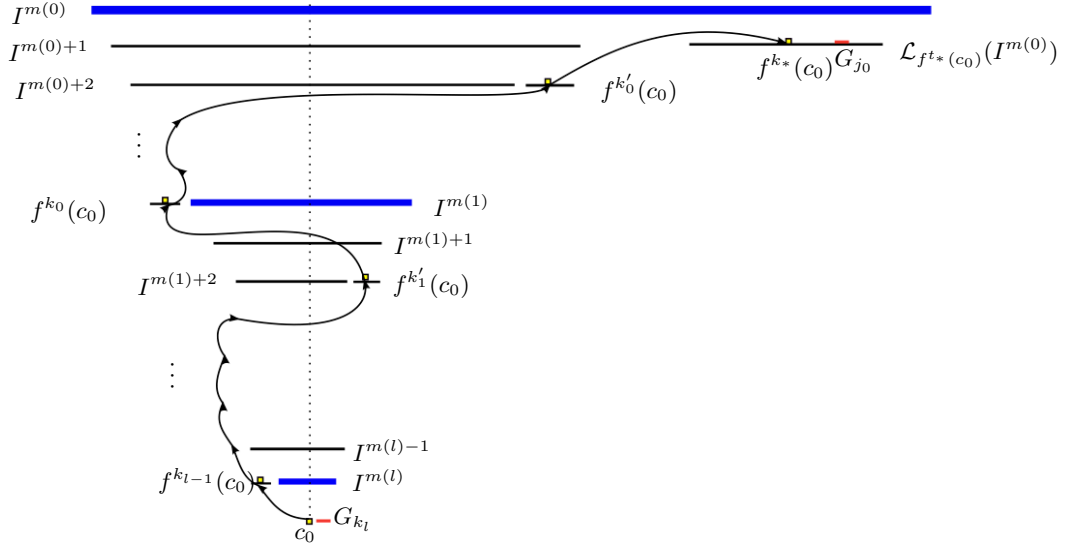
 \caption{Finding $k_*$.\label{fig:t_0}}
 \end{figure}

\medskip
 \noindent\textbf{Combinatorial Remark 1.} Note that by the definition of $j_0$ there exists at most one interval $G_j$ with $j_0\leq j\leq s_{i-1}$ with $G_j\subset I_i^\infty.$
 
\medskip

We start by pulling back along the first part of the chain, from time $s_{i-1}$ to time $j_0$. The loss of angle along this segment of the chain will be controlled by Proposition~\ref{combinatorial depth angle control}, so we need find a way to split the chain into smaller segments, each of them with bounded combinatorial depth. The argument is complicated by the fact that the chain may enter a terminating component of the landing map to $I_i$ that contains a critical point and therefore has infinite combinatorial depth. Let us study the combinatorial depth of the chain $\{G_j\}_{j=j_0}^{s_{i-1}}$ with respect to the interval $I_{i-1}$. By construction, there exists at most one interval $G_j$  with  $j_0< j\leq s_{i-1}$ and $G_j\subset I_{i-1}^{\infty}$. Since $I_i\subset I^\infty_{i-1}$, for each $c\in\crit(I_{i-1}; \{G_j\}_{j=j_0}^{s_{i-1}}),$ 
there exists at most one interval $G_{j}$ with $j_0< j<s_{i-1}$ such that $G_j\subset\hat{\mathcal L}_c(I_{i})$. Let $t_1$ be maximal so that $j_0< t_1<s_{i-1}$ and $\hat{\mathcal{L}}_{G_{t_1}}(I_{i})$ contains a critical point. 
If $t_1$ is not defined,let  $t_1=j_0.$ Let $t_1'> t_1$ be minimal such that $G_{t_1'}\subset I_{i-1}$ and let $\mathbb G= \{G_j\}_{j=t_1'}^{s_{i-1}}.$ There exist two possibilities either
\begin{itemize}
\item[\textit{(i)}] $k(I_{i-1}, \mathbb G)$ is finite or
\item[\textit{(ii)}] $k(I_{i-1}, \mathbb G)= \infty$.
\end{itemize}
First, we estimate the combinatorial depth of the chain $\mathbb G$ if $(i)$ holds.

\medskip

\noindent\textit{Claim 1.}  \textit{If $k(I_{i-1},\mathbb G)$ is finite there exists $a=a(\nu)>0$ such that $k(I_{i-1},\mathbb G)< a$.}

Assume $k(I_{i-1}, \mathbb G)>a$, where $a$ will be chosen later. By the definition of combinatorial depth, there exists $c'\in\crit(I_{i-1}; \mathbb G)$ such that $k(I_{i-1}, \{G_j\}_{j=n_{c'}}^{m_{c'}})>a/b$. Let $a_*$ be the integer part of  $a/b$ and let $c'\in\crit(I_{i-1}; \mathbb G)$ be so that $m_{c'}$ is maximal and $k(I_{i-1}, \{G_j\}_{j=n_{c'}}^{m_{c'}})>a_*.$ By the definition of combinatorial depth,  there exists $G_j\subset \mathcal C^{a_*}(\hat{ \mathcal L}_{c'}(I_{i-1}))$. Assume $j'$ is maximal with $t_{1}'\leq j'\leq s_{i-1}$ and   $G_{j'}\subset \mathcal C^{a_*}(\hat{ \mathcal L}_{c'}(I_{i-1}))$; since $j'>t_1$, $\mathcal{L}_{c'}(I_i)\subset \mathcal{C}^{a_*}(\hat{ \mathcal L}_{c'}(I_{i-1}))$. Since $k(I_{i-1},\mathbb G)$ is finite none of the intervals  $ \mathcal C^n(\hat{\mathcal L}_{c'}(I_{i-1}))$ is terminating for $1\leq n\leq a_* $. By Lemma \ref{KSSLemma91} there exists $\rho>0$ such that for any $1<n<a_*$
$$(1+2\rho) \mathcal C^{n+1}(\hat{\mathcal L}_{c'}(I_{i-1}))\subset \mathcal C^n(\hat{\mathcal L}_{c'}(I_{i-1})).$$
If $c'=c_0,$ and $a$ is sufficiently big, so $a_*$ is sufficiently big, this would imply $|I_{i-1}|/  |I_{i+1}|>\nu$. So let us assume $c'\neq c_0$. For each $n<a_*$ let $p_n$ be minimal such that $f^{p_n}(c_0)\subset  \mathcal C^{n}(\hat{\mathcal L}_{c'}(I_{i-1}))$; ${p_n}$ is the first entry time of $c_0$ to $ \mathcal C^{n}(\hat{\mathcal L}_{c'}(I_{i-1}))$.  Let $C_n^n$ and $C_{n+1}^n$  be the pullbacks of $ \mathcal C^{n}(\hat{\mathcal L}_{c'}(I_{i-1}))$ and $ \mathcal C^{n+1}(\hat{\mathcal L}_{c'}(I_{i-1}))$  containing $c_0$ along the orbit $c_0, f(c_0),\dots, f^{p_n}(c_0)$, respectively. By Theorem \ref{real Koebe} (1) there exists $\rho'>0$ such that for any $n<a_*$
$$(1+2\rho')C^n_{n+1}\subset C^n_n.$$
Note that  $a_n\leq a_{n+1}$ so  $$C^{n+1}_{n+1}\subset  C^n_{n+1} \subset C^n_n.$$
If $a$ is big enough, so  $a_*$ is big enough, this implies that $\Gamma(I_{i})\subset \LL_{c_0}(\mathcal{L}_{c'}(I_i))$ is deep-inside $I_{i-1}$, which contradicts $|I_{i-1}|/|I_{i+1}|<\nu$. \endpfclaim
\medskip

\noindent \textbf{Remark.} From the proof of Claim 1 we get the
following. Assume that none of the intervals $ \mathcal
C^n(\hat{\mathcal L}_{c'}(I_{i-1}))$ is terminating for $1\leq n\leq
n_0$ and that there exists $G_j\in \mathbb G$ with $G_{j}\subset
\mathcal C^{n_0}(\hat{\mathcal L}_{c'}(I_{i-1}))$. Then  $n_0 <a$. Now
we pull back from time $s_{i-1}$ to time $j_0$.

\medskip

\noindent\textit{Claim 2. There exist $\lambda_0\in(0,1)$ and an interval $K_0\supset G_{j_0}$ that is well-inside $I^{m(0)}=I_{i-1}$, so that
$$U_{j_0}\subset D_{\lambda_0\theta}(K_0).$$} 
By Theorem \ref{real geometry}(a) we know that there exists $\rho>0$
so that $I_{i-1}$ is $\rho$-nice. Since $|I_{i-2}|/|I_{i+1}|<\nu,$
Corollary \ref{cor:delta free} implies there exists $\delta>0$ such that $I_{i-1}$ is $\delta$-free. Shrinking $\delta$, if necessary, we can assume that $I_{i-1}$ is $\delta$-nice and $\delta$-free. Let $ t<s_{i-1}$ be maximal so that $G_{t}\subset I_{i-1}$. Since $I_{i-1}$ is $\delta$-free we can apply Lemma \ref{chain pullback} to find  $\lambda'>0$ so that 
$$U_t\subset D_{\lambda'\theta}(\LL_{G_{t}}(I_{i-1})),$$
with $(1+2\delta)\LL_{G_{t}}(I_{i-1})\subset I_{i-1}$. There are two cases to consider:

\medskip
\noindent\textit{Case (i): $k(I_{i-1}, \mathbb G)$ is finite.}

Since we are assuming $|I_{i-1}|/|I_{i+1}|<\nu$, Claim 1 tells us that $k(I_{i-1}, \mathbb G)< a$. In this case Proposition~\ref{combinatorial depth angle control} allows us to control the loss of angle from time $t$ to time $t_1'.$ Since $G_{t_1'}$ is the first entry of $G_{t_1}$ to $I_{i-1}$ we can control the loss of angle from time $t_1'$ to time $t_1$ using Corollary \ref{chain space}; this means that we can find an interval $J_1$  and constants $\lambda_1,\delta_1>0$ such that  $$U_{t_1}\subset D_{\lambda_1\theta} (J_1)$$
with $G_{t_1}\subset J_1\subset (1+2\delta_1)J_1\subset \hat{\mathcal
  L}_{G_{t_1}}(I_{i-1})$. Let $j_0\leq t_2< t_1$ be maximal so that
$\hat{\mathcal{L}}_{G_{t_2}}(I_{i})$ contains a critical point. If
$t_2$ is not defined, the proposition follows after one application of
Corollary \ref{chain space}. If $t_2$ is defined let $t_2'$ with
$t_2<t_2'\leq t_1$ be minimal such that $G_{t_2'}\subset I_{i-1}$. If
we let $\mathbb G= \{G_j\}_{j=t_2'}^{t_1},$ either Case (i) or Case (ii)
holds. If Case (i) holds, we repeat the previous argument.
See below for the argument in Case (ii).
The proposition follows after repeating this argument
at most
once for each point
$c\in\crit(I_{i-1}; \{G_j\}_{j=j_0}^{t_1-1}),$ so at most $b-1$ times.

\medskip

\noindent\textit{Case (ii), $k(I_{i-1},\mathbb G)= \infty$.} 
By the definition of combinatorial depth,
there exist $n\in \N$, $t_1'\leq j\leq s_{i-1}$ and a point
$c\in\crit(I_{i-1}; \mathbb G),$  so that $G_j\subset \mathcal
C^n(\hat{\mathcal L}_{c}(I_{i-1}))$, where $ \mathcal
C^n(\hat{\mathcal L}_{c}(I_{i-1}))$ is a periodic interval. For each
$c\in\crit(I_{i-1}; \mathbb G)$, pick $n_c$ minimal and $m_c$ maximal
so that $t_1'\leq n_{c}\leq m_c \leq s_{i-1}$ and $G_{n_c},
G_{m_c}\subset\hat{ \mathcal L}_c(I_{i-1})$. 
Let $c\in \crit(I_{i-1};\mathbb G)$ be such that $m_c$ is maximal.
From the remark at the end of Claim 1 we know the following. There exists $m< a$ such that $C^{m}(\hat{\mathcal L}_{c}(I_{i-1}))$ is periodic, otherwise $|I_{i_1}|/|I_{i+1}|>\nu$. Let $m_0$ be minimal such that $C^{m_0}(\hat{\mathcal L}_{c}(I_{i-1}))$ is a periodic interval.  Let  $\mathcal P=\mathcal C^{m_0}(\hat{\mathcal L}_{c}(I_{i-1}))$. By the choice of $c$ there exists $G_j\subset \mathcal P$. Even more, since $\mathcal P$ is a periodic interval and 
$G_{k_l}$ is the only element of the chain contained in $I_{i}$, there exists only one interval
$G_{\hat{t}_1}\subset\mathcal P$ with $j_0\leq \hat t_1<s_{i-1}$. By
the definition of $\hat t_1$ and the choice of $c$ we can use
Proposition~\ref{combinatorial depth angle control} and Corollary
\ref{chain space} to control the loss of angle from time $s_{i-1}$ to
time $\hat t_1$. Thus, there exists an interval $J_1$ and constants $\lambda_1,\delta_1>0$ with $G_{\hat t_1}\subset J_1\subset (1+2\delta_1)J_1\subset \hat{\mathcal L}_{c}(I_{i-1})$ such that  $$U_{\hat t_1}\subset D_{\lambda_1\theta} (J_1).$$
Pick $j_0\leq t_2<\hat t_1$ maximal such that $G_{t_2}\subset \mathcal{C}^n(\mathcal L_{c'}(I_{i-1}))$ for some $c'\in \crit(I_{i-1}; \{G_j\}_{j=j_0}^{\hat t_1})$ and $\mathcal{C}^n(\mathcal L_{c'}(I_{i-1}))$ terminating; if $t_2$ is not defined let $t_2=j_0$. Let $\hat t_2<\hat t_1$ be maximal such that $G_{\hat t_2}\subset \hat{\mathcal L}_{c'}(I_{i-1})$. Proposition~\ref{combinatorial depth angle control} and Corollary \ref{chain space} allow us to control the loss of angle from time $\hat t_1$ to time $t_2$. The proof follows repeating the previous argument at most once for each critical point of $f$, so at most $b$ times .
\endpfclaim

\medskip
Let us continue with the proof of the proposition. Now that we have
pulled back to  time $j_0$ we pull back to the desired time $s_i'=j_0-k_*$. 
By the choice of $j_0$ we can restrict all but a bounded number of these pullbacks
to monotone branches of return maps to intervals in the principal nest.
See Figure \ref{fig:to m-hat} for an idea of how this will be done. Let us first define the following times; let $\tilde m(1)<\tilde m(2)<\dots<\tilde m(k)$ be the levels of the principal nest such that there exists a turning point $c$ of $R_{I^{\tilde m(i)-1}}$ such that $R_{I^{\tilde m(i)-1}}(c)\notin I^{\tilde m(i)},$ $i=1,2,\dots k$. Set $\tilde m(0)=0$, so that $I^{\tilde m(0)}=I^0=I_{i-1}$. Notice that the terminating interval $\tilde I$ is equal to $I^{\tilde m(k)}$ and that the return time of $c_0$ to $\tilde I$ is $r$.
Let $j_0'<j_0$ be minimal so that for all $t$, $j_0'\leq t< j_0$, we have that if $G_t\subset I^{\tilde m(0)}$, then $G_{t}\subset I^{\tilde m(0)+1}\setminus I^{\tilde m(1)}$. Define $j_1<j_0'$ to be maximal so that $G_{j_1}\subset I^{\tilde m(1)}$. Similarly, define $j_1'>j_2>j_2'>\dots j_k$.
 \begin{figure}[htb] \centering \def\svgwidth{300pt}
 \scriptsize
 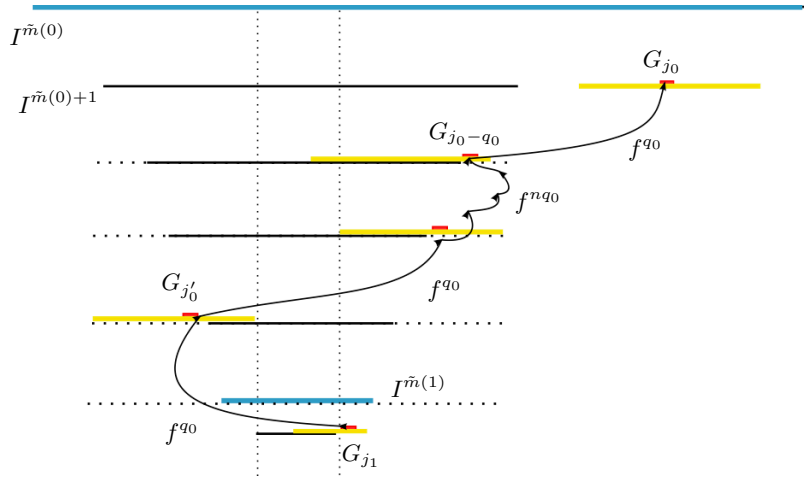
 \caption{Pulling back from time $j_0$ to $j_1$. \label{fig:to m-hat}}
 \end{figure}
If $\tilde m(1)=\tilde m(0)+1$ or $\tilde m (0)+2$, then by applying Corollary \ref{chain space} once or twice there exist $\lambda_1\in(0,1)$, and an interval $K_1$ that is well-inside $I^{\tilde m(1)}$ such that 
$$U_{j_1}\subset D_{\lambda_1\theta}(K_1),$$
so suppose that $\tilde m(1)>\tilde m(0)+2$.
Let $q_0$ be such that $f^{q_0}|I^{\tilde m(0)+1}=R_{I^{\tilde m(0)}}|I^{\tilde m(0)+1}$.
Since $|I^{\tilde m(0)+1}|/|I^{\tilde m(0)+2}|<\nu$, and $I^{\tilde m(0)+1}$ is well-inside $I^{\tilde m(0)}$, Lemma \ref{bounds1} implies that there exists $\kappa_0>0$ and $C$ such that for any critical point $c$ of $R_{I^{\tilde m(0)+1}}|I^{\tilde m(0)+2}$,
$$|f^{q_0}(c)-c|\geq\kappa_0 |I^{\tilde m(0)}|,\; \text{and} \; |D(f^{q_0})(x)|\leq C\ \text{for\ all\ }x\in I^{\tilde m(0)+2}.$$
Since $f^{q_0}$ is a first return map, by Corollary \ref{chain space}, there exist $\lambda'\in(0,1)$ and an interval $K'$ which is well-inside $I^{\tilde m(0)+1}$ such that
$$U_{j_0-q_0}\subset  D_{\lambda'\theta}(K').$$
We have that $G_{j_0-2q_0}\subset I^{\tilde m(0)+2}\setminus I^{\tilde
  m(1)}$. Let $c_1$ be the turning point of $f^{q_0}$ on the boundary
of the monotone branch of $f^{q_0}$ containing $G_{j_0-2q_0}$. Let
$\mathcal{Y}$ denote this monotone branch. We will now pull back to
get the base of the Poincar\'e disk inside of $\mathcal{Y}$. 
Since $c_1\in I^{\tilde m(1)}$, $G_{j_0-q_0}$ does not contain $f^{q_0}(c_1)$. Let $c'$ be the critical point of $f$ with the property that $f^k(c_1)=c'$ with $0\leq k< q_0$ maximal. Let $\alpha$ be the fixed point of $f^{q_0}$ closest to $c_0$ and let  $\alpha'\in \comp_{c'}f^{q_0-k}(I^{\tilde m(0)+1})$ be the preimage of $\alpha$ that lies on the opposite side of $c'$ to  $G_{j_0-2q_0+k}$. By Lemma \ref{bounds1}, $|c_1-\alpha|$ is comparable to $|I^{\tilde m(0)}|,$ so there exists an interval $J_0\subset (c_1,\alpha)$ such that $|J_0|>\delta_0|I^{\tilde m(0)}|$ and $J_0$ does not contain any critical values of $f^{q_0}$. Let $J'_0$ be the pullback of $J_0$ in $(c',\alpha')$,
and let $J_1'$ be the pullback by $f^{-(q_0-k)}$ of the connected component of  $I^{m(0)+1}\setminus I^{m(0)+2}$ on the same side of $c_1$ as $G_{j_0-q_0}$, so that
the convex hull of $J_0'\cup J_1'$ contains $G_{j_0-2q_0+k}$.
There exists $\delta'_0>0$ such that $|J_{0}'|, |J_1'|>\delta'_0|\mathcal{L}_{c'}(I^{\tilde m(0)+1})|$.
Let $K''$ be the convex hull of $J_0'\cup J_1'$. There exists $\delta''>0$ and an interval $K_0'$ which is well-inside $K''$ such that $|K_0'|$ is comparable to $|\LL_{c'}(I^{m(0)})|$ and  $(1-2\delta'')K_0'\setminus(1-2\delta'')^{-1}K_0'$ contains no critical values of $f^{k}$.
Since $G_{j_0-2q_0+k}$ does not contain $c'$
by Lemma \ref{different branch angle control}, there exists $\lambda'_0\in(0,1)$
$$U_{j_0-2q_0+k}\subset D_{\lambda_0'\theta}(K_0').$$
Since $K_0'$ is $\delta''$-free we can apply by Lemma \ref{chain pullback} to find $\lambda_1\in(0,1)$ and an interval $K_1$ that is well-inside $\mathcal{Y}$, 
such that $G_{j_0-2q_0}\subset K_1$ and
$$U_{j_0-2q_0}\subset D_{\lambda_1\theta}(K_1).$$
Now, we use an argument similar to the one in the proof of Claim 6 of
Proposition~\ref{good start or good deep} to pull back to a Poincar\'e disk based in $I^{\tilde m(1)}$.
By Proposition~\ref{monotone pullbacks} and at most one application of Corollary \ref{chain space}, we have that there exists $\lambda_2\in(0,1)$ such that 
$$U_{j_0'}\subset D_{\lambda_2\theta}(\mathcal Y'),$$
where $\mathcal{Y}'$ is an outermost monotone branch of $f^{q_0}|I^{\tilde m(0)+1}$.
So we
have that there exists $\lambda_3\in(0,1)$ and an interval $K_3$ that is well-inside $I^{\tilde m(1)}$ such that 
$$U_{j_1}\subset  D_{\lambda_3}(K_3).$$
Let $q_1$ be the return time of $I^{\tilde m(1)+1}$ to $I^{\tilde m(1)}$. By Corollary \ref{chain space}, we have that there exists $\lambda_4\in (0,1)$ and an interval $K_4$ which is well-inside $I^{\tilde m(1)+1}$   such that 
$$U_{j_1-q_1}\subset  D_{\lambda_4\theta}(K_4).$$
Observe that $G_{j_1-q_1}\subset I^{m(1)+1}\setminus I^{m(1)+2}$.

\medskip

Since $I^{m(0)}=I_{i-1}$ is $\delta$-nice, there exists $\delta_1>0$ such that $I^{m(1)+1}$  is $\delta_1$-nice. 
If $j\in\{0,\dots, k-1\}$ is such that there exists $j'\in\{0,\dots,l-1\}$ such that $\tilde m(j)=m(j')$, then by Lemma \ref{KSSLemma91} there exists $\delta_1'>0$ such that $I^{\tilde m(j)+1}$ is $\delta_1'$-nice.
If this is not the case, then $\mathcal{L}_{c_0}(I^{\tilde m(j)})$ is
well-inside $I^{\tilde m(j)}$, since there exist two turning points
$c_1, c_2$, critical points of $f^{q_0}$, such that $\mathcal{L}_{c_0}(I^{\tilde m(j)})\subset
(c_1,c_2)$, and by Lemma \ref{bounds1}, the interval $(c_1,c_2)$ is
well-inside $I^{\tilde{m}(j)}$. Shrinking $\delta_1$ if necessary,
$I^{\tilde m(j)+1}$ is $\delta_1$-nice for all $j=0,\dots, k$. Since
$|I_{i-1}|/|I_i|<\nu$, $k$ is bounded. We complete the proof the
proposition by repeating the previous argument, to pull back from time $\tilde m(0)$ to $\tilde m(1)$, $k-1$ times. 

\medskip
 \noindent\textbf{Combinatorial Remark 2.} Note that there exist at most two intervals $G_j$ with $s_i'\leq j\leq s_{i-1}$ and $G_j\subset I_i^\infty$, and one of them is $G_{s_i'}$. Thus, there exists at most one interval $G_j$ with  $s_i'\leq j\leq s_{i-1}$ with $G_j\subset I_{i+1}.$ 
\end{pf}

In the next proposition, we will deal with the first step of the pullback when $I_{i-1}$ is terminating, see Figure \ref{diagram2}. 

\begin{prop}\label{terminating angle control}
For each $\nu\geq 1$ there exist $\mu=\mu(\nu)\in(0,1)$ and $\delta>0,$ and for each $\theta\in(0,\pi)$ there exists $\varepsilon>0$ such that the following holds. Assume that $I_{i-1}$ is terminating with $|I_{i-1}|<\varepsilon$ and that
$$\frac{|I_{i-2}|}{|I_{i+1}|}<\nu.$$
Let $\{G_{j}\}_{j=0}^{s}$ be a disjoint chain with $G_{s}=I_{n},
G_{0}\subset I_{n}$ for some $n>i$ and
$G_{0}\cap\omega(c_{0})\neq\emptyset$. Suppose  $s_{i-1}\leq s$ is
such that at least three of the intervals $G_{0}, G_{1},\dots
G_{s_{i-1}}$ are contained in $I_{i+1}$. Then there exists
$r_{i-1}<s_{i-1}$ with $G_{r_{i-1}}\subset \mathcal{R}(I_{i-1})$ such
that if $$U_{s_{i-1}}= D_{\theta}(I_{i-1})\cap\mathbb{C}_{G_{s_{i-1}}}\,and\;U_{j}=\mathrm{Comp}_{G_{j}}f^{-(s_{i-1}-j)}(U_{s_{i-1}}),$$
then there exists an interval $K$ so that
$$U_{r_{i-1}}\subset D_{\mu\theta}(K),$$ 
and $(1+2\delta)K\subset \mathcal{R}(I_{i-1})$. Moreover, there are at most three times $j$ with $r_{i-1}\leq j\leq s_{i-1}$ and $G_{j}\subset I_{i+1}$. 
\end{prop}

\begin{pf}
Let $r$ be the period of $I_{i-1}^\infty$. Since $I_{i-1}$ is
terminating, $I_{i-1}^{\infty}$ is a periodic interval, and
$\mathcal{R}(I_{i-1})$ is bounded between $\alpha_{i-1}$ and
$\tau(\alpha_{i-1})$, where $\alpha_{i-1}$ is the orientation
reversing fixed point of $R_{I_{i-1}}$ closest to $c_{0}$. To be
definite let $Y_{-1}\neq Y_{0}$ be the other component of
$I_{i-1}^{\infty}\setminus (f^r|I_{i-1}^\infty)^{-1}(\alpha)$ that contains $\alpha_{i-1}$ in its boundary, let $Y_{-a}$ be the component containing $\beta$, the fixed point of $R_{I_{i-1}^\infty}$ contained in $\partial I_{i-1}^\infty$, and $Y_{a}$ be the component containing $\tau(\beta)$ in its boundary. For simplicity, we will denote $\alpha_{i-1}$ by $\alpha$.

\medskip
\noindent\textit{Step 1: Pulling back to $I_{i-1}^{\infty}.$} 
We define a time $s_{i-1}^\infty$ with $0\leq s_{i-1}^\infty<s_{i-1}$ so that  the following holds. There exist  $\lambda_0\in(0,1)$  and an interval $K_0$ well-inside $I_{i-1}^\infty$  so that
$$U_{s_{i-1}^{\infty}} \subset  D_{\lambda_0\theta}(K_0).$$
Since $I_{i-1}$ is terminating, only one component of its return domain intersects $\omega(c_{0})$, so $G_{s_{i-1}}\subset I_{i-1}^{\infty}$. 
Let us first pull back by one iterate of $f^r$.  Since $|I_{i-2}|/|I_{i+1}|<\nu$ we can apply 
Corollary \ref{cor:delta free} to find $\rho>0$ so that $I_{i-1}$ is $\rho$-free. 
Notice that since $I_{i-1}$ is terminating, all critical points of $f^r|I_{i-1}^1$ are contained in $I_{i-1}^{\infty}$. Let $0< j\leq r $ be minimal so that $\mathcal{L}_{f^{r-j}(c_0)}(I_{i-1})$ contains a turning point $c'$ of $f$  and $c'\notin G_{s_{i-1}-j}$. Then $G_{s_{i-1}-j}$ is contained in an interval bounded on one side by a preimage $\hat \alpha$ of $\alpha$ (such that $f^{r-j}(\hat \alpha)=\alpha$) and on the other either by $c'$ or by a boundary point of $\mathcal{L}_{f^{r-j}(c_0)}(I_{i-1})$. Let $\alpha'$ be the symmetric point to $\hat \alpha$ with respect to $c'$.
Since $I_{i-1}$ is $\rho$-free we can apply Lemma \ref{bounds1}(2) to conclude that $|c'-\alpha'|$ is comparable to $|\mathcal{L}_{f^{r-j}(c_0)}(I_{i-1})|$. As in the previous lemma,
this means that there exist an interval $K'$, $\lambda'\in(0,1)$  and $\delta'>0$
such that $$U_{s_{i-1}-j}\subset D_{\lambda'\theta}(K'),$$
where $K'$ is an interval that is bounded on one side by a boundary point $x$ of $\mathcal{L}_{f^{r-j}(c_0)}(I_{i-1})$ and on the other by a point $x'\in (c'-\delta'|\mathcal{L}_{f^{r-j}(c_0)}(I_{i-1})|,c'+\delta'|\mathcal{L}_{f^{r-j}(c_0)}(I_{i-1})|)$, on the opposite side of $c'$ as $x$. Additionally, we can assume $K'$ has the property that $(1+\delta')K'\setminus (1+\delta')^{-1}K'$ is disjoint from $\cup \{f(c),f^2(c),\dots f^r(c)\}$, where the union is taken over all critical points $c$ of $f$.
Replacing $\mathcal{L}_{f^{r-j}(c_0)}(I_{i-1})$ by $K'$ and repeating this process at most once for each critical point of $f$, we have that there exists an interval
$K_1\subset I_{i-1}^1$, $\lambda_1\in(0,1)$ and $\delta_1>0$
such that 
$$U_{s_{i-1}-r}\subset D_{\lambda_1\theta}(K_1),$$
where $K_1$ is of one of the following types.
\begin{enumerate}[(a)]
\item Assume $G_{s_{i-1}-r}$ is contained in one of the two outermost monotone branches of $f^r|I_{i-1}^1$. Then $K_1$ is bounded between a boundary point $x_1$ of $\mathcal{L}_{c_0}(I_{i-1})$ and a point $x_2\in (c_1-\delta_1|I_{i-1}^\infty|,c_1+\delta_1|I_{i-1}^\infty|)$ on the opposite side of $c_1$ as $x_1$ where $c_1$ is an outermost turning point of $f^r$;
\item Assume $G_{s_{i-1}-r}$ is contained in $I^\infty_{i-1}$ minus the two outermost monotone branches of $f^r|I_{i-1}^1$. Then there exist turning points $c_1<c_2$ of $f^r$ and points $x_1\in(c_1-\delta_1|I_{i-1}^\infty|,c_1)\subset I_{i-1}^\infty$ and $x_2\in (c_2,c'+\delta_1|I_{i-1}^\infty|)\subset I_{i-1}^\infty$ such that $K_1=(x_1,x_2)$ is well-inside $I_{i-1}^\infty$;
\item  Assume that  $G_{s_{i-1}-r}$ contains one of the turning points
  from the two 
outermost monotone branches of $f^r|I_{i-1}^1$. Then, $K_1$ is bounded on one side by a boundary point $x_1$ of $\mathcal{L}_{c_0}(I_{i-1})$ and on the other side by a point $x_2\in (c_{0}-\delta_1|I_{i-1}^\infty, c_{0}+\delta_1|I_{i-1}^\infty|)$ on the opposite side of $c_0$ as $x_1$.
\end{enumerate}

Suppose first that (a) holds. Let $\mathcal{Y}$ denote the outermost monotone branch of $f^r|I_{i-1}^1$ containing  $G_{s_{i-1}-r}$ and let $c$ be the turning point on its boundary. If $\mathcal{Y}$ is the orientation reversing branch, pulling back under one more iterate of $f^r$, as in the proof of Claim 6 of Proposition~\ref{good start or good deep}, we obtain an interval $K_2$ that is well-inside $I_{i-1}^\infty$, and $\lambda_2\in(0,1)$ such that 
$$U_{s_{i-1}-2r}\subset D_{\lambda_2\theta}(K_2),$$
and $s_{i-1}^\infty=s_{i-1}-2r$.
We know that $s_{i-1}-2r>0$, since at least two of $G_{0}, G_{1},\dots G_{s_{i-1}-1}$ are contained in $I_{i+1}$. Even more, the only interval in $G_{s_{i-1}^\infty}, G_{s_{i-1}^\infty+1},\dots,G_{s_{i-1}-1}$ that could be contained in $I_{i+1}$ is $G_{s_{i-1}^\infty}$.
Now suppose that $\mathcal{Y}$ is the orientation preserving branch. Let $q\geq 1$ be maximal so that 
$G_{s_{i-1}-qr}\subset \mathcal{Y}$ and for $t, 1\leq t\leq q$, $G_{s_{i-1}-tr}\subset \mathcal{Y}$. If $q=1$, then either 
\begin{itemize}
\item there exists an interval $K_2$ which is well-inside $I_{i-1}^\infty$, and $\lambda_2\in(0,1)$ such that 
$$U_{s_{i-1}-2r}(z)\subset D_{\lambda_2\theta}(K_2),$$
and $s_{i-1}^\infty=s_{i-1}-2r$ or
\item $G_{s_{i-1}-2r}$ is contained in the outermost orientation reversing branch of $f^r|I_{i-1}^1$, $s_{i-1}-3r>0$ and, as above, 
we have that there exists an interval $K_2$ which is well-inside $I_{i-1}^\infty$, and $\lambda_2\in(0,1)$ such that 
$$U_{s_{i-1}-3r}(z)\subset D_{\lambda_2\theta}(K_2).$$
\end{itemize}
So suppose that $q>1$. Let $c_{-1}$ be the preimage of $c$ under $f^{r}$ contained in $\mathcal Y$. By Lemma \ref{bounds1}, the interval $(c_{-1},c)$ has length comparable to $I_{i-1}$, and by Lemmas \ref{SLemma55} and Lemma \ref{fundamental domains} there exists a fundamental domain for $f^r|\mathcal Y$
contained in $I_{i-1}^1\setminus I_{i-1}^\infty$ of size comparable to $I_{i-1}$. So
by Proposition~\ref{monotone pullbacks}, there exists $\lambda_2>0$ such that
$$U_{s_{i-1}-qr}\subset D_{\lambda_2\theta}(\mathcal Y).$$
Arguing as in the case when $q=1$, we have that $s_{i-1}^\infty= s_{i-1}-(q+1)r$ or 
$s_{i-1}^\infty= s_{i-1}-(q+2)r$, and there exist an interval $K_3$ which is well-inside $I_{i-1}^\infty$ and $\lambda_3\in(0,1)$ such that 
$$U_{s_{i-1}^{\infty}} \subset D_{\lambda_3\theta}(K_3).$$
Now, suppose that either (b) or (c) holds.
If (b) holds, then $s_{i-1}-r=s_{i-1}^\infty$ and we are done with this step.
If (c) holds, applying the argument at the start of Step 1 to pull back once under  $f^r$, 
we have that either (a) or (b) holds, since the chain $\{G_j\}_{j=0}^{s_{i-1}}$ is disjoint. In any case, we have that there exists an interval $K_0$ which is well-inside
$I_{i-1}^\infty$, and $\lambda_0\in(0,1)$ such that 
$$U_{s_{i-1}^\infty}\subset D_{\lambda_0\theta}(K_0).$$

\noindent \textbf{Combinatorial Remark 3.} Note that there exist at most two intervals $G_j\subset \mathcal R(I_{i-1})$ with $s_{i-1}^\infty\leq j\leq s_{i-1}$, so there exists at most one interval $G_j\subset I_{i+1}$ with $s_{i-1}^\infty\leq j\leq s_{i-1}$. 

\medskip 
Since $c$ is a critical point of $f^r$, we have that $|c-\alpha|$ is comparable to $|I_{i-1}^\infty|$, so we can choose the previous constants so that there exists $\delta_0>0$ such that $(1+\delta_0)K_0\setminus (1+2\delta_0)^{-1}K_0$ does not contain
$\alpha$.
From now on, we will consider pullbacks by $f^{r}|I_{i-1}^{\infty}$.

\medskip
\noindent\textit{Step 2: Pulling back to a $\hat \rho$-nice interval.}

We define a time $t$, with $0<t<s_{i-1}^\infty$ with the following properties. There exist a $\hat \rho$-nice interval $X$, for $\hat \rho>0$ universal, and $\lambda''>0$ so that $G_{t}\subset X$ and $U_{t}\subset D_{\lambda''\theta}(K'')$, for an interval $K''$ that is well-inside $X$. We will make use of  the intervals $Y_j$ defined in Section \ref{sec:enhanced nest} to find $X$. Remember $Y_0=\mathcal R(I_{i-1})$; $Y_{-a}$ contains $\beta$, the fixed  point of $f^r$ in the boundary of $I_{i-1}^\infty$; and $Y_a$ contains $\tau(\beta)$.
 
\medskip
\noindent\textit{Step 2(a).} First, suppose that $\alpha\notin
K_{0},$ then pulling back once under $R_{I_{i-1}^{\infty}}$, we
can find an interval $K_{1}'$, and a
$\lambda_{1}'\in(0,1)$ such that $U_{s_{i-1}^\infty-r}\subset D_{\lambda_{1}'\theta}(K_{1}')$ where $K_{1}'$ is well-inside $Y_{j}$ for some $j$. If $j=0$ we set $r_{i-1}=s_{i-1}^\infty-r$  and we are done with the proof of this proposition in this case. Note that if $j=0$ there are at most two times $j$ with $r_{i-1}\leq j\leq s_{i-1}$ and $G_{j}\subset I_{i+1}$. In the remainder of this step, we will assume that $j\neq 0$. Note that this implies that $s_{i-1}^\infty-r> 0$. 

\medskip
\noindent\textit{Case 1:} If $j\notin\{a, -a,1\}$, then   by Lemma~\ref{rho nice landing domains}
we know that $Y_{j}$ is $\hat\rho$-nice, so we set $X=Y_{j}$ and $t=s_{i-1}^{\infty}-r$. 
Note that in this case there exists no $G_j$ with $G_j\subset \mathcal R(I_{i-1})$ and $ t\leq j\leq s_{i-1}^\infty-r.$

\medskip
\noindent\textit{Case 2:} If $j=a$, then $s_{i-1}^{\infty}-2r>0$ and
after one more pull back by $f^r$ we can find an interval $K_{2}'$ and a $\lambda_{2}'\in(0,1)$ such that $U_{s_{i-1}^{\infty}-2r}\subset D_{\lambda_{2}'\theta}(K_{2}')$ where $K_{2}'$ is well-inside $Y_{j'}$. Since $K_{1}'\subset Y_{a} $, we have that $j'\notin\{a,-a\}.$ If $j'\notin\{0,1\}$ we set $t= s_{i-1}^{\infty}-2r>0$  and $X=Y_{j'}$ and go to Step 3. Note that if $j'\notin\{0,1\}$ there exists no  $G_j\subset \mathcal R(I_{i-1})$ with $ t\leq j\leq s_{i-1}^\infty-r.$ If $j'=0,$ we set  $r_{i-1}= s_{i-1}^{\infty}-2r$ and we are done with the proof of the proposition. If $j'=-1$ we go to Case 4.

\medskip
\noindent\textit{Case 3:} If $j=-a$, let $q$ be such that
$qr=s_{i-1}^{\infty}-r$, and let $q'\leq q$ be minimal so that for
$k$, $q'\leq k\leq q$, $G_{kr}\subset Y_{-a}$. Since $f^r|_{Y_{-a}}$
is monotone we can apply Proposition~\ref{monotone pullbacks} to
control the loss of angle until time $q'r$. Then after pulling back
once more, we can find an interval $K_{2}'$ and a
$\lambda_{2}'\in(0,1)$ such that $U_{(q'-1)r}(z)\subset
D_{\lambda_{2}'\theta}(K_{2}')$. Let $Y_{j'}$ be the puzzle piece
containing $K_2'$, then $K_2'$ is well-inside $Y_{j'}$. If $j'=0$ we set $r_{i-1}=(q'-1)r$ and we are done with the proof of the proposition. If $j'\notin\{-a, 1\}$, then $Y_{j'}$ is $\delta'$-nice for some $\delta'>0$ so we set $X=Y_{j'}$ and $t=(q'-1)r$. Note that if $j'\notin\{-a,0, 1\}$ there exists no $G_j\subset \mathcal R(I_{i-1})$ with $ t\leq j\leq s_{i-1}^\infty-r.$ If $j'=a$, go to Case 2, and if $j'=1$ go to Case 4.

\medskip
\noindent\textit{Case 4:} Suppose that $j=-1$. If $G_{s_{i-1}^\infty-2r}$ is contained $Y_0$, we set $r_{i-1}= s_{i-1}^\infty-2r$ and we are done with the proof of the proposition; if $G_{s_{i-1}^\infty-2r}\subset Y_{j'}$ with $j'\notin\{a,-a\}$ we set $t=s_{i-1}^\infty-2r$ and $X=Y_{j'}$; if $G_{s_{i-1}^\infty-r}\subset Y_{j'}$ with $j'\in\{a,-a\},$ we complete the proof of this step using the appropriate previous case. Observe that if we go from Case 2 to Case 4 and back to Case 2 or from Case 3 to Case 4 and back to Case 3, then by Lemma \ref{rho nice landing domains} we can set $X=\mathcal L_{G_t}(Y_{a})$ or $X=\mathcal L_{G_t}(Y_{-a})$, for a suitable time $t$. 

\medskip

\noindent\textbf{Combinatorial Remark 4.} If after applying Step 2(a) we have not proved the proposition, then there exists no $G_j\subset \mathcal R(I_{i-1})$ with $t\leq j\leq s_{i-1}^{\infty}-r$. This means that there exists at most one $G_j\subset I_{i+1}$  with $t\leq j \leq s_{i-1}$.

\medskip
\noindent\textit{Step 2(b).} Suppose that $\alpha\in(1+\delta_0)K_{0}.$
Recall that $K_0\supset G_{s_{i-1}^{\infty}}$. First, as in Step 1 of this proposition, pulling back by one iterate of $f^r$
we have that there exists an interval $K_1'$ which is well-inside $I_{i-1}^\infty$ and $\lambda_1'\in(0,1)$ such that 
$$U_{s_{i-1}^\infty-r}\subset D_{\lambda_1'\theta}(K_1').$$
If $G_{s_{i-1}^{\infty}-r}$ is not contained in $\mathcal{\tilde{Y}}_{\alpha}$ we are in Step 2(a). 
Otherwise, $G_{s_{i-1}^\infty-r}\subset \mathcal{\tilde Y}_\alpha$ and we can ensure that 
the pullback of $K_1'$ under $f^{r}|\mathcal{\tilde{Y}}_{\alpha}$ is  
well-inside $\mathcal{\tilde{Y}}_{\alpha}$.  Let $p$ be such that $pr=s_{i-1}^{\infty}-r$, and let $p'\leq p$ be minimal so that for $j$, $p'\leq j\leq p$, $G_{jr}\subset \mathcal{\tilde{Y}}_{\alpha}$. We have that $\mathcal{\tilde{Y}}_{\alpha}$ is bounded by two turning points, so by Lemma \ref{bounds1} we can find fundamental domains for the map from this monotone branch to itself of size comparable to $|I_{i-1}^{\infty}|$. By Proposition~\ref{monotone pullbacks} we can control the loss of angle until time $p'r$.
Observe that for $p'\leq k\leq p$ only $G_{pr}$ could be in $I_{i+1}$.
Let $j_{1}$ be such that $G_{p'r-r}\subset Y_{j_{1}}$. Since $G_{p'r-r}$ is not contained in 
$Y_0\cup Y_{-1}$ we are in Step (2a).

\medskip
\noindent\textbf{Combinatorial Remark 5.} If after applying Steps 2(a) and 2(b) we have not proved the proposition, then there exists no $G_j\subset \mathcal R(I_{i-1})$ with $t\leq j\leq s_{i-1}^{\infty}-r$. This means that between time $s_{i-1}$ and time $t$, we can have at most one visit to $I_{i+1}$. Even more, if $\mathcal R(I_{i-1})$ is non-periodic, there exists at most two intervals $G_j\subset I_{i}$ with $t\leq j\leq s_{i-1}^{\infty}.$

\medskip
\noindent\textit{Step 3: Pulling back to $\mathcal{R}(I_{i-1})$.}

Let $r_{i-1}<t$ be maximal with $G_{r_{i-1}}\subset \mathcal{R}(I_{i-1}).$
Let $t'$, $r_{i-1}<t'<t$, be minimal with the property that $G_{t'}\subset X$. Since  $G_{r_{i-1}}$ is the only element of the chain contained in $\mathcal{R}(I_{i-1})$ we know that for each $c\in\crit(X; \{G_j\}_{j=r_{i-1}}^{t})$ there exists at most one interval $G_{j}$ such that $G_j\subset\hat{\mathcal L}_c(\mathcal{R}(I_{i-1}))$. Let $t_1<t$ be maximal so that  $G_{t_1}\subset\hat{\mathcal L}_c(\mathcal{R}(I_{i-1}))$. If $t_1$ is not defined let  $t_1=r_{i-1}$. Let $t_1'> t_1$ be minimal such that $G_{t_1'}\subset X$ and let $\mathbb G= \{G_j\}_{j=t_1'}^{t}$. There are two cases to consider:
\begin{itemize}
\item[(i)] $k(X, \mathbb G)< \infty$ or
\item[(ii)] $k(X,\mathbb G)= \infty$.
\end{itemize}
\medskip
For each $c\in\crit(X; \mathbb G)$ pick $t_1'\leq n_{c}< t$ minimal such that $G_{n_c}\subset \mathcal L_c(X)$  and $t_1'\leq m_c \leq t$ maximal such that $G_{m_c}\subset \mathcal L_c(X)$.

\medskip
\noindent\textit{Claim:}  \textit{If $k(X,\mathbb G)$ is finite, there exists $a=a(\nu)>0$ such that $k(X,\mathbb G)< a$.}

The proof of this claim is analogous to the proof of Claim 1 of Proposition~\ref{folding lemma}; one only needs to substitute $I_{i-1}$ by $X$ and $I_i$ by $\mathcal R(I_{i-1})$.
Following the proof we get that if $a$ is large enough and $c'=c_0$, then $|\mathcal{R}(I_{i-1})|/  |I_{i+1}|>\nu$. If $a$ is large enough and $c'\neq c_0$ we get that $|\mathcal{R}(I_{i-1})|/|\Gamma(\mathcal{R}(I_{i-1}))|>\nu$ 
\endpfclaim
 
\medskip

\noindent\textit{Case (i).} Since we are assuming $|\mathcal{R}(I_{i-1})|/|I_{i+1}|<\nu$ the claim tells us that $k(X, \mathbb G)< a$. In this case the proposition follows arguing as in Case (i) of Proposition~\ref{folding lemma}, substituting $j_0$ by $r_{i-1}$ and $I_{i-1}$ by $X$ in the proof.

\medskip
\noindent\textit{Case (ii).} In this case the proposition follows arguing as in Case (ii)  of Proposition~\ref{folding lemma}, substituting $j_0$ by $r_{i-1}$, $s_{i-1}$ by $t$, $I_{i-1}$ by $X$ and $I_i$ by $\mathcal R(I_{i-1})$ in the proof.

\medskip
From Combinatorial Remark 5 we get that at most three intervals in the chain $\{G_{j}\}_{j=r_{i-1}}^{s_{i-1}}$ are contained in $I_{i+1}$ and the proposition follows. Observe that if exactly three of these intervals are contained in $I_{i-1}$, $G_{s_{i-1}}$ and $G_{r_{i-1}}$ are both contained in $I_{i+1}$. Moreover, if $\mathcal R(I_{i-1})$ is non-periodic there exists at most three intervals $G_j\subset I_{i}$ with $r_{i-1}\leq j\leq s_{i-1}^{\infty}$. If there are exactly three of these intervals contained in $I_i$, $G_{s_{i-1}}$ and $G_{r_{i-1}}$ are both contained in $I_{i}$.
\end{pf}

\noindent \textbf{Combinatorial Remark 6.} It is important to observe that in the proof of Proposition~\ref{terminating angle control} we do not have control on the number of visits to $I_i$ if $\mathcal R(I_{i-1})$ is periodic. However,  if $I_i\neq \mathcal R(I_{i-1})$ we can prove the proposition assuming that at least three of the intervals $G_0,\ldots, G_{s_{i-1}}$ are contained in $I_i$ (instead of in $I_{i+1}$). 

\medskip
The proof of the previous proposition also shows the following.

\begin{cor}\label{cor:terminating angle control}
In the setting of the previous proposition, 
suppose that $\hat I$ is a terminating interval with
$I_{i-1}\subset \hat I\subset I_{i-2}$. Let $\hat K$ be an interval $\delta'$-well-inside $\hat I$ for some $\delta'>0$ with $ G_{s_{i-1}}\subset \hat K$. Let
$$U_{s_{i-1}}= D_{\theta}(\hat K)\cap\mathbb{C}_{G_{s_{i-1}}}\,and\;U_{j}=\mathrm{Comp}_{G_{j}}f^{-(s_{i-1}-j)}(U_{s_{i-1}}).$$
Then there exist an interval $K$ and $\delta>0$ such that
$$U_{r_{i-1}}\subset D_{\mu\theta}(K),$$ 
and $(1+2\delta)K\subset \mathcal{R}(I_{i-1})$. Moreover, there are at most three  times $j$ with $r_{i-1}\leq j\leq s_{i-1}$ and $G_{j}\subset I_{i+1}$.
\end{cor}

Corollary \ref{cor:terminating angle control}  is the result we need to apply after Proposition~\ref{folding lemma} in Figure \ref{diagram1}. In this application the interval  $\hat I$ will satisfy $I_{i}\subset \hat I\subset I_{i-1},$ and the time $s_{i-1}$ will be substituted by the time $s_i'$  obtained from Proposition~\ref{folding lemma} (rather than the time $s_i$). In the next section we will show that $s_i'$ satisfies the assumption needed in Proposition~\ref{terminating angle control}. 

The next proposition completes the results needed to complete the pullbacks  described in Figures \ref{diagram1} and \ref{diagram2}. Note that in Figure \ref{diagram1} we use Proposition~\ref{prop:from renormalization to landing} for $I_i$ instead of $I_{i-1}$, so the indexes have to be changed accordingly. 

\begin{prop}\label{prop:from renormalization to landing}
For each $\nu>0$ 
there exists $\mu=\mu(\nu)\in(0,1)$ 
and for each $\theta\in (0,\pi)$ 
there exists $\varepsilon>0$  such that the following holds.
Let $I_{i-1}$ be  terminating with  $|I_{i-1}|<\varepsilon$. 
Assume that $\mathcal{R}(I_{i-1})$ is non-terminating and that
$$\frac{|I_{i-2}|}{|I_{i+1}|}<\nu.$$
Assume $n>i$ and consider a disjoint chain $\{G_{j}\}_{j=0}^{s}$ with $G_{s}=I_{n},$ $G_{0}\subset I_{n}$, $G_0\cap\omega(c_0)\neq\emptyset$ and $G_{r_{i-1}}\subset \mathcal{R}(I_{i-1})$ for  some $0<r_{i-1}\leq s$. Let $s_{i}$ be maximal so that $G_{s_{i}}\subset I_{i}$ with $0\leq s_i< r_{i-1}$ and define  $$U_{r_{i-1}}=D_{\theta}(\mathcal{R}(I_{i-1}))\cap\mathbb{C}_{G_{r_{i-1}}}\;\;and\;\;U_{j}=\mathrm{Comp}_{G_{j}}f^{-(r_{i-1}-j)}(U_{r_{i-1}}).$$
Then 
$$U_{s_{i}}\subset D_{\mu\theta}(I_{i}).$$
Moreover, there are at most two times $j$ with $G_{j}\subset I_{i}$ and $s_{i}\leq j\leq r_{i-1}.$

\end{prop}

\begin{pf}
Let $r$ be the period of $I_{i-1}^\infty$.
By Proposition~\ref{delta free}, $\mathcal{R}(I_{i-1})$ is $\rho$-free. 
Let $s'<r_{i-1}$ be maximal so that $G_{s'}\subset \mathcal{R}(I_{i-1})$. Let $\hat I:= \mathcal{L}_{G_{s'}}(\mathcal{R}(I_{i-1}))$. By Lemma \ref{chain pullback}, pulling back once under the return map to ${\mathcal{R}(I_{i-1})}$ we have that
$$U_{s'} \subset D_{\lambda\theta} (\hat I),$$
for $\lambda \in (0,1)$ and $G_{s'}\subset \hat I \subset \mathcal{R}(I_{i-1})$.
If $\hat I=I_{i}$ we are done, so we suppose that this is not the case.  As in Corollary \ref{first landing} we define $V_0$ to be the component of $f^{-2r}(\mathcal{R}(I_{i-1}))$ that contains $\alpha$ in its boundary, $W_0$ to be the component of $f^{-2r}(\mathcal{R}(I_{i-1}))$ that contains $\tau(\alpha)$ in its boundary, $V_1$ be the component of $f^{-4r}(\mathcal{R}(I_{i-1}))$ that contains $\alpha$ in its boundary, and $W_1$ be the component of $f^{-4r}(\mathcal{R}(I_{i-1}))$ that contains $\tau(\alpha)$ in its boundary.

\medskip
\noindent \textit{Case 1:} By Corollary \ref{first landing}, if $\hat I$ is different from $V_0$ and from $W_0$, then there exists $\rho'>0$ such that $\hat I$ is  $\rho'$-nice and $\rho'$-free.  
Let $s''$ be minimal such that $G_{s''}\subset \hat{I}$ with  $s_{i}< s''\leq s'$. 
For each critical point $c$, there is at most one $j_{c}$ such that $G_{j_{c}}\subset \hat{\mathcal{L}}_{c}(I_{i})$ with $s'' \leq j_c< s'$. Let $j_{c}$ be the maximal such  $G_{j_{c}}\subset \hat{\mathcal{L}}_{c}(I_{i})$ for $s'' \leq j_c<s'$, where the maximum is taken over all $c\in \crit(f)$. If $j_{c}$ is not defined let $j_{c}=s''$. If  $j_{c}>s''$ let $s_{1}$ be minimal with $G_{s_{1}}\subset\hat{I}$ and $j_{c}<s_{1}\leq s'$. Then if $I_{i}$ is non-terminating,
$k(\hat{I},\{G_{j}\}_{j=s_{1}}^{s'})$ is bounded  by a constant $K(\nu)$ by Proposition~\ref{lem61}.  If $I_{i}$ is terminating, we split the chain into at most $b$ parts, one for each critical point in $\crit(\hat I, \{G_{j}\}_{j=s_{1}}^{s'})$, on which we can control the loss of angle using Corollary
\ref{chain space} and Proposition~\ref{combinatorial depth angle control} as in Case (ii) of Proposition~\ref{folding lemma}. In any case, we can find  $\lambda_1\in (0,1)$ and an interval $K_1$ which is well-inside $\hat I$ such that 
$$U_{s''}\subset D_{\lambda_1\theta}(K_1).$$
%
%
%
%
By Corollary \ref{chain space} there exists $\mu\in(0,1)$ so that 
\begin{equation}\label{s_{i}'''}
U_{s_{i}}\subset D_{\mu\theta}(I_{i})
\end{equation}

\medskip
\noindent\textit{Case 2:}  If $\hat I=W_0$ let $p$ be maximal such that $G_p\subset \mathcal{R}(I_{i-1})$ and $s_i\leq p<s'$. Our goal in this case is to get inside a $\rho$-nice interval, so that we can make use of the argument of the previous case. Applying Lemma \ref{chain pullback}  we can find $\lambda'>0$ such that
$$U_p\subset D_{\lambda'\theta}(J),$$
where $J=\mathcal L_{G_p}(W_0)$. If $p=s_i$ we are done with the proof of the proposition, so let us assume this is not the case. If $J$ is not contained in $V_0$ or $W_0$, then Corollary \ref{first landing} implies that there exists $\rho'>0$ such that $J\subset \mathcal L_{G_p}(\mathcal{R}(I_{i-1}))$  for some $\mathcal L_{G_p}(\mathcal{R}(I_{i-1}))$ which is $\rho'$-nice and $\rho'$-free. In this case we can proceed as in Case 1 to get (\ref{s_{i}'''}). If $J\subset V_0$, then the fact that $V_0$ has size comparable to the size of $\mathcal{R}(I_{i-1})$ implies there exists $\rho_1>0$ such that $(1+2\rho_1)J\subset V_0$.  Applying Theorem \ref{real Koebe} (2) we can find $\rho_1'>0$ so that the following holds. Given $x\in J\cap\omega(c_0)$, $(1+2\rho_1')\mathcal L_x(J)\subset \mathcal L_x(V_0)\subset J$. Which means that $J$ is $\rho_1'$-nice. Shrinking  $\rho_1'$ if necessary we can assume that $J$ is $\rho_1'$-free. Now we can proceed as in Case 1 to get (\ref{s_{i}'''}). If $J\subset W_0$, then a similar argument implies the existence of  $\rho_1'>0$ so that the following holds. Given $x\in J\cap\omega(c_0)$, then $(1+2\rho_1')\mathcal L_x(J)\subset \mathcal L_x(W_0)=J$. Which means that $J$ is $\rho_1'$-nice. Again, shrinking  $\rho_1'$ if necessary we can assume that $J$ is $\rho_1'$-free. 

\medskip
\noindent\textit{Case 3:} Assume $\hat I=V_0.$ If $G_{s'}$ is at most the second return of $G_{s_{i}}$ under the first return map to $\mathcal{R}(I_{i-1})$, then we can get (\ref{s_{i}'''}) applying Lemma \ref{chain pullback}  at most twice, so let us assume this is not the case. In this case we will pull back along the outermost orientation preserving monotone branch to end up in a $\rho$-nice interval or in $W_0$ as in the previous cases. We do not make use of Proposition~\ref{monotone pullbacks}, since we do not have fundamental domains on both sides of this branch. Instead, we use the free space to give us intervals at the boundary of $V_0$ whose pullbacks shrink exponentially. We will use the same notation as in the proof of $(3)$ in Proposition~\ref{delta free}. Let $\hat{V}_0 =\Comp_\alpha f^{-r}(V_0) \cup V_{0}$,
$\hat{V}_1= \comp_{\alpha}f^{-r}(V_1)\cup V_1$ and  $\hat{V}_{j+1}= \comp_\alpha f^{-r}(\hat{V}_j)$ for $j\in \N$. From the proof of Proposition~\ref{delta free} we know that $f^{2r}|_{V_1}$ is monotone and that there exists $\delta_1>0$ such that $(1+2\delta_1)\hat V_2\subset \hat V_1.$ Since $G_{s_i}\subset I_i$ there exists $G_j$ with $s_{i}\leq j< s'$ such that  $G_j\subset \mathcal L_x(\mathcal{R}(I_{i-1}))$, for some $\mathcal L_x(\mathcal{R}(I_{i-1}))$ different from $V_0$. Let $q'$ be minimal such that $s_{i}\leq q'<s'$, $G_{q'}$ is not contained in $V_1$ and for each $q'<j<s'$ with $G_j\subset \mathcal{R}(I_{i-1})$ we have that $G_j\subset V_1$. Set $q=s'-2r$, note that $s'-2r>0$ by assumption. By Lemma \ref{chain pullback} we can  find $\lambda_1\in(0,1)$ so that$$U_q\subset D_{\lambda_1\theta}(V_2).$$

If $q'=q,$ under one more pull back we are in one of the previous
cases so assume $q>q'$. Since $\hat V_2$ is $\rho'$ well-inside $\hat
V_1$ and $\mathcal{R}(I_{i-1})$ is $\delta$-free, taking
$\delta''=\min\{\rho',\delta\}$ we get that $V_2$ is $\delta''$
well-inside the monotone branch of $f^r$ containing $\alpha$. Let
$V_{n_{0}+N}$ be as in  Case 3 of the proof of Proposition~\ref{delta
  free}. Then, $V_{n_{0}+N}\cap\omega(c_{0})=\emptyset$ and the size
of $V_{n_0+N}$ is comparable to the size of $V_0$. Let $V$ be the
preimage of $V_{n_{0}+N}$ under $f^{2r}|_{V_1}$. If there is a turning
point $c$ of  $f^{2r}|_{V_{0}}$ let $\tilde{V}$ be the preimage of
$V_{n_{0}+N}$ under $f^{2r}$ that is symmetric to $V$ with respect to
$c$. If $f^{2r}|_{V_{0}}$ is monotone let $V'$ be the preimage of
$V_{n_0+N}$ under $f^{2r}$ that contains $\tau(\alpha)$ and let
$\tilde V$ be the preimage of $V'$ under $f^{2r}|_{V_0}.$  We can
apply Lemma \ref{bounds1} to $f^r|_{I^{\infty}_{i-1}},$ so we get real
bounds for $f^{2r}|_{V_0}$; these bounds along with the definition of
$V_{n_{0}+N}$ imply that  the sizes of $V$ and $\tilde V$ are
comparable to the size of $V_0$. By construction, the pullbacks of
$\tilde V$ under $f^{2r}|_{V_1}$ are disjoint.  The pullbacks of $V$
under $f^{2r}$ are not disjoint but they will shrink exponentially. To
see this, define $F_{q}=V$ and $F_{j}=(f^{2r}|V_{1})^{-1}(F_{j+1})$
for $q'\leq j<q$. By Theorem \ref{real Koebe} (2), there exists
$\eta\in(0,1)$ universal such that for each $j$, $j=q',\dots, q-1$,
$|F_{j+1}|<\eta|F_{j}|$. Hence $\sum_{j=q'}^{q}|F_{j}|$  is bounded by
a constant depending on $\varepsilon$. We can apply the same proof
used in Proposition~\ref{monotone pullbacks} to pull back from time $q$ to time $q'$. In this case the fact that the pullbacks of $V$ are not disjoint will not matter since we have exponential contraction of the chain $\{F_j\}_{j=q'}^{q}$ and this is enough to control the loss of angle. Thus, we can find $\lambda_2\in(0,1)$  so that $$U_{q'}\subset D_{\lambda_2\theta}(V_2).$$
Let $t_{1}<q'$  be maximal such that $G_{t_{1}}\subset
\mathcal{R}(I_{i-1})$,. If $G_{t_{1}}\cap V_{0}=\emptyset$ we can
apply Lemma \ref{chain pullback}, since $V_0$ is $\rho'$-free , to get
to Case 2 or Case 3. If $\hat G_{t_{1}}\subset V_{0},$ then after one
more pull back we end up either in one of the previous cases or inside a landing domain to $V_{0}$ that is well-inside $V_{0}$, which will be $\delta''$-nice, $\delta''$-free, for some $\delta''>0$, and argue as in Case 1.
\end{pf}

\noindent\textbf{Combinatorial Remark 7.} If in the previous proposition there exist two intervals $G_j\subset I_i,$ we must have that $G_{s_i}$ and $G_{r_{i-1}}$ are such intervals. 

\subsection{The non-renormalizable case}

\begin{prop}[Non-renormalizable case]\label{pullback to previous level}
For each $\nu>0$ there exists $\mu=\mu(\nu)\in(0,1)$ 
and  for each $\theta\in (0,\pi)$ there exists $\varepsilon>0$  such that the following holds.  
Suppose that $I_{i-1}$ and $I_{i}$ are non-terminating,  $|I_{i-1}|<\varepsilon$ and
$$\frac{|I_{i-2}|}{|I_{i+1}|}<\nu.$$
Assume $n>i$ and consider a disjoint chain $\{G_{j}\}_{j=0}^{s}$ with $G_{s}=I_{n},$ $G_{0}\subset I_{n},$ $G_{s_{i-1}}\subset I_{i-1}$ for  $0<s_{i-1}\leq s$. Let $p_{i-1}$ be so that $I_{i}$ is a pullback of $I_{i-1}$ under $f^{p_{i-1}}$. Suppose   $x\in G_{0}\cap\omega(c_{0})$ is so that $f^{s_{i-1}}(x)\in I_{i-1}$ for some $p_{i-1}< s_{i-1}\leq s $ and that  $f^{s'_{i}}(x)\in I_{i}$  for some $0\leq s'_{i} < s_{i-1} -p_{i-1}$. Assume that
$$\#\{k; \, s'_{i} \leq k\leq s_{i-1}\mbox{ and  } f^k(x)\in I_{i}\} \le  4.$$

Let $$U_{s_{i-1}}=D_{\theta}(I_{i-1})\cap\mathbb{C}_{G_{s_{i-1}}}\;and\;U_{j}=\mathrm{Comp}_{G_{j}}f^{-(s_{i-1}-j)}(U_{s_{i-1}}).$$
Then 
$$U_{s'_{i}}\subset D_{\mu\theta}(I_{i}).$$
where $\mu=\mu(\nu, \delta)$. 

\end{prop}
\begin{pf}
By Corollary \ref{cor:delta free} there exists $\delta>0$ such that
$I_{i-1}$ is $\delta$-free. Since $f^{s'_{i}}(x) \in I_{i}$, $f^{p_{i-1}+s'_{i}}(x)\in I_{i-1}$, $ f^{s_{i-1}}(x)\in I_{i-1}$ and $s'_{i} < s_{i-1} -p_{i-1}$ we can decompose $f^{s_{i-1}}$ in the following way $f^{s_{i-1}}(x) =R^{q}_{I_{i-1}}\circ f^{p_{i-1}}\circ f^{s_{i}'}(x)$ 
with  $q\ge 1$. Let $f^{s'}=R^{q-1}_{I_{i-1}}\circ f^{p_{i-1}}\circ f^{s'_{i}}(x).$ By Proposition~\ref{lem:N-modalpb} we know that there exists $\lambda >0$ such that 
$$U_{s'} \subset D_{\lambda\theta} (\mathcal{L}_{f^{s'}(x)}(I_{i-1}) )$$
\noindent with $G_{s'}\subset \mathcal{L}_{f^{s'}(x)}(I_{i-1}) \subset (1+2\delta) \mathcal{L}_{f^{s'}(x)}(I_{i-1})\subset I_{i-1} $. 
Let $s''$ be maximal with $s_{i}'+p_{i-1} \leq s''+p_{i-1} \leq s'$
and $G_{s''}\subset I_{i}$. If $s'= s'' +p_{i-1}$, do nothing. If
$s''+p_{i-1} < s'$ let $s''+p_{i-1} = k_0 < k_1< \dots < k_m= s'$ be
such that $G_{k_j}\in I_{i-1}$. By maximality of $s''$  and Lemma
\ref{KSS83} there are at most three intervals $G_{k_j}$ such that $G_{k_j}\subset I_{i}$.
Let $k_{j_0}$ be maximal such that $j_0\in \{0, \ldots m\}$ and
$G_{k_{j_0}}\in I_{i}.$  Proposition~\ref{combinatorial depth angle
  control} allows to control the loss of angle when we pull back from time $s'$ to time $k_{j_0+1}$ and Corollary \ref{chain space} controls the loss of angle pulling back from time $k_{j_0+1}$ to $k_{j_0}$. Repeating this procedure at most three times, we have that 
$$U_{s'' +p_{i-1}}\subset D_{\lambda'\mu\theta}(K'),$$
where $K'$ is well-inside $I_{i-1}$ and  $\mu$ depends only on $\hat{k}(I_{i-1},I_{i})$, which by Lemma \ref{lem61} is bounded by a constant that depends on $\nu$. If $I_{i-2}$ is non-terminating, then $I_{i-1}$ is $\delta$-externally free, $\delta$-internally free, and $I_{i}$ is a pullback of $I_{i-1}$ of bounded order (depending only on $b$), so, since $|I_{i-1}|/|I_{i}|<\nu$, we can apply Proposition~\ref{lem:N-modalpb}  and find $\mu'>0$ such that
$$U_{s''}\subset D_{\mu'\theta}(I_{i}).$$
If $I_{i-2}$ is terminating, then as in the Case 3 of the proof of the previous proposition, there are puzzle pieces $U$ and $V$ attached to $\alpha$ and $\tau(\alpha)$ in $\mathcal{R}(I_{i-2})$, that are disjoint from $\omega(c_{0})$ and with size comparable to $|\mathcal{R}(I_{i-2})|$, pulling these back to $I_{i-1}$ we obtain puzzle pieces $U_{1}$, $V_{1}$ attached to $\partial I_{i-1}$, with size comparable to $|I_{i-1}|$ that are disjoint from $\omega(c_{0})$. So arguing as in Proposition~\ref{lem:N-modalpb}, we get that 
$$U_{s''}\subset D_{\mu'\theta}(I_{i}).$$
If $s''=s'_{i}$ we are done. In case $s'' >s'_{i}$. Since $I_{i-1}$ is $\delta$-free Corollary \ref{cor:free-space} implies that $I_{i}$ is $\delta'$-free, for some $\delta'>0$. This and the fact that  $\#\{k; \, s'_{i} \leq k\leq s_{i-1}\mbox{ and } f^k(x)\in I_{i}\}< 4$ allows us to prove the proposition applying Lemma \ref{chain pullback} at most twice. 
\end{pf}

\subsection{Large scaling factors}\label{large scaling factors}

In the remainder of this section we show that as soon as one has a
large scaling factor, the corresponding first return map extends to an analytic
box mapping or a qr box mapping, see Definitions~\ref{defn:box
  mapping} and \ref{defn:quasi-box mapping}. The results are stated
for maps in $\mathcal A_{\underline b}^3,$ but if the maps are in $\mathcal A_{\underline b}$ the qr box mappings are box mappings.


\begin{prop}\label{big space gpl}
There exist $\varepsilon>0$ and $\nu_1\geq 1$ such that the following holds.
Suppose that $I_{i-2}$ and $I_{i-1}$ are non-terminating and $|I_{i-2}|<\varepsilon$. If $$\frac{|I_{i-2}|}{|I_{i-1}|}>\nu_{1},$$ then
\begin{enumerate}[(1)]
\item the first return map to $\Gamma(I_{i-1}^1)$ restricted to
  components that intersect $\omega(c_0)$ extends to a qc quasi-box mapping with range $D_{\pi/2}(\Gamma(I^1_{i-1}))$; 
\item if $I_{i}$ is non-terminating, then the first return map to
  $I_{i}$ restricted to the components of the domain intersecting
  $\omega(c_{0})$ extends to a qr box mapping with range $D_{\pi/2}(I_{i})$. 
\end{enumerate}
\end{prop}

\begin{pf}
Let us prove the first statement. To simplify the notation let
$J=\Gamma(I_{i-1}^1)$. By Theorem \ref{real geometry} (a) there exists
$\rho>0$ so that $I_{i-1}$ is $\rho$-externally free. Moreover, if
$c_0$ is even we know that $I_{i-1}\setminus (1+2\rho)^{-1}I_{i-1}$ is
disjoint from $\omega(c_0)$. Since  $J$ is a pullback of $I_{i-1}$ of
bounded order there exists $\rho'>0$ so that $J$ is $\rho'$-externally
free; if  $c_0$ is even we additionally have that $J\setminus
(1+2\rho')^{-1}J$ is disjoint from $\omega(c_0)$.  Fix $C>0$ from Lemma \ref{big scaling factors real bounds}, then $J$ is $C$-nice.
Pick $x\in J\cap\omega(c_0)$ and let $s$ be the return time of $x$ to $J$.  Since $J$ is $\rho'$-externally free the return time of $x$ to $(1+2\rho')J$ is $s$.
Let us first consider the case when $c_0$ is even.  Let $V=D_{\pi/2}(J)$ and $U=\comp_{x} f^{-s}(V).$ Since, in this case, $J\setminus (1+2\rho')^{-1}J$ is disjoint from $\omega(c_0)$ we can apply Lemma \ref{chain pullback} to find  $\lambda\in (0,1)$, depending only on $\rho$, so that $U\subset D_{\lambda\pi/2}(\LL_x(J)).$
Taking $C$ be big enough we get that
 $$D_{\lambda\pi/2}(\LL_x(J))\subset D_{\pi/2}((1+2C)\mathcal{L}_{x}(J))\subset D_{\pi/2}(J).$$
In this case the first part of the proposition follows taking $\nu>1/\varepsilon$, where $\varepsilon$ is the constant associated to $C$ in Lemma \ref{big scaling factors real bounds}. Let us now assume that $c_0$ is odd; then all the critical points in $\omega(c_0)$ are odd. Since $J$ is $\rho'$-externally free, $\tilde I:=(1+\rho')J$ is $\rho''$-free for some $\rho''>0$ that depends only on $\rho$. Note that since $J$ is $\rho'$-externally free $\LL_x(\tilde I)\subset J$. Let $V=D_{\pi/2}(\tilde{I)}$ and $U=\comp_{x} f^{-s}(V).$ Since $\tilde I$ is $\rho''$-free we can apply Lemma \ref{chain pullback} to find  $\lambda'\in (0,1)$, depending on $\rho$, so that $U\subset D_{\lambda'\pi/2}(\LL_x(\tilde I))$. Since all the critical points in $\omega(c_0)$ are odd $f^s|_{\LL_x(J)}$ and $f^s|_{\LL_x(\hat I)}$ (the return maps to $J$ and $\tilde I$) are monotone and there exists $\kappa>0$, depending on $\rho$, the number of critical points of $f$ and their orders, so that $ |\mathcal{L}_{x}(J)|>\kappa|\LL_x(\tilde I)|.$ This means that there exists $C'>0$ depending on $\kappa$ and $C$, with $C'\to \infty$ as $C\to \infty$, so that 
$$(1+2C')\LL_x(\tilde I)\subset J,$$
since $(1+2C)\LL_x(J)\subset J$. Taking $C$ big enough we get
$$D_{\lambda'\pi/2}(\LL_x(\tilde I)) \subset D_{\pi/2}((1+2C')\LL_x(\tilde I))\subset D_{\pi/2}(J).$$
In this case the first part of the proposition follows taking $\nu>1/\varepsilon$, where $\varepsilon$ is the constant associated to $C$ in Lemma \ref{big scaling factors real bounds}. 

\medskip
Let us now prove the second statement. By Theorem \ref{real geometry} (a),  $I_{i}$ is $\rho$-externally free. Moreover, if $c_0$ is even we know that $I_{i}\setminus (1+2\rho)^{-1}I_{i}$ is disjoint from $\omega(c_0)$. Fix $x\in\omega(c_0)\cap I_{i}$. We know the return time of $x$ to $I_{i}$ is the same as its return time to $(1+2\rho)I_i$; let this common return time be equal to $q$.  Let $\tilde{F}_{n}$ be defined as in the proof of Lemma \ref{big scaling factors real bounds}. From Lemma \ref{big scaling factors real bounds} we know that there exists $C>0$ so that $I_i$ is $C-$strongly nice and that $$\LL_x(I_i)\subset (1+2C)\LL_x(I_{i})\subset \LL_x(\Gamma^T\mathcal B(\tilde F_n))\subset I_i.$$
Note that $(1+C)\mathcal{L}_{x}(I_i)$ is compactly contained in
$\LL_x(\Gamma^{T}\mathcal{B}(\tilde{F}_{n}))$. Recall that
given $x,y\in I_i$ so that $\LL_x(I_i)$ is different from
$\LL_y(I_i)$, then $\LL_x(\Gamma^{T}\mathcal{B}(\tilde{F}_{n}))$ is
disjoint from $\LL_y(\Gamma^{T}\mathcal{B}(\tilde{F}_{n})).$
Assume $c_0$ is even and let $V=D_{\pi/2}(I_i)$ and $U=\comp_{x} f^{-q}(V).$ Since, in this case, $I_i\setminus (1+2\rho)^{-1}I_i$ is disjoint from $\omega(c_0)$ we can apply Lemma \ref{chain pullback} to find  $\lambda_1\in (0,1)$ so that $U\subset D_{\lambda_1\pi/2}(\LL_x(I_i)).$ Taking $C$ sufficiently large we get that $$D_{\lambda_1\pi/2}(\LL_x(I_i))\subset D_{\pi/2}((1+C)\mathcal{L}_{x}(I_i)),$$
 so the
proposition follows in this case. If $c_0$ is odd we let
$\tilde{I}=(1+\rho)I_{i}$. By definition, $\tilde{I}$ is $\rho_1$-free
for some $\rho_1>0$. Let  $V=D_{\pi/2}(\tilde I)$ and $U=\comp_{x}
f^{-q}(V).$  By Lemma \ref{chain pullback} there exists a constant
$\lambda_2\in (0,1)$ so that $U\subset D_{\lambda_2\pi/2}(\LL_x(\tilde
I)).$ We know that $ |\mathcal{L}_{x}(I_i)|>\kappa|\LL_x(\tilde I)|.$ This means that there exists $C'>0$ depending on $\kappa$ and $C$, with $C'\to \infty$ as $C\to \infty$, so that 
$$(1+2C')\LL_x(\tilde{I})\subset \LL_x(\Gamma^T\mathcal B(\tilde F_n)).$$
Making $C$ big enough we get
$$D_{\lambda_2\pi/2}(\LL_x(\tilde I))\subset D_{\pi/2}((1+C')\mathcal{L}_{x}(\tilde I)).$$
Since $(1+C')\mathcal{L}_{x}(\tilde I)$ is compactly contained in $\LL_x(\Gamma^{T}\mathcal{B}(\tilde{F}_{n}))$ the proposition follows.
\end{pf}

\medskip

\begin{cor}\label{big bounds boxmapping}
There exist $\varepsilon>0$ and $\nu_1\geq 1$ such that the following holds. Suppose that $I_{i-2}$ and $I_{i-1}$ are non-terminating, $|I_{i-2}|<\varepsilon,$ and $$\frac{|I_{i-2}|}{|I_{i-1}|}>\nu_{1}.$$ Let $\hat I_i=\cup \hat{ \mathcal L}_c(I_i),$ where the union is taken over all $c\in \crit(f)\cap\omega(c_0).$ Then, the first return map to $\hat I_i$ restricted to the domains intersecting $\omega(c_0)$ extends to a qr box mapping $F:\mathcal U\to \mathcal V$, where $F|_{U}$ is at most unicritical for all $U\in \mathcal U.$ Furthermore, there exists a universal constant $\rho>0$ so that each component in $\mathcal U$ has $\rho$-bounded geometry.
\end{cor}

\begin{pf}
Let $\nu_1\geq 1$ from Proposition \ref {big space gpl}. Observe that since $I_i$ is $\rho$-externally free, taking $\nu_1$
larger if necessary, we can guarantee the following. For
each $c\in\crit(f)\cap\omega(c_0), c\neq c_0$,
let $r_c>0$ be minimal so that $f^{r_c}(c)\in I_i$,
$V_{c_0}=D_{\pi/2}(I_i)$ and $V_c=\comp_c f^{-r_{c}}(V_0)$.
Then for any two distinct $c,c'\in\crit(f)\cap\omega(c_0)$, we have that
$V_c\cap V_{c'}=\emptyset.$
Thus the qr box mapping given by Proposition~\ref{big space gpl}
induces a qr box mapping $F\colon \mathcal U\rightarrow\mathcal V$
where $\mathcal V$ is a neighbourhood of $\crit(f)\cap\omega(c_0)$, where the map $F|_{U}$ has at most one critical point for each $U\in \mathcal U$.
Recall that there exists a constant $\delta>0$ such that
either $I_i$ is $\delta$-free or all the critical points in
$\omega(c_0)$ are odd. By Corollary~\ref{cor:diffeo lower bounds}
and
Lemma~\ref{lem:z^l lower bound},
there exists a constant
$\theta'=\theta'(\delta)\in(\pi/2,\pi)$ such that if $I_i$ is
$\delta$-free, then
each component $U$ of $\mathcal U$ 
contains $D_{\theta'}(U\cap\mathbb R)$.
Whether $I_i$ is $\delta$-free or not, since $F$ is a quasiregular map
with bounded degree and bounded qc distortion, there exists $\delta'>0$
such that each component of $U$ has $\delta'$-bounded geometry.

\end{pf}


\begin{lem}\label{big space pl}
There exist constants $\nu_1\geq 1$ and $\varepsilon>0$ such that the
following holds. Suppose that $I_{i-1}$ is terminating, 
$|I_{i-1}|<\varepsilon$ and $$|I_{i-1}|/|I_{i-1}^{\infty}|>\nu_1.$$
Then the first return map to $I_{i-1}^{\infty}$ extends to a qr box mapping with range $D_{\pi/2}(I^{2}_{i-1}).$
\end{lem}
\begin{pf}
Let $r$ be the return time of $c_0$ to $I_{i-1}$. By Lemma \ref{bounds1} we know that given $\delta_1>0$ there exists $\kappa>0$ so that if $|I_{i-1}^2|>\delta_1|I_{i-1}|$, then $|I^\infty_{i-1}|>\kappa |I_{i-1}|$. This means that taking $\nu_1$ small enough there exists $C'>0$, $C'\to \infty$ as $\nu_1\to \infty$ so that $(1+2C')I^2_{i-1}\subset I_{i-1}.$ Since $I_{i-1}$ is terminating, $c_0$ is an even critical point, either $I_{i-1}^1$ is deep-inside $I_{i-1}$ or $I_{i-1}^2$ is deep-inside $I_{i-1}^1$. In any case we get that $I_{i-1}^2$ is deep-inside $I_{i-1}^1$ and $I_{i-1}^3$ is deep-inside $I_{i-1}^2$. In particular $I^2_{i-1}$ is $\delta$-free for some $\delta>0$. Note that we can fix $\delta$ so that  $I^2_{i-1}$ is $\delta$-free for all $\nu_1>C_0$ for some $C_0>0$. Let $V=D_{\pi/2}(I^2_{i-1})$ and $U=\comp_{c_0} f^{-r}(V)$. By Lemma~\ref{chain pullback} there exists $\lambda\in(0,1)$, depending on $\delta$, so that 
$$U\subset D_{\lambda\pi/2}(I^3_{i-1}).$$
Making $\nu_1$ big enough we can guarantee that $(1+2C'')I^3_{i-1}\subset I^2_{i-1}$ for some $C''=C''(\nu_1)>0$ so that 
$$ D_{\lambda\pi/2}(I^3_{i-1})\subset  D_{\pi/2}((1+C'')(I^3_{i-1}))$$
note that $ (1+C'')(I^3_{i-1})$ is compactly contained in $I_{i-1}^2,$ so the proposition follows.
\end{pf}

A direct consequence of the previous lemma we have:

\begin{cor}\label{cor:big space pl}

There exist constants $\nu_1\geq 1$ and $\varepsilon>0$ such that the following holds. Suppose that $I_{i-1}$ is non-terminating, $I_i$ is terminating, $|I_{i-1}|<\varepsilon$ and $$|I_{i-1}|/|I_{i}|>\nu_1.$$
Then the first return map to $I_{i}^{\infty}$ extends to a qr box mapping with range $D_{\pi/2}(I^{2}_{i}).$
\end{cor}
\begin{pf}
Since $I_i$ is terminating $c_0$ is even and $I_i$ and $I_i\setminus (1+\delta)^{-1}I_i\cap\omega(c_0)$ is empty. Since $I_{i-1}$ is non-terminating we can apply
Corollary \ref{nice} to conclude $I^1_i$ is deep-inside $I_i$. Hence there exists $C>0$ large so that $(1+2C)I_i^1\setminus I_i$ is disjoint from $\omega(c_0),$ and the proof follows as in the previous lemma.
\end{pf}

Combining the results from this section we obtain the following:

\begin{prop}\label{big space term}
There exist constants $\nu_1\geq 1$ and $\varepsilon>0$ such that the following holds. Suppose that $I_{i-1}$ is terminating, $|I_{i-1}|<\varepsilon$ and $$|I_{i-1}|/|I_{i}|>\nu_1.$$  Then we have one of the following:
\begin{enumerate}
\item the first return map to $I_{i-1}^{\infty}$ extends to a qr box mapping with range $D_{\pi/2}(I^{2}_{i-1})$;
\item $I_i$ is terminating and the first return map to
  $I_{i}^{\infty}$ extends to a qr box mapping with range $D_{\pi/2}(I^{2}_{i});$
\item $I_i$ and $I_{i+1}$ are non-terminating and the first return map
  to $I_{i+1}$ restricted to the components of the domain intersecting
  $\omega(c_{0})$ extends to a qr box mapping with range $D_{\pi/2}(I_{i+1});$
\item $I_i$ is non-terminating, $I_{i+1}$ is terminating and the first
  return map to $I_{i+1}^{\infty}$ extends to a qc  box mapping with range $D_{\pi/2}(I^{2}_{i+1}).$
\end{enumerate}
\end{prop}

\begin{pf}
Note that since $I_{i-1}$ is terminating, $c_0$ is even. 
If $I_{i-1}^\infty$ is small compared to $I_{i-1}$, by Lemma \ref{big space pl}, \textit{(1)} holds. 
So assume the size of $I_{i-1}^\infty$ is comparable to the size of $I_{i-1}$.
By Lemma \ref{bounds1}, the size of $\mathcal R(I_{i-1})$ is comparable to the size of $I_{i-1}$. In particular, $\mathcal R(I_{i-1})$ is not periodic, for otherwise we would have that $|I_{i}|$ is comparable to $|I_{i-1}|$. 
Since $I_{i}$ is deep-inside $\mathcal R(I_{i-1})$,
Corollary \ref{nice} implies that $I_{i}$ is $C$-nice for some $C>0$, with $C\to \infty$ as $\nu_1\to \infty$. 
If $I_{i}$ is terminating, we can take $C$ so big that
$(1+2C')I_{i}^\infty\subset I_{i}$ with $C'\rightarrow\infty$ as
$C\rightarrow\infty$, and \textit{(2)} follows from Lemma \ref{big
  space pl}. If $I_i$ is non-terminating,
by Theorem \ref{real geometry} (d), $I_{i+1}$ is $C''$-nice and $C''$-externally free, with $C''\to \infty$ as $C\to \infty$. If $I_{i+1}$ is terminating,  (4) follows from Corollary \ref{cor:big space pl}. If $I_{i+1}$ is non- terminating we argue as follows. Let $x\in\omega(c_0)\cap I_{i+1},$ with return time to $I_{i+1}$ equal to $r$.  Let $V=D_{\pi/2}(I_{i+1})$ and $U=\comp_xf^{-r}(V).$ By \ref{chain pullback} we know there exists $\lambda\in(0,1)$ so that $U\subset D_{\lambda\pi/2}(\LL_x(I_{I_{i+1}}))$. Making $C''$ sufficiently large we can guarantee $$D_{\lambda\pi/2}(\LL_x(I_{i+1}))\subset D_{\pi/2}((1+C'')\mathcal{L}_{x}(I_{i+1}))\subset V.$$ 
Since the intervals $(1+C'')\mathcal{L}_{x}(I_{i+1})$ are pairwise disjoint \textit{(3)} holds.
\end{pf}

The following result will be used in Section \ref{sec:box mapping} to prove complex bounds for renormalizable maps.

\section{Extension to a qc quasi-box mapping}\label{sec:ext to quasi-box mapping}
The results in this section apply to maps $f\in\mathcal A_{\underline
  b}$ and  to asymptotic holomorphic extensions of maps  $f\in\mathcal
A^3_{\underline b}$, as in Subsection~\ref{sub-sec:asympholoS}. We
recall that we assume that the Standing Assumptions on page
\pageref{standing assumptions} hold.


 In this section, we will prove that we can associate to $f$ a qc quasi-box mapping (see Theorem~\ref{quasi-box mapping}). In particular,  if $f$ is a real analytic map the qc quasi-box mapping is quasi-box mapping  (Definition~\ref{defn:quasi-box mapping})

\subsection{Some inductive statements used in the construction of a quasi-box mapping}
Let $\{I_i\}_{i\ge 0}$ be the generalized enhanced nest, defined in 
Section~\ref{sec:enhanced nest}. In the following, we will consider points $x\in I_n\cap\omega(c_0)$ and 
their corresponding first return domains $\LL_x(I_n).$ Our aim is to show that the first return map to $I_n$ restricted to the union of these components can be extended to a qc quasi-box mapping with range 
$D_\theta(I_{n-M})\cap \C_{I_n}$, for  some $M>0$ universal.

\medskip
Recall that, if  the interval $I_{n-j-1}$ is non-terminating there exists $p_{n-j-1}>0$ so that  $I_{n-j}$ is a pullback of $I_{n-j-1}$ by $f^{p_{n-j-1}}$.

\begin{prop}\label{prop:s-seq}
For every $\nu\geq 1$ there exists $\mu=\mu(\nu)$ and for each $\theta\in(0,\pi)$ there exists $\varepsilon>0$ such that the following holds. Assume  $n>5$ and $|I_n|<\varepsilon$. 
 Let $J=\LL_x(I_n)$ for a point $x\in I_n\cap \omega(c_0)$, and let
 $s>0$ be its first return time to $I_n$. Let $\{G_{j}\}_{j=0}^{s}$ be the chain with $G_{s}=I_{n}$ and $G_{0}=J$. Fix $M\in\mathbb{N}$, $4<M<n$. 
Assume that for $j=0,1,\dots,M$
$$\frac{|I_{n-j-1}|}{|I_{n-j}|}<\nu.$$
Then there exist integers 
$$s_{n-4}<\dots<s_{n-M-1}<s_{n-M}:=s$$
so that the following holds: 
\begin{enumerate}
\item $G_{s_{n-i}}\subset I_{n-i};$ 
\item if $I_{n-i}$ and $I_{n-i+1}$ are periodic, $G_{s_{n-i}}\subset I_{n-i+1}$; 
\item if $I_{n-i}$ is non-periodic, there exist at most six intervals $G_{j}\subset I_{n-i}$ with $s_{n-i}\leq j\leq s;$ 
\item if $I_{n-i}$ is periodic, there exist at most five intervals $G_{j}\subset I_{n-i+1}$ with $s_{n-i+1}\leq j\leq s.$ 
\end{enumerate}
For each $i\in\{4,\dots, M-1\}$ we let $$U_{s_{n-i-1}}= D_{\theta}(I_{n-i-1})\cap\mathbb{C}_{G_{s_{n-i-1}}}\,and\, \,U_{j}=\mathrm{Comp}_{G_{j}}f^{-(s_{n-i-1}-j)}(U_{s_{n-i-1}}),$$ for $0\leq j\leq s_{n-i-1},$ then
\begin{enumerate}[(i)]
\item if $I_{n-i-1}$ and $I_{n-i}$ are non-terminating,
$$U_{s_{n-i}}\subset D_{\mu\theta}(\mathcal{L}_{G_{s_{n-i}}}(I_{n-i}));$$

\item if $I_{n-i-1}$ is non-terminating, $I_{n-i}$ is terminating and $i>5$, then there exists an interval $K$ such that
$$U_{s_{n-i+1}}\subset D_{\mu\theta}(K),$$
where $G_{s_{n-i+1}}\subset K$ and $K$ is well-inside $I_{n-i+1}$;

\medskip
\item if $I_{n-i-1}$ is terminating, then there exists an interval $K$ such that
$$U_{s_{n-i}}\subset D_{\mu\theta}(K),$$
where $G_{s_{n-i}}\subset K$, and $K$ is well-inside $I_{n-i}$.
\end{enumerate}
\end{prop}
\begin{pf} 
The proof of this proposition  will be divided into steps. First we will construct the sequence of times $s_{n-i}$ and then we will do the pullbacks.

\medskip
\noindent\textit{Step 1: definition of $s_{n-i}.$} We state how we
will define the time $s_{n-i}$ if $s_{n-i-1}$ is defined. 

\medskip
\noindent\textit{Case (a). Assume $I_{n-i-1}$ and $I_{n-i}$ are both non-terminating and there are at least four intervals $G_j\subset I_{n-i}$ with $0\leq j\leq s_{n-i-1}$.} 
By Lemma \ref{KSS83} (1), we can choose  $s_{n-i}'<s_{n-i-1}$ maximal
with $G_{s_{n-i}}\subset I_{n-i}$ and $s_{n-i-1}-s_{n-i}'>p_{n-i-1}.$
By maximality of $s_{n-i}'$ and Lemma \ref{KSS83} there are at most
three times $j$, $s_{n-i}'\leq j\leq s_{n-i-1}$ such that
$G_{j}\subset I_{n-i}$ and none of the intervals $G_j$ with
$s_{n-i}\leq j\leq s_{n-i-1}$ are contained in $I_{n-i+1}.$ Since
there are at least four intervals $G_j\subset I_{n-i}$ with $0\leq
j\leq s_{n-i-1}$ we can define  $s_{n-i}$ to be maximal with
$s_{n-i}<s_{n-i}'$ and $G_{s_{n-i}}\subset I_{n-i}$. Observe that there are at most four intervals $G_{j}\subset I_{n-i}$ with $s_{n-i}\leq j\leq s_{n-i-1}$ and one of them is $G_{s_{n-i}}$.

\medskip
\noindent\textit{Case (b). Assume $I_{n-i-1}$ is non-terminating,
  $I_{n-i}$ is terminating and there are at least three intervals $G_j\subset I_{n-i}$ with $0\leq j\leq s_{n-i-1}$.}  Let $r$ be the return time of $c_0$ to $I_{n-i}$. Since there are at least three intervals $G_j\subset I_{n-i}$ with $0\leq j\leq s_{n-i-1}$, $s_{i-1}>2r$ and we can define $s_{n-i}=s_{n-i}'$, where $s_{n-i}'$ is as in Proposition~\ref{folding lemma}.  Observe that with this definition there are at most two intervals $G_{j}\subset I_{n-i}$ with $s_{n-i}\leq j\leq s_{n-i-1},$ and one of them is $G_{s_{n-i}}$.

\medskip

\noindent\textit{Case (c). Assume $I_{n-i-1}$ is terminating, $I_{n-i}\neq \mathcal R(I_{n-i-1})$ and  there are at least four intervals $G_j\subset I_{n-i}$ with $0\leq j\leq s_{n-i-1}$. }  
From Combinatorial Remark 6 we can define $r_{n-i-1}<s_{n-i-1}$ as in
Proposition~\ref{terminating angle control} so that
there are at most three intervals $G_{j}\subset I_{n-i}$
with $r_{n-i-1}\leq j\leq s_{n-i-1}.$ Thus we can define
$s_{n-i}<r_{n-i}$ to be maximal such that $G_{s_{n-i}}\subset
I_{n-i}$. Even more, if there are four intervals $G_j\subset I_{n-i}$ with $s_{n-i}\leq j\leq s_{n-i-1}$ we must have that $G_{s_{n-i-1}}$, $G_{r_{n-i-1}}$ and $G_{s_{n-i}}$ are in $I_{n-i}$. 
\medskip

\noindent\textit{Case (d). Assume $I_{n-i-1}$ is terminating, $I_{n-i}=\mathcal R(I_{n-i-1})$ and there are at least three intervals $G_j\subset I_{n-i+1}$ with $0\leq j\leq s_{n-i-1}$. }  In this case, we define $r_{n-i-1}<s_{n-i-1}$ as in Proposition~\ref{terminating angle control}. If  $I_{n-i+1}\neq \mathcal R(I_{n-i})$  or  $I_{n-i+1}= \mathcal R(I_{n-i})$  and $G_{r_{n-i-1}}\subset I_{n-i+1}$ let $s_{n-i}=r_{n-i-1}$. Observe that with this definition there are at most three intervals $G_{j}\subset I_{n-i+1}$ with $s_{n-i}\leq j\leq s_{n-i-1}.$ If $I_{n-i+1}= \mathcal R(I_{n-i})$ and $G_{r_{n-i-1}}\nsubseteq I_{n-i+1}$, then by Proposition~\ref{terminating angle control} there are at most two intervals $G_{j}\subset I_{n-i+1}$ with $s_{n-i}\leq j\leq s_{n-i-1}.$ If we let $q$ be the return time of $c_0$ to $I_{n-i}$ we must have that $G_{r_{n-i-1}-q}\subset I_{n-i+1}$. So in this case we define $s_{n-i}=r_{n-i-1}-q$. Observe that with this definition there are at most three intervals $G_{j}\subset I_{n-i+1}$ with $s_{n-i}\leq j\leq s_{n-i-1}.$ Even more, if there exist three intervals $G_j\subset I_{n-i+1}$ with $s_{n-i}\leq j\leq s_{n-i-1}$ we must have that $G_{s_{n-i-1}}$ and $G_{s_{n-i}}$ are in $I_{n-i+1}$.  

\medskip
\noindent\textit{Step 2: construction of the sequence $s_{n-4}<\dots <s_{n-M}$.}
We will construct the sequence of times $s_{n-i}$ inductively. To do so we need to check that once $s_{n-i-1}$ has been defined we still have enough visits to $I_{n-i}$ or $I_{n-i+1}$, depending on the case, to define $s_{n-i}$ as in Step 1. Note that if $i+1>4$  there exist at least 16 intervals $G_j\subset I_{n-i}$ and that there are at least 8 intervals $G_j\subset I_{n-i+1}$. Let us prove the basis of the induction. Since $M>4$ there exist at least 16 intervals $G_j\subset I_{n-M+1}$ and at least 8 of them are contained in $I_{n-M+2}.$ So we can define $s_{n-M+1}$ using the corresponding case from Step 1. Observe that if $I_{n-M+1}$ is different from $\mathcal R(I_{n-M})$ there exist at most four intervals $G_j\subset I_{n-M+1}$ with $s_{n-M+1}\leq j\leq s.$ This shows the induction step in this case. If $I_{n-M+1}\neq \mathcal R(I_{n-M})$, then by Step 1 we know that there exist at most three intervals $G_j\subset I_{n-M+2}$ with $s_{n-M+1}\leq j.$ Since at least 8 intervals $G_j\subset I_{n-M+2}$ we can define $s_{n-M+2}<s_{n-M+1}$ maximal with $G_{s_{n-M+2}}\subset I_{n-M+2}$. In this case there exist at most four intervals $G_j\subset I_{n-M+2}$ with $s_{n-M+2}\leq j\leq s$ and the induction step follows. The statement of the proposition says the following: if $I_{n-i}$ is non-periodic, there exist at most six intervals $G_{j}\subset I_{n-+1i}$ with $s_{n-i+1}\leq j\leq s;$ if $I_{n-i}$ is periodic, there exist at most six intervals $G_{j}\subset I_{n-i+1}$ with $s_{n-i}\leq j\leq s.$ It is easy to see that the bound is in fact six, and not four like the basis of the induction suggests, by constructing $s_{n-M+2}$ and $s_{n-M+3}$; we will leave this to the reader. 
\medskip

Assume $s_{n-i-1}< s_{n-i-2}<\ldots< s_{n-M}$ are defined, and that they definition has been made according the corresponding case in Step .1  Now we construct the time $s_{n-i}$, for $i>M-1$. Note that if we want to define $s_{n-i}$ we must have that $i+1>4.$

\medskip
\noindent\textit{Case(a'): assume $I_{n-i-1}$ and $I_{n-i}$ are non-terminating.}
Since $I_{n-i-1}$ is non-terminating the induction hypothesis tells us that there exist at most six intervals $G_j\subset I_{n-i-1}$ with $s_{n-i-1}\leq j\leq s$; one of them is $G_s$. 
Thus, there exist at most three intervals $G_j\subset I_{n-i}$ with $s_{n-i-1}\leq j\leq s.$ 
Since $G_0,G_{s}\subset I_{n-i+1}$, Lemma~\ref{KSS83} (2) tells us that there exist at least sixteen intervals $G_{j}\subset I_{n-i}$ with $0\leq j\leq s$. So there are at least thirteen intervals $G_{j}\subset I_{n-i}$ with $0\leq j\leq s_{n-i-1},$ and we can define $s_{n-i}$ as in Case $(a).$ Observe that with this definition we still have that there exist at most six intervals $G_{j}\subset I_{n-i}$ with $s_{n-i}\leq j \leq s$.

\medskip
\noindent\textit{Case(b'): assume $I_{n-i-1}$  is non-terminating and $I_{n-i}$ is terminating and $i>5$.} 
By the same argument given in the previous case, we know that there exist at most three intervals $G_j\subset I_{n-i}$ with $s_{n-i-1}\leq j\leq s.$ Since $i+1>4$  there exist at least 16 intervals $G_j\subset I_{n-i},$ so  at least 13 intervals $G_j\subset I_{n-i}$ for $0\leq j\leq s_{n-i-1}.$  This means that we can define $s_{n-i}$ as in Case $(b)$. Observe that with this definition we still have that there exist at most six intervals $G_{j}\subset I_{n-i}$ with $s_{n-i}\leq j \leq s$.

\medskip
\noindent\textit{Case(c'): Assume $I_{n-i-1}$ is terminating $I_{n-i}\neq \mathcal R(I_{n-i-1})$}.
We have two cases, depending on if  $I_{n-i-1}$ is or not periodic. 

\medskip
\noindent\textit{Case(c'-a): Assume $I_{n-i-1}$ is non-periodic.} By
the induction hypothesis and the argument given in Case (b'), we know that  at least 13 intervals $G_j\subset I_{n-i}$ for $0\leq j\leq s_{n-i-1}.$  This means that we can define $s_{n-i}$ as in Case $(c)$. Observe that with this definition we still have that there exist at most six intervals $G_{j}\subset I_{n-i}$ with $s_{n-i}\leq j \leq s$. 

\medskip
\noindent\textit{Case(c'-b): Assume $I_{n-i-1}=\mathcal R(I_{n-i-2})$.} Then, by Case $(d)$ we know that  $s_{n-i-1}=r_{n-i-2}$ is as in Proposition~\ref{terminating angle control}. By the induction hypothesis there exist at most five intervals $G_j\subset I_{n-i}$ with $s_{n-i-1}\leq j \leq s$. Since  $i+1>4$  there exist at least 16 intervals $G_j\subset I_{n-i},$ so at least 12 intervals $G_j\subset I_{n-i}$ for $0\leq j\leq s_{n-i-1}.$  This means that we can define $s_{n-i}$ as in Case $(c)$.

\medskip
\noindent\textit{Case(d'): assume $I_{n-i-1}$ is terminating  and $I_{n-i} = R(I_{n-i-1})$.}
As in the previous case we need to consider two cases, depending if $I_{n-i-1}$ is or not periodic. 

\medskip
\noindent\textit{Case(d'-a): Assume $I_{n-i-1}$ is non-periodic.}
Arguing as in the Case $(c'-a)$ we get that there exist at most three intervals $G_j\subset I_{n-i}$ with $s_{n-i-1}\leq j\leq s,$ and at most two of them are in $I_{n-i+1}$. Since $i+1>4$  there exist at least 8 intervals $G_j\subset I_{n-i+1},$ so there exist at least 6 intervals $G_j\subset I_{n-i+1}$ for $0\leq j\leq s_{n-i-1}.$ This means that we can define $s_{n-i}$ as in Case $(d)$. Observe that with this definition we still have that there exist at most five intervals $G_{j}\subset I_{n-i+1}$ with $s_{n-i}\leq j \leq s$.

\medskip
\noindent\textit{Case(d'-b): Assume $I_{n-i-1}=\mathcal R(I_{n-i-2})$.} By the induction hypothesis, we know that there exist at most five intervals $G_j\subset I_{n-i}$ with $s_{n-i-1}\leq j\leq s,$ and at most three of them are contained in $I_{n-i+1}$.  Since $i+1>4$ there exist at least 8 intervals $G_j\subset I_{n-i+1},$ so there exist at least 5 intervals $G_j\subset I_{n-i+1}$ for $0\leq j\leq s_{n-i-1}.$ This means that we can define $s_{n-i}$ as in Case $(d)$. Observe that with this definition we still have that there exist at most five intervals $G_{j}\subset I_{n-i+1}$ with $s_{n-i}\leq j \leq s$.

\medskip
So the sequence $s_{n-4}\ldots, <s_{n-M}=s$ is well defined. \medskip

\noindent\textit{Step 3: pulling back.} Now, we will make use of the
way the times $s_{n-i}$ have been defined to prove $(i),(ii)$ and
$(iii).$ For this part of the proof it will be useful to see Figures
\ref{diagram1} and \ref{diagram2}. Fix $i\in\{4,\dots, M-1\}$, and
let $$U_{s_{n-i-1}}=
D_{\theta}(I_{n-i-1})\cap\mathbb{C}_{G_{s_{n-i-1}}}\,and\,
\,U_{j}=\mathrm{Comp}_{G_{j}}f^{-(s_{n-i-1}-j)}(U_{s_{n-i-1}}),$$
where $0\leq j\leq s_{n-i-1}.$ 

\medskip
\noindent\textit{Case(a').} By the definition of $s_{n-i}$ we can apply Proposition~\ref{pullback to previous level} to prove $(i).$

\medskip
\noindent\textit{Case(b').} In this case we will prove that $(ii)$ holds. For this we need $s_{n-i-1}, s_{n-i}$ and $s_{n-i+1}$ to be defined, so we assume $i>5$. By definition of $s_{n-i}$ we can apply Proposition~\ref{folding lemma} to find an interval $K'$ well-inside $\hat I$ and $\lambda'\in(0,1)$ so that $$U_{s_{n-i}}\subset D_{\lambda'\theta}(K').$$
Since $I_{n-i}$ is terminating $s_{n-i+1}\leq r_{n-i}$ and we can apply Corollary \ref{cor:terminating angle control} to find $\lambda''>0$ and an interval $K''$ well-inside $\mathcal R(I_{n-i})$ so that 
$$U_{r_{n-i}}\subset D_{\lambda''\theta}(K'').$$
If $I_{n-i+1}=\mathcal R(I_{n-i})$, then $s_{n-i+1}=r_{n-i}$ and  $(ii)$ holds. If $s_{n-i+1}<r_{n-i}$ we apply Proposition  \ref{prop:from renormalization to landing}  to get $(ii)$.

\medskip
\noindent\textit{Case(c').}  if $I_{n-i-1}$ is terminating  and
$I_{n-i}\neq \mathcal R(I_{n-i-1})$. Then we can apply Proposition
\ref{terminating angle control} to pull back from time $s_{n-i-1}$ to
time $r_{n-i-1},$ and  Proposition~\ref{prop:from renormalization to
  landing} to pull back from time $r_{n-i-1}$ to time $s_{n-i}$, so $(iii)$ holds.

\medskip
\noindent\textit{Case(d').}  By the definition of $s_{n-i}$, $(iii)$ holds from Proposition~\ref{terminating angle control} and at most one application of Corollary \ref{chain space}. 
\end{pf} 

From now on, $s_{n-4}<\dots < s_{n-M}$ will denote the sequence defined in the above proposition and we will assume that $n$ is large.

\begin{lem}\label{lem62} 
For each $\nu\geq 1$  there exists $C>0$ and  for each $\theta\in (0,\pi)$ there exists $\varepsilon>0$ such that the following holds. Suppose that $|I_{n}|<\varepsilon$ 
and fix $M\in\mathbb{N},$ with  $4< M<n$. 
Assume that for $j=0,1,\dots,M$
$$\frac{|I_{n-j}|}{|I_{n-j+1}|}<\nu.$$ 
Given $x\in I_n\cap\omega(c_{0})$ consider the chain $\{G_{j}\}_{j=0}^{s}$ with $G_{s}=I_{n}$ and $G_{0}=\LL_x(I_n)$. For $i\in\{4,\dots, M-1\}$ let  $$w\in\mathrm{Comp}_{G_{0}}f^{-s_{n-i-1}}(V_{n-i-1})\cap\mathbb{C}_{G_{0}}\;where\;V_{n-i-1}:=D_{\theta}(I_{n-i-1})\cap\mathbb{C}_{G_{s_{n-i-1}}}.$$
Then, if $\ell$ is the order of $c_0$ the following holds:
\begin{enumerate}

\item if $I_{n-i-1}$ and $I_{n-i}$ are non-terminating,
either
$$w_{s_{n-i}}\in D_{\theta}(I_{n-i})\;or\;w\in D_{\theta'}(I_{n-i}),$$
where $\theta'=C\theta\Big(\frac{|I_{n-i}|}{|I_{n}|}\Big)^{\ell-1};$

\item if $I_{n-i-1}$ is non-terminating and $I_{n-i}$ is terminating and $i>5$, then either
$$w_{s_{n-i+1}}\in D_{\theta}(I_{n-i+1})\;or\;w\in D_{\theta'}(I_{n-i+1}),$$
where $\theta'=C\theta\Big(\frac{|I_{n-i+1}|}{|I_{n}|}\Big)^{\ell-1};$

\item if $I_{n-i-1}$ is terminating,  
then either
$$w_{s_{n-i}}\in D_{\theta}(I_{n-i})\;or\;w\in D_{\theta'}(I_{n-i}),$$
where $\theta'=C\theta\Big(\frac{|I_{n-i}|}{|I_{n}|}\Big)^{\ell-1}.$






\end{enumerate}
\end{lem}
\begin{pf}

We will first show the real bounds that will be needed in the
proof. By Corollary \ref{cor:delta free},  there exists $\delta>0$ so
that $I_{n-M}, \ldots, I_{n-1}$ are $\delta$-free.  Since we are
assuming $|I_{n-1}|$ is comparable to $|I_{n+1}|$, Corollary
\ref{cor:delta free} or Proposition~\ref{delta free}, depending on the
definition of $I_n$, imply that $I_n$ is $\delta$-free except if
$I_{n}=\mathcal R(I_{n-1})$ and it is periodic; in this case the we
need  $|I_{n-1}|$ to be comparable to $|I_{n+2}|$ to conclude that
$I_n$ is $\delta$-free, which we do not have. In this particular case,
we will use a slightly different argument involving the interval
$I_{n-1}$ and its pullbacks, instead of the intervals $G_j$. To do so
we will define a new chain as follows. If $I_{n}=\mathcal R(I_{n-1})$
and it is periodic let $\{\hat G_j\}^{s}_{j=0}$ be the chain with
$\hat G_s=I_{n-1}$ and $\hat G_j=\comp_{G_j}f^{-1}(\hat G_{j+1})$ for $0\leq j <s-1.$ Note that the intersection multiplicity of this chain is bounded by two, since the chain $\{G_j\}^{s-1}_{j=0}$ is disjoint and $I_n$ has period two inside $I_{n-1}$.  
By Corollary \ref{cor:free-space}, there exists $\hat \delta>0$ so that if $G_s=I_n$ is $\delta$-free, then $G_j$ is $\hat \delta$-free for all $0\leq j\leq s$, and if $\hat G_s$ has been defined, then $\hat G_j$ is $\hat \delta$-free for all $0\leq j\leq s$, since $\hat G_s=I_{n-1}$ is $\delta$-free. Finally, by Theorem \ref{real geometry} we know that if $I_{n-i-1}$ is non-terminating, then it is $\hat\rho$-nice for some $\hat\rho>0$.

\medskip
\noindent\textit{Step 1: Pulling back to an interval well-inside an interval in the enhanced nest.}
 Let $U_{j}=\mathrm{Comp}_{G_{j}}f^{-(s_{n-i-1}-j)}(V_{n-i-1})$ for $0\leq j\leq s_{n-i-1}.$ We will study each of the cases from the statement separately. 
\begin{enumerate}

\item If $I_{n-i-1}$ and $I_{n-i}$ are non-terminating.  By Proposition~\ref{prop:s-seq} there exists $\mu\in (0,1)$ and an interval $K$ well-inside $I_{n-i}$ such that $U_{s_{n-i}}\subset D_{\mu\theta}(K),$ where $K=\mathcal{L}_{ G_{s_{n-i} } }(I_{n-i})$ and $G_{s_{n-i}}\subset K$. In this case, let  $I^{f}_{n-i}$ be the pullback of $I_{n-i-1}$ under $f^{p_{n-i-1}-1}$ containing $f(I_{n-i})$. In this case define $i_0=n-i$ and $s^*=s_{n-i}$.

\item If $I_{n-i-1}$ is non-terminating and $I_{n-i}$ is terminating and $i>5$.  By Proposition~\ref{prop:s-seq} there exists $\mu\in (0,1)$ and an interval $K$ well-inside $I_{n-i+1}$  such that
$U_{s_{n-i+1}}\subset D_{\mu\theta}(K),$ with $G_{s_{n-i+1}}\subset K$.
Assume the return time of $c_{0}$ to $\mathcal{R}(I_{n-i})$ is equal to $q$ and let
$I^{f}_{n-i+1}$ be the pullback of $\mathcal{R}(I_{n-i})$ under $f^{q-1}$ containing $f(I_{n-i+1})$. In this case define $i_0=n-i+1$ and $s^*=s_{n-i+1}$.

\item If $I_{n-i-1}$ is terminating we have two cases.
\begin{enumerate}[(a)]
\item If $I_{n-i}\neq \mathcal R(I_{n-i-1}).$ By Proposition~\ref{prop:s-seq} there exists $\mu\in (0,1)$ and an interval $K$ well -inside $I_{n-i}$ such that $U_{s_{n-i}}\subset D_{\mu\theta}(K),$ with $G_{s_{n-i}}\subset K$. Assume the return time of $c_{0}$ to $\mathcal{R}(I_{n-i-1})$ is equal to $q$ and let
$I^{f}_{n-i}$ be the pullback of $\mathcal{R}(I_{n-i-1})$ under $f^{q-1}$ containing $f(I_{n-i})$.
In this case define $i_0=n-i$ and $s^*=s_{n-i}$.

\item If $I_{n-i}=\mathcal R(I_{n-i-1}).$ By Proposition~\ref{prop:s-seq} there exists $\mu\in (0,1)$ and an interval $K$ well -inside $I_{n-i}$ such that $U_{s_{n-i}}\subset D_{\mu\theta}(K),$ with $G_{s_{n-i}}\subset K$. Assume the return time of $c_{0}$ to $I_{n-i}$ is equal to $q$ and let $I^{f}_{n-i}$ be the pullback of $I_{n-i}$ under $f^{q-1}$ containing $f(I_{n-i})$. In this case define $i_0=n-i$ and $s^*=s_{n-i}$.
\end{enumerate}
\end{enumerate}

\medskip
Observe that, independently of the case, the definition of $i_0$ and $s^*$ give the following:
$$U_{s^*}\subset D_{\mu\theta}(K),$$

and $G_{s^*}\subset K\subset I_{n-i_0}$.
We always have that $$\comp_{c_0}f^{-1}(I^f_{i_0})=I_{n-i_0}.$$

\medskip
We know that $I_{n-i_{0}-1}$ is $\delta$-free.  If $I_{n-i_{0}-1}$ is terminating, Proposition~\ref{delta free} tells us that there exists $\delta'>0$ so that $\mathcal{R}(I_{n-i_{0}-1})$ is $\delta'$-free. So shrinking $\delta$, if necessary, we can assume that if $I_{n-i_0-1}$ is terminating, $I_{n-i_0-1}$  and  $\mathcal{R}(I_{n-i_{0}-1})$ are $\delta$-free. By definition,  $I^f_{n-i_0}$ is a pullback of bounded order of a $\delta$-free interval so, regardless of the definition of $i_0$, the interval $I_{n-i_{0}}^{f}$ is $\delta''$-free for some $\delta''>0$. In the remainder of the proof $C_{k}$, for $k\in\N$, stands for a positive constant that depends only on the previous constants. Observe that, in any case, if  $w_{s^*}\in D_{\theta}(I_{n-i_{0}})$ we are done with the proof of the proposition, so let us assume $w_{s^*}\notin D_{\theta}(I_{n-i_{0}})$. Then $$w_{s^*}\in D_{\mu\theta}(K)\setminus D_{\theta}(I_{n-i_{0}}).$$  By Lemma \ref{Poincare disks} we have that there exists $C_1>0$ so that
$$w_{s^*}\in D_{C_{1}\frac{|G_{s^*}|}{|I_{n-i_{0}}|}\theta}(G_{s^*}).$$

\noindent\textit{Step 2: Pulling back to the start.}

If $I_{n-i_0}$ is non-periodic, then by  Proposition~\ref{prop:s-seq}  we know that $f^{s-s^*}\colon G_{s^*}\rightarrow G_{s}=I_{n}$ is at most a sixth iterate of the first return map $R_{I_{n-i_{0}}}.$ So there exists $\nu_0>0$ so that $|G_{s^*}|/|I_{n}|\geq\nu_{0}$. Indeed, if $|G_{s^*}|/|I_{n}|$ is small, then the derivative of $R_{I_{n-i_{0}}}$ would be big at some point. Then, since $I_{n-i_{0}}$ is $\delta$-free Theorem \ref{real Koebe} would imply that one of the components of the domain of $R_{I_{n-i_{0}}}$ is small compared to $|I_{n-i_{0}}|.$ Theorem \ref{real geometry} or Lemma \ref{bounds1}, depending on whether $I_{n-i_{0}}$ is non-terminating or terminating, allows us to conclude that $I_{n-i_{0}+1}$ is small compared to $I_{n-i_{0}}$, which is a contradiction.  
If $I_{n-i_0}$ is periodic, then by  Proposition~\ref{prop:s-seq}  we
know that $f^{s-s^*}\colon G_{s^*}\rightarrow G_{s}=I_{n}$ is at most a
fifth iterate of the first return map $R_{I_{n-i_{0}+1}},$ 
so we conclude that $|G_{s^*}|/|I_{n}|\geq\nu_{0}$ using the same
argument as before. Similarly, if $I_{n}=\mathcal R(I_{n-1})$ and it
is periodic, we have that $|\hat G_{s^*}|/|I_{n-1}|\geq\nu_{0}$.

\medskip 
Assume $I_{n}\neq \mathcal R(I_{n-1})$. Then, since $G_{s^*}$ is $\delta$-free and the intervals $G_{0},\dots, G_{s^*}$ are disjoint  we have that 
$$w_{1}\in D_{C_{2}\frac{|G_{s^*}|}{|I_{n-i_{0}}|}\theta}(G_{1}).$$
Because $|G_{s^*}|/|I_{n}|\geq \kappa_{0}$, we have 
\begin{equation}\label{eqn:start92}
w_{1}\in D_{C_{3}\frac{|I_{n}|}{|I_{n-i_{0}}|}\theta}(G_{1}).
\end{equation}
By Lemma \ref{Poincare disks} we have that 
$$w_{1}\in D_{C_{4}\frac{|I_{n}|}{|I_{n-i_{0}}|}\frac{|I^{f}_{n-i_{0}}|}{|G_{1}|}\theta}(I^{f}_{n-i_{0}}).$$
Since $I_{n-i_{0}}^{f}$ is $\delta''$-free, Lemma \ref{lem42} gives us that 
$$w\in D_{C_{5}\frac{|I_{n}|}{|I_{n-i_{0}}|}\frac{|I^{f}_{n-i_{0}}|}{|G_{1}|}\theta}(I_{n-i_{0}}).$$
Note that $G_{0}\subset I_{n}$, and therefore $|I_{n}|^{\ell}/|G_{1}|\geq C_{6}.$  Finally, since $I_{n-i_{0}}$ is $\delta$-free we have that $|I^{f}_{n-i_{0}}|\geq C_{7}|I_{n-i_{0}}|^{\ell}.$ Hence 
$$w\in D_{C_{8}\frac{|I_{n-i_{0}}|^{\ell-1}}{|I_{n}|^{\ell-1}}\theta}(I_{n-i_{0}}).$$

Assume $I_{n}=\mathcal R(I_{n-1})$ and it is periodic. Then using the same argument as above, substituting $G_j$ by $\hat G_j$ and $I_n$ by $I_{n-1}$ we have the following:

$$w\in D_{C_{8}\frac{|I_{n-i_{0}}|^{\ell-1}}{|I_{n-1}|^{\ell-1}}\theta}(I_{n-i_{0}}).$$
Now, since $|I_n|> \nu |I_{n-1}|$ there exists $C_9\in(0,1)$ so that 
$$w\in D_{C_{9}\frac{|I_{n-i_{0}}|^{\ell-1}}{|I_{n}|^{\ell-1}}\theta}(I_{n-i_{0}}).$$
Taking $C$ as the minimum between $C_8$ and $C_9$ the proposition follows.

\end{pf}

\subsection{The construction of a qc quasi-box mapping}

\begin{thm}[The existence of a quasi-box mapping]\label{quasi-box mapping}
For each $\nu\geq1$ there exists $M>0$, and for each  $\theta\in (0,\pi/2]$
there exists $\varepsilon>0$ so that, if $| I_{n-M} |<\varepsilon$ the following holds. If for all $j\in \{n-M, \ldots n+1\}$ we have that $$|I_{j-1}|/|I_{j}|\leq \nu,$$  then the first return mapping to $I_{n}$ extends to a quasi-box mapping with range $D_{\theta}(I_{n-M}).$
Furthermore, there exists a universal constant $\nu_0\geq 1$ so that, if $|I_{j-1}|/|I_{j}|>\nu_{0}$, for $j$ sufficiently large, then we obtain a quasi-box mapping extending the return map to a puzzle piece $\hat I$, with $I_{j+2} \subset \hat I\subset I_{j-1}$.
\end{thm}
\noindent\textbf{Remark.} 
See the corollary following the proof of the theorem for a more precise statement
in the non-renormalizable case.

Let us recall that the domains of the qc quasi-box mapping may intersect each other and that they do not have to be compactly contained in the range.
\begin{pf}
We will divide the proof of the theorem in steps. In this proposition
when we talk about first return maps to an interval $I$ containing
$c_0$, we will be referring to the first return map to $I$ restricted
to the components of $R_I$ that intersect $\omega(c_0)$. 
In what follows, assume that all intervals $I_j$ from the
  generalized enhanced nest have $j$ large enough so that
  $|I_{j}|<\varepsilon$,
for all $\varepsilon$ in this section and in 
Section~\ref{sec:pullback along enhanced nest}.

\medskip
\noindent\textit{Step 1.} 
Let $\nu\geq 1$ and let $C$ be the constant associated to $\nu$ by 
Lemma~\ref{lem62}.
Choose an integer $M$, with $4<M$ so that for all $n$ sufficiently large $$D_{\mu^3\theta}(I_{n-4})\subset D_{\theta}(I_{n-M+1}),$$ where $\mu\in(0,1)$ is the constant from Proposition~\ref{prop:s-seq}, and so that 
for  each $i=n-M,\dots,n-1$ either
\begin{equation}
C(|I_{i}|/|I_{n}|)^{\ell-1}\geq 1 \mbox{ or }
D_{C \left(\frac{|I_{i}|}{|I_{n}|}\right)^{\ell-1}\theta}(I_{i})\subset D_{\theta}(I_{n-M+1}),\label{eq:inclu}\end{equation} 
where $\ell$ is the order of $c_0$. To see that such $M$ exists notice that because of real bounds there exists $\rho>0$
so that $(1+2\rho)I_{j}\subset I_{j-1}$, for every $j\in \N.$ So we can take $M$ large enough, so that $(M+1)\rho>1/C$. 

Assume that for all $j\in \{n-M, \ldots n+1\},$ 
$$|I_{j-1}|/|I_{j}|\leq \nu.$$ 
Let $x\in\omega(c_0)\cap I_n$ and let $s>0$ be minimal so that
$f^s(x)\in I_n$. Let $\{G_j\}_{j=0}^s$ be the chain with $G_s=I_n$ and
$G_0=\mathcal{L}_x(I_n)$.
Let $s_{n-4}<\dots<s_{n-M-1}<s_{n-M}:=s $ be the sequence defined in Proposition~\ref{prop:s-seq}. Let $$w\in\mathrm{Comp}_{G_{0}}f^{-s_{n-M}}(V_{n-M})\cap\mathbb{C}_{G_{0}}\;\text{where}\;V_{n-M}:=D_{\theta}(I_{n-M})\cap\mathbb{C}_{G_{s_{n-i-1}}}$$
and define $w_{s_{n-i}}= f^{s_{n-i}}(w)$ for $i\in\{4,\dots, M-1\}$. By Lemma \ref{lem62} we know that one of the following holds:

\begin{enumerate}

\item If $I_{n-M}$ and $I_{n-M+1}$ are non-terminating,
either
$$w_{s_{n-M+1}}\in D_{\theta}(I_{n-M+1})\;or\;w\in D_{\theta'}(I_{n-M+1}),$$
where $\theta'=C\theta\Big(\frac{|I_{n-M+1}|}{|I_{n}|}\Big)^{\ell-1}$.
\item if $I_{n-M}$ is non-terminating and $I_{n-M+1}$ is terminating, then either
$$w_{s_{n-M+2}}\in D_{\theta}(I_{n-M+2})\;or\;w\in D_{\theta'}(I_{n-M+2}),$$
where $\theta'=C\theta\Big(\frac{|I_{n-M+2}|}{|I_{n}|}\Big)^{\ell-1}$;

\item If $I_{n-M}$ is terminating either 
$$w_{s_{n-M+1}}\in D_{\theta}(I_{n-M+1})\;or\;w\in D_{\theta'}(I_{n-M+1}),$$
where $\theta'=C\theta\Big(\frac{|I_{n-M+1}|}{|I_{n}|}\Big)^{\ell-1}$;
\end{enumerate}
By the choice of $M$, if  the second assertion holds, on the corresponding case, we get that $w\subset D_\theta(I_{n-M})\cap \C_{I_n},$ and the theorem holds. If the first assertion holds, we repeat the argument. If after at most $M-3$ repetitions we have not proved the theorem we get 
$$w_{s_{n-4}}\in D_{\theta}(I_{n-4}), \, \, \text{where} \, \, \theta'=C\theta\Big(\frac{|I_{n-4}|}{|I_{n}|}\Big)^{\ell-1}.$$

We know $|I_{i-2}|/|I_{i+1}|<\nu_0$ for $i\in \{n-4,\ldots, n-1\}$.
So we can follow the corresponding proofs in Section \ref{sec:pullback
  along enhanced nest} to control the loss of angle and obtain 
\begin{equation}\label{eqn:endseq}
w\in D_{\mu^3\theta}(I_{n-4}).
\end{equation}
For example, see the proof of Theorem~\ref{thm:complex bounds infinitely ren}.

\medskip

\noindent\textit{Step 2.} 
Now we show that there exists a constant $\nu_{0}\geq 1$, such that if $|I_{j-1}|/|I_{j}|>\nu_{0},$ for some $j\in \{n-M, \ldots n+1\}$
then we obtain a complex extension for a return map to a puzzle piece $\hat I$, with $$I_{j+2} \subset \hat I\subset I_{j-1}.$$ Let $\nu_{0}$ be the maximum of the constants $\nu_1$  from Section \ref{large scaling factors}. Assume that  $|I_{j-1}|/|I_{j}|>\nu_{0}$ for $j\in \{n-M, \ldots, n+1\}$. 

\medskip
\noindent\textit{Case (a): If $I_{j-1}$ and $I_j$ are
  non-terminating.} By Proposition~\ref{big space gpl} we know that
the return
mapping to $I_{j+1}$ extends to a qr box
mapping with range $D_{\pi/2}(I_{j+1}).$

\medskip 
\noindent\textit{Case (b):  If $I_{j-1}$ is non-terminating  and $I_j$
  is terminating}. By Corollary \ref{cor:big space pl} the first
return map to $I^\infty_j$ extends to a qr box mapping with range $D_{\pi/2}(I_{j}^2).$

\medskip 
\noindent\textit{Case (c):  If $I_{j-1}$ is terminating.} By
Proposition~\ref{big space term} one of the following holds:  the
first return map to $I^\infty_{j-1}$ extends to a qr box mapping with
range $D_{\pi/2}(I_{j-1}^2)$; $I_j$ is terminating and the first
return map to $I^\infty_{j}$ extends to a qr box mapping with
range $D_{\pi/2}(I_{j}^2)$; $I_j$ and $I_{j+1}$ are non-terminating
and the first return map to $I_{j+2}$ extends to a qc
complex box mapping with range $D_{\pi/2}(I_{j+2})$; $I_{j}$
is non-terminating  and $I_{j+1}$ is terminating the first return map
to $I^\infty_{j+1}$ extends to a qr box mapping with range $D_{\pi/2}(I_{j+1}^2).$

that $|I_{j-1}|/|I_{j}|>\nu_{0},$ we can construct a qc
quasi-box mapping as in the statement of the theorem. 


\end{pf}



\begin{cor}[Quasi-box mapping in the non-renormalizable
  case]\label{cor:non-ren quasi-box}
Suppose that $f$ is non-renormalizable. Given $\nu\geq 1,$ there
exists $M>0$ with the following properties. For $\theta\in(0,\pi/2],$
if $n\in\mathbb N$ is sufficiently large:
\begin{enumerate}
\item If $\nu\geq \nu_0$, where $\nu_0$ is the constant from
  Proposition~\ref{big space gpl}, and $|I_{k}|/|I_{k+1}|>\nu_0$ for
  some $k\in\{n,n-1,\dots,n-M\}$ then: the first return mapping to
  $I_{k+2}$ extends to a qr box mapping $F\colon \mathcal U\to
  \mathcal V$ with range $\mathcal V=D_{\pi/2}(I_{k+2})$ and $\mathcal U\subset
  D_{\theta}(I_{k+2})$.
\item If for all $k\in\{n,n-1,\dots,n-M\}$ we have that
  $|I_{k}|/|I_{k+1}|\leq \nu$, then the first return mapping to $I_{n}$ extends to a qc quasi-box
  mapping $F\colon \mathcal U\rightarrow \mathcal V$ with range
  $\mathcal V=D_{\theta}(I_{n-M})\cap\mathbb{C}_{I_{n}}$ and $\mathcal U\subset
  D_{\theta}(I_{n-M+1})$.
\end{enumerate}
\end{cor}

This concludes the construction of the qc quasi-box mapping.

\section{Box mappings and complex bounds}\label{sec:box mapping}
In this section, we make use of the qc quasi-box mappings from
Theorem~\ref{quasi-box mapping} to construct
qr box mappings (this notion was introduced in Subsection~\ref{sub-sec:kqcboxS}).
In particular, when $f$ is real analytic the box mappings we obtain are holomorphic.

As usual, for a map $f\in \mathcal A_{\underline b}^3$  we abuse
notation and denote by $f$ its asymptotically holomorphic extension of
order 3, given in Subsection~\ref{sub-sec:asympholoS}.
 

\subsection{Finitely renormalizable maps}
We follow similar arguments to those of Levin-van Strien in \cite{LS-inflection}, which were also used in \cite{Shen}. However, we augment them to show that we can obtain complex bounds (i.e.\ bounds for the moduli of relevant annuli)
for the qr box mappings. 
Throughout this subsection we assume that $f\colon M\rightarrow M$ is at most finitely many times renormalizable.

\medskip
As an intermediate step in our construction of qr box mappings we will
make use of  \emph{smooth box mappings}, which are maps of the form
$F_S\colon \mathcal{U}_S\rightarrow \mathcal{V}_S,$ which have all
the properties of a complex box mapping (see page \pageref{def:box mapping}) except that we only require $F$ to be $C^3$.

\medskip
Let $N>0$ be maximal so that
none of the intervals $I_{N-i}$ with $i=0,1,2,3$ is terminating.
For ease of exposition, assume $N=0,$ and that all of the intervals $I_{0}\supset I_{1}\supset I_{2}\supset\dots$ of the enhanced nest about $c_0$ are non-terminating. 
\medskip

If $I$ is a nice interval containing $c_0$, we let 
$$\hat I=\bigcup_{c\in\crit(f)\cap\omega(c_0)}\mathcal{\hat L}_{c}(I)$$
and let $\mathrm{D}(\hat I)$ denote the union of all first return domains to $\hat I$ that intersect $\omega(c_0)$. We denote the first return map to $\hat I$ by $$R_{\hat I}\colon D(\hat I)\rightarrow \hat I.$$ 
Our aim is to prove the following:
\begin{thm}[Complex bounds in the non-renormalizable case]
\label{thm:complex bounds non-ren}
There exists $\delta>0$ such that the following holds.
Suppose that $$R_{\hat I_0}\colon \mathrm{D}(\hat I_0)\rightarrow \hat I_0$$
is non-renormalizable (at $c_0$). Then for all $m$ sufficiently large,
there exists a qr box mapping $F\colon \mathcal U\rightarrow \mathcal
V$ extending $R_{\hat I_m}$ with the property that $\mathcal V$ is
$\delta$-nice, $\delta$-free and $\mathcal U$ has $\delta$-bounded
geometry.
\end{thm}
The proof of the theorem above will occupy the following proposition and the next
two subsections.

\begin{prop}[cf. \cite{LS-inflection} pages 425-427]\label{prop:smoothtoanalytic}
Given $\nu\geq 1$, there exists $\varepsilon>0$, $m\in\mathbb N$, $\eta\geq 1$ and $\theta\in(0,\pi)$ such that the following holds.
Suppose that $$R_{\hat I_0}\colon \mathrm{D}(\hat I_0)\rightarrow \hat I_0$$
is non-renormalizable at $c_0$. Then given $n\in \N$ with $|I_{n-m}|<\varepsilon$ and so that $|I_{j-1}|/|I_{j}|\leq\nu,$ for all $j\in\{n-m +1,\ldots n+1\}$ the real return mapping $R_{\hat I_n}$ extends to a 
qr box mapping 
$F\colon \mathcal U\rightarrow \mathcal V$ with the following properties: the map $F|_{U}$ is at most unicritical for all $U\in \mathcal U;$ if $V$ is a  component of $\mathcal V$, then $\diam (V)\leq \eta|V\cap\mathbb R|$ and $D_{\theta}(V\cap\mathbb R)\subset V$.

\end{prop}


\begin{pf}
Let $M$ be the constant associated to  $\nu$ by Corollary~\ref{cor:non-ren quasi-box}.
We will choose the constant $m>M$ in the course of the proof;
$M$ will be independent of $n$.
How small $|I_{n-m}|$ needs to be will be determined by the proof.
As usual, let $p_n\in\mathbb N$ be so that 
$I_{n+1}=\comp_{c_0}f^{-p_n}(I_n)$.
By Theorem~\ref{real geometry},
the intervals $I_n$ are $\delta>0$ externally and internally free, for
some $\delta=\delta(\nu),$ and each component of
the internal free space is of size comparable to $|I_n|$.

\medskip

\noindent\textit{Claim: Given $\delta_1>0$, there exists $\varepsilon_1>0$ and $\theta_1\in(\pi/2,\pi)$ such that the following holds.
Suppose $I$ is a $\delta_1$-externally free nice interval, with $I\setminus (1+2\delta_1)^{-1}I$ disjoint form $\omega(c_0)$ and $|I|<\varepsilon_1.$ Let $V=D_{\theta_1}(I)$ and consider $G_0^1=\LL_x(I)$ and $G_0^2=\LL_y(I)$ 
disjoint, with $x,y\in \omega(c_0).$ Let 
$U_0^1$ and $U_0^2$ be the 
components of the landing map to
$V$ that contain $G_0^1$ and $G_0^2$, respectively. Then $U_0^1$ and $U_0^2,$  are disjoint.}

\medskip

\noindent
\textit{Proof of claim:}
Let $F$ be the interval that defines the external free space of $I$. Consider the chains $\{G_j^i\}_{j=0}^{s_i}$ associated to  $G_0^1$ and $G_0^2,$ so $G^i_{s_i}=I,$ and let $\{H_j^i\}_{j=0}^{s_i}$ be the corresponding chains associated to the pullbacks of $F$, for $i=1,2.$ Let $U_j^i=\comp_{G^i_j}f^{s_i-j}(V)$ for $j=0,\ldots s_i.$
If the order the chains $\{G_j^i\}_{j=0}^{s_i}$ is zero, then provided that $\theta_1\in(0,\pi)$ is chosen sufficiently close to $\pi,$ and $\varepsilon_1$ is small enough,
Corollary \ref{cor:diffeo lower bounds} implies that
$U^1_0$ and $U^2_0$ are contained in the geometric circles based on their real traces, so the claim follows. Otherwise, let $s<\min\{s_1,s_2\}$ be maximal so that either $G^1_s$ or $G^2_s$
contains a critical point. 
By Lemma~\ref{chain pullback},
we can choose $\theta_1\in(\pi/2, \pi)$ so that
$U_s^1\cap U_s^2=\emptyset.$ 
Observe that this can be done provided $\varepsilon_1$ is sufficiently small.
From this it follows that 
$U_0^1\cap U_0^2=\emptyset$.
\endpfclaim

\medskip

Let $\theta_1=\theta_1(\delta)$ and $\varepsilon_1(\delta)$
be the constants given by the claim, and assume $|I_{n-m}|<\varepsilon_1.$

\medskip

\noindent\textit{Step 1: Construction of a smooth box mapping at level
  $I_{n-m}$.}
  
  \textit{There exists $\theta_2\in(\theta_1,\pi)$, depending
only on  $\delta$
and a smooth box mapping 
$F_S\colon \mathcal{U}_S\rightarrow\mathcal{V}_S$ with the following properties:}
\textit{
\begin{itemize}
\item The map $F_s$ extends $R_{\hat I_{n-m}}\colon \mathrm{D}(\hat I_{n-m})\rightarrow \hat I_{n-m}$ in a neighbourhood of $\mathrm{D}(\hat I_{n-m}),$  it is asymptotically holomorphic and has the same critical points as $R_{\hat I_{n-m}};$ 
\item for each $U\in \mathcal U_S,$ $F|_U$ is at most unicritical;
\item for each $U\in \mathcal U_S,$ $D_{\theta_2}(U\cap\mathbb R)\subset U$ and $F_S=f^r$ on $D_{\theta_2}(U\cap\mathbb R)$, where $r$ is the return time of $U$ to $\mathcal V$ under $f$.
\end{itemize}}

\medskip

\noindent\textit{Step 1a: Obtaining a complex extension such that the
  components of its domain are disjoint.}
We begin by constructing an intermediate smooth mapping 
$\tilde{F}_S\colon \tilde{\mathcal{U}}_S\rightarrow\tilde{\mathcal{V}}_S$ 
where $\tilde{\mathcal{V}}_S$ is a neighbourhood of
$\crit(f)\cap\omega(c_0)$
and the components of $\tilde{\mathcal{U}}_S$ are 
disjoint, but not necessarily contained in $\tilde{\mathcal{V}}_S$.

Let $V_{c_0}=D_{\theta_1}(I_{n-m}).$ For each
$c\in\crit(f)\cap\omega(c_0)$,
let $V_c=\hat{\mathcal{L}}_c(V_{c_0}),$ and 
let $$\tilde{\mathcal{V}}_S=\bigcup_{c\in\crit(f)\cap\omega(c_0)}V_c.$$
Observe that $\tilde{\mathcal V}_S\cap\mathbb R=\hat I_{n-m}$.
For each $x\in\omega(c_0)\cap \hat I_{n-m}$,
let $\tilde U(x)=\mathcal{L}_x(\tilde{\mathcal V}_s)$ and let $$\tilde{\mathcal{U}}_S=\bigcup_{x\in\omega(c_0)\cap\hat
  I_{n-m}}\tilde U_x.$$
Notice that if $x\in \omega(c_0)\cap\hat{\mathcal L}_c(I_{n-m})$
for $c\neq c_0$, then $\tilde U(x)=\tilde V_c.$
The claim implies that the components of $\tilde{\mathcal{V}}_S$ are
pairwise disjoint, and that the components of $\tilde{\mathcal{U}}_S$
together with the components of $\tilde{\mathcal V}_S\setminus V_{c_0}$
are pairwise disjoint.
By making $\varepsilon_1$ smaller, if necessary, Lemma~\ref{lem42},
Corollary~\ref{cor:GSSProp2}, Lemma~\ref{lem:z^l lower bound} and
Corollary~\ref{cor:diffeo lower bounds}, imply that there exist
$\theta_2\in(\theta_1,\pi)$ and $\theta_2'\in(0,\pi)$,
depending only on $\delta$ and $\underline b$,
 with the following property.
For each component $U$ of  $\tilde{\mathcal U}_S$ or of $\tilde{\mathcal V}_S$,
 $$D_{\theta_2}(U\cap\mathbb R)\subset U\subset  D_{\theta_2'}(U\cap\mathbb R).$$
Define a mapping $\tilde F_S\colon \tilde{\mathcal
  U}_S\rightarrow\tilde{\mathcal V}_S$ as follows, for each component
$U$ of $\tilde{\mathcal U}_S$, let $s>0$ denote its first entry time
to $\tilde{\mathcal V}_S$ under $f$, so that $f^s(U)$ is a component
$V$ of $\tilde{\mathcal V}_S.$ Set $\tilde F_S|U=f^s|U$.

\medskip
\noindent\textit{Step 1b. Obtaining a smooth box mapping.}
For each critical point $c\in\omega(c_0),$ let $A_c$ be topological disk that properly contains 
$\tilde V_c$. We can choose these disks so that $\{A_c:c\in\omega(c_0)\cap\crit(f)\}$ 
is a collection of domains with pairwise disjoint closures such that
given $x\in\omega(c_0)\cap \hat I_{n-m}$,
if $\tilde U(x)$ intersects $A_c$ and $\tilde U(x)\neq \tilde V_c$, then $\tilde U(x)$ is compactly contained in $A_c$. 
For each such $\tilde U(x)$ let $B(x)$ be a topological disk such that
$\tilde U(x)\Subset B(x)\Subset A_c$. 
Moreover, choose the domains $B(x)$ so that they have pairwise disjoint closures. 
For each $\tilde U(x)$, let $c$ be such that $\tilde V_c=\tilde{F}_S(\tilde U(x))$ 
and extend $\tilde F_S\colon \tilde U(x)\rightarrow \tilde V_c$ to a smooth map $F_S|B(x)\rightarrow A_c$. Let 
$$\mathcal{V}_S=\bigcup_{c\in\omega(c_0)\cap\crit(f)} A_c\; \mathrm{\; \; and\;\; } 
\mathcal{U}_S=\Big(\bigcup_{x\in\omega(c_0)} B(x)\Big)\cup \Big(\bigcup_{c\in\omega(c_0)\cap\crit(f)\setminus\{c_0\}}A_c\Big).$$
Then $$F_S\colon \mathcal{U}_S\rightarrow \mathcal{V}_S$$ is a smooth box
mapping with the required properties. This concludes Step 1.

\medskip

From now on, $\theta_2$ is the constant given by Step 1.

\medskip

\noindent\textit{Step 2: Choosing $m$.}
We will choose $m$ so that, roughly, the domain of the quasi-box mapping 
given by Corollary~\ref{cor:non-ren quasi-box}
that extends the return map to $I_{n}$
is contained in the part of the plane,
close to the real line, where $F_S$ agrees with 
the corresponding iterate of $f$.
 
 From Corollary~\ref{cor:non-ren quasi-box} we know there exists $M\in \N$, depending on $\nu$, so that the first return map to $I_n$ extends to a quasi-box mapping with range $V_Q=D_{\pi/2}(I_{n-M})\cap\mathbb C_{I_{n}}$.
For each critical point $c\in\crit(f)\cap\omega(c_0)$, $c\neq c_0$,
let $r_c>0$ be minimal so that $f^{r_c}(c)\in I_{n}$.
Set $\hat V_{c_0}=V_Q$ and $\hat V_c=\comp_c f^{-r_c}(\hat V_{c_0})$.
Since $|I_{j-1}|/|I_j|<\nu,$ for
  all $n-1\leq j\leq n-M+1$, there exists
$\theta_3\in(0,\pi)$ so that
$$\hat V_c\subset D_{\theta_3}(\hat{\mathcal{L}}_c(I_{n-M}))\cap\mathbb C_{\hat{\mathcal{L}}_c(I_{n})}.$$ To see this, 
choose any $\hat V_c$ with  $c\neq c_0$.
There exists $x\in \omega(c_0)\cap I_{n}$ such that if $s>0$ is
minimal so that $f^s(x)\in I_{n}$, then there exists $0<r_x<s$ so that
$f^{r_x}(x)\in \mathcal{L}_c(I_{n})$.
Observe that the landing times of $f^{r_x}(x)$ and $c$ to $I_{n}$
under $f$ are the same; they are both $r_c=s-r_x$.
Let $\{G_j\}_{j=0}^{s}$ be the chain with $G_{s}=I_{n}$
and $G_0=\mathcal{L}_x(I_{n}),$ and
let 
$$s_{n-4}<\dots<s_{n-M}:=s$$ 
be the sequence given by Proposition~\ref{prop:s-seq}.
Let $w\in \comp_xf^{-s}D_{\pi/2}(I_{n-M})$ and $w_{i}=f^i(w)$ for
  $i\in \N.$
Either there exists $s_{n-j}$ so that $s_{n-j+1}< r_x < s_{n-j}$ or $0<r_x<s_{n-4}$.
Let us assume $s_{n-j+1}<r_x< s_{n-j}.$ By Proposition~\ref{prop:s-seq},
there exists $\mu\in(0,1)$ such that
$w_{s_{n-j}}\in D_{\mu\pi/2}(\mathcal L_{G_{s_{n-j}}}(I_{n-j}))\subset I_{n-j}$.
By Lemma ~\ref{lem61}, we know
$k(I_{n-j},\{G_{j}\}_{j=r_x}^{s_{n-j}})$ is bounded. This and the fact
that $I_{n-j}$ is $\delta$-nice, allow us to apply
Proposition~\ref{combinatorial depth angle control} and  Corollary \ref{chain space} to find
$\theta_3\in (0,\pi)$ so that $w_{r_x}\in D_{\theta_3}(\hat{\mathcal{L}}_c(I_{n-j}))$. The proof when $0<r_x<s_{n-4}$ is analogous to this one.

Now choose $m$ so that for each critical point
$c\in\crit(f)\cap\omega(c_0)$
$$D_{\theta_3}(\hat{\mathcal{L}}_c(I_{n-M+1}))\subset 
D_{\theta_2}(\hat{\mathcal{L}}_c(I_{n-m})).$$

\medskip

\noindent\textit{Step 3. Intersecting the smooth box mapping and the
  quasi-box mapping.}
Let $$\mathcal V_Q=\bigcup_{c\in\crit(f)\cap\omega(c_0)} \hat V_c,$$ where
the $\hat V_c$ are as constructed at the beginning of Step 2.
We construct a quasi-box mapping $F\colon \mathcal U_Q\rightarrow\mathcal
V_Q$
extending the real return map to $I_{n}$
as follows. 
Since the hypotheses from the second part of 
Corollary~\ref{cor:non-ren quasi-box} hold, we obtain a qc quasi-box mapping $\hat F_Q\colon \hat{\mathcal U}_Q\rightarrow \hat V_{c_0}$ extending the return map to $I_{n}.$

Now let $\mathcal U_Q=\hat{\mathcal U}_Q$. 
For each $U\in\mathcal U_Q$ choose $x\in U\cap\omega(c_0)$ and let
$s>0$ be minimal so that $f^s(x)\in\mathcal V_Q$. Set $F_Q|U=f^s|U$.

We get a smooth box mapping 
$F_{S,n}\colon \mathcal{U}_{S,n}\rightarrow \mathcal{V}_{S,n}$
extending $R_{\hat I_{n}}$ as follows.
There exists $k_{n}$ such that
$F_S^{k_{n}}|_{I_{n}}=f^{p_{n-1}+p_{n-2}+\dots+p_{n-m}}|_{I_{n}}.$
Let $V_{S,n}(c_0)=\comp_{c_0}F_S^{-k_{n}}(A(c_0)).$ For each critical point $c\in\crit(f)\cap\omega(c_0),$ let $i_c>0$ 
be minimal so that $F_S^{i_c}(c)\in I_{n}.$ Let $V_{S,n}(c)= \comp_c(F_S^{-i_c}(V_{S,n}(c_0)))$ and 
 $$\mathcal{V}_{S,n}=\bigcup_{c\in\crit(f)\cap\omega(c_0)}V_{S,n}(c).$$
If $x\in\omega(c_0)\cap\mathcal{V}_{S,n}$,
let $k\geq 1$ be minimal so that $F_S^k(x)\in \mathcal{V}_{S,n}.$ Let $V=\comp_{F_S^k(x)}\mathcal{V}_{S,n}$,
and set $U_{S,n}(x)=\comp_x(F_S^{-k}(V))$.

Finally, for each critical point $c\in\omega(c_0)$,
let $V_c=V_{S,n}(c)\cap\hat V_c$ and let $\mathcal{V}$ be the union of the domains $V_c$. 
For each $x\in\omega(c_0)\cap\mathcal{V},$ let 
$\hat U(x)$ be the component of $\hat{\mathcal U}_Q$ that contains $x$ and set
$U(x)=U_{S,n}(x)\cap\hat U(x).$ Let
$\mathcal{U}$ be the union of all domains $U(x)$.
Define $F|U(x)=\hat F|U(x)$. See Figure \ref{map intersection}. Then
$F\colon \mathcal{U}\rightarrow\mathcal{V}$ is a qr box mapping
that extends the real return map to $I_{n}$.

Since 
for all $j$, with $ 1\leq j\leq m+1$, the interval $I_{n-j}$ is $\delta$-free,
there exists $\theta\in(0,\pi),$
depending on $\theta_2$ and $\delta$,
so that  for any component $V$ of $\mathcal V$,
$V\supset D_{\theta}(V\cap\mathbb R).$
Hence there exists $\eta>0$ such that $\diam (V)<\eta|V \cap \mathbb R|$.

\begin{figure}
\resizebox{8cm}{!}{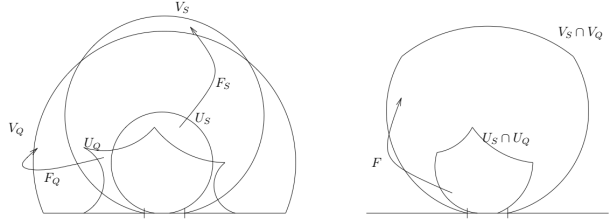}
\caption{Intersecting a smooth box mapping with a quasi-box mapping.}\label{map intersection}
\end{figure}
\end{pf}

\begin{cor}\label{box sequence}
Suppose that $R_{\hat I_0}:D(\hat I_0)\rightarrow \hat I_0$ is
non-renormalizable at $c_0$.
Then for all $n$ sufficiently large, there exists a qr box mapping 
$F_n\colon \mathcal U_n\rightarrow\mathcal V_n$ extending $R_{\hat I_n}.$
\end{cor}

\begin{pf}
Let $\nu_1\geq 1$ be the constant coming  from Corollary \ref{big bounds boxmapping}. Let $M$ be the constant associated to $\nu_1$ by Proposition  \ref{prop:smoothtoanalytic}. There are two cases: either there exists $N\in \N$ so that $|I_{j}|/|I_{j+1}|\leq \nu_1$ for all $j>N,$ or there exist infinitely many integers $j$ such that $|I_{j}|/|I_{j+1}|> \nu_1.$  Using either Proposition \ref{prop:smoothtoanalytic}, or Corollary \ref{big bounds boxmapping}, we construct $F:\mathcal U \to \mathcal V$  a qr box mapping extending $R_{\hat I_n,}$ for an integer $n$ sufficiently large. 

Let $k_1\in \N$ be such that $F^{k_1}|_{I_{n+1}}=f^{p_{n}}|_{I_{n+1}}.$ Let $V(c_0)$ be  the component of $\mathcal V$ that contains $c_0$, and define
$$V_1(c_0)=\comp _{I_{n+1}}F^{-k_1}(V(c_0)).$$
We can obtain a qr box mapping $F_1\colon \mathcal{U}_1\rightarrow\mathcal{V}_1$
that extends $R_{\hat I_{n+1}}\colon \mathrm{D}(\hat I_{n+1})\rightarrow \hat I_{n+1}$
using, for example, the argument used in the proof of Step 3 in Proposition \ref{prop:smoothtoanalytic}. Similarly, we  can define a sequence of qr box mappings
$F_i\colon \mathcal{U}_i\rightarrow\mathcal{V}_i$ each extending $R_{\hat I_{n+i}}$.
\end{pf}
We will bound the quasiconformal distortion of the maps $F_i$ in the next subsection.

To prove complex bounds, as stated in Theorem~\ref{thm:complex bounds non-ren},
it suffices to prove certain geometric bounds known as upper and lower
bounds in \cite{KSS}, see page~\pageref{page:upperlower}.
This is done in Subsections~\ref{sub-sec:upperbounds} and \ref{sub-sec:lower bounds}.

\subsection{Upper Bounds for finitely renormalizable maps}\label{sub-sec:upperbounds}

To prove Proposition~\ref{prop:upper bounds2}, we follow the approach used in the
proof of \cite[Proposition~8.3]{KSS}.

It will be convenient to point out the following modification
to the proof in \cite{KSS} now. In \cite{KSS}, an interval $I$ is
said to be in $\mathcal{T_\xi}$ if $I$ is $\xi$-nice and
$(1+2\xi)I\setminus (1+2\xi)^{-1}I$ is disjoint from the post-critical
set of $f$. We modify this definition to take into account odd
critical points: if $c_0$ is even, then we will say that an interval
$I$ \textit{is in} $\mathcal T_{\xi}$ \label{Txi} if $I$ is $\xi$-nice and $\xi$-free; if
$c_0$ is odd, then we say that $I$ \textit{is in} $\mathcal T_{\xi}$
if $I$ is $\xi$-nice and $\xi$-externally free. It is worth remarking
that this implies that $I$ is $\xi'$-strongly nice, with
$\xi'\rightarrow\infty$ as $\xi\rightarrow\infty$.

Let $F:\mathcal U\rightarrow\mathcal V$ be a qr box mapping given by 
Proposition~\ref{prop:smoothtoanalytic} or by Corollary \ref{big bounds boxmapping}
and let $V_{c_0}$ be the component of $\mathcal V$ that contains $c_0$.
For ease of exposition, assume $F$ extends the return map to $\hat
I_0=\cup_{c\in\crit(f)\cap\omega(c_0)}\hat{\mathcal
  L}_{c}(V_{c_0}\cap\mathbb R)$, and let $\I_0=V_{c_0}$. Assuming that 
$\I_{n}$ is defined, let $\I_{n+1}=\comp_{c_0}f^{-p_n}(\I_n)$,
where $p_n$ is defined as on page \pageref{p-n}.
Notice that for all $n\in \N\cup\{0\},$ $I_n=\I_n\cap\mathbb{R}$.
We define the complex puzzle pieces $\mathcal{B}(\I)$,
$\mathcal{A}(\I)$ and $\Gamma(\I)$, associated to a domain $\I$, exactly as we did for real puzzle pieces.

We will need the following version of 
Proposition~\ref{combinatorial depth angle control},
which gives us information about the loss of angle 
even when the combinatorial depth $k(I,\mathbb G)$ is large, 
see definition on page \pageref{comb depth symbols}.

\begin{prop}\cite[Proposition 11.2]{KSS}\label{prop:KSS prop 11.2}
For each $\delta>0,$ there exists $\mu\in(0,1)$ and $\delta'>0$,
and for each $\theta\in(0,\pi),$ there exists $\varepsilon>0$, such
that the following holds.
Let $I$ be a $\delta$-nice interval with $|I|<\varepsilon,$ and let $\mathbb
G=\{G_i\}_{i=0}^s$
be a disjoint chain with $G_0, G_s$ nice intervals and
$G_0\cap\omega(c_0)\neq \emptyset$. Let $\hat G_s$ be an interval with
$G_s\subset \hat G_s\subset (1+2\delta)\hat G_s\subset I.$ Let
$V=D_{\theta}(\hat G_s)$, and write $U_i=\comp_{G_i}f^{-(s-i)}(V),$
$i=0,1,\dots, s.$ Then there exists an inteval $\hat I\supset G_0$
with $(1+2\delta')\hat I\subset I$ and such that
$$U_0\subset D_{\mu^{k(I,\mathbb G)}\theta}(\hat I).$$
\end{prop}
\begin{pf}
The proof of this result proceeds by induction as in
\cite{KSS}.
To obtain the proof when $N=0$, we use the proof of
the $N=0$ case of
Proposition~\ref{combinatorial depth angle control}.
If $c$ is even, which we can and will assume whenever $c_0$ is even,
\cite[Lemma 11.2]{KSS} holds by 
Lemma~\ref{KSSLemma93}. 
When all critical points are odd,
we do not need \cite[Lemma 11.2]{KSS}
to obtain that certain intervals are $\hat \rho$-nice, since we can
use
Proposition~\ref{lem:N-modalpb}
instead of Corollary~\ref{chain space}.
The remainder of the 
proof of Proposition~11.2 in \cite{KSS} can be repeated verbatim,
where we use
Lemma~\ref{KSSLemma97} instead of \cite[Lemma 9.7]{KSS},
Theorem~\ref{real Koebe} (1) instead of \cite[Lemma 9.2]{KSS},
Lemma~\ref{KSSLemma94} instead of \cite[Lemma 9.4]{KSS},
Proposition~\ref{good start or good deep} instead of \cite[Lemma 11.1]{KSS},
and Lemma~\ref{chain pullback} and its Corollaries,
Corollary~\ref{free space well inside} and Corollary~\ref{chain space},
instead of Lemmas 13.5 and 13.6 of \cite{KSS}.
\end{pf}


The following proposition is the main tool that is
needed
to combine the bounded and big geometry cases.
Without control on the geometry, there is no way to control the combinatorial depth between consecutive levels of the enhanced nest, but we are
still able to control the geometry at deeper levels,
see the last part of the following proposition.
\label{pageref:depth}

\begin{prop}\cite[Proposition 11.3]{KSS}\label{prop:KSS 11.3}
For each $\delta>0$ and $N\geq 0$ there exist $\mu(N,\delta)\in(0,1)$
and $C\in(0,1),$ and for each $\theta\in(0,\pi)$  there exists $\varepsilon>0$ so that 
the following holds. Let $I\owns c_0$ be an interval in $\mathcal
T_{\delta}$ with $|I|<\varepsilon$.
Let $J\owns c_0$  be an (at most) $N$-modal pullback of 
$I$, and let $t$ be so that $J=\comp_{c_0}(f^{-t}(I))\cap\mathbb R.$ 
Assume $x\in J\cap \omega(c_0)$ is so that $f^s(x)\in J$ for some $s\geq t.$ Let 
$$\nu=\#\{0\leq j\leq s-t: f^j(x)\in J\}.$$
Let $s_0=0<s_1<\ldots <s_\nu$ be the times for which $s_j\leq s-t$ and
$f^{s_i}(x)\in J.$ Consider the chain $\mathbb{G}:=\{G_i\}_{i=0}^s$
defined by 
$G_s=J$, and $G_i\ni f^i(x).$
Let 
$U_s=D_{\theta}(I)\cap\mathbb{C}_{G_{s}}$ and $U_{i}=\mathrm{Comp}_{G_{i}}f^{-(s-i)}(U_s).$ Then, 
$$U_{0}\subset D_{\theta'}(J),$$ 
where $$\theta'=\min\left[ \mu^{\hat{k}(I,J)} \left( \Pi^{\nu-1}_{j=0} C\rho_j\right) \cdot \theta, \theta_0\right],$$
and ${\hat{k}(I,J)}$ is defined on page \pageref{comb depth symbols 2} and $\rho_j>0$ is so that $$(1+2\rho_j)\LL_{f^{s_i}(x)}(J)\subset J$$ and $\theta_0$ is defined in Lemma \ref{Poincare disks}.

Moreover, there exists a universal (large) constant 
$\xi>0$ and $\hat \nu\in \N$, which depends on 
$\delta$ and $N$, such that if $J$ is $\xi$-nice and 
$\xi$-externally free, $\nu\geq \hat\nu$ and 
for each $c\in \crit(f)\cap \omega({c_0}),$ different form $c_0$
$$\#\{0\leq j\leq s: f^j(x)\in \LL_c(J)\}\geq \nu_0,$$
then 

$$\theta'=\min[\theta, \theta_0].$$
\end{prop}

\begin{pf}
The proof of this is the same as the proof of \cite[Proposition
11.3]{KSS} using Proposition~\ref{prop:KSS prop 11.2} instead of
\cite[Proposition 11.2]{KSS},
Lemma~\ref{KSSLemma96} instead of \cite[Lemma 9.6]{KSS},
Lemma~\ref{Poincare disks} instead of \cite[Lemma 13.4]{KSS}, and
Lemmas~\ref{chain pullback}~and ~\ref{lem:N-modalpb} instead of Lemma 13.5.
\end{pf}

The following result is the key initial estimate needed to prove the
Upper Bounds:

\begin{prop}\cite[Theorem 11.1]{KSS}
There exist $\theta\in(0,\pi)$ and $n_0\in\mathbb N$
so that for all $n$ sufficiently large, $\I_n\subset D_{\theta}(I_{n-n_0}).$
\end{prop}
\begin{pf}
We can repeat the proof Theorem 11.1 of \cite{KSS} after making the
following substitutions:
we use Lemma~\ref{Poincare disks} instead of \cite[Lemma 13.4]{KSS},
Proposition~\ref{prop:KSS 11.3} instead of 
\cite[Proposition 11.3]{KSS},
Theorem~\ref{real geometry} instead of \cite[Proposition 8.1]{KSS}
Theorem~\ref{quasi-box mapping} instead of
\cite[Proposition 11.4]{KSS}, and
Lemma~\ref{KSS83} instead of \cite[Lemma 8.3]{KSS}.
\end{pf}

\begin{prop}\cite[Proposition 11.5]{KSS}\label{prop:KSS prop 11.5}
There exists $\theta\in(0,\pi)$ such that for all $n$ sufficiently
large and for each
$x\in\omega(c_0)\cap I_n$,
$$\mathcal{L}_x(\I_n)\subset D_{\theta}(I_n).$$
\end{prop}
This proposition follows immediately from the 
following two lemmas whose proofs can be copied
from \cite{KSS} making the same substitutions that we have already
pointed out:
\begin{lem}\cite[Lemma 11.4]{KSS}
There exist $N$ and $\theta_1\in(0,\pi)$ such that for each
$n$ sufficiently large there exists $m\in\{n,n+3\}$
such that for each $x\in \omega(c_0)\cap I_m$
one has $$\comp_x\mathrm{Dom}(R^N_{\I_m})\subset D_{\theta_1}(I_m).$$
\end{lem}

\begin{lem}\cite[Lemma 11.5]{KSS}
There exists a constant $\mu\in(0,1)$ such that for each $N\geq 2,$
each $\theta\in(0,\pi),$ and each $n$ sufficiently large, if
$$\cup_{x\in\omega(c_0)\cap
  I_n}\comp_x\mathrm{Dom}(R^N_{\I_n})\subset D_{\theta}(I_n),$$
 then
$$\cup_{x\in\omega(c_0)\cap
  I_{n+1}}\comp_x\mathrm{Dom}(R^{N-1}_{\I_{n+1}})\subset D_{\mu\theta}(I_{n+1}).$$
\end{lem}

Using this, we can bound the quasiconformal distortion as we pullback
through the enhanced nest.

\begin{lem}
There exists a constant $\eta_1>0$ such that the following holds.
Let $\hat I_0$ be a nice
neighbourhood of $\crit(f)$ such that $R_{\hat I_0}$ extends to a
qr box mapping $F\colon \mathcal U\rightarrow\mathcal V$.
Let $\I_0=\Comp_{c_0}(\mathcal V)$. Assuming that $\I_{n}$ has been
constructed let $\I_{n+1}=\Comp_{c_0}F^{-p_n}(\I_n),$ where $p_n$ 
is as in the construction of the enhanced nest
 (see page~\pageref{iterate for pullback of enhanced nest}).
Let $s=p_n+p_{n-1}+\dots+p_0$. Then $F^s\colon \I_{n+1}\rightarrow\I_0$ is
$1+\eta_1\mu(I_0)^{1/2}$-quasiregular.
\end{lem}
\begin{pf}
By Lemma~\ref{lem:quasiregular pullbacks} and
Propositions ~\ref{prop:KSS prop 11.5} and 
there exists a constant $\eta>0$ such that
each $F^{p_n}\colon \I_{n+1}\rightarrow\I_{n}$ is $(1+\eta\mu(I_{n})^{1/2})$-quasiregular.
Since, by Theorem~\ref{real geometry}, there exists a universal constant $\rho>0$ such that each $I_n$ is
$\rho$-nice, the lengths of the intervals $I_n$ decays exponentially. 
Thus there exists a constant $\eta_1\geq 0$ such that
$F^s\colon \I_{n+1}\rightarrow\I_0$ is $1+C\eta_1\mu(I_0)^{1/2}$-quasiregular.
\end{pf}

Finally with Proposition~\ref{prop:KSS prop 11.5} in hand we can repeat the
proof of the Upper Bounds,
Proposition~8.3 from \cite{KSS}.
\begin{prop}[Upper Bounds for the enhanced nest] 
\label{prop:upper bounds}
There exists a constant $\eta>0$ such that for all $n$ sufficiently
large the following hold.
\begin{itemize}
\item $\mathrm{diam}(\I_{n})\leq \eta|I_{n}|$, and;
\item there exists a topological disk $\Omega\supset \I_{n}$ such that
$(\Omega\setminus\I_{n})\cap\omega(c_0)=\emptyset$ and
$$\mod(\Omega\setminus\I_{n})>1/\eta.$$
\end{itemize}
\end{prop}
\begin{pf}
By construction there exists a positive integer $\nu$ and
by Proposition~\ref{prop:KSS prop 11.5}, there exists  $\theta'\in(0,\pi)$
such that
$f^{\nu}\colon \mathcal{B}(\I_n)\rightarrow \I_n$ is a proper map with
bounded degree,
$f^{\nu}(\mathcal{A}(\I_n))=\mathcal{L}_{f^\nu(c)}(\I_n),$ and
$\mathcal{L}_{f^{\nu}(c)}(\I_n)\subset D_{\theta'}(I_n)$.

Then, since  by Theorem~\ref{real geometry}, there exists $\rho>0$ 
so that $I_n$ is $\rho$-nice, $\rho$-externally free, the pullback of $I_n$ to 
$\mathcal{B}(I_n)$  is of bounded degree, and we have that either
$c_0$ is of even order and $I_n$ is $\rho$-free or
$c_0$ is of odd order, 
there exists $\lambda\in(0,1)$ such that
$$\mathcal{A}(\I_n)\subset D_{\lambda\theta}(\mathcal{B}(I_n)),$$
and
$$\frac{\diam(\mathcal{A}(\I_n))}{|\mathcal{A}(I_n)|}\leq
\eta(\theta)\max\Bigg\{1,\Big(\frac{\diam(\I_n)}{|I_n|}\Big)^{1/2}\Bigg\}.$$
To prove this estimate, we let $r>0$ be the return time of
$\mathcal{L}_{f^{\nu}(c_0)}(\I_n)$ to $\I_n$,
and decompose $f^{\nu+r}\colon \mathcal A(\I_n)\rightarrow\I_n$ into
$f^{\nu+r-1}\colon f(\mathcal A(\I_n))\rightarrow\I_n$ and $f\colon \mathcal
A(\I_n)\rightarrow f(\mathcal A(\I_n)),$ apply 
Lemma~\ref{lem:quasiregular pullbacks}, and use the fact that
there is a critical point $c_0\in\mathcal{A}(\I_n)\subset\I_n$ of degree $d\geq 2$.
Since $\mathcal{B}(I_n)$ is $\rho'$-nice and $\mathcal{A}(I_n)$ is
$\rho'$-externally free 
for some $\rho'>0$ given by Corollaries~\ref{nice}~and~\ref{cor:free-space},
it follows that there exist
$\xi>0$, $\theta'\in(0,\pi/2)$ and a
topological disk $\mathcal{A}(\I_n)',$ $\mathcal{A}(\I_n)\Subset
\mathcal{A}(\I_n)'\subset D_{\theta'}(\mathcal B(I_n))$,
so that
$\mod(\mathcal{A}(\I_n)'\setminus \mathcal{A}(\I_n))>\xi$ and $(\mathcal{A}(\I_n)'\setminus\mathcal{A}(\I_n))\cap\omega(c_0)=\emptyset$.
Let $\tilde \nu\in\mathbb N$ be so that
$f^{\tilde \nu}\colon \mathcal{BA}(\I_n)\rightarrow\mathcal A(\I_n)$.
Let $F_n=\Comp_{c_0}f^{-\tilde \nu}(\mathcal B(I_n))$.
Then $f^{\tilde \nu}|F_n$ is a proper map of bounded degree,
and either $F_n$ is $\delta'$-free for some
$\delta'>0$ or $c_0$ is an odd critical point.
Likewise, for each $i,$ $i=1,2,\dots, T$, let $t_i\in\mathbb N$, be so
that
$f^{t_i}\colon \Gamma^{i}\mathcal{BA}(\I_n)\rightarrow\Gamma^{i-1}\mathcal{BA}(\I_n)$.
Inductively construct a sequence of domains $F^i_n$, $i=0,\dots, T$,
setting $F^0_n=F_n$ and $F^i_n=\Comp_{c_0}f^{-t_i}(F^{i-1}_n)$.
Each $F_n^i$ is a pullback of $F^{i-1}_n$ of bounded degree and 
either $F^i_n$ is $\delta''$-free for some $\delta''>0$ or $c_0$ is an
odd critical point.
Set $\I'_{n+1}=\comp_{c_0}(f^{-(\tilde \nu+t_1+\dots t_T
}(\mathcal{A}'(\I_n))).$
Observe that $(\mathbf I'_{n+1}\setminus \mathbf
I_{n+1})\cap\omega(c_0)=\emptyset.$

By Lemma~\ref{lem:quasiregular pullbacks}, we have
that $f^{\tilde \nu+t_1+\dots+t_T}|\I'_{n+1}$ is a $1+c\mu(I_{n-1})^{1/2}$-quasiregular
mapping with bounded degree.
It follows that there exists $c'\geq 1$ so that
$\mod(\mathbf I'_{n+1}\setminus \mathbf I_{n+1})>\xi/c'.$
As before,
$$\frac{\diam(\I_{n+1})}{|I_{n+1}|}\leq \eta(\theta)\max\Bigg\{1,\Big(\frac{\diam(\I_n)}{|I_n|}\Big)^{1/2}\Bigg\}.$$
Since $\diam(\mathbf I_0)<\eta|I_0|,$
it follows that
$\diam(\I_n)/|I_n|$ is bounded from above.

\end{pf}

It is worth noticing that from
Proposition~\ref{prop:upper bounds2}, Proposition \ref{prop:KSS prop 11.5} and the definitions of the
operators $\mathcal A$ and $\mathcal B$ we have:

\begin{cor}
There exists $\varepsilon>0$ and a universal constant
$\delta>0$ such that, for all $n>0$ for which $|I_n|<\varepsilon$ the
following holds.
The puzzle piece $\I_n$ is $\delta$-nice and $\delta$-free. Hence the return mapping to
$\I_n$ is $\delta$-extendible.
Even more, there exists an interval $\tilde I_n\supset (1+2\delta)I_n$ 
$\tilde I_n\supset (1+2\delta)I_n$ with $(\tilde I_n\setminus
I_n)\cap\omega(c_0)$ so that $D_{\theta'}(\tilde I_n)\supset \I_n$.
\end{cor}

\subsection{Lower Bounds for finitely renormalizable maps}\label{sub-sec:lower bounds}

Lower Bounds for the geometry of puzzle pieces in the enhanced nest follows from the next two lemmas.
The first lemma is Lemma 10.1 of \cite{KSS}; however, we cannot repeat the proof given there, since $f$ does not have to be analytic.
\medskip

Let 
$$\eta_n=\inf_{x\in\omega(c_0)\cap I_n}\frac{d(x,\partial \I_n)}{|I_n|}.$$

\begin{lem}\label{lem:KSS 10.1}\cite[Lemma 10.1]{KSS}
There exists $\varepsilon>0$ so that for all $n>0$ for which $|I_n|<\varepsilon$ the following holds:
\begin{itemize}
\item[1.]There exists a constant $\delta>0$ such that $\eta_{n+1}>\delta\eta_n.$
\item[2.]There exist $\kappa>0$, $\nu>1$ such that if $|I_{n}|/|I_{n+1}|>\nu$, then
$$\eta_{n+2}\geq\min(\kappa,2\eta_{n+1}).$$
\end{itemize}
\end{lem}
\begin{pf}
Note that we can decompose the map $f^{p_n}\colon  I_{n+1}\rightarrow I_n$
into a bounded number of maps, each of the form $z^d$ followed by a
diffeomorphism with bounded distortion.

To prove the first estimate, observe that $\eta_{n+1}\geq
\delta\eta_n$ is equivalent to
$$\frac{\inf_{x\in\omega(c_0)\cap I_n}d(x,\partial \I_n)}{\inf_{x\in\omega(c_0)\cap I_{n+1}}d(x,\partial \I_{n+1})}\leq \frac{1}{\delta}\frac{|I_n|}{|I_{n+1}|}.$$
For any $x\in\omega(c_0)\cap I_n$, since $f^{p_n}$ has bounded degree and,
by Theorem~\ref{real geometry},
$I_n$ has external free space,
it follows from Theorem~\ref{real Koebe} that
there exists a constant $c'>1$ such that 
$$|(f^{p_n})'(x)|\leq c'\frac{|I_n|}{|I_{n+1}|}.$$

By Lemma~\ref{lem:quasiregular pullbacks},
there exists $\eta>0$, such that $f^{p_n}|\I_{n+1}$
is $(1+\eta|\mu(I_{n})|^{1/2})$-quasiregular.
So by the Stoilow Factorization Theorem we can express
$f^{p_n}|\I_{n+1}=h\circ g,$ where $g\colon \I_{n+1}\rightarrow \I_{n+1}$ 
is $(1+\eta|\mu(I_{n})|^{1/2})$-quasiconformal and
$h\colon \I_{n+1}\rightarrow\I_{n}$ is holomorphic. 
Let $\lambda_n$ be an affine mapping that scales $\I_n$ to unit size.
Then the mapping
$H(z)=\lambda_nf^{p_n}(\frac{z}{\lambda_{n+1}})$
is a $(1+\eta|\mu(I_{n})|^{1/2})$-quasiregular mapping between domains
of unit size. As $n\rightarrow\infty$ this mapping converges to a
holomorphic mapping.
Arguing by contradiction, it follows that 
there exist constants $c>1$ and $\kappa_1>0$ such that
$$|(f^{p_n})'(z)|\leq c\frac{|I_n|}{|I_{n+1}|}$$
for $z\in\mathbb{C}$ with $\dist(z,\omega(c_0)\cap
I_{n+1})\leq\kappa_1|I_{n+1}|$, and the first estimate follows
immediately.

Similarly, from the Koebe Distortion Theorem and the fact that $f^{p_n}|\I_{n+1}$ is
$(1+\eta|\mu(I_{n})|^{1/2})$-quasiregular, we have that
there exists $\kappa_2>0$ such that for any $z\in B(c_0,2\kappa_2|I_{n+1}|),$
$$|(f^{p_n})'(z)|\leq \frac{1}{2}\frac{|I_n|}{|I_{n+1}|}.$$
Let $\nu>1$ be large. 
Notice that when $|I_{n}|/|I_{n+1}|>\nu$,
then by Lemma~\ref{lem:small return domains} 
there exists $\nu_1>0$,
so that for any $x\in\omega(c_0)\cap I_n$ we have
$|I_{n+1}|/|\mathcal L_x(I_{n+1})|>\nu_1.$
From the definitions of the operators $\mathcal {A,B}$ and
Theorem~\ref{real Koebe}~(2),
there exists a constant $\eta=\eta(\varepsilon)\rightarrow 0$ as
$\varepsilon\rightarrow 0$, so that
$\omega(c_0)\cap\mathcal{BA}(I_{n+1})\subset\mathcal{AA}(I_{n+1})$
is contained in a $\eta|\mathcal{BA}(I_{n+1})|$-neighbourhood of the
critical point, and the second estimate in the
lemma follows.
\end{pf}

The following is Lemma 10.2 of \cite{KSS}. 
The proof differs from that given in \cite{KSS} 
because $f$ is not analytic and we have to treat 
the big geometry case separately from the bounded geometry case.
\begin{lem}\label{lem:KSS 10.2}\cite[Lemma 10.2]{KSS}
There exists  $k_0\in \N$ and a constants $\gamma>0$  and $\varepsilon>0$ such
that for all $n>0$ with $|I_n|<\varepsilon$ and
for all $x\in\omega(c_0)\cap I_n$,
$$B(x,\gamma|\comp_x\mathrm{Dom}(R^{k_0}_{I_n})|)\subset\I_n.$$
\end{lem}
\begin{pf}
Let $N=p_{n-1}+p_{n-2}+\dots+p_0$. 
Observe that $f^N(\I_n)=\I_0$.
By Lemma 8.3 of \cite{KSS}, we have that $2p_{n-1}\geq N$ and $r(I_n)\geq N/6$.
Let $x\in\omega(c_0)\cap I_n$ and let $W=\hat{\mathcal{L}}_{f^N(x)}(\I_n)$. 
Let $U=\comp_x(f^{-N}(W))\cap\mathbb{R}.$
Then $U=\comp_x(\mathrm{Dom}(R^{k_0}_{I_n}))$ for some $k_0$, and
since $r(I_n)\geq N/6$ we have that  $k_0\leq 6$ and that the order of $f^N\colon U\rightarrow W$ is bounded by $6b$.

If either $c_0$ is even or $|I_{n-3}|/|I_{n-2}|$ is bounded from
above, then by Theorem~\ref{real geometry} $I_n$ is $\rho$-free
for some $\rho>0$, and we can repeat
the proof of Lemma 10.2 of \cite{KSS}.
So, suppose that $|I_{n-3}|/|I_{n-2}|$ is big and
$c_0$ is of odd order. In this case, $I_{n}$ is
$C$-nice and $C$-externally free with $C$ big.
Let $0<C'<C$ be so that
$\tilde I:=(1+2C')I_{n}$ is $\delta$-externally
free.
Let $\tilde W=\hat{\mathcal{L}}_{f^N(x)}(\tilde I)$,
and $\tilde U=\comp_x(f^{-N}((1+2\delta)W))\cap\mathbb{R}.$
Since $C$ is big, we can choose $C'$ so that $\delta>0$ does not
depend on $f$ and so that $\tilde W\supset (1+2\delta)W$.
Then $\tilde W$ is $\rho''$-free for some $\rho''>0$.
By Proposition~\ref{prop:smoothtoanalytic}, there exists $\sigma\in (0,\pi/2)$
so that $D_{\pi-\sigma}(\tilde W)\subset D_{\pi-\sigma}(I_0)\subset\I_0$.
So Lemmas \ref{lem:quasiregular pullbacks}, \ref{lem42} and \ref{lem:z^l lower bound}, 
imply that $\I_n\supset\comp_{x}f^{-N}(D_{\pi-\sigma}(\tilde W))\supset D_{\theta''}(\tilde U)$,
where $\theta''\in(0,\pi/2)$.
\end{pf}

\begin{prop}

There exist beau constants $\eta>0$ and $\varepsilon>0$ such that for all $n>0$ with $|I_n|<\varepsilon$

$$B(c_0,\eta|I_n|)\subset\I_n.$$
\end{prop}
\begin{pf}
The proof follows the proof of Proposition 10.1 of \cite{KSS} with
some minor adjustments to deal with the case when $c_0$ is odd. 
Indeed, there exists a constant $\varepsilon'>0$ such that
if $0<\eta_n<\varepsilon'$,
then $|I_{n+1}|/|I_{n}|<\varepsilon,$ where $\varepsilon$ is the
constant from the second estimate in Lemma~\ref{lem:KSS 10.1}.
To see this, observe that if $\eta_{n}$ is very small,
then by Lemma~\ref{lem:KSS 10.2}
there exists $x\in\omega(c_0)\cap I_n$ such that
$|U|/|I_n|$ is very small where 
$U\owns x$ is a component of the domain of $R^{k_0}_{I_n}.$
Since $k_0$ is bounded ($\leq 6$, see the proof of Lemma~\ref{lem:KSS 10.2}), this implies that there exists a return domain $J$ to $I_n$,
$J\cap\omega(c_0)\neq\emptyset,$ and such that $|J|/|I_n|$ is very
small. It follows from Theorem~\ref{real geometry} (e)
that $|I_{n+1}|$ is very small compared to $|I_n|$. 
Thus it follows from
the second statement of Lemma~\ref{lem:KSS 10.1} that when $\eta_n<\varepsilon'$,
$\eta_{n+2}\geq\min\{\kappa,2\eta_{n+1}\}$. By the first statement of
Lemma~\ref{lem:KSS 10.1}, we have
$\eta_{n+2}\geq\min(\kappa',2\eta_{n+1})$ for all $n$, where $\kappa'=\min(\delta^2\varepsilon',\kappa)$.
Since $\eta_0,\eta_1$ are bounded away from 0 the result follows. 
\end{pf}

Thus we have concluded the proof of complex bounds in the non-renormalizable case, i.e.\ 
Theorem~\ref{complex bounds}  (in the non-renormalizable case) and Theorem~\ref{thm:complex bounds non-ren}. 
Let us now turn to the infinitely renormalizable case.

\subsection{Infinitely renormalizable maps}

To  construct a qc polynomial-like extension from a qc quasi-box mapping for infinitely renormalizable maps,
we use a more geometric approach.
Let $L\colon \mathbb{R}_+\to \mathbb{R}_+$.  A map $h\colon X\to Y$ between two metric spaces is called an {\em L-quasi-isometry} if for any $\varepsilon > 0$,
$\dist(h(x), h(y))\leq \max\{L(\varepsilon) \dist(x, y),\varepsilon\},$ for all  $x, y \in X$.
Quasiconformal maps are quasi-isometries with respect to the hyperbolic metric:
\begin{lem}[Lemma 2.2, \cite{ALdM}]
For every $\kappa \geq 1$ there exists $L_\kappa \colon  \mathbb{R}_+ \rightarrow \mathbb{R}_+$  such that if $h \colon  S \rightarrow \tilde S$ is a $\kappa$-qc map between two hyperbolic Riemann surfaces, then $h$ is a $L_\kappa$ quasi-isometry in the hyperbolic metric. Furthermore, for every $\varepsilon>0$, $\lim_{\kappa\rightarrow 1} L_\kappa(\varepsilon)=1$. 
\end{lem}

\begin{lem}[cf. \cite{LY} Lemma 2.4, \cite{McMullen-3mflds} Proposition 4.10]\label{lem:LY}
For every $\varepsilon,\eta>0$, there exists $\kappa_0>1$ such that if $1\leq \kappa\leq \kappa_0$, there exists $\delta>0$ such that the following holds. 
Let $U\subset V$ be two real-symmetric topological disks and $f\colon U\rightarrow V$ be a
real-symmetric qc branched covering with non-escaping critical points and compact
Julia set. Suppose that $U$ contains an $\eta$-neighbourhood,
$\mathcal N$ of the
$J(f)$ in the hyperbolic metric on $U$, and that $f|\mathcal N$ is
$\kappa$-quasiregular. Then there are real-symmetric topological disks $U\supset V'\supset U'\supset K(f)$ such that the restriction $f\colon U'\rightarrow V'$ is polynomial-like. Moreover, if $\mathrm{mod}(U\setminus K(f))\geq \varepsilon>0$, then $\mathrm{mod}(V'\setminus U')\geq\delta(\varepsilon,\ell,\eta)>0$, where $\ell$ is the degree of $f$.
\end{lem}
\begin{pf}

The proof of this lemma goes as the proof of Proposition 4.10 of
\cite{McMullen-3mflds}, we will provide a sketch of the argument.
Notice that since $f$ is a real-symmetric mapping between
real-symmetric domains $U$ and $V$, $K(f)$  and the
$\eta$-neighbourhood of $K(f)$ in the hyperbolic metric on $U$ are
both real-symmetric.
Observe that we can assume that $\mod(V\setminus K(f))$ is bounded
from above. 
This implies that $\mathcal N\setminus K(f)$ is an annulus with modulus bounded
from below, so that its core curve $\gamma$ is a $\kappa'$-quasicircle.
Let $\mathcal N_{\gamma}$ denote the region bounded by $\gamma$.
Then, by the Schwarz Reflection Principle,
$f|\mathcal N_{\gamma}$ extends to a $\kappa''$-quasiregular mapping
$\tilde f\colon U\rightarrow V$ that agrees with $f$ on $\mathcal N_{\gamma}$.
Consider the hyperbolic Riemann surfaces $\mathcal N_0\subset U_0\subset V_0$
obtained by doubling $\mathcal N_\gamma\setminus K(f), U\setminus K(f)$ and $V\setminus K(f)$,
respectively, across their ends corresponding to $K(f)$.
Then $\mathcal N_0\subset U_0\subset V_0$ are annuli with the same
core geodesic $\gamma$,
and $\tilde f|U\setminus K(f)$ extends to a symmetric mapping
$F\colon U_0\rightarrow V_0$ that sends $\gamma$ to itself.

Let us show that $\|DF(z)\|>1$ (in the hyperbolic metric) for $z$
sufficiently close to $\gamma$. 
Take $\tau>0$ and $z\in \gamma$, consider the balls
$B_{U_0}(z,\tau)$, $B_{V_0}(z,\tau)$  in terms of the hyperbolic
metrics on $U_0$ and $V_0$, respectively. 
Provided
we take $\tau>0$ sufficiently small, we get that the distance of the
ball
$B_{U_0}(z,\tau)$  to $V_0\setminus U_0$ is strictly positive.
Then $F(B_{U_0}(z,\tau))\supset B_{U_0}(F(z),\tau)$ provided $\kappa$ and $\tau$ are
sufficiently small. Indeed, $F|\mathcal N_0$ is a composition of a $\kappa$-qc map
 and a
conformal isometry from the hyperbolic metric on $U_0$ to the
hyperbolic metric on $V_0$. Therefore $F(B_{U_0}(z,\tau))\supset
B_{V_0}(F(z),r\tau)$ where $r<1$ can be chosen arbitrarily close to
one provided $\kappa$ is sufficiently close to one. Moreover, the
contraction of the inclusion $U_0\rightarrow V_0$  on
$B_{U_0}(z,\tau)$ is bounded from above by a constant $t<1$ since the
distance of the ball $B_{U_0}(z,\tau)$  to $V_0\setminus U_0$ is
strictly positive. It follows
that $F(B_{U_0}(z,\tau))\supset B_{U_0}(F(z),\kappa\tau)$ for some $\kappa >1$. 

Let $V'_0\subset V_0$ be the $\tau$-collar neighbourhood of $\gamma$,
and let $U'_0=F^{-1}(V'_0)$. Set $U'=(U'_0\cap U)\cup K(f)$ and
$V'=f(U')$.
Then $f\colon U'\rightarrow V'$ is a polynomial-like map and
$\mod(V'\setminus U')$ is bounded away from zero. 

\end{pf}

\begin{thm}[Complex bounds in the infinitely renormalizable case]\label{thm:complex bounds infinitely ren}
There exist $\delta>0$ and $C>0$ such that the following
holds. Suppose that $f\colon M\rightarrow M$ is infinitely renormalizable at $c_0$,
with periodic intervals $M\supset J_1\supset J_2\supset\dots $ with
periods $1<j_1<j_2<\dots$.
Suppose that $i$ is sufficiently large. Then there exists a $\kappa(V_i)$-qc
polynomial-like mapping $F_i\colon U_i\rightarrow V_i$ which extends
$f^{j_i}\colon J_i\rightarrow J_i$ such that $U_i, V_i$ are
real-symmetric,  
 $\mod(V_i\setminus U_i)>\delta$, $U_i$ has $\delta$-bounded geometry and $\mathrm{diam}(U_i)<C|J_i|$.
\end{thm}
\begin{pf}
Let $I_n$ be the smallest puzzle piece in
the generalized enhanced nest such that $I_n\supset J_i$. Then $I_n$ is a
terminating interval.
Note that $J_i=I_n$ 
if and only if $I_{n-1}$ is terminating and $J_i=\mathcal R(I_{n-1})$.

\medskip

\noindent\textit{Step 1: Moduli bounds.}
\textit{We will show that if $i$ is sufficiently big,
the return mapping to $J_i$ 
extends to a quasi-box mapping $F_Q\colon U_Q\rightarrow V_Q$,
and that $\mod(V_Q\setminus K(F_Q))$ is universally bounded from below.}

We will
show that $\mod(V_Q\setminus K(F_Q))$ is bounded from below
by proving that the filled Julia set of the
quasi-box mapping, $K(F_Q)$, is contained in a Poincar\'e disk, $D_{\theta}(K),$ where $K$ is an interval well-inside $U_Q\cap\mathbb R,$ with $|K|$
comparable to $|I_n^\infty|.$

We will divide the proof in two cases, depending on the geometry of
the generalized enhanced nest. 
To do so, we define a constant $\nu>0$ that will
separate  the bounded geometry case from the big geometry one.

\medskip
\noindent\textit{Claim: there exists $\nu>0$ so that the following holds:} 
\begin{itemize}
\item[(a)] If $|I_{n}|/|I^\infty_{n}|>\nu,$ then the first return mapping to $I^\infty_{n}$ extends to a polynomial-like map $F:U\to V$ with $V=D_{\pi/2}(I^{j}_{n})$ with $ |I^j_{n}|\asymp  |I_{n}^\infty|,$ for some $j\in \N.$ 
 Furthermore, there exists an interval $K$ well-inside $I^{j+1}_n$, with $|K|\asymp |I_{n}^\infty|,$ so that $K(F)\subset D_{\pi/2}(K).$
\item[(b)] If $I_{n-1}$ is non-terminating and $|I_{n-1}|/|I_{n}|>\nu$ then the first return mapping to $I_{n}$ extends to a qr box mapping $F:U\to V$ with $V=D_{\pi/2}(I^{j}_{n})$ for some $j\in \N.$  Furthermore, there exists an interval $K$ well-inside $I^{j+1}_n$, with $|K|\asymp |I_{n}^\infty|$ so that $K(F)\subset D_{\pi/2}(K).$
\item[(c)] If $I_{n-1}$ is terminating and
  $|\mathcal{R}(I_{n-1})|/|I_n|>\nu,$ then the first return mapping to
  $I_n^\infty$ extends to a qr box mapping with range
  $D_{\pi/2}(I^{j}_{n})$ for some $j\in \N.$  Furthermore, there exists an interval $K$ well-inside $I^{j+1}_n$, with $|K|\asymp |I_{n}^\infty|$ so that $K(F)\subset D_{\pi/2}(K).$

\item[(d)] Assume $I_{n-1}$ non-terminating and $|I_{n-1}|/|I_n|>\nu.$ Then the return map to $I_{n-1}^\infty$ extends to a polynomial like map $F:U\to V$ with $V\cap \mathbb R= I_{n-1}^j$, for some $j\in\N$, with the property that $ |I^j_{n-1}|\asymp  |I_{n-1}^\infty|.$

 \end{itemize}

\noindent \textit{Proof of Claim:}
We will prove (a), the other cases follow using similar arguments. By the proof of Lemma \ref{big space pl}, we know that $I_n^2$ is $C$-nice and $C$-free for some $C>0,$ with $C\to \infty$ as $\nu\to \infty$. The Koebe Principle, Theorem \ref{real Koebe} and Lemma \ref{bounds1} imply that there exists  $j\in \N$ so that $|I^j_{n}|/|I_n^{j+1}|>\nu_1,$ where $\nu_1>0$ is defined in Lemma \ref{big space pl}, with $|I_n^j|$ is comparable to $|I_n^{\infty}|.$ The claim follows using a similar argument as the one used in Lemma \ref{big space pl}, except for the fact that we pull back two times under the return map to $I^\infty_n,$ instead of one.

\medskip

\medskip
Cases (a), (b) and (c) reduce the proof to the following: either,
\begin{itemize}
\item $I_{n-1}$ is terminating and $|I_{n}^\infty|$ is comparable to
$|\mathcal{R}(I_{n-1})|$, or
\item
$I_{n-1}$ is non-terminating and $|I_{n-1}|$ is comparable to $|I_n^\infty|$.
\end{itemize}

\medskip
Let $M\in\mathbb N$ be the constant associated to $\nu$ 
given by Theorem~\ref{quasi-box mapping}.

\medskip

\noindent\textit{Case A: assume $J_i\subsetneq I_n$.}

\medskip

\noindent\textit{Case A.1: Bounded geometry.} 
Suppose that  for all
$0\leq j\leq M,$ 
$|I_{n-j}|/|I_{n-j+1}|\leq \nu.$
Then $|I_{n+1}|\asymp |I_{n-M}|$.
Let $F_Q\colon U_Q \rightarrow V_Q$ be the qc quasi-box mapping associated
to the return mapping to $I_n$
given by Theorem~\ref{quasi-box mapping},
where $V_Q=D_{\pi/2}(I_{n-M})\cap\mathbb C_{I_{n}}$ and 
$U_Q\subset D_{\pi/2}(I_{n-M+1})$. 

Now, we will show the moduli bounds between $V_Q$ and $K(F_Q)$
  hold.
 We will choose a constant $k_0\in\mathbb N$ in the course of the proof.
Let $\mathbb G=\{G_l\}_{l=0}^{k_0 j_i}$ with $G_{k_0 j_i}=J_i$ and $G_0=J_i$.
To begin, we take $k_0$ at least $4(4)+1.$
Let us recall that in Proposition~\ref{prop:s-seq},
given a point $x\in I_n\cap \omega(c_0)$ with return time to $I_n$
equal to $s$ 
we constructed
a sequence of times $0\leq s_{n-4}<s_{n-5}<\dots<s_{n-M}=s$ 
for which the pullbacks of certain Poincar\'e disks are well controlled.
In cases (a), (b) and (c) of the proof of
Proposition~\ref{prop:s-seq}, assuming that $s_{n-i-1}$ has been
defined for $0\leq i\leq M-1$, we can define $s_{n-i}$ provided that there are at
least four intervals $G_l\subset I_{n-i}$ with $0\leq l<s_{n-i-1}$. In
case (d),  
to be able to define  $s_{n-i}$ we need at least three of the intervals
 $G_l$ with $0\leq l<s_{n-i-1}$ to be contained in $I_{n-i+1}$. 
Now, we apply the argument form Proposition~\ref{prop:s-seq} 
to the chain $\mathbb G$. Observe that, since $k_0-1\geq 16,$ we can repeat
Step 2 of Proposition~\ref{prop:s-seq}, to carry on the construction
of the sequence $s_i$,
to obtain $s_{n-1}<s_{n-2}<s_{n-3}<s_{n-4}<\dots<s_{n-M}=k_0j_i$. 
Then, using Lemma~\ref{lem62}, as in the proof of
Theorem~\ref{quasi-box mapping}, we find $\theta_1\in(0,\pi),$ so that
$$\comp_{c_0}f^{-(k_0j_i-s_{n-1})}(V_Q)\subset D_{\theta_1}(I_{n-1}).$$

Suppose that $I_{n-1}$ is non-terminating. 
Using Proposition~\ref{folding lemma} and the choice of $k_0$
we define $0\leq
s_n<s_{n-1}$, and obtain an interval $K_2$ and an angle $\theta_2\in(0,\pi)$
so that 
$$\comp_{c_0}f^{-(k_0j_i-s_{n})}(V_Q)\subset D_{\theta_2}(K_2),$$
where $K_2$ is well-inside 
the largest terminating interval $I^{i_0}_{n-1}$ in the principal nest
$I^0_{n-1}\supset I^1_{n-1}\supset I^2_{n-1}\supset\dots$. Under these circumstances, Lemma \ref{KSS82} and the definition of successor imply that $I_{n+1}$ is a pullback of $I_n$ of bounded order; $\mathcal{BA}(I_{n})$ is a pullback of 
$I_{n}$ with order bounded from above by $2b^2$,
and each successor $\Gamma^{i+1}(\mathcal{BA}(I_{n}))$ is a pullback of
order at most $2b-1$ of $\Gamma^{i}(\mathcal{BA}(I_{n})).$
Thus
$I_{n+1}$ is a pullback of $I_{n}$ with order at most
$2b^2+5b(2b-1)$. For each $i\geq i_0,$ the pullback $I^{i+1}_{n-1}$ of
$I^i_{n-1}$ has order $b,$ since $I^{i_0}_{n-1}$ is
terminating. Therefore, increasing $k_0$ by $2b+5(2b-1)$ we can pull
back $\Comp_{c_0}f^{-(k_0j_i-s_{n-1})}(V_Q)$ by at most $15b-5$
iterates of $f^{j_i}$ and obtain an interval $K$  that is well-inside
$I_n$ and $\theta\in(0,\pi)$ such that $K(F_Q)\subset D_{\theta}(K).$ 
Furthermore, since $J_i\subset K$ and $|I_n|\asymp |I_{n+1}|$ we get that $|K|\asymp |I_{n}^\infty|.$

Now suppose that $I_{n-1}$ is terminating. Then $I_n=\mathcal{L}_{c_0}(\mathcal
R(I_{n-1}))$. Since $J_i\subsetneq I_n$, we have that
$\mathcal{L}_{c_0}(\mathcal R(I_n))\neq \mathcal R(I_n)$.
Recall that $\comp_{c_0}f^{-(k_0j_i-s_{n-1})}(V_Q)\subset D_{\theta_1}(I_{n-1}).$
By Proposition~\ref{terminating angle control},
using three returns to
$I_{n},$ see Combinatorial Remark 6,
there exist $r_{n-1}<s_{n-1}$
and $\theta_2\in(0,\pi)$ such that
$\comp_{c_0}f^{-(k_0j_i-r_{n-1})}(V_Q)\subset D_{\theta_2}(\mathcal{R}(I_{n-1}))$.


Now, by Proposition~\ref{prop:from renormalization to landing},
Combinatorial Remark 7,
and Lemma~\ref{chain pullback}, increasing $k_0$, if necessary, and
pulling back by at most two more iterates of 
$f^{j_i}$, we have that there exist $s_{n}<r_{n-1}$, an interval $K$ well-inside $I_{n}$ 
and $\theta\in(0,\pi)$ such that 
$\comp_{c_0}f^{-(k_0j_i-s_{n})}(V_Q)\subset D_{\theta}(K).$ Thus
$K(F_Q)\subset D_{\theta}(K)\cap U_Q.$
Furthermore, since $J_i\subset K$
and $|I_n|\asymp |I_{n+1}|$ we get that $|K|\asymp |I_{n}^\infty|.$

\medskip

\noindent\textit{Case A.2: Unbounded geometry.}
Suppose that there exists $j, 0\leq j\leq M$ such that 
$|I_{n-j}|/|I_{n-j+1}|>\nu.$
Let $j_0, 0\leq j_0\leq M,$ be minimal so that $|I_{n-j_0}|/|I_{n-j_0+1}|>\nu.$

\medskip
\noindent \textit{Case A.2.1: Assume $j_0\geq 2$.}

If $I_{n-j_0}$ is non-terminating, then it is easy to modify the argument in
the bounded geometry case:
If $I_{n-j_0}$ and $I_{n-j_0+1}$ are both non-terminating,
then by  Case (a) in Theorem \ref{quasi-box mapping}, the first return mapping to
$I_{n-j_0+2}$
extends to a quasiregular box mapping with range $D_{\pi/2}(I_{n-j_0+2})$.
Since 
for all $0\leq i\leq j_0-1$, we have that $|I_{n-i}|/|I_{n-i+1}|<\nu$,
we argue as in
Theorem~\ref{quasi-box mapping} to see that there exists a 
qc quasi-box mapping $F\colon U_Q\rightarrow V_Q$ 
that extends the return mapping to $I_n$, and we can repeat the
argument from the bounded geometry case to find $\theta\in(0,\pi)$
and an interval $K,$ well-inside $I_n,$ so that $K(F_Q)\subset
D_{\theta}(K)$. Furthermore, since $J_i\subset K$ and $|I_n|\asymp |I_{n+1}|$ we get that $|K|\asymp |I_{n}^\infty|.$
Similarly, if $I_{n-j_0}$ is non-terminating, and $I_{n-j_0+1}$ is terminating,
then by Case (b) in Theorem \ref{quasi-box mapping}, the return mapping from $I_{n-j_0+1}^1$ to $I_{n-j_0+1}$ extends
to a qc polynomial-like mapping with range $D_{\pi/2}(I^{2}_{n-j_0+1})$, and we can repeat the proof of the bounded geometry case to conclude the proof in this case. So if
Again, since for all $0\leq i\leq j_0-1$, we have that $|I_{n-i}|/|I_{n-i+1}|<\nu$,
we can repeat the arguments in the bounded geometry case.

Now, assume that $I_{n-j_0}$ is terminating.
We need to consider in which case of Proposition~\ref{big space term}
we are in.
Let us first treat the cases that immediately reduce to the
previous one: 
\begin{itemize}
\item If we are in Case (2) of 
Proposition~\ref{big space term}, then the first return mapping to
$I_{n-j_0+1}^{\infty}$
extends to a qr box mapping with range $D_{\pi/2}(I_{n-j_0+1}^2)$.
\item If we are in Case (3) of Proposition~\ref{big space term},
then the first return mapping to $I_{n-j_0+2}$ extends to a qr box
mapping with range $D_{\pi/2}(I_{n-j_0+2})$.
\item If we are in Case (4) of Proposition~\ref{big space term},
then the first return mapping to $I_{n-j_0+2}^{\infty}$ extends to a
qr box mapping with range $D_{\pi/2}(I_{n-j_0+2}^2)$.
\end{itemize}
So suppose that we are in Case (1) of  
Proposition~\ref{big space term},
which means that $|I_{n-j_0}^\infty|$ is small compared to
$|I_{n-j_0}|$, 
so the first return mapping to $I_{n-j_0}^\infty$ extends to
a qr box mapping with range $D_{\pi/2}(I_{n-j_0}^2)$.
Notice that we must have that 
$|I_{n-j_0+1}| \asymp|\mathcal R(I_{n-j_0})|$ for otherwise,
we would have that $I_{n-j_0+1}$ is $C$-nice for some $C>0,$ large, 
which contradicts the choice of $j_0$.

So we have that 
$|I_{n-j_0+1}| \asymp|\mathcal R(I_{n-j_0})|,$ and
  $|I_{n-j_0}^\infty|$ is much smaller than $|I_{n-j_0}^2|.$
Observe that if
if $|I_{n-j_0+1}| \asymp|\mathcal R(I_{n-j_0})|$ and
$|I_{n-j_0}|\asymp |I_{n-j_0}^2|$, then by Lemma~\ref{bounds1},
$|I_{n-j_0}^\infty|\asymp |I_{n-j_0}|$.
Now we can repeat the
argument of the bounded geometry case.
This concludes the argument when $j_0\geq 2$.
\medskip

\noindent\textit{Case A.2.2: Assume $j_0=1$.} 
Then $|I_{n}|$ is small compared to $|I_{n-1}|.$ Observe that by the
assumptions made after the proof of the claim we must have that
$I_{n-1}$
is terminating and $|\mathcal{R}(I_{n-1})|$ is comparable to
$|I_n^\infty|.$ 
These imply that
$|\mathcal{R}(I_{n-1})|$ is much smaller than $|I_{n-1}|$. 
This fact, along with Lemma~\ref{bounds1},
 implies 
that $|I_{n-1}|$ is much bigger than $|I_{n-1}^2|.$ Then, we have that there exists a qc polynomial-like
mapping 
$F\colon U\rightarrow V$ that 
extends the return mapping to $I_{n-1}^\infty$ with $\diam(V)$
comparable to $|I_{n+1}|$.
Now we can repeat the argument in the bounded geometry case to
pullback one step to complete the proof.

\medskip
\noindent\textit{Case A.2.3: Assume $j_0=0$.}

There are two cases to consider depending on whether $I_{n-1}$ is
terminating, either:
\begin{itemize}
\item  $I_{n-1}$ is non-terminating and $|I_{n-1}|\asymp |I_n|$
and $|I_n|/|I_{n+1}|>\nu$ or
\item $I_{n-1}$ is terminating, $|I_{n}^\infty|\asymp
  |\mathcal{R}(I_{n-1})|$, and $|I_{n}|/|I_{n+1}|>\nu.$ 

\end{itemize}

\medskip

Assume $I_{n-1}$ is non-terminating.
Since $|I_n|$ is comparable to $|I_{n-1}|,$
we can use the $j_0\geq 1$ cases to obtain a qc
quasi-box mapping 
$F'_Q\colon U'_Q\rightarrow V'_Q$ with
$V'_Q=D_{\pi/2}(I_{n-j'_0})\cap\mathbb C_{I_{n-1}}$
for some $j_0'\geq 1,$
that extends the return map to $I_{n-1}.$
We obtain a quasi-box mapping $F_Q\colon U_Q\rightarrow V_Q$ extending the return map to $I_n,$
with $V_Q=D_{\pi/2}(I_{n-j'_0})\cap\mathbb C_{I_{n}}$ 
and $U_Q=\comp_{c_0}f^{-j_i}(V_Q)$ since $F_Q$
is an iterate of $F_Q'$.
Using the
argument in the cases for $j_0\geq 1$, we have that there exist
an interval $K_1$ well-inside
$I_{n-1}$ 
and
$\theta_1\in(0,\pi)$ such that 
$$K(F_Q)\subset D_{\theta_1}(K_1).$$

As in the bounded geometry case, we are going to use
the argument for
Proposition~\ref{folding lemma}, to find $s'_n>0$,
 $\theta_2\in(0,\pi)$ and $K_2\subset I_{n-1}^{i_0}$ such that $\comp_{c_0} f^{-(k_0j_i-s'_n)}(V_Q)\subset D_{\theta_2}(K_2),$
where $i_0\in\mathbb N$ is minimal with $I_{n-1}^{i_0}$ terminating.
However, we cannot apply Proposition~\ref{folding lemma} 
directly since we do not have that $I_{n-2}$ is comparable to
$I_{n+1}.$
First, notice that since $c_0$ is even and $I_{n-1}$ is
non-terminating, Theorem~\ref{real geometry} implies that there
exists $\rho>0$ such that $I_{n-1}$ is $\rho$-free.
Second, observe that the chain $\mathbb{G}$ is bound to enter a 
terminating component of the landing domain to $I_{n-1}$ containing
critical point of $f$, so $k(I_{n-1},\mathbb G)=\infty.$ 
To use the same argument as in Case
(ii) of Proposition~\ref{folding lemma}, which deals with the
situation when a chain has infinite combinatorial depth. Recall we are
assuming $|I_{n-1}|$ is comparable to $|I_n^\infty|$. This means there
exists
$m\in\mathbb N,$ depending only
on the ratio $|I_{n-1}|/|I_n^\infty|$,
so that $\mathcal C^m(\hat{\mathcal L}_c(I_{n-1}))$ is a
periodic interval. So the proof of
Proposition~\ref{folding lemma}  Case (ii) goes through and we obtain  $\comp_{c_0} f^{-(k_0j_i-s'_n)}(V_Q)\subset D_{\theta_2}(K_2),$ with $s_n',$ $\theta_2$ and $K_2$ as above.
Finally, we pullback a bounded number of times by $f^{j_i},$
to find $\theta\in(0,\pi)$ and $K$ well-inside $I_n$ so that
$$K(F_Q)\subset D_{\theta}(K).$$
Recall, we are working under the assumption that $|I_{n-1}|$ 
is comparable to $|I_n^\infty|.$ Since $I_n^\infty\subset K,$
we get that $|K|\asymp |I_n^\infty|$.

\medskip
Assume now that $I_{n-1}$ is terminating and $|I_{n}^\infty|\asymp
  |\mathcal{R}(I_{n-1})|.$ Hence, $|I_n|$ is comparable to $|I_{n-1}^\infty|$.
 If $|I_{n-1}|$ is comparable to $|I_{n-1}^\infty|,$
then from the cases when $j_0\geq 1$, 
we have that the return mapping to $I_{n-1}$ extends to a quasi-box
mapping $F'_Q:U'_Q\rightarrow V'_Q$,
with $|V'_Q\cap\mathbb R|$ comparable to $|I_{n-1}^\infty|$.
Moreover we have that there exists an interval 
$K$ well-inside of $I_{n-1}$ and $\theta\in(0,\pi/2)$ so that
$K(F'_Q)\subset D_{\theta}(K)$.
On the other hand, if
$|I_{n-1}|/|I_{n-1}^\infty|>\nu$, then by part (a) of the claim, applied to $I_{n-1}$,
there is a qc polynomial-like mapping
$F\colon U\rightarrow V$
that extends the return mapping to $I_{n-1}^\infty$ with
$|V\cap\mathbb R| \asymp |I_{n-1}^\infty|$. 
As in the case when $I_{n-1}$ is non-terminating,
we can obtain a quasi-box mapping $F_Q\colon U_Q\rightarrow V_Q$
that extends the return map
$f^{j_i}\colon I_n^1\rightarrow I_n$ 
as an iterate of $F'_Q:U_Q'\rightarrow V'_Q$ or of 
$F:U\rightarrow V$.
Moreover,  there exists $\rho_1>0$ so that $I_n^j$ is $\rho_1$-free.
The proof in the bounded geometry case is analogous, one just needs to consider $V\cap \mathbb R=I^{2}_{n-1}.$
We will explain how to use the arguments for
Proposition~\ref{terminating angle control} and
Proposition~\ref{prop:from renormalization to landing} to pullback through one more step in the enhanced nest. Recall we are assuming $I_n$ is terminating and
$|I_{n+1}|$ is much smaller than $|I_n|.$

Since $V\cap\mathbb R$ is $\rho_1$-free,
 we can use the argument from
Step 1 of Proposition~\ref{terminating angle control} (compare 
Corollary~\ref{cor:terminating angle control})
to find a $k_1$, $j_i\leq k_1\leq 4j_i$ (see Combinatorial Remark 3), an interval $K_1$ well-inside
$I_{n-1}^{\infty}$
and $\theta_1\in(0,\pi)$ so that

$$\comp_{c_0} f^{-(k_0j_i-k_1)}(V_Q)\subset D_{\theta_1}(K_1).$$

Now we argue as in Steps 2 and 3 of Proposition~\ref{terminating angle
control}. It is important to notice that since
$\mathcal C(\mathcal R(I_{n-1}))\subset\mathcal{L}_{c_0}\mathcal
R(I_{n-1})=I_n$ is terminating, 
$k(X,\mathbb G)=\infty$, where $X$ is defined as in
Step 2 of Proposition~\ref{terminating angle
control}. Because of this, we do not need to rely on having bounded
geometry at a deeper level in order to control the loss of angle,
since we do have control of the combinatorial depth in the
parts of the
chain between critical points,
see Case (ii) of the proof of  Proposition~\ref{folding lemma}.
Thus we 
obtain $k_2,$ $j_i\leq k_2\leq 7j_i,$ see Combinatorial Remark 6,
an interval $K_2$, well-inside $\mathcal R(I_{n-1})$ and 
and $\theta_2\in(0,\pi)$ such that
$$\comp_{c_0} f^{-(k_0j_i-k_2)}(V_Q)\subset D_{\theta_2}(K_2).$$

We are assuming $\mathcal{R}(I_{n-1})$ is not a periodic interval and that
$|I_n^\infty|$ is comparable to $|\mathcal R(I_{n-1})|.$
We use the same argument as the one used in Proposition~\ref{delta free} to prove that $\mathcal{R}(I_{n-1})$
is $\rho'$-free for some $\rho'>0.$ If $\mathcal{R}(I_{n-1})$ is not periodic and
terminating, then Proposition~\ref{delta free} (1) implies that there
exists $\rho'>0$ so that
$\mathcal{R}(I_{n-1})$
is $\rho'$-free. Now, assume $\mathcal{R}(I_{n-1})$
is non-terminating. 
Observe that
$$I_{n}^\infty\subset \Gamma^2(\mathcal{R}(I_{n-1}))\subset
I_{n-1}^{\infty},$$ and compare this with Equation~(\ref{eqn:free3}).
The proof of Proposition~\ref{delta
    free} (3) applies to this situation. Since $|I_{n}^\infty|$ is comparable to
  $|I_{n-1}^\infty|$,  there exists $\rho'>0$ so that $\mathcal{R}(I_{n-1})$
is $\rho'$-free.
Finally, we repeat the argument used to prove
Proposition~\ref{prop:from renormalization to landing} to obtain
$k_3$, $0\leq k_3\leq 9j_i$ and $\theta_3\in(0,\pi)$ 
so that $$\comp_{c_0}f^{-(k_0j_i-k_3)}(V_Q)\subset D_{\theta_3}(I_n).$$
Observe that, since $\mathcal R(I_{n-1})$ is $\rho'$-free there exists $\rho''>0$ so that $I_n$ is $\rho''$-free. So we can apply Lemma~\ref{chain pullback},
to obtain $\theta\in(0,\pi)$ and an interval $K$ well-inside $I_n$ so that
 $$\comp_{c_0}f^{-(k_3+1)j_i}(V_Q)\subset D_{\theta}(K),$$
with $I_n^\infty\subset K.$ Hence, $|K|\asymp |I_n^\infty|.$ Since, $K(F_Q)\subset D_{\theta}(K)\cap U_Q$ the result follows.

\medskip

\noindent\textit{Case B: Assume $J_i= I_n$.}

In this case,  $I_n=\mathcal{R}(I_{n-1})=J_i$ is a periodic interval.
Let $T_i\supset J_i$ be the
largest interval so that $f^{j_i}$ has no critical points in
$T_i\setminus J_i$, and let $S_i=\comp_{c_0}f^{-j_i}(T_i)$.
By Lemma~\ref{bounds1},
since, $\partial J_i=\{\alpha,\tau(\alpha)\},$ and $j_i=2j_{i-1}$
where $\alpha$ is the orientation reversing fixed point closest to
$c_0$ of $f^{j_{i-1}}$,
there exists $\eta>0$
such that $(1+2\eta)J_i\subset S_i\subset(1+\eta)S_i\subset T_i$,
compare \cite[Lemma 8.7]{Shen}. In the case when $I_n=J_i$,
the interval $K$ will always be $S_i$.

\medskip
\noindent\textit{Case B.1: Bounded geometry.}

First suppose that for no $j,$ $0\leq j\leq M+1$
that $|I_{n-j}|/|I_{n-j+1}|>\nu.$ 
Following the proof of the Proposition~\ref{quasi-box mapping},
one sees that it applies just as well to the return
mapping $f^{j_i}\colon S_i\rightarrow T_i$ as it does to
$f^{j_i}\colon J_i\rightarrow J_i$.
This means we can construct 
$F_Q\colon U_Q\rightarrow V_Q$ a qc quasi-box mapping
that extends the return mapping to $T_i$, where
$V_Q=D_{\pi/2}(I_{n-M-1})\cap\mathbb C_{T_i}$. 
Notice that the chain 
$\{G_i\}^{j_i}_{i=0}$ with $G_{j_i}=T_i$ and
$G_0=\Comp_{c_0}f^{-j_i}(T_i)=S_i$ has intersection multiplicity two.

Let us show that there exists $\theta_1\in(0,\pi/2)$ so that
$K(F_Q)\subset D_{\theta_1}(S_i)$.
To see this, we consider two cases separately.
First, suppose that for every $j,$ $0\leq j\leq M$,
$I_{n-j-1}$ is terminating, and
$I_{n-j}$ is of period two under $R_{I_{n-j-1}}$.
In this case, we have that each interval
$I_{n-j}=J_{i-j},$ a periodic interval.
We have that
$f^{j_{i-M}}|_{S_{i-M}}:S_{i-M}\rightarrow T_{i-M}$,
and the chain,
$\{\tilde G_j\}_{j=0}^{j_{i-M}}$ with
$\tilde G_{j_{i-M}}=I_{i-M-1}$ and
$\tilde G_0=I_{i-M-1}$ has order $2b$.
Thus by the argument in Step 1 of the proof of
Proposition~\ref{terminating angle control},
we have that there exists $\theta\in(0,\pi)$
such that
$\comp_{c_0}f^{-{j_{i-M}}}(V_Q)\subset
D_{\theta}(S_{i-M}).$
Repeating this argument for the return map to each
$J_{i-j}$, $0\leq j\leq M,$ we have that there exists $k_0,0<k_0\leq M$,
and $\theta_1\in(0,\pi/2)$,
such that
$$\comp_{c_0} f^{-j_i k_0}(V_Q)\subset D_{\theta_1}(S_i).$$
Thus $K(F_Q)\subset D_{\theta_1}(S_i).$
Alternatively, there exists $i_0,$ $0< i_0\leq M$ minimal so that
$I_{n-i_0}$ is terminating, $I_{n-i_0}$ does not have period two under
$R_{I_{n-i_0-1}}$ and for all $j<i_0$, $I_{n-j}$ does have period two under
$R_{I_{n-j-1}}.$ Let $F_Q':U_Q'\rightarrow V_Q'$ denote the quasi-box mapping that extends the return mapping to $I_{n-i_0}$.
By Case A of this proof,
we have that there exists $\theta'\in(0,\pi/2)$,
and an interval $K'$ with $|K'|\asymp |I_{n-i_0}|$ such that
such that $K(F_Q')\subset D_{\theta'}(K').$
Now we can argue as in the case when all of the returns to
$I_{n-j}$ were of period two under $R_{I_{n-j-1}}$ to see that
there exists an angle $\theta_1\in(0,\pi/2)$ such that
$K(F_Q)\subset D_{\theta_1}(S_i)$.
Thus $K(F_Q)\subset D_{\theta_1}(S_i)\cap D_{\pi/2}(I_{n-M})$.

\medskip
\noindent\textit{Case B.2: Unbounded geometry.}

Assume there exists $j,$ $0\leq j\leq M+1$ such that
$|I_{n-j}|/|I_{n-j+1}|>\nu$. Let $j_0\geq 0$ be minimal so that $|I_{n-j_0}|/|I_{n-j_0+1}|>\nu$.

\medskip

\noindent\textit{Case B.2.1: Suppose $j_0\geq 2$.} 
There exists
$i_0,$ $0< i_0\leq j_0-1$ minimal so that
$I_{n-i_0}$ is terminating, $I_{n-i_0}$ does not have period two under
$R_{I_{n-i_0-1}}$ and for all $j<i_0$, $I_{n-j}$ does have period two under
$R_{I_{n-j-1}}.$ Let $F_Q':U_Q'\rightarrow V_Q'$ denote the quasi-box
mapping that extends the return mapping to $I_{n-i_0}$, given by Case
A,
and repeat the proof of Case B.1 to obtain that
there exists
a $\theta\in(0,\pi)$ so that
$$K(F_Q)\subset D_{\theta}(S_i).$$

\medskip

\noindent\textit{Case B.2.2: Suppose $j_0=1,$ or $j_0=0$ and $|I_{n-1}|$ is much bigger than
$|I_n|.$}
Assume $j_0=1.$  By Lemma~\ref{bounds1}, we know that
$|\mathcal{R}(I_{n-1})|$ is comparable to $|I_{n-1}^\infty|.$
Furthermore, since we are assuming $|I^{\infty}_n|$ 
is comparable to $|\mathcal R(I_{n-1})|,$
then by Lemma~\ref{bounds1},
we have that $|I_{n-1}^2|$ is much smaller than $|I_{n-1}|$.
By the claim,we can find a polynomial-like mapping 
$F\colon U\rightarrow V$ that extends the return mapping to
$I_{n-1}^\infty$ so that
$V\cap\mathbb R=I_{n-1}^{j}$ and $|I_{n-1}^{j}| \asymp |I_{n-1}^\infty|.$
We obtain a quasi box mapping that extends the
mapping $f^{j_i}:J_i\rightarrow J_i$, 
by setting $V_Q=V\cap\mathbb C_{T_i}$, 
$U_Q=\comp_{c_0}f^{-j_i}(V_Q)$
and defining $F_Q|U_Q=F^2$.
Since the boundary points of $T_i$ are critical values of $f^{j_i}$
and $|T_i|$ is comparable to $I_{n-1}^{i_0}$,
by the argument from Step 1 of the proof of Proposition~\ref{terminating angle control}
there exists $\theta\in(0,\pi)$ such that 
$$\comp_{c_0}f^{-j_i}(V_Q)\subset D_{\theta}(S_i).$$
Notice that this argument does not depend on $|I_n|$ and $|I_{n+1}|$
being comparable, so we also use it to cover the case when $j_0=0$
and $|I_{n-1}|$ is much bigger than
$|I_n|.$

\medskip

\noindent\textit{Case B.2.3: Suppose $j_0=0$ and $|I_{n-1}|$ is
  comparable to $|I_n|.$}
This implies that $|S_{i-1}|$ is comparable to $|S_i|.$
We use the argument from either the bounded geometry case or the case when $j_0\geq 2$
with $n-1$ in place of $n$
to obtain 
a qc quasi-box mapping $F'_Q\colon U_Q'\rightarrow V_Q'$ that extends the
return mapping to $I_{n-1}$, and
an angle $\theta_1\in(0,\pi)$ such that 
$$K(F_Q')\subset D_{\theta_1}(S_{i-1}).$$
We obtain $F_Q\colon U_Q\rightarrow V_Q$ extending the
first return mapping to $J_i$ as the second iterate of $F_Q',$
just as in Case B.2.2.
Now, since $|S_{i-1}|$ is comparable to $|S_i|,$
by the argument from Step 1 of Proposition~\ref{terminating angle control}
there exists $\theta_2\in(0,\pi)$
so that
$$K(F_Q)\subset \comp_{c_0}f^{-j_i}(V_Q)\subset D_{\theta_2}(S_i).$$ 

\medskip
\noindent\textit{Step 2: Controlling the dilatation.}
Since there exist constants $C_1$ and $\rho''>0$ such that
$$\diam J_i\leq \diam D_{\theta}(K) <C_1|J_i|,$$
and $(K\setminus(1+2\rho''')^{-1}K)$ is disjoint from
the post-critical set,
there exists a constant $C_2>0$ such that 
the dilatation of 
$f^{j_i}|D_{\theta}(K)$ is bounded from above by
$C_2\sum_{k=0}^{j_i}|\mathcal{L}_{f^k(c_0)}(J_i)|^2$.

Now, by Lemma~\ref{lem:LY}, we have that there exists a 
qc polynomial-like mapping
$F_i\colon U_i\rightarrow V_i$ that extends the return mapping to $J_i$
such that $\mod(V_i\setminus U_i)$ is bounded from below.

\medskip
\noindent\textit{Step 3: $U_i$ has bounded geometry at $c_0$.}
Since $\mod(V_i\setminus K(F_i))$ is bounded from below and $V_i\supset
K(F_i)\supset J_i$, it follows, for example from \cite[Proposition 4.8]{McMullen-3mflds},
that there exists
$\rho>0$ so that the ball with radius $\rho$ centered at $F_i(
c_0 )$ is contained in $V_i.$
If the biggest $\rho$ for which is true is small, then $V_i$ comes very near $J_i$ but still contains $J_i$, which means that
the modulus of $V_i\setminus K(F_i)$ is necessarily small.
Then since $\mod(V_i\setminus U_i)$ is bounded from below, it follows that
there also exists a $\rho'(\rho)>0$ so that the ball with radius $\rho'$ centered at $c_0$ is contained in $U_i$.


\end{pf}

\end{document}